\documentclass[11pt,reqno]{amsart}

\usepackage[english]{babel}
\usepackage{amsmath,amssymb,amsthm,mathtools,mathrsfs,microtype}
\usepackage{enumitem,booktabs,array,longtable,needspace,fancyvrb,listings}
\usepackage{xcolor}
\usepackage{tikz}
\usetikzlibrary{arrows.meta}
\usepackage{cite}
\usepackage[hypertexnames=false]{hyperref}
\allowdisplaybreaks
\setlist{topsep=0.4em,itemsep=0.25em,parsep=0.1em}

\newcommand{\E}{\mathbb E}

\newtheorem{theorem}{Theorem}[section]
\newtheorem{proposition}[theorem]{Proposition}
\newtheorem{lemma}[theorem]{Lemma}
\newtheorem{corollary}[theorem]{Corollary}
\theoremstyle{remark}
\newtheorem{remark}[theorem]{Remark}
\numberwithin{equation}{section}

\hypersetup{
  colorlinks=true,
  linkcolor=black,
  citecolor=blue,
  urlcolor=blue,
  pdftitle={On the Gardner transition in the Ising pure p-spin glass II},
  pdfauthor={Yuxin Zhou},
  bookmarksdepth=3
}

\makeatletter
\let\reftagform@=\tagform@
\def\tagform@#1{\maketag@@@{(\ignorespaces\textcolor{magenta}{#1}\unskip\@@italiccorr)}}
\renewcommand{\eqref}[1]{\textup{\reftagform@{\ref{#1}}}}
\makeatother

\title{On the Gardner Transition in the Ising Pure $p$-Spin Glass II}
\author{Yuxin Zhou}
\date{}

\begin{document}
\begin{abstract}
In \cite{zhou}, we identified, for every $p\geq3$, a unique first
critical inverse temperature $\beta_1^p$, proved that the Parisi measure
is replica symmetric (RS) for $0<\beta\leq\beta_1^p$, and proved that it is
one-step replica symmetry breaking (1-RSB) on a nonempty interval immediately
above $\beta_1^p$.  In this sequel, we determine the rest of the phase
diagram.  There is a unique second critical inverse temperature
$\beta_2^p>\beta_1^p$.  The Parisi measure is 1-RSB for $\beta_1^p<\beta\leq\beta_2^p$, while for
$\beta>\beta_2^p$ there are $0<q_\beta<q'_\beta<1$ such that
\[
 \operatorname{supp}\mu_\beta=\{0\}\cup[q_\beta,q'_\beta].
\]
On the interior of the interval the measure has a smooth density, and hence it is full
replica symmetry breaking (FRSB).  Together
with \cite{zhou}, this proves the two transitions predicted by Gardner
for the Ising pure $p$-spin glass. 
\end{abstract}

\maketitle

\setcounter{tocdepth}{1}
\tableofcontents
\clearpage

\section{Introduction and main results}
\label{sec:introduction}

The Ising pure $p$-spin glass is one of the basic mean field spin glass
models.  We consider the model without external field and with
$p\geq3$.  For background on the model and the Parisi theory, see
\cite{mezardparisivirasoro}.  For $N\geq1$, set
$\Sigma_N=\{-1,+1\}^N$.  For 
$\sigma=(\sigma_1,\ldots,\sigma_N)\in\Sigma_N$, the Hamiltonian at finite temperature $\beta$ is
\begin{equation}
 H_N(\sigma)=\frac{\beta}{N^{(p-1)/2}}
 \sum_{1\leq i_1,\ldots,i_p\leq N}
 g_{i_1,\ldots,i_p}\sigma_{i_1}\cdots\sigma_{i_p},
 \label{eq:Hamiltonian}
\end{equation}
where the $g_{i_1,\ldots,i_p}$ are independent standard Gaussian random
variables and $\beta>0$ is the inverse temperature.  For two
configurations, let
\[
 R_{1,2}=\frac1N\sum_{i=1}^N\sigma_i^1\sigma_i^2,
\]
which is their normalized overlap.  The Gaussian field is centered and
has covariance
\begin{equation}
 \E H_N(\sigma^1)H_N(\sigma^2)=N\xi(R_{1,2}),
 \qquad \xi(u)=\beta^2u^p.
 \label{eq:xi}
\end{equation}

The limiting free energy is given by the Parisi variational formula
\cite{talagrandparisi,panchenko}.  For odd $p$, we use Panchenko's
extension.  The Parisi functional is defined on
the probability measures on $[0,1]$ and has a unique minimizer
\cite{auffingerchenunique}.  This minimizer is called the Parisi measure at temperature $\beta$
and is denoted by $ \mu_\beta$.  The structure of its support is crucial for characterizing the mean field spin glass models.

We say that the model is replica symmetric at $\beta$, abbreviated RS, if $\mu_\beta$ is a
Dirac measure. We say that the model is  $k$-step replica symmetry
breaking at $\beta$, abbreviated $k$-RSB, if the support of $\mu_\beta$ consists of $k+1$ atoms. We say that the Parisi measure is full replica symmetry breaking at $\beta$,
abbreviated FRSB, if the support of $\mu_\beta$ contains a nontrivial interval. By \cite[Theorem~2(ii)]{auffingerchen}, its distribution function is
then infinitely differentiable on the interior of that interval, so
the measure has a smooth density there. In the present model without external field, the origin is always in the support of $\mu_\beta$
\cite[Theorem~1(i)]{auffingerchen}. Therefore RS means
$\mu_\beta=\delta_0$, while 1-RSB means
\[
  \mu_\beta=m\delta_0+(1-m)\delta_q, \text{ for some } 0<m,q<1.
\]

Gardner predicted two transitions for this model \cite{gardner}.
As $\beta$ increases, the Parisi measure should first leave the RS phase,
then remain 1-RSB for an interval, and finally enter the FRSB phase.
Chen \cite{Chen} identified the high-temperature threshold through the
free energy.  In the present normalization, Theorem~1 of \cite{zhou}
identifies the exact RS boundary and proves that the Parisi measure is
1-RSB on a nonempty interval immediately above it.  It also provides a
computational method to locate the first critical temperature.  In
particular, this gives the first rigorous
finite-temperature example of 1-RSB in a mean field Ising spin glass.

The present paper is a direct sequel to \cite{zhou}.   
The two regions left undetermined in \cite[Figure~1]{zhou} were the
remainder of the conjectured 1-RSB phase and the FRSB phase.  We retain the notation and the framework
of \cite{zhou}.  We determine
both regions here and thereby confirm Gardner's prediction. 

\begin{theorem}
\label{thm:main}
For every $p\geq3$, there are unique numbers
\[
 0<\beta_1^p<\beta_2^p<\infty
\]
such that the model is:
\begin{enumerate}
\item  RS, i.e. $\mu_\beta=\delta_0$ if and only if
      $0<\beta\leq\beta_1^p$.
\item   1-RSB, i.e. $\mu_\beta=m \delta_0 +(1-m)\delta_q$
      if and only if
      $ \beta_1^p<\beta\leq\beta_2^p.$
\item  FRSB if and only if $\beta>\beta_2^p$.  More precisely,  there are $0<q_\beta<q'_\beta<1$ such that
      $ \operatorname{supp}\mu_\beta=\{0\}\cup[q_\beta,q'_\beta]$.
\end{enumerate}
\end{theorem}


\begin{figure}[htbp]
\centering
\begin{tikzpicture}[
  x=1cm,y=1cm,
  >=Stealth,
  every node/.style={font=\small}
]
  \draw[->,line width=.75pt] (0,0)--(12.4,0)
       node[below right] {$\beta$};
  \draw[line width=1.15pt] (0,0)--(4.0,0);
  \draw[line width=1.15pt] (4.0,0)--(8.1,0);
  \draw[line width=1.15pt] (8.1,0)--(11.9,0);
  \fill (0,0) circle (1.8pt);
  \fill (4.0,0) circle (1.8pt);
  \fill (8.1,0) circle (1.8pt);
  \node[above] at (2.0,.08) {\textbf{Replica symmetric}};
  \node[above] at (6.05,.08) {\textbf{1-RSB}};
  \node[above] at (10.0,.08) {\textbf{FRSB}};
  \node[below] at (0,-.02) {$0$};
  \node[below] at (4.0,-.02) {$\beta_1^p$};
  \node[below] at (8.1,-.02) {$\beta_2^p$};
\end{tikzpicture}
\caption{Phase transitions of the Parisi measure $ \mu_\beta$ with respect
to $\beta$.  $\mu_\beta$ is RS for
$0<\beta\leq\beta_1^p$, 1-RSB for
$\beta_1^p<\beta\leq\beta_2^p$, and FRSB for
$\beta>\beta_2^p$.  The exact RS boundary and a nonempty 1-RSB interval
to its right were established in \cite[Theorem~1]{zhou}. The present
sequel identifies $\beta_2^p$ and completes
the phase diagram.  }
\label{fig:phase-diagram}
\end{figure}
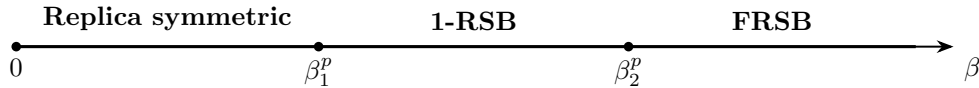

\begin{figure}[htbp]
\centering
\begin{tikzpicture}[
  x=1cm,y=1cm,
  >=Stealth,
  every node/.style={font=\scriptsize}
]
\begin{scope}[xshift=0cm]
  \draw[->,line width=.65pt] (0,0)--(4.0,0) node[below] {$u$};
  \draw[->,line width=.65pt] (0,0)--(0,2.45)
       node[above] {$\alpha_\beta(u)$};
  \node[above] at (1.9,2.48)
       {$\beta_1^p<\beta<\beta_2^p$};
  \node[left] at (0,2) {$1$};
  \node[left] at (0,.82) {$m$};
  \node[below] at (2.1,0) {$q$};
  \node[below] at (3.55,0) {$1$};
  \draw[line width=1.05pt] (0,.82)--(2.1,.82);
  \draw[line width=1.05pt] (2.1,2)--(3.55,2);
  \draw[densely dashed] (2.1,0)--(2.1,2);
  \draw[densely dashed] (3.55,0)--(3.55,2);
  \fill (0,.82) circle (1.5pt);
  \draw[fill=white,line width=.65pt] (2.1,.82) circle (1.7pt);
  \fill (2.1,2) circle (1.5pt);
  \fill (3.55,2) circle (1.5pt);
\end{scope}

\begin{scope}[xshift=5.1cm]
  \draw[->,line width=.65pt] (0,0)--(4.0,0) node[below] {$u$};
  \draw[->,line width=.65pt] (0,0)--(0,2.45)
       node[above] {$\alpha_\beta(u)$};
  \node[above] at (1.9,2.48) {$\beta=\beta_2^p$};
  \node[left] at (0,2) {$1$};
  \node[left] at (0,.82) {$m_2^p$};
  \node[below] at (2.1,0) {$q_2^p$};
  \node[below] at (3.55,0) {$1$};
  \draw[line width=1.05pt] (0,.82)--(2.1,.82);
  \draw[line width=1.05pt] (2.1,2)--(3.55,2);
  \draw[densely dashed] (2.1,0)--(2.1,2);
  \draw[densely dashed] (3.55,0)--(3.55,2);
  \fill (0,.82) circle (1.5pt);
  \draw[fill=white,line width=.65pt] (2.1,.82) circle (1.7pt);
  \fill (2.1,2) circle (1.5pt);
  \fill (3.55,2) circle (1.5pt);
  \node[text=red,font=\large] at (2.1,1.45) {$\times$};
\end{scope}

\begin{scope}[xshift=10.2cm]
  \fill[gray!12] (1.05,0) rectangle (2.65,2.12);
  \draw[->,line width=.65pt] (0,0)--(4.0,0) node[below] {$u$};
  \draw[->,line width=.65pt] (0,0)--(0,2.45)
       node[above] {$\alpha_\beta(u)$};
  \node[above] at (1.9,2.48) {$\beta>\beta_2^p$};
  \node[left] at (0,2) {$1$};
  \node[below] at (1.05,0) {$q_\beta$};
  \node[below] at (2.65,0) {$q'_\beta$};
  \node[below] at (3.55,0) {$1$};
  \fill (0,.42) circle (1.5pt);
  \draw[line width=1.05pt] (0,.42)--(1.05,.42);
  \draw[densely dotted,gray,line width=1pt] (1.05,.42)--(1.05,.66);
  \draw[line width=1.15pt]
       (1.05,.66) .. controls (1.45,.80) and (2.18,1.52) .. (2.65,1.69);
  \draw[densely dotted,gray,line width=1pt] (2.65,1.69)--(2.65,2);
  \draw[line width=1.05pt] (2.65,2)--(3.55,2);
  \draw[densely dashed] (1.05,0)--(1.05,.42);
  \draw[densely dashed] (2.65,0)--(2.65,1.69);
  \draw[line width=2.1pt] (1.05,.04)--(2.65,.04);
  \fill (1.05,.04) circle (1.5pt);
  \fill (2.65,.04) circle (1.5pt);
\end{scope}
\end{tikzpicture}
\caption{The distribution of the Parisi measure across the second
transition.  From left to right: $\mu_\beta$ is 1RSB for
$\beta_1^p<\beta<\beta_2^p$;  1-RSB marginally  at
$\beta=\beta_2^p$; and FRSB for
$\beta>\beta_2^p$.  The first two panels give the exact distribution
functions.  In the last panel,
$\operatorname{supp}\mu_\beta=\{0\}\cup[q_\beta,q'_\beta]$ with
$0<q_\beta<q'_\beta<1$.  Thus the distribution function jumps at the
isolated origin, is constant on $(0,q_\beta)$, is infinitely
differentiable on $(q_\beta,q'_\beta)$ by
\cite[Theorem~2(ii)]{auffingerchen}, and is constant after
$q'_\beta$. }
\label{fig:second-transition-cdfs}
\end{figure}
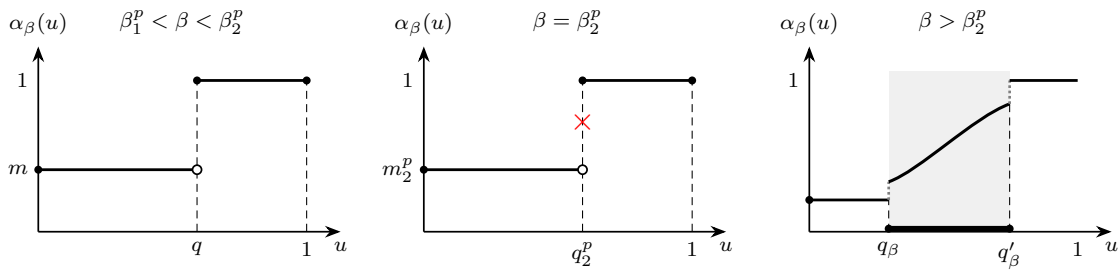

For a direct comparison with the distribution-function picture in
\cite[Figure~2]{zhou}, write
\[
 \alpha_{\beta}(u)= \mu_\beta([0,u]),\qquad 0\leq u\leq1.
\]
Figure~\ref{fig:second-transition-cdfs} describes the second transition.

The first assertion of Theorem~\ref{thm:main} and the existence of a
nonempty 1-RSB interval above $\beta_1^p$ are
\cite[Theorem~1]{zhou}.  The first new step is to determine the entire
1-RSB interval. We apply the criterion introduced in
\cite[Theorem~3]{zhou}, where a function $f_\mu$ with respect to $\mu$
is used to determine whether $\mu$ is the Parisi measure.  For a measure
$\mu$ with $\operatorname{supp}\mu=\{0,q\}$, we prove the strict
sign of $f_\mu$ on both sides of its positive atom whenever
$\Gamma_\mu'(q)<1$.  At equality, we prove that $\Gamma_\mu'(q)$
crosses $1$ with the increase of temperature $\beta$.
This gives the unique endpoint $\beta_2^p$ and excludes both RS and
1-RSB above it.

The second new step is to exclude a bounded gap between two positive
support points of a Parisi measure.  Once the Parisi measure is known to
be neither RS nor 1-RSB, its support has two distinct positive points;
the no-gap theorem forces the interval between them into the support, which implies the exact form
$\{0\}\cup[q_\beta,q'_\beta]$.  This is the FRSB conclusion.
Recent works of Lopatto \cite{lopatto} and Chen \cite{chenfrsb} establish
the FRSB phases for the Sherrington--Kirkpatrick model.
The present proof concerns the different pure $p$-spin geometry with
$p\geq3$.

The remainder of the paper is organized as follows. Section~\ref{sec:second-transition} derives the equations determining the second transition. Section~\ref{sec:sign} determines the complete 1-RSB interval. Section~\ref{sec:positive-gaps} rules out positive internal gaps and completes the proof of Theorem~\ref{thm:main}.


\subsection*{Acknowledgments}
The proofs in the appendix were drafted by large language models and have not yet received their final authorial revision.  The author will carefully verify, revise, and rewrite these proofs in a subsequent version and assumes full responsibility for their correctness and presentation.

\section{The second transition of the Parisi measure}
\label{sec:second-transition}

We first recall how the first critical point was located in \cite{zhou}.
With
\[
 \theta(q)=q\xi'(q)-\xi(q),\qquad Y_q=\sqrt{\xi'(q)},
\]
the pair $(\beta_1^p,q_1^p)$ is the unique solution in
$(0,\infty)\times(0,1)$ of
\begin{equation}
\left\{
\begin{aligned}
 &\frac{\E[\cosh(Y_qg)\log\cosh(Y_qg)]}
        {\E\cosh(Y_qg)}
  -\frac12\xi'(q)-\frac12\theta(q)=0,\\
 &\frac{\E[\tanh^2(Y_qg)\cosh(Y_qg)]}
        {\E\cosh(Y_qg)}-q=0.
\end{aligned}
\right.
\label{eq:first-boundary}
\end{equation}
Theorem~1 in \cite{zhou} identifies \(\beta_1^p\) as the exact boundary of the RS regime and shows that the 1-RSB phase emerges immediately beyond it. This result serves as the starting point of the present sequel.

The first transition occurs when the two atoms of the 1RSB measure degenerates
at $m=1$.  At the second transition, the two atoms remain distinct and
the loss of the 1-RSB phase is detected by the equality
$\Gamma_\mu'(q)=1$.  We now make this characterization precise.

We use the notation of \cite{auffingerchen,zhou}.  For a probability
measure $\mu$ on $[0,1]$, let $\alpha_\mu(u)=\mu([0,u])$ and let
$\Phi_\mu$ solve
\begin{equation}
 \partial_u\Phi_\mu(x,u)
 =-\frac{\xi''(u)}2\left(
 \partial_{xx}\Phi_\mu(x,u)
 +\alpha_\mu(u)(\partial_x\Phi_\mu(x,u))^2\right),
 \qquad \Phi_\mu(x,1)=\log\cosh x.
 \label{eq:pde}
\end{equation}
If $B$ is a standard Brownian motion, set
\begin{equation*}
 M(u)=B(\xi'(u)),\qquad
 W_\mu(u)=\int_{[0,u]}
 \{\Phi_\mu(M(u),u)-\Phi_\mu(M(t),t)\}\,d\mu(t),
 \label{eq:MW}
\end{equation*}

The exponential weight is normalized:
\begin{equation}
 \E e^{W_\mu(u)}=1,\qquad 0\leq u\leq1.
 \label{eq:weight-normalization}
\end{equation}
For a finite-step measure this follows by iterated Gaussian conditioning
in the Cole--Hopf representation; the general case follows by step
approximation.  This is also \cite[Proposition~2]{auffingerchen}.

Now define
\begin{equation}
 \Gamma_\mu(u)=
 \E\left[(\partial_x\Phi_\mu(M(u),u))^2e^{W_\mu(u)}\right].
 \label{eq:Gamma}
\end{equation}
The basic differentiation formula is
\begin{equation}
 \Gamma_\mu'(u)=\xi''(u)
 \E\left[(\partial_{xx}\Phi_\mu(M(u),u))^2e^{W_\mu(u)}\right].
 \label{eq:Gamma-prime}
\end{equation}
Following \cite{zhou}, define
\begin{equation}
 F_\mu(u)=\Gamma_\mu(u)-u,
 \qquad
 f_\mu(u)=\frac12\int_0^u\xi''(t)F_\mu(t)\,dt.
 \label{eq:f}
\end{equation}
The variational criterion \cite[Theorem~3]{zhou} says that $\mu= \mu_\beta$ if and only if
\begin{equation}
 f_\mu(u)\leq0 \text{ for } u \in [0,1],
 \qquad
 f_\mu(u)=0 \text{ for } u\in\operatorname{supp}\mu.
 \label{eq:criterion}
\end{equation}

We now take
\begin{equation*}
 \mu=m\delta_0+(1-m)\delta_q,\qquad 0<m,q<1,
 \label{eq:mu}
\end{equation*}
and write $Y_{u-t}=\sqrt{\xi'(u)-\xi'(t)}$ for $u>t\geq0$ and
$Y_u=\sqrt{\xi'(u)}$.  The corresponding Parisi functional can be written as
\begin{align}
 \mathcal P(\mu)
 ={}&\log2+\frac1m\log\E\cosh^m(Y_qg)
 +\frac12\{\xi'(1)-\xi'(q)\}\notag\\
 &-\frac12\{q\xi'(q)-\xi(q)\}m
 -\frac12\int_q^1u\xi''(u)\,du,
 \label{eq:restricted-P}
\end{align}
where $g$ is a standard Gaussian random variable.  Differentiating in $m$ and
$q$ gives the two stationary equations from \cite[(15)]{zhou}:
\begin{equation}
\left\{
\begin{aligned}
 C_\beta^1(m,q)
 :={}&-\frac1{m^2}\log\E\cosh^m(Y_qg)
 +\frac1m\frac{\E[\cosh^m(Y_qg)\log\cosh(Y_qg)]}
 {\E\cosh^m(Y_qg)}
 -\frac12\{q\xi'(q)-\xi(q)\}=0,\\
 D_\beta^1(m,q)
 :={}&\frac{\E[\tanh^2(Y_qg)\cosh^m(Y_qg)]}
 {\E\cosh^m(Y_qg)}-q=0.
\end{aligned}\right.
\label{eq:stationary}
\end{equation}
Indeed,
\begin{equation*}
 \frac{\partial}{\partial m}\mathcal P(\mu)=C_\beta^1(m,q)=f_\mu(q),
 \qquad
 \frac{\partial}{\partial q}\mathcal P(\mu)
 =-\frac{1-m}{2}\xi''(q)D_\beta^1(m,q).
 \label{eq:stationary-derivatives}
\end{equation*}
Thus \eqref{eq:stationary} is exactly
\begin{equation}
 f_\mu(q)=0,\qquad \Gamma_\mu(q)=q.
 \label{eq:contacts}
\end{equation}

At the right endpoint of the support $u=q$, \eqref{eq:Gamma-prime} takes the
explicit form
\begin{equation}
 \Gamma_\mu'(q)=\xi''(q)
 \frac{\E\cosh^{m-4}(Y_qg)}{\E\cosh^m(Y_qg)}.
 \label{eq:second-marginal}
\end{equation}
The second transition is then obtained by adding
$\Gamma_\mu'(q)=1$ to the two stationary equations as follows:

\begin{proposition}
\label{prop:second-boundary}
For every $p\geq3$, the system
\begin{equation}
\left\{
\begin{aligned}
 C_\beta^1(m,q)&=0\\
 D_\beta^1(m,q)&=0\\
 \beta^2p(p-1)q^{p-2}
 \frac{\E\cosh^{m-4}(Y_qg)}{\E\cosh^m(Y_qg)}&=1
\end{aligned}
\right.,
\label{eq:second-boundary}
\end{equation}
has a unique solution $(\beta,m,q)\in(0,\infty)\times(0,1)^2$.
We denote it by $(\beta_2^p,m_2^p,q_2^p)$.
\end{proposition}

Let us explain the role of the last equation.  By
\eqref{eq:contacts}, if a 1-RSB measure $\mu$ is the Parisi measure, then it satisfies
$f_\mu(q)=f_\mu'(q)=0$.  At this point, differentiation of
\eqref{eq:f} gives
\begin{equation}
 f_\mu''(q)=\frac{\xi''(q)}2\{\Gamma_\mu'(q)-1\}.
 \label{eq:second-curvature}
\end{equation}
If a measure $\mu$ is the Parisi measure, then $f_\mu\leq0$, so $q$ is a
local maximum and necessarily $\Gamma_\mu'(q)\leq1$.  The proof below
shows that, for a 1-RSB measure, strict inequality is
sufficient: it implies $f_\mu<0$ away from $0$ and $q$, and hence the
measure is the Parisi measure.

At equality, the stationary solution is marginal.  We prove that the
nearby stationary solutions form a locally unique curve and that,
along this curve,
\[
 \frac d{d\beta^2}\Gamma_\mu'(q)>0.
\]
Thus $\Gamma_\mu'(q)<1$ immediately below the marginal point and
$\Gamma_\mu'(q)>1$ immediately above it.  This identifies the unique
solution of \eqref{eq:second-boundary} as the right endpoint of the
1-RSB phase.  The support argument in
Section~\ref{sec:positive-gaps} turns the exclusion of RS and 1-RSB
above this endpoint into FRSB.

\section{Proof of 1-RSB Phase: $\beta^p_1<\beta \leq \beta^p_2$}\label{sec:sign}
\subsection{The sign of $f_\mu$}

In this section, we establish the 1-RSB phase for
$\beta_1^p<\beta\leq\beta_2^p$ by using $f_\mu$ introduced in
\cite[Theorem~3]{zhou}.
The next proposition is the only input needed on $(0,q)$.  Its full
proof is included in Appendix~\ref{app:contact}.

\begin{proposition}
\label{prop:left-computation}
Assume \eqref{eq:contacts} and $\Gamma_\mu'(q)<1$.  If
$u\in(0,q)$ and
\begin{equation}
 \Gamma_\mu(u)=u\Gamma_\mu'(u),
 \label{eq:left-contact}
\end{equation}
then
\[
 \Gamma_\mu''(u)<0.
\]
\end{proposition}

\begin{proposition}
\label{prop:left-sign}
Under the assumptions of Proposition~\ref{prop:left-computation},
\[
 f_\mu(u)<0,\qquad 0<u<q.
\]
\end{proposition}

\begin{proof}
The Cole--Hopf calculation in
Appendix~\ref{ft:sec:slope-coordinate} gives
$0\leq\partial_{xx}\Phi_\mu\leq1$ on $(0,q)$.  On $[q,1]$,
$\partial_{xx}\Phi_\mu(x,u)=\cosh^{-2}x$, so the same bound holds.
Together with
\eqref{eq:weight-normalization} and \eqref{eq:Gamma-prime}, this gives
$0\leq\Gamma_\mu'(u)\leq\xi''(u)$.  Moreover, $M(0)=0$ and
$\Phi_\mu(\cdot,0)$ is even, so $\Gamma_\mu(0)=0$.  Since $p\geq3$,
\begin{equation}
 \lim_{u\downarrow0}\frac{\Gamma_\mu(u)}u=0.
 \label{eq:left-limit}
\end{equation}
Direct differentiation gives
\begin{equation}
 \frac d{du}\frac{\Gamma_\mu(u)}u
 =\frac{u\Gamma_\mu'(u)-\Gamma_\mu(u)}{u^2}.
 \label{eq:ratio-first}
\end{equation}
At a critical point of $\Gamma_\mu(u)/u$, a second differentiation
and Proposition~\ref{prop:left-computation} give
\begin{equation}
 \frac {d^2}{du^2}\frac{\Gamma_\mu(u)}u
 =\frac{\Gamma_\mu''(u)}u<0.
 \label{eq:ratio-second}
\end{equation}
Thus every critical point is a strict maximum, and there can be at
most one critical point.  On the other hand,
\begin{equation}
 \frac{\Gamma_\mu(q)}q=1,
 \qquad
 \left.\frac d{du}\frac{\Gamma_\mu(u)}u\right|_{u=q-}
 =\frac{\Gamma_\mu'(q)-1}{q}<0.
 \label{eq:ratio-end}
\end{equation}
Together with \eqref{eq:left-limit}, this says that
$\Gamma_\mu(u)/u$ crosses $1$ once: it is below $1$ first and above
$1$ afterwards.  Consequently $F_\mu$ is first negative and then
positive.  By \eqref{eq:f}, $f_\mu$ first decreases and then
increases.  Since $f_\mu(0)=f_\mu(q)=0$, it is strictly negative on
$(0,q)$.
\end{proof}

We next prove the sign on the other side of $q$.  

\begin{proposition}
\label{prop:right-sign}
Assume \eqref{eq:contacts} and $\Gamma_\mu'(q)\leq1$.  Then
\begin{equation}
 \Gamma_\mu''(u)<0,\qquad
 \Gamma_\mu(u)<u,
 \qquad q<u\leq1.
 \label{eq:right-conclusion}
\end{equation}
In particular,
\[
 f_\mu(u)<0,\qquad q<u\leq1.
\]
\end{proposition}

Combining Propositions~\ref{prop:left-sign} and
\ref{prop:right-sign} with \eqref{eq:contacts}, we have our desired conclusion: If for $\mu=m\delta_0+(1-m)\delta_q$ satisfying that \eqref{eq:stationary} and $\Gamma_\mu'(q)<1$, then
\[
 f_\mu(u)\leq0 \text{ for } 0\leq u\leq1,
\text{ and }
 f_\mu(u)=0 \text{ for } u\in\{0,q\}.
\]
Consequently, $\mu$ is the Parisi measure.

We then prove Proposition \ref{prop:right-sign} as follows:
\begin{proof}[Proof of Proposition \ref{prop:right-sign}]
At $q$, \eqref{eq:Gamma-prime} and the endpoint form of the Parisi PDE
give
\begin{equation}
 \Gamma_\mu'(q)=\xi''(q)
 \frac{\E\cosh^{m-4}(Y_qg)}{\E\cosh^m(Y_qg)}.
 \label{eq:endpoint-replicon}
\end{equation}
For fixed $\xi'(q)$, differentiation in $m$ shows that
$ \frac{\E[\tanh^2(Y_qg)\cosh^m(Y_qg)]}{\E\cosh^m(Y_qg)}$ increases with $m$
and $ \frac{\E\cosh^{m-4}(Y_qg)}{\E\cosh^m(Y_qg)}$ decreases with $m$.
These signs are immediate covariances: $\tanh^2 x$ and
$\log\cosh x$ increase together as $|x|$ increases, whereas
$\cosh^{-4}x$ decreases as $\log\cosh x$ increases.  Since
$0<m\leq1$, the relations above imply
\begin{equation}
 \frac{\E\cosh^{m-4}(Y_qg)}
 {\E[\tanh^2(Y_qg)\cosh^m(Y_qg)]} \geq
 \frac{e^{-\xi'(q)/2}\E\cosh^{-3}(Y_qg)}
 {1-e^{-\xi'(q)/2}\E\cosh^{-1}(Y_qg)}.
 \label{eq:m-one}
\end{equation}
Gaussian differentiation gives
\begin{align}
 \frac d{d\xi'(q)}
 \left\{1-e^{-\xi'(q)/2}\E\cosh^{-1}(Y_qg)\right\}
 &=e^{-\xi'(q)/2}\E\cosh^{-3}(Y_qg),
 \label{eq:heat-one}\\
 \frac d{d\xi'(q)}
 \left\{e^{-\xi'(q)/2}\E\cosh^{-3}(Y_qg)\right\}
 &=4e^{-\xi'(q)/2}\E\cosh^{-3}(Y_qg)
   -6e^{-\xi'(q)/2}\E\cosh^{-5}(Y_qg)\notag\\
 &>-2e^{-\xi'(q)/2}\E\cosh^{-3}(Y_qg).
 \label{eq:heat-two}
\end{align}
Integrating \eqref{eq:heat-one} and using \eqref{eq:heat-two} yields
\begin{equation}
 1-e^{-\xi'(q)/2}\E\cosh^{-1}(Y_qg)
 <\frac{e^{2\xi'(q)}-1}{2}
 e^{-\xi'(q)/2}\E\cosh^{-3}(Y_qg).
 \label{eq:heat-integrated}
\end{equation}
Using $D_\beta^1(m,q)=0$, \eqref{eq:endpoint-replicon},
\eqref{eq:m-one}, and \eqref{eq:heat-integrated}, we obtain
\begin{equation}
 \Gamma_\mu'(q)>
 \frac{2(p-1)\xi'(q)}{e^{2\xi'(q)}-1}.
 \label{eq:endpoint-lower}
\end{equation}
Therefore $\Gamma_\mu'(q)\leq1$ implies
\begin{equation}
 e^{2\xi'(q)}-1>2(p-1)\xi'(q).
 \label{eq:endpoint-stability}
\end{equation}

For $q\leq u\leq1$, differentiation of \eqref{eq:Gamma-prime} gives
\begin{align}
 \Gamma_\mu''(u)
 ={}&\xi'''(u)
 \E[\cosh^{-4}(M(u))e^{W_\mu(u)}]\notag\\
 &+\xi''(u)^2\left\{
 4\E[\cosh^{-4}(M(u))e^{W_\mu(u)}]
 -6\E[\cosh^{-6}(M(u))e^{W_\mu(u)}]\right\}.
 \label{eq:right-second}
\end{align}
The Gaussian integration-by-parts calculation is given at the beginning
of Appendix~\ref{app:right-ratio}.
We claim that
\begin{equation}
 \frac{\E[\cosh^{-6}(M(u))e^{W_\mu(u)}]}
 {\E[\cosh^{-4}(M(u))e^{W_\mu(u)}]}
 >\frac{3\xi'(u)+1}{4\xi'(u)+1}.
 \label{eq:right-ratio}
\end{equation}
The complete one-dimensional calculation is given in
Appendix~\ref{app:right-ratio}.

Substitution of \eqref{eq:right-ratio} into
\eqref{eq:right-second}, together with
\[
 \xi''(u)=\frac{(p-1)\xi'(u)}u,\qquad
 \xi'''(u)=\frac{(p-1)(p-2)\xi'(u)}{u^2},
\]
gives
\begin{eqnarray}
 \Gamma_\mu''(u)
 <\frac{(p-1)\xi'(u)}{u^2}
 \E[\cosh^{-4}(M(u))e^{W_\mu(u)}] \cdot\left\{p-2-
 \frac{2(p-1)\xi'(u)(\xi'(u)+1)}{4\xi'(u)+1}\right\}.
 \label{eq:right-upper}
\end{eqnarray}
It remains to check the sign of the braces.  First,
\eqref{eq:endpoint-stability} implies
\begin{equation}
 2(p-1)\xi'(q)(\xi'(q)+1)
 >(p-2)(4\xi'(q)+1).
 \label{eq:quadratic-at-q}
\end{equation}
Indeed, if \eqref{eq:quadratic-at-q} failed, then
\[
 (p-1)\{1+2\xi'(q)-2\xi'(q)^2\}\geq4\xi'(q)+1.
\]
The expression in braces would be positive, and therefore
\[
 2(p-1)\xi'(q)
 \geq\frac{2\xi'(q)(4\xi'(q)+1)}
 {1+2\xi'(q)-2\xi'(q)^2}
 >e^{2\xi'(q)}-1,
\]
contrary to \eqref{eq:endpoint-stability}.  
 After multiplication by the positive denominator, the last strict inequality
is
\begin{equation}
 1+4\xi'(q)+6\xi'(q)^2
 -\{1+2\xi'(q)-2\xi'(q)^2\}e^{2\xi'(q)}>0.
 \label{eq:elementary}
\end{equation}
The left side vanishes at $0$, and its derivative is
\[
 4\{1+3\xi'(q)-(1-\xi'(q)^2)e^{2\xi'(q)}\}>0.
\]
For $\xi'(q)\geq1$ this is immediate.  For
$0<\xi'(q)<1$, the derivative of
$ \log\frac{1+3\xi'(q)}{1-\xi'(q)^2}-2\xi'(q)$ is 
\[
 \frac{(1-2\xi'(q))^2+\xi'(q)^2+6\xi'(q)^3}
 {(1+3\xi'(q))(1-\xi'(q)^2)}>0.
\]
This proves \eqref{eq:elementary}.

Finally,
\[
 2(p-1)\xi'(u)(\xi'(u)+1)-(p-2)(4\xi'(u)+1)
\]
is a quadratic with positive leading coefficient, negative value at
$\xi'(u)=0$, and exactly one positive root.  It is positive at
$\xi'(q)$ by \eqref{eq:quadratic-at-q}, and hence remains positive for
$u\geq q$.  Equation \eqref{eq:right-upper} now yields
$\Gamma_\mu''(u)<0$.  Since $\Gamma_\mu(q)=q$ and
$\Gamma_\mu'(q)\leq1$, it follows that $\Gamma_\mu(u)<u$ for $u>q$.
Thus $f_\mu'(u)<0$ there.  The equality $f_\mu(q)=0$ completes the
proof.
\end{proof}

\subsection{The marginal crossing}

In this section we consider the solutions of \eqref{eq:stationary} at which
$\Gamma_\mu'(q)=1$.  Its complete proof is included in
Appendix~\ref{app:marginal}.

\begin{proposition}
\label{prop:transversality}
Suppose that \eqref{eq:stationary} holds at
$(\beta,m,q)\in(0,\infty)\times(0,1)^2$ and
$\Gamma_\mu'(q)=1$.  The nearby solutions of
\eqref{eq:stationary} form a locally unique $C^1$ curve, and $\beta^2$
is a valid coordinate on this curve.  Along the curve,
\begin{equation}
 \frac d{d\beta^2}\Gamma_\mu'(q)>0
 \label{eq:transversality}
\end{equation}
at the marginal point.
\end{proposition}

\subsection{Completion of the 1-RSB phase}

We now pass from the local results to the complete 1-RSB interval.  The
Parisi functional is jointly continuous in $\beta$ and $\mu$, and its
minimizer is unique.  Consequently, if $\beta$ converges, then the
corresponding Parisi measures converge weakly: every subsequential limit
minimizes the limiting functional and therefore equals its unique
minimizer.

We first show that the strict 1-RSB phase is open.  Suppose that at
$\beta$ the Parisi measure is
$m\delta_0+(1-m)\delta_q$, where $0<m,q<1$ and
$\Gamma_\mu'(q)<1$.  For nearby inverse temperatures, minimize the
Parisi functional over the compact family
\begin{equation}
 \{m\delta_0+(1-m)\delta_q:0\leq m,q\leq1\}.
 \label{eq:two-atom-family}
\end{equation}
Every limit of these restricted minimizers is a restricted minimizer at
$\beta$.  At $\beta$ the original Parisi measure belongs to
\eqref{eq:two-atom-family}, so the restricted minimum equals the global
minimum.  Hence every restricted minimizer at $\beta$ is the Parisi
measure and, by uniqueness, equals the original measure.  Since its
two-atom representation is interior and unique, weak convergence also
gives convergence of $m$ and $q$.  The nearby restricted minimizers are
therefore in the interior of \eqref{eq:two-atom-family}, satisfy
\eqref{eq:stationary}, and still have $\Gamma_\mu'(q)<1$.
Propositions~\ref{prop:left-sign} and \ref{prop:right-sign}, together
with \eqref{eq:criterion}, show that they are the Parisi measures.  This
proves the claimed openness.

By \cite[Theorem~1]{zhou}, the Parisi measure is 1-RSB on an interval
immediately to the right of $\beta_1^p$.  At every point of this
interval, $q$ is a local maximum of $f_\mu$, so
\eqref{eq:second-curvature} gives
$\Gamma_\mu'(q)\leq1$.  Equality at an interior point would, by the
local argument below, make that point a right endpoint of the 1-RSB
phase, contradicting the same interval immediately to its right.
Therefore $\Gamma_\mu'(q)<1$ after possibly shortening the interval.
Let
$\beta_2^p$ be the right endpoint of
the largest such interval issuing from $\beta_1^p$.  This endpoint is
finite: \cite[Theorem~3 and Example~1]{auffingerchen} show that for all
sufficiently large $\beta$ the pure $p$-spin Parisi measure is neither
RS nor 1-RSB.

We next examine any finite boundary point of the strict 1-RSB phase
above $\beta_1^p$.  Let strict 1-RSB Parisi measures converge to that
point.  After passing to a subsequence, their two parameters converge
to $m,q\in[0,1]$, and the preceding continuity gives the limiting
Parisi measure $m\delta_0+(1-m)\delta_q$.  Neither $m=1$ nor $q=0$ is
possible, since either case gives $\delta_0$ above the exact RS
boundary.  If $m=0$, the limit is $\delta_q$; since $0$ belongs to the
support of every Parisi measure in zero external field
\cite[Theorem~1(i)]{auffingerchen}, this forces $q=0$, already
excluded.  Finally, $q=1$ is impossible.  Indeed, passing to the limit
in $D_\beta^1(m,q)=0$ gives
\[
 1=\frac{\E[\tanh^2(Y_1g)\cosh^m(Y_1g)]}
          {\E\cosh^m(Y_1g)}<1,
\]
a contradiction.
Thus $0<m,q<1$, and \eqref{eq:stationary} passes to the limit.  Since
the limiting measure is Parisi, $f_\mu\leq0$ and $q$ is a local maximum
of $f_\mu$.  Equation~\eqref{eq:second-curvature} therefore gives
$\Gamma_\mu'(q)\leq1$.  Strict inequality would put the boundary point
in the open strict 1-RSB phase.  Hence every such boundary point
satisfies
\begin{equation}
 \Gamma_\mu'(q)=1.
 \label{eq:boundary-marginal}
\end{equation}

Conversely, take any solution of \eqref{eq:second-boundary} in
$(0,\infty)\times(0,1)^2$.  Proposition~\ref{prop:transversality}
shows that the nearby stationary solutions form a unique curve
parametrized by $\beta^2$, with $\Gamma_\mu'(q)<1$ immediately below
the given point and $\Gamma_\mu'(q)>1$ immediately above it.  The
strict-sign result makes all the solutions immediately below it Parisi
measures.  Their limit is therefore the Parisi measure at the marginal
point.

There can be no 1-RSB Parisi measures immediately above this point.
Otherwise, choose such inverse temperatures decreasing to the marginal
one.  Continuity of the unique Parisi minimizer gives weak convergence
to the marginal measure.  Since its two-atom representation is interior
and unique, the corresponding parameters converge to $m$ and $q$.
Local uniqueness then puts them on the same stationary curve, where
$\Gamma_\mu'(q)>1$, contrary to \eqref{eq:second-curvature} and
$f_\mu\leq0$.  Thus every solution of
\eqref{eq:second-boundary} is locally a right endpoint of the strict
1-RSB phase.  It lies above $\beta_1^p$, because its local stationary
curve would otherwise give a 1-RSB Parisi measure below the exact RS
boundary.

There cannot be two components of the strict 1-RSB phase.  Indeed, the
phase contains an interval immediately to the right of $\beta_1^p$.
If a later component existed, its finite left endpoint would, by the
compactness argument above, satisfy \eqref{eq:second-boundary}.  But
every such point has strict 1-RSB measures immediately below it and no
1-RSB measures immediately above it, contradicting that it is the left
endpoint of a component.  Therefore the strict phase is exactly
\begin{equation}
 \beta_1^p<\beta<\beta_2^p.
 \label{eq:strict-phase}
\end{equation}

At $\beta_2^p$, compactness and \eqref{eq:boundary-marginal} give an
interior 1-RSB Parisi measure satisfying
\eqref{eq:second-boundary}.  Every other solution of that system is the
same right endpoint.  At this inverse temperature every such solution
is the Parisi measure, so uniqueness of the Parisi minimizer and of the
representation $m\delta_0+(1-m)\delta_q$, for $0<m,q<1$, give the same
$m$ and $q$.  This proves Proposition~\ref{prop:second-boundary}.

Finally, for $\beta>\beta_2^p$ the Parisi measure is not RS by
\cite[Theorem~1]{zhou}.  If it were 1-RSB, then
\eqref{eq:second-curvature} would give $\Gamma_\mu'(q)\leq1$.  Strict
inequality would put $\beta$ in \eqref{eq:strict-phase}, while equality
would give another solution of \eqref{eq:second-boundary}.  Both are
impossible.  Hence the Parisi measure is 1-RSB precisely for
$\beta_1^p<\beta\leq\beta_2^p$ and is neither RS nor 1-RSB above
$\beta_2^p$.

\section{Proof of FRSB Phase: $\beta > \beta^p_2$}\label{sec:positive-gaps}

We now establish the FRSB phase above $\beta_2^p$ and show that $\operatorname{supp}\mu_\beta=\{0\} \cup [q_\beta,q'_\beta]$, for some $0<q_\beta<q'_\beta<1$.  We prove by contradiction that if a measure leaves a positive gap from its support, then it cannot satisfy all the boundary conditions for the Parisi measure.

\begin{theorem}\label{thm:no-positive-gap}
Let $\mu$ be a probability measure such that
\[
 0<q<q'<1,\qquad q,q'\in\operatorname{supp}\mu,
 \qquad \operatorname{supp}\mu\cap(q,q')=\varnothing,
\]
then the following four conditions
\begin{equation}
 \Gamma_\mu(q)=q,\qquad \Gamma_\mu(q')=q',\qquad
 \Gamma_\mu'(q)\leq1,\qquad \Gamma_\mu'(q')\leq1.
 \label{eq:endpoint-conditions}
\end{equation}
cannot hold simultaneously.  More precisely, the endpoint
conditions would imply
\begin{equation}
 \Gamma_\mu(s)-s\Gamma_\mu'(s)>0,\qquad q<s<q'.
 \label{eq:strict-N-physical}
\end{equation}
Consequently,
\[
 \Gamma_\mu(s)<s,\qquad
 \Gamma_\mu(s)=s\Gamma_\mu'(s)
\]
cannot hold simultaneously for $q<s<q'$.  
\end{theorem}

\begin{remark} The assumption $q>0$ is essential.  The proof does not exclude a gap whose
left endpoint is zero, and is therefore compatible with an isolated support
point at the origin.
\end{remark}

We need the following proposition in the proof of
Theorem~\ref{thm:no-positive-gap}.

\begin{proposition}
\label{prop:contact-fold}
Suppose that
\[
    0<q<q'<1,\qquad q,q'\in\operatorname{supp}\mu,
    \qquad \operatorname{supp}\mu\cap(q,q')=\varnothing.
\]
If $s\in(q,q')$ satisfies
\[
    \Gamma_\mu(s)=s\Gamma_\mu'(s),
    \qquad
    \Gamma_\mu''(s)=0,
\]
then
\[
    \Gamma_\mu'''(s)>0.
\]
\end{proposition}

\begin{proof}
Appendix~\ref{app:analytic-details} derives the required differentiated
identity and reduces the conclusion to
\eqref{eq:app-final-safe-bound} on the region
\eqref{eq:app-complete-feasible-region}. Appendix
\ref{app:exact-certification} proves this inequality.  In the notation
of Appendix~\ref{app:analytic-details}, direct differentiation gives
\[
 \Gamma_\mu''(s)=\{\xi''(s)\}^2I(\xi'(s)).
\]
At the fold $I(\xi'(s))=0$, and therefore
\[
 \Gamma_\mu'''(s)=\{\xi''(s)\}^3I'(\xi'(s))>0.
\]
The approximation
argument at the end of Appendix~\ref{app:analytic-details} extends the
result from finite-step measures to arbitrary order parameters.
\end{proof}

\begin{lemma}
\label{lem:one-crossing}
Under the assumptions of Theorem~\ref{thm:no-positive-gap},
\[
    \Gamma_\mu(s)-s\Gamma_\mu'(s)>0,
    \qquad q<s<q'.
\]
\end{lemma}

\begin{proof}
Appendices~\ref{app:ig:one-crossing}--\ref{app:ig:homotopy} prove
\eqref{app:ig:N1-positive}.  For $t=\xi'(s)$, the definitions
\eqref{app:ig:NI} and the identity
$s\xi''(s)=(p-1)\xi'(s)$ give
\[
 N_1(\xi'(s))
 =\Gamma_\mu(s)
  -(p-1)\xi'(s)\frac{\Gamma_\mu'(s)}{\xi''(s)}
 =\Gamma_\mu(s)-s\Gamma_\mu'(s).
\]
The strict positivity in \eqref{app:ig:N1-positive} proves the lemma.
\end{proof}

\begin{proof}[Proof of Theorem~\ref{thm:no-positive-gap}]
Lemma~\ref{lem:one-crossing} gives
\[
    \Gamma_\mu(s)-s\Gamma_\mu'(s)>0,
    \qquad q<s<q'.
\]
Consequently,
\[
    \frac{d}{ds}\left(\frac{\Gamma_\mu(s)}{s}\right)
    =
    \frac{s\Gamma_\mu'(s)-\Gamma_\mu(s)}{s^2}
    <0.
\]
Thus $\Gamma_\mu(s)/s$ is strictly decreasing on $(q,q')$. However,
the endpoint conditions give
\[
    \frac{\Gamma_\mu(q)}q
    =
    \frac{\Gamma_\mu(q')}{q'}
    =1,
\]
which is impossible. This proves
Theorem~\ref{thm:no-positive-gap}.
\end{proof}

\begin{corollary}\label{cor:positive-support-convex}
If $ \mu_\beta$ is the Parisi measure and
$0<q<q'<1$ belong to $\operatorname{supp} \mu_\beta$, then
\[
 [q,q']\subseteq\operatorname{supp} \mu_\beta.
\]
\end{corollary}

\begin{proof}
Suppose that some point of $(q,q')$ is outside the support.  Since the
support is closed, the connected component of its complement containing
that point is an open interval $(a,b)\subset(q,q')$ whose endpoints
belong to $\operatorname{supp} \mu_\beta$.  The optimality criterion
\eqref{eq:criterion} gives $f_{ \mu_\beta}(a)=f_{ \mu_\beta}(b)=0$ and
$f_{ \mu_\beta}\leq0$.  Since both endpoints lie in $(0,1)$,
$f_{ \mu_\beta}'(a)=f_{ \mu_\beta}'(b)=0$, and hence
\[
 \Gamma_{ \mu_\beta}(a)=a,\qquad \Gamma_{ \mu_\beta}(b)=b.
\]
Since these are one-sided local maxima and
$\Gamma_{\mu_\beta}(a)=a$, $\Gamma_{\mu_\beta}(b)=b$,
\[
 f_{\mu_\beta}''(a+)
 =\frac{\xi''(a)}2\{\Gamma_{\mu_\beta}'(a+)-1\}\leq0,
 \qquad
 f_{\mu_\beta}''(b-)
 =\frac{\xi''(b)}2\{\Gamma_{\mu_\beta}'(b-)-1\}\leq0.
\]
Therefore
\[
 \Gamma_{\mu_\beta}'(a)\leq1,
 \qquad
 \Gamma_{\mu_\beta}'(b)\leq1.
\]
Theorem~\ref{thm:no-positive-gap} excludes $(a,b)$.
\end{proof}

\begin{proof}[Completion of the proof of Theorem~\ref{thm:main}]
The first transition is \cite[Theorem~1]{zhou}.  The argument of
Section~\ref{sec:sign} proves that the Parisi measure is 1-RSB
precisely for $\beta_1^p<\beta\leq\beta_2^p$, and that for
$\beta>\beta_2^p$ it is neither RS nor 1-RSB.

Fix $\beta>\beta_2^p$.  In zero external field, zero belongs to the
support of every Parisi measure \cite[Theorem~1(i)]{auffingerchen}.  Also,
one cannot have $1\in\operatorname{supp} \mu_\beta$.  Indeed,
the definition of $\Gamma_{ \mu_\beta}$ gives
\[
 \Gamma_{ \mu_\beta}(1)
 =\E[\tanh^2(M(1))e^{W_{ \mu_\beta}(1)}]
 <\E e^{W_{ \mu_\beta}(1)}=1.
\]
By continuity, $f_{ \mu_\beta}'<0$ immediately to the left of $1$.  If
$1$ belonged to the support, then $f_{ \mu_\beta}(1)=0$, so
$f_{ \mu_\beta}>0$ immediately to its left, contrary to
\eqref{eq:criterion}.

We next show that the origin is isolated.  At $\beta=\beta_2^p$,
equation \eqref{eq:second-marginal} and the marginal equality give
\[
 1=\xi''(q_2^p)
 \frac{\E\cosh^{m_2^p-4}(Y_{q_2^p}g)}
      {\E\cosh^{m_2^p}(Y_{q_2^p}g)}.
\]
The ratio is strictly smaller than $1$, since it is the expectation of
$\operatorname{sech}^4(Y_{q_2^p}g)$ under the probability measure
with density proportional to $\cosh^{m_2^p}(Y_{q_2^p}g)$ and
$Y_{q_2^p}g$ is nondegenerate.  Hence $\xi''(q_2^p)>1$, and so
$\xi''(1)>1$ at $\beta_2^p$.  The same inequality holds for every
$\beta>\beta_2^p$.  Since $\xi''(0)=0$, there is a point in $(0,1)$ at
which $\xi''=1$.  Theorem~1(ii) of \cite{auffingerchen} therefore gives
a neighborhood of the origin containing no positive support point.
Thus $0$ is isolated in $\operatorname{supp}\mu_\beta$.

If the support had at most one positive point, it would therefore be
either $\{0\}$ or $\{0,q\}$ with $0<q<1$.  The first case is RS and the
second is 1-RSB, both impossible above $\beta_2^p$.  Hence there are two
distinct positive support points.  Because the origin is isolated and
the support is compact, its positive part has a smallest and a largest
point.  Denote them by $q_\beta$ and $q'_\beta$.  The preceding paragraph
and the exclusion of $1$ give
\[
 0<q_\beta<q'_\beta<1.
\]
Corollary~\ref{cor:positive-support-convex} yields
\[
 [q_\beta,q'_\beta]\subseteq\operatorname{supp}\mu_\beta.
\]
The reverse inclusion away from zero follows from the definitions of
$q_\beta$ and $q'_\beta$.  Consequently,
\[
 \operatorname{supp}\mu_\beta=\{0\}\cup[q_\beta,q'_\beta].
\]
Finally, \cite[Theorem~2(ii)]{auffingerchen} shows that the distribution
function is infinitely differentiable on $(q_\beta,q'_\beta)$.
Equivalently, $\mu_\beta$ has a smooth density there.  This does not
rule out atoms at either endpoint of the positive interval.

For $\beta\leq\beta_2^p$, the Parisi measure is RS or 1-RSB and its
support is finite, so it is not FRSB.  This completes all three parts of
the theorem.
\end{proof}

\appendix

\section{The estimate to the right of the positive atom}
\label{app:right-ratio}

For $q<u\leq1$, the explicit form of $\Phi_\mu$ on $[q,1]$ and the
definition of $W_\mu$ give, conditionally on $M(q)=z$ and $M(u)=x$,
\[
 e^{W_\mu(u)}
 =\frac{e^{-\{\xi'(u)-\xi'(q)\}/2}}
 {\E\cosh^m(Y_qg)}\cosh x\cosh^{m-1}z.
\]
Consequently,
\[
 \E[\cosh^{-4}(M(u))e^{W_\mu(u)}]
 =\frac{e^{-\{\xi'(u)-\xi'(q)\}/2}}
        {\E\cosh^m(Y_qg)}
   \E[\cosh^{-3}(M(u))\cosh^{m-1}(M(q))].
\]
Differentiating the Gaussian variance $\xi'(u)-\xi'(q)$ and integrating
by parts gives
\[
\begin{aligned}
 \frac d{du}\E[\cosh^{-4}(M(u))e^{W_\mu(u)}]
 =\xi''(u)\bigl\{
 &4\E[\cosh^{-4}(M(u))e^{W_\mu(u)}]\\
 &-6\E[\cosh^{-6}(M(u))e^{W_\mu(u)}]\bigr\}.
\end{aligned}
\]
Indeed,
\[
 \frac12\left\{\frac{d^2}{dx^2}\cosh^{-3}x-\cosh^{-3}x\right\}
 =4\cosh^{-3}x-6\cosh^{-5}x.
\]
Differentiating \eqref{eq:Gamma-prime} proves \eqref{eq:right-second}
for $u>q$, and the identity at $u=q$ follows by continuity.

We now prove \eqref{eq:right-ratio}.  The quotient on its left-hand side
equals
\begin{equation}
 \frac{\displaystyle
 \int_{\mathbb R}\cosh^{-5}x
 \int_{\mathbb R}
 \exp\left\{-\frac{z^2}{2\xi'(q)}
 -\frac{(x-z)^2}{2\{\xi'(u)-\xi'(q)\}}\right\}
 \cosh^{m-1}z\,dz\,dx}
 {\displaystyle
 \int_{\mathbb R}\cosh^{-3}x
 \int_{\mathbb R}
 \exp\left\{-\frac{z^2}{2\xi'(q)}
 -\frac{(x-z)^2}{2\{\xi'(u)-\xi'(q)\}}\right\}
 \cosh^{m-1}z\,dz\,dx}.
 \label{eq:right-integrals}
\end{equation}
Normalize the even function
\begin{equation}
 x\longmapsto
 \int_{\mathbb R}
 \exp\left\{-\frac{z^2}{2\xi'(q)}
 -\frac{(x-z)^2}{2\{\xi'(u)-\xi'(q)\}}\right\}
 \cosh^{m-1}z\,dz
 \label{eq:right-density}
\end{equation}
to have integral one.  We first prove that the second derivative of its
negative logarithm is at least $1/\xi'(u)$.

Differentiating twice in $x$ gives
\begin{align}
 &-\frac{d^2}{dx^2}\log
 \int_{\mathbb R}
 \exp\left\{-\frac{z^2}{2\xi'(q)}
 -\frac{(x-z)^2}{2\{\xi'(u)-\xi'(q)\}}\right\}
 \cosh^{m-1}z\,dz \notag\\
 &\qquad=\frac1{\xi'(u)-\xi'(q)}
 -\frac{\operatorname{Var}(z\mid x)}
 {\{\xi'(u)-\xi'(q)\}^2}.
 \label{eq:right-density-second}
\end{align}
The negative logarithm of the conditional density of $z$ has second
derivative
\[
 \frac1{\xi'(q)}+\frac1{\xi'(u)-\xi'(q)}
 +(1-m)\cosh^{-2}z
 \geq
 \frac1{\xi'(q)}+\frac1{\xi'(u)-\xi'(q)}.
\]
The following direct integration by parts gives the required variance
bound.  The preceding lower bound on the second derivative implies
\begin{align*}
 &\left\{\frac1{\xi'(q)}
 +\frac1{\xi'(u)-\xi'(q)}\right\}
 \operatorname{Var}(z\mid x)\\
 &\quad\leq
 \E\left[(z-\E[z\mid x])
 \frac d{dz}\bigl(-\log \text{conditional density}\bigr)(z)
 \,\middle|\,x\right]=1.
\end{align*}
Here the term obtained by evaluating the derivative at
$\E[z\mid x]$ vanishes, and the last equality is integration by parts.
Therefore
\[
 \operatorname{Var}(z\mid x)
 \leq\frac{\xi'(q)\{\xi'(u)-\xi'(q)\}}{\xi'(u)},
\]
and \eqref{eq:right-density-second} is at least $1/\xi'(u)$.
When $u=q$, the limiting density is proportional to
\[
 \exp\left\{-\frac{x^2}{2\xi'(q)}\right\}\cosh^{m-1}x,
\]
whose negative logarithm has second derivative
$1/\xi'(q)+(1-m)\cosh^{-2}x\geq1/\xi'(q)$.  Thus the same conclusion
holds at $u=q$.

Since \eqref{eq:right-density} is even, the derivative of its negative
logarithm is at least $x/\xi'(u)$ for $x>0$ and at most
$x/\xi'(u)$ for $x<0$.  In the following display, expectation is taken
with respect to the normalized form of \eqref{eq:right-density}.
Using
\[
 \frac d{dx}\{\tanh x\cosh^{-3}x\}
 =4\cosh^{-5}x-3\cosh^{-3}x,
\]
integration by parts gives
\begin{align*}
 4\E\cosh^{-5}x-3\E\cosh^{-3}x
 &\geq\frac1{\xi'(u)}\E[x\tanh x\cosh^{-3}x]\\
 &>\frac1{\xi'(u)}\E[\tanh^2x\cosh^{-3}x]\\
 &=\frac1{\xi'(u)}
   \{\E\cosh^{-3}x-\E\cosh^{-5}x\}.
\end{align*}
The strict inequality follows from
$x\tanh x>\tanh^2x$ for $x\ne0$.  Rearranging and using
\eqref{eq:right-integrals}, together with its limiting form when $u=q$,
proves \eqref{eq:right-ratio}.

\section{The contact-curvature computation}
\label{app:contact}

This appendix proves Proposition~\ref{prop:left-computation}.  The
notation of the main text remains in force.  Temporary quantities are
introduced only inside the calculation.  Displayed equation numbers in
this appendix are local to the calculation.

\subsection{Clock change and the global contact branch}\label{ft:sec:clock}

From the definition above,

\[
 \Gamma_\mu(t)=\mathbb E\left[
   \Phi_{\mu,x}(M(t),t)^2e^{W_\mu(t)}\right],
 \qquad M(t)=\mathcal B_{\xi'(t)}.
\]

Set

\[
 s=\xi'(u),\qquad v=\xi'(q),\qquad
 Q_v(s)=\Gamma_\mu((\xi')^{-1}(s)).
\]

For \(0<u<q\), the cumulative mass in the Parisi PDE is the constant \(m\).
Writing \(r=v-s\), the logarithmic heat transform gives, up to an additive constant independent
of \(x\),

\[
 \Phi_\mu(x,u)=\frac{1}{m}\log F_r(x),\qquad
 F_r(x)=\mathbb E\cosh^m(x+\sqrt r Z).
\]

Here \(Z\sim N(0,1)\), and

\[
 \phi_s(x)=\frac{1}{\sqrt{2\pi s}}e^{-x^2/(2s)}
\]

is the centered Gaussian density of variance \(s\).

Only the atom at zero enters \(W_\mu(u)\), and hence

\[
 e^{W_\mu(u)}=\frac{F_r(M(u))}{F_v(0)},\qquad
 F_v(0)=\int_{\mathbb R}\phi_s(x)F_r(x)\,dx.          \tag{1.0}
\]

If

\[
 B={1\over m}(\log F_r)',\qquad C=B',
\]

then the definition gives

\[
 Q_v(s)={\int\phi_sF_rB^2\over F_v(0)}.               \tag{1.0a}
\]

We record the derivative because it is easy to lose track of which
quantity is held fixed.  Here $v$ is fixed, so $r=v-s$.  Put

\[
 \rho_s(x)={\phi_s(x)F_{v-s}(x)\over F_v(0)}.
\]

The two heat equations for $\phi_s$ and $F_{v-s}$, followed by two
integrations by parts, give, for every smooth $f=f(s,x)$,

\[
 {d\over ds}\int\rho_s f
 =\int\rho_s\left(f_s+{1\over2}f_{xx}+mBf_x\right).   \tag{1.0b}
\]

Along the same fixed-$v$ path,

\[
 B_s=-{1\over2}B_{xx}-mBB_x.
\]

Taking $f=B^2$ in (1.0b), and writing $C=B_x$, makes every
term cancel except $C^2$.  Therefore

\[
 Q_v'(s)={\int\phi_sF_rC^2\over F_v(0)}.              \tag{1.0c}
\]

This proves the derivative identity directly in the two-atom setting.

Set \(P=p-1\). Since \(\xi'(u)\) is proportional to \(u^P\),

\[
 \Gamma_\mu(u)=u\Gamma_\mu'(u)
 \iff Q_v(s)=PsQ_v'(s),                                      \tag{1.1}
\]

and, at a contact,

\[
 \Gamma_\mu''(u)=\left({Ps\over u}\right)^2 I_\alpha(s),
 \qquad
 I_\alpha=Q_v''+{\alpha\over s}Q_v',
 \qquad \alpha={P-1\over P}.                              \tag{1.2}
\]

By (1.0a), up to the positive partition-function factor \(F_v(0)\), the
fixed-\(r\) contact
numerator is

\[
 \mathcal N(s,r)=\int_{\mathbb R}\phi_s(x)F_r(x)
       \{B(x)^2-PsC(x)^2\}\,dx.                         \tag{1.3}
\]

For every fixed \(r\), \(\mathcal N(\,\cdot\,,r)\) has exactly one zero
and crosses it upward.
Here is the short proof.  Put

\[
 I(x)=\int_0^xFC^2,\qquad \mathcal M(x)={F(x)B(x)^2\over xI(x)},
 \qquad t=mxB,\qquad a={xC\over B},\qquad
 z=-{C_x\over2C}.
\]

In the \(B\) coordinate, \(h(B(x))=F(x)C(x)\) is strictly decreasing, since
the slope cone \(z\ge mB\) gives \((FC)'=FC(mB-2z)<0\).  Hence

\[
 e={Bh(B)\over\int_0^Bh}\in(0,1),
 \qquad x(\log \mathcal M)'=t-1+2a-ae.                \tag{1.4}
\]

Thus \(\mathcal M'>0\) when \(t\ge 1\). When \(t\le 1\),
the elementary product inequality and monotonicity of
\(B\) give

\[
 \mathcal M(x)\le {e^t-1\over t}\le e-1<2\le P.      \tag{1.5}
\]

Since \(\mathcal M(0+)=1\) and \(\mathcal M(\infty)=\infty\),
\(\mathcal M-P\) changes sign exactly once,
from minus to plus.  After \(X=x^2\), two integrations by parts turn (1.3)
into a positive multiple of

\[
 \zeta\int_0^\infty e^{-\zeta X}K_P(X)\,dX,
 \qquad \zeta={1\over2s},
 \qquad K_P(x^2)=2I(x)\{\mathcal M(x)-P\}.           \tag{1.6}
\]

The elementary one-change lemma for the Laplace kernel proves the claimed
uniqueness and upward crossing.

Let that root be \(s=\sigma(r)\) and put

\[
 V(r)=r+\sigma(r).                                      \tag{1.7}
\]

Implicit differentiation, using \(\mathcal N_s>0\), gives

\[
 V'={\mathcal N_s-\mathcal N_r\over\mathcal N_s}.       \tag{1.8}
\]

For clarity, let \(Z_v>0\) denote the omitted partition-function factor and
let \(\mathscr L=\partial_s-\partial_r\), the derivative along a line on
which the total horizon \(v=s+r\) is fixed.  Then

\[
 \mathcal N=Z_v\{Q_v-PsQ_v'\},\qquad
 \mathscr L\mathcal N=-Z_vPsI_\alpha.                 \tag{1.9}
\]

Consequently

\[
 V'=-{Z_vPs\over \mathcal N_s}I_\alpha.               \tag{1.10}
\]

At a critical point of \(V\), one has both (1.1) and \(I_\alpha=0\).  Moreover,
there \(\sigma'=-1\), and a second implicit differentiation gives

\[
 V''=-{\mathscr L^2\mathcal N\over\mathcal N_s},\qquad
 \mathscr L^2\mathcal N=-Z_vP\{I_\alpha+sI_\alpha'\}.
\]

It follows that

\[
 \boxed{V''>0\iff I_\alpha'(s)>0.}                     \tag{1.11}
\]

It remains to prove the local implication on the right of (1.11).

\subsection{Exact standardized local identity}

At a simultaneous contact and shifted fold put \(y=x/\sqrt{s}\) and let

\[
 d\pi(y)={e^{-y^2/2}F_r(\sqrt s\,y)\,dy
 \over\int_{\mathbb R}e^{-z^2/2}F_r(\sqrt s\,z)\,dz}.
\]

Set

\[
 z=-{C_x\over2C},\qquad
 J=-m-{1\over2}C_{BB},\qquad
 N={x\over s}+z-mB,
\]

where subscripts \(B\) mean differentiation in the inverse-slope coordinate.
The preserved slope invariant gives the substantive cone

\[
 J,J_B\ge0,\qquad 0\le CJ_B\le3zJ.                    \tag{2.0}
\]

The slope-coordinate maximum-principle argument establishing this cone is
given in Subsection~\ref{ft:sec:slope-coordinate}.  It is included because the last inequality in (2.0)
is both important and not visually obvious.

Use the standardized fields

\[
 a=msC,\quad b=m\sqrt s B,\quad u=\sqrt s\,z,
 \quad j=sCJ,\quad n=\sqrt s\,N,
 \quad \gamma=s^{3/2}C(3zJ-CJ_B).
\]

and define the probability law

\[
 d\nu={a^2\,d\pi\over\int a^2\,d\pi}.
\]

Every unadorned expectation \(E\) below is with respect to \(\nu\).  Primes
below mean differentiation in \(y\).  Direct substitution, together with
(2.0), gives

\[
 a'=-2au,\quad b'=a,\quad u'=a+j,\quad
 j'=uj-\gamma,\quad n'=1+j,\quad 0\le\gamma\le3uj.     \tag{2.1}
\]

Put

\[
 c={P-1\over2P},\qquad
 \mathfrak p=a-2u^2-c,\qquad
 \mathfrak q=j-u(n+u)+c.                                \tag{2.2}
\]

The contact--fold identities are

\[
 E_\nu\mathfrak p=E_\nu\mathfrak q=0,\qquad
 E_\nu(b/a)^2=P.                                       \tag{2.3}
\]

Indeed, contact gives the last equality.  At the shifted fold the local
shifted kernel equals \(-\mathfrak p/s\), giving
\(E_\nu\mathfrak p=0\); integration by parts against the score
\(-(n+3u)\) of \(\nu\) then gives \(E_\nu\mathfrak q=0\).

Let \(\omega\) be the positive-half-line density of \(\nu\).  The centered Stein
tails are boundary-normalized by

\[
 \omega(y)v(y)=\int_0^y\omega(t)\mathfrak p(t)\,dt
              =-\int_y^\infty\omega(t)\mathfrak p(t)\,dt,
\]

\[
 \omega(y)\delta(y)=\int_0^y\omega(t)\mathfrak q(t)\,dt
              =-\int_y^\infty\omega(t)\mathfrak q(t)\,dt. \tag{2.3a}
\]

Thus \((\omega v)'=\omega\mathfrak p\) and
\((\omega\delta)'=\omega\mathfrak q\); all boundary terms below vanish.
Since both \(\mathfrak p\) and \(\mathfrak q\) are strictly decreasing and centered,
\(v,\delta>0\) on the open positive half-line.

For completeness, the fixed-horizon differentiation leading to the next
identity is recorded explicitly.  Write

\[
 \rho_{s,r}(x)={\phi_s(x)F_r(x)\over F_v(0)}.
\]

On the positive half-line use \(B\) as the coordinate and define

\[
 w(B)=2\rho_{s,r}(x(B))C(x(B)),\qquad
 \Phi_{\rm L}=2z^2-mC.
\]

The factor \(2\) accounts for evenness.  Since \(dB=C\,dx\),

\[
 Q_v'(s)=\int_0^1w\,dB,
 \qquad Q_v''(s)=2\int_0^1w\Phi_{\rm L}\,dB.
\]

Thus \(w\,dB\) is exactly the unnormalized \(C^2\)-biased law in the
slope coordinate.  Put

\[
 \Psi=\Phi_{\rm L}+{c\over s},\qquad
 \mathcal J=\int_0^1w\Psi\,dB,
\]

The exact identities are

\[
 I_\alpha=2\mathcal J,
\]

\[
 \Psi_s=CJ\Psi+2Cz(3zJ-CJ_B)
       -{c\over s}\left(CJ+{1\over s}\right).        \tag{2.3b}
\]

Differentiating \(\mathcal J\), including the fixed-slope derivative of
\(w\), and integrating its centered kernel using (2.3a) gives the fully
differentiated shifted identity

\[
 D:={s^2I_\alpha'(s)\over2Q_v'(s)}
 =E_\nu\{v[n(1+j)+uj]+2\delta\gamma-c(1+j)\}.          \tag{2.4}
\]

This includes the term \(-cEj\); omitting it is the normalization error in
the tempting tail-only argument.

Put

\[
 A=a(0),\quad k={j(0)\over A},\quad
 X=A-a,\quad Z=a+u^2-A,\quad d=EX,\quad e=EZ,          \tag{2.5}
\]

and

\[
 F=3X+2Z,\quad L=n^2+Z,\quad
 G=\mathfrak q(0)-\mathfrak q,\quad
 R_\gamma=j(0)-j+{Z\over2}.                           \tag{2.6}
\]

Because \(\mathfrak p=(A-c)-F\) and
\(\mathfrak q=\mathfrak q(0)-G\),

\[
 {1\over2}L'=n(1+j)+uj,\qquad R_\gamma'=\gamma.
\]

Two exact integrations by parts in (2.4) therefore give

\[
\boxed{
 D={1\over2}\operatorname{Cov}(F,L)
   +2\operatorname{Cov}(G,R_\gamma)-c(1+Ej).}          \tag{2.7}
\]

\subsection{Profile-shape lemmas}

Write

\[
 \kappa={j\over a},\qquad \theta={\gamma\over uj},
 \qquad \rho={n\over u}.
\]

The following three monotonicities hold on the positive half-line:

\[
 \kappa'\ge0,\qquad \theta'\ge0,\qquad \rho'>0.         \tag{3.1}
\]

The first follows at once from

\[
 \kappa'={3uj-\gamma\over a}\ge0.
\]

For the second, in inverse-slope coordinates put

\[
 R={CJ_B\over zJ}=3-\theta.
\]

Here $t$ is the remaining heat time and $\mathcal A$ is a smooth drift
coefficient; its explicit form is irrelevant to the maximum-principle sign.

It satisfies

\[
 R_t={C^2\over2}R_{BB}+\mathcal A R_B
 +Cm(R-3)(R-4)+CJ(R-2)(R-6).                           \tag{3.2}
\]

After differentiating in \(B\), the inhomogeneous source is

\[
 z\{-2m(R-3)(R-4)+J(R-2)^2(R-6)\}<0,\qquad 0\le R\le3. \tag{3.3}
\]

The terminal data are zero, the \(B=1\) derivative is negative, and the
center singularity is the regular radial five-dimensional operator after
writing \(R_B=BT\).  The forward maximum principle gives \(R_B\le 0\), hence
\(\theta_B\ge 0\).

For \(\rho\), the curvature field obeys the maximum-principle inequality
\(z-xz_x\ge 0\).  Exact substitution then gives

\[
 u(1+j)-n(a+j)>0,\qquad
 \rho'={u(1+j)-n(a+j)\over u^2}>0.                    \tag{3.4}
\]

The center value of \(\theta\) has the quantitative bound

\[
 \boxed{\theta_0\ge {k\over1+k}.}                      \tag{3.5}
\]

Here a subscript $0$ means evaluation at $B=0$, and a dot means
differentiation in $t$.  At the center the same equation (3.2), including
the singular drift limit, gives

\[
 \dot R_0\le C_0m\{(R_0-3)(R_0-4)+k(R_0-2)(R_0-6)\},
\]

\[
 \dot k=-{C_0m\over2}k(1+k)(4-R_0).                    \tag{3.6}
\]

At the threshold $R_0=2+1/(1+k)$, the reaction term in the preceding
inequality is $-2k/(1+k)$.  Equation (3.6) says that $k$ decreases,
so this threshold increases.  A first upward contact is impossible.
Therefore $R_0<2+1/(1+k)$, which is (3.5).

The derivatives with respect to \(X\) now show that all four functions in
(2.6) are increasing convex functions of \(X\), vanishing at zero.  For
example,

\[
 Z_X=\kappa,\qquad (R_\gamma)_X={\theta\kappa\over2},
\]

\[
 L_X={\rho\over a}+\kappa(\rho+1),                    \tag{3.7}
\]

and the corresponding formula for \(G_X\) is

\[
 G_X={1\over2a}+1+\kappa+{\theta\kappa\over2}
       +{(1+\kappa)\rho\over2}.                       \tag{3.8}
\]

If \(H_1,H_2\) are increasing convex functions of \(X\) and vanish at zero,
size bias and the opposite-monotonicity covariance identity give

\[
 \operatorname{Cov}(H_1,H_2)
 \ge\beta EH_1EH_2,\qquad
 \beta={\operatorname{Var}X\over(EX)^2}.              \tag{3.9}
\]

To see this directly, write \(H_i=Xh_i\) with \(h_i\) increasing; apply that
covariance identity under the \(X^2\)-biased law, then use that this law stochastically
dominates the \(X\)-biased law.

\subsection{Density, moments, and the scalar lower bound}

Let

\[
 \rho_0={1+Ak\over A(1+k)}.                             \tag{4.1}
\]

The exact pushforward density identity, (3.4), and convexity of \(Z\) imply
that \(W=\sqrt{X}\) has density

\[
 f_W(w)=(A-w^2)^{(1+\rho_0)/2}h(w),\qquad h\ 
 \hbox{nonincreasing}.                                  \tag{4.2}
\]

This law is a mixture of Beta-type laws restricted to initial intervals
\([0,t]\) (equivalently, with the upper tail truncated).  For
$0<t<\sqrt A$, let $E_t$ denote expectation under the probability density
proportional to $(A-w^2)^{(1+\rho_0)/2}$ on $0\le w\le t$.  The
truncated second-moment curve is convex: if \(z=t^2\),
\(m(t)=E_tX\), and \(H(t)=E_tX^2\), then

\[
 {d\over dt}{dH\over dm}
 ={2tE_t[(X-z)^2]\over(z-m)^2}>0.                       \tag{4.3}
\]

The endpoint tangent and the decreasing-density moment inequality give

\[
 {d\over A}\le {1\over4+\rho_0},                       \tag{4.4}
\]

\[
 \beta\ge\max\left\{{4\over5},
 { (\rho_0+7)\delta-1\over(\rho_0+6)\delta^2}-1\right\},
 \qquad \delta={d\over A}.                             \tag{4.5}
\]

The fold means are

\[
 3d+2e=A-c,\qquad EF=A-c,\qquad EG=Ak+c.                \tag{4.6}
\]

From (3.4)--(3.5),

\[
 EL\ge\rho_0^2d+(1+\rho_0^2)e,                         \tag{4.7}
\]

\[
 R_\gamma\ge {1\over2}{k\over1+k}Z,\qquad
 Ej\le Ak+{e\over2(1+k)}.                              \tag{4.8}
\]

Substitution into (2.7) gives the rigorous scalar lower bound

\[
\boxed{
\begin{aligned}
D\ge \mathcal D(A,k,\delta;c):={}&
 {\beta\over2}(A-c)\{\rho_0^2d+(1+\rho_0^2)e\}\\
&+\beta(Ak+c){k\over1+k}e
-c\left\{1+Ak+{e\over2(1+k)}\right\},
\end{aligned}}                                         \tag{4.9}
\]

where \(d=A \delta\), \(e=(A-c-3d)/2\), and \(\beta\) is replaced by the right
side of (4.5).

There is one more exact profile constraint.  Integrating the differential
inequality for \(\kappa\) and using the lower-truncated \(t^{1/2}\) mixture of
the \(a/A\) density gives

\[
 e\le\Phi(k)d,\qquad
 \Phi(k)={5k(1+k)\over2+3k}.                            \tag{4.10}
\]

Thus

\[
 \delta\ge\delta_\Phi(A,k;c)
 ={A-c\over A\{3+2\Phi(k)\}}.                         \tag{4.11}
\]

Equations (4.4), (4.9), and (4.11) reduce local transversality to a
three-variable rational inequality.

\subsection{Exact parameter-domain reduction}

Contact alone gives the strict size bound

\[
 E_\nu a>1-{1\over P},                                 \tag{5.1}
\]

by the one-dimensional integration-by-parts argument.  Hence \(A>1-1/P\).

For fixed center data \((A,k)\), compare with the terminal profile

\[
 u_T=L\tanh(Ly),\quad a_T=A\operatorname{sech}^2(Ly),
 \quad b_T={A\over L}\tanh(Ly),\quad j_T=ka_T,
 \qquad L^2=A(1+k).
\]

Let

\[
 d\pi_T(y)\propto e^{-y^2/2}\cosh(Ly)^{1/(1+k)}\,dy.
\]

Since \((j/a)'\ge 0\), the quantity \(u^2+(1+k)a\) is increasing.  ODE
comparison gives

\[
 u\ge u_T,\qquad a\le a_T,\qquad b\le b_T.             \tag{5.2}
\]

A decreasing-input integral inequality, followed by the monotone density
ratio \(\mathrm d\pi/\mathrm d\pi_T\), yields

\[
 \boxed{P\le P_T(A,k)},                                 \tag{5.3}
\]

where, with \(\widehat m=1/(1+k)\) and \(h=A(1+k)\),

\[
 P_T(A,k)=
 {J_{-\widehat m}(h)-J_{2-\widehat m}(h)
  \over hJ_{4-\widehat m}(h)},
 \qquad
 J_r(h)=\int_{\mathbb R}e^{-z^2/(2h)}\operatorname{sech}^r z\,dz.
                                                               \tag{5.4}
\]

For fixed \(\widehat m=1/(1+k)\), the fixed-remaining-time primitive
argument used to prove (1.6) applies to this terminal profile: as a function of
\(A\), \(P_T(A,k)-P\) has exactly one zero and crosses it upward.

Indeed, at remaining time zero,
$F=\cosh^{\widehat m}$, $B=\tanh$, and
$C=\operatorname{sech}^2$.  The contact numerator in (1.3) is exactly

\[
 \int_{\mathbb R}\phi_h(z)\cosh^{\widehat m}z
 \{\tanh^2z-Ph\operatorname{sech}^4z\}\,dz
 ={hJ_{4-\widehat m}(h)\over\sqrt{2\pi h}}
 \{P_T(A,k)-P\}.
\]

The prefactor is positive and $A=\widehat m h$.  Thus the upward
one-crossing in $h$ proved in (1.6) is exactly an upward one-crossing in
$A$.  Hence

\[
 P_T(A_*,k)<P
 \quad\Longrightarrow\quad
 \bigl(P_T(A,k)\ge P\bigr)\Longrightarrow A>A_* .     \tag{5.4a}
\]

This is the transfer direction used in every exclusion below.

The exact terminal certificates give the following implications whenever
\(A<4\):

\[
\begin{array}{c|c}
P&\hbox{remaining possible center region}\\ \hline
2& {8\over25}\le k\le {541\over360},\quad
 A\ge {259\over200}+{9\over5}k,\\[1mm]
3& {1\over4}\le k\le {17\over30},\quad
 A\ge {23\over10}+3k,\\[1mm]
4& \hbox{empty},\\
P\ge5&\hbox{empty}.
\end{array}                                             \tag{5.5}
\]

For completeness, the exact interval partition behind (5.5) is:

Here ``controlled below $L(k)$'' means that the verifier proves
$P_T(A,k)<P$ for every $1-1/P\le A\le L(k)$.  ``Admits no
$A\le4$'' means that it proves the same strict inequality for every
$1-1/P\le A\le4$.  In view of (5.4a), these are exclusions of whole
continuous regions, not checks along their boundaries.

\begin{itemize}
\item For \(P=2\), \(\widehat m \in [10/13,1]\) is controlled below
  \(259/200+(8/5)k\); \(\widehat m \in [360/901,10/13]\) is controlled below
  \(259/200+(9/5)k\); and \(\widehat m \in [1/5,360/901]\) admits no \(A\le 4\).
\item For \(P=3,4,5\), the low-\(k\) affine lines are respectively
  \(23/10+3k\), \(3+(7/2)k\), and \(18/5+4k\); the complementary ranges
  \(\widehat m \in [1/5,30/47]\), \([1/5,7/9]\), and \([1/5,10/11]\) admit no
  \(A\le 4\).
\item On \(\widehat m \in [1/20,1/5]\), an exact boundary certificate gives
  \(P_T(4,k)<2\); this stronger \(P=2\) bound is reused for every \(P\ge 2\).
  For \(\widehat m\le 1/20\),

  put $h=4/\widehat m$.  Since
  $\cosh^{\widehat m}z\tanh^2z\le e^{\widehat m|z|}$,
  completing the square gives

\end{itemize}
\[
\begin{aligned}
J_{-\widehat m}(h)-J_{2-\widehat m}(h)
&\le 2e^{2\widehat m}\sqrt{\frac{8\pi}{\widehat m}},\\
hJ_{4-\widehat m}(h)
&\ge \frac8{\widehat m}e^{-\widehat m/8}
       \operatorname{sech}^4(1).
\end{aligned}
\]

  In the second line we only integrate over $[-1,1]$.  Dividing the two
  bounds gives

\[
 P_T(4,k)\le e^{17\widehat m/8}\cosh^4(1)
              {\sqrt{8\pi\widehat m}\over4}<2.         \tag{5.6}
\]

  For the last strict inequality, first note that

\[
 \log{9\over8}=\int_8^9{dt\over t}>{1\over9}>{17\over160}.
\]

  Also, after the terms (1+1/2), the ratios of consecutive terms in the
  series for (cosh(1)) are at most (1/30).  Therefore

\[
 \cosh(1)<{3\over2}+{1/24\over1-1/30}
 ={179\over116}<{31\over20}.
\]

  Finally,

\[
 {22\over7}-\pi=\int_0^1{t^4(1-t)^4\over1+t^2}\,dt>0.
\]

  Since $\widehat m\le1/20$, it remains only to use the rational
  comparisons

\[
 {44\over35}<\left({9\over8}\right)^2,
 \qquad
 \left({9\over8}\right)^2\left({31\over20}\right)^4<8.
\]
  Fixed-\(\widehat m\) terminal one-crossing then excludes every \(A<4\).

For clarity, derive the rational center bound used in the low-\(k\)
exclusions.  Put \(x=E_\nu a=A(1-\delta)\) and
\(U=E_\nu u^2=(x-c)/2\).  The odd one-dimensional variance estimate and the
lower-truncated \(T^{1/2}\) chord bound give

\[
 Ak\left(1+{5\delta\over4}\right)
 \ge {4\over3}(U-c)+{1\over3}(1-\delta){4A\over1+3A}.
\]

The right side is \((1-\delta)J_1(A)-2c\), where
\(J_1(A)=2A(A+1)/(1+3A)\).  Rearranging gives, whenever the numerator is
positive,

\[
 \delta\ge
 {J_1(A)-2c-Ak\over J_1(A)+(5/4)Ak},\qquad
 J_1(A)={2A(A+1)\over1+3A},                            \tag{5.7}
\]

together with (4.4).  Exact polynomial-basis conversion proves that the lower
bound in (5.7) exceeds the upper bound in (4.4) on the excluded boxes.

More explicitly, this strict comparison is certified on the following
strips, always from the stated affine lower boundary up to $A=4$:

\[
\begin{array}{c|c|c}
P&k\text{-interval}&\text{affine lower boundary}\\ \hline
2&[0,3/10]&259/200+(8/5)k\\
2&[3/10,8/25]&259/200+(9/5)k\\
3&[0,1/4]&23/10+3k\\
4&[0,2/7]&3+(7/2)k\\
P\ge5&[0,1/10]&18/5+4k.
\end{array}                                             \tag{5.8}
\]

On each strip the exact inequality is

\[
 {J_1(A)-2c-Ak\over J_1(A)+(5/4)Ak}
 >{1\over4+\rho_0}.
\]

The left side is a lower bound for $\delta$ by (5.7), while the right
side is an upper bound by (4.4), so the strip is impossible.  For
$P\ge5$, the cleared difference decreases with $c$, and the verifier
uses the worst endpoint $c=1/2$.  To check a polynomial on a box, the
program first maps the box to $[0,1]^2$.  If
$\sum c_{ij}s^it^j$ has degrees at most $n_1,n_2$, its coefficient in
the positive basis $s^u(1-s)^{n_1-u}t^v(1-t)^{n_2-v}$, after the usual
binomial normalization, is

\[
 \sum_{i\le u,\,j\le v}c_{ij}
 {\binom ui\over\binom{n_1}i}
 {\binom vj\over\binom{n_2}j}.
\]

Every one of these rational numbers is positive in the printed verifier.
This proves positivity on the whole box and explains exactly how the
table (5.5) follows from the terminal and density certificates.

The terminal certificates are analytic interval proofs, not quadrature.
The regional dyadic verifier uses a degree-six interpolation majorant with a
signed seventh-derivative remainder; the \(A=4\) verifier uses the separate
degree-ten positive-series majorant.  Gaussian tails are enclosed by the
even/odd continued-fraction convergents, and the denominator uses the
Fourier-product lower bound.  Every final comparison is an outward-rounded
integer comparison.

\subsection{Exact positivity of the scalar lower bound}

Three exact continuous-box certificates finish (4.9):

\begin{enumerate}
\item \(A\ge 4\), simultaneously for every \(1/4\le c\le 1/2\) and every \(k\ge 0\).
   Use \(h=1/A \in [0,1/4]\) and \(\widehat m=1/(1+k) \in [0,1]\) to compactify the
   domain.
\item The compact \(P=2\) box in (5.5).
\item The compact \(P=3\) box in (5.5).

\end{enumerate}
The verifier maps each box affinely to a unit cube and evaluates the
rational right side of (4.9) using signed \texttt{\detokenize{__int128}} intervals at scale
\(2^45\).  Addition, subtraction, multiplication, and division are rounded
outward and checked for overflow.  A box that straddles the switch in
(4.5) uses the lower branch \(\beta=4/5\) only after the coefficient of
\(\beta\) has been certified positive.  In the large-\(A\) program the
quantity called \texttt{out} is
\(\widehat m h^2\mathcal D(A,k,\delta;c)\); this has the same sign as
\(\mathcal D\) for finite \(A\) and \(k\).  The boundary values \(h=0\)
and \(\widehat m=0\) only certify the corresponding limits.  In the two
compact programs \texttt{out} is \(\mathcal D\) itself.  In all three
programs the unit coordinate called \texttt{r} runs through the entire
interval for \(\delta\) between (4.11) and (4.4).  Adaptive subdivision
only reduces interval dependency; no floating-point sign is accepted.

The reproduced outputs are

\begin{Verbatim}[fontsize=\small,frame=single]
PASS exact large-A corrected-D certificate
visited=333280 leaves=160375 empty=7289 maxdepth=20
smallest accepted scaled lower=1.25560859487450216e-08

P2 PASS exact dyadic D scalar certificate
visited=129874 leaves=65961 maxdepth=15
smallest accepted lower=6.60175999200873775e-07

P3 PASS exact dyadic D scalar certificate
visited=6934 leaves=4491 maxdepth=5
smallest accepted lower=4.06296681774165336e-05
\end{Verbatim}

Thus \(D>0\) at every simultaneous contact--shifted fold.  By (2.4),

\[
 \boxed{I_\alpha'(s)>0.}                                \tag{6.1}
\]

\subsection{Endpoint stability selects the good branch}\label{ft:sec:endpoint}

The endpoint assumptions become

\[
 Q_v(v)=q,\qquad PvQ_v'(v)=q\Gamma_\mu'(q)<q=Q_v(v).    \tag{7.1}
\]

Since the fixed-remaining-time contact numerator crosses upward, (7.1)
implies

\[
 v>\sigma(0)=V(0).                                     \tag{7.2}
\]

Every critical point of \(V\) is a strict minimum by (1.11) and (6.1).
Hence \(V'\) can cross zero only from negative to positive.  If
\(V(r)=v>V(0)\), the mean-value theorem gives a point before \(r\) where
\(V'>0\); a later nonpositive derivative would force a forbidden
positive-to-negative crossing.  Therefore \(V'(r)>0\) at every contact of
the stable horizon.

From (1.10),

\[
 V'(r)>0\iff I_\alpha(s)<0.                             \tag{7.3}
\]

Finally (1.2) gives

\[
 \Gamma_\mu''(u)=\left({Ps\over u}\right)^2I_\alpha(s)<0.
\]

This proves Proposition~\ref{prop:left-computation}.

\subsection{The slope-coordinate inequalities}\label{ft:sec:slope-coordinate}

This subsection proves the differential inequalities used in the main
argument. They look less mysterious if one remembers that the slope, not
the spatial variable, is being used as the coordinate.

Fix $0<m<1$, let $Z$ be standard normal, and, for $r\ge0$, set

\[
 U(r,x)={1\over m}\log \mathbb E\cosh^m(x+\sqrt r Z),
 \qquad B=U_x(r,x),\qquad C=U_{xx}(r,x).
\]

Under the probability weight proportional to
$\cosh^m(x+\sqrt r Z)$, differentiation gives
\[
 B=\mathbb E[\tanh(x+\sqrt rZ)],\qquad
 C=\mathbb E[\operatorname{sech}^2(x+\sqrt rZ)]
   +m\operatorname{Var}(\tanh(x+\sqrt rZ)).
\]
The expectation and variance in this line are taken under that weight.
Since $\operatorname{sech}^2=1-\tanh^2$,
\[
 C=1-(1-m)\mathbb E[\tanh^2(x+\sqrt rZ)]
      -m\{\mathbb E[\tanh(x+\sqrt rZ)]\}^2\leq1.
\]
Thus $-1<B<1$, $0<C\leq1$, and $x\mapsto B$ increases from $-1$ to $1$.
In this subsection every $B$-derivative is taken after writing
$C=C(r,B)$. Define

\[
 K=-{C_B\over2B}-m,\qquad
 J=K+BK_B=-m-{1\over2}C_{BB},\qquad
 z=-{1\over2}C_B=(m+K)B.                              \tag{A.1}
\]

The quotient $C_B/B$ extends smoothly through $B=0$, by evenness.
We shall prove

\[
 K\ge0,\qquad J\ge0,\qquad J_B\ge0,
 \qquad 0\le CJ_B\le3zJ.                              \tag{A.2}
\]

Only the last inequality is counterintuitive. It says that the growth of
$J$ in the slope coordinate is controlled by the decay of the curvature.

The forward heat time $r$ satisfies

\[
 U_r={1\over2}(U_{xx}+mU_x^2).
\]

Differentiating at fixed $B$ gives

\[
 x_r=KB,\qquad C_r=-C^2J.                              \tag{A.3}
\]

Two further $B$-derivatives give

\[
 J_r={C^2\over2}J_{BB}+2CC_BJ_B
 +(C_B^2-2C(J+m))J,                                   \tag{A.4}
\]

\[
 \begin{split}
 (J_B)_r={}&{C^2\over2}(J_B)_{BB}+3CC_B(J_B)_B\\
 &+(3C_B^2-6C(J+m)-2CJ)J_B-6C_B(J+m)J.
 \end{split}                                          \tag{A.5}
\]

We use the following elementary comparison principle. Suppose a function
is continuous on a closed bounded strip, is smooth in its interior, and
satisfies
\[
 f_r=d f_{BB}+b f_B+cf+g,\qquad d\ge0,
\]
with bounded coefficients, nonnegative initial and boundary values, and
$g\ge0$. Then $f\ge0$. Indeed, replace $f$ by
$e^{-\lambda r}f$, where $\lambda>\sup c$; a negative interior minimum
contradicts the equation. Applying this argument to $-f$ gives the
reversed statement.

We next check the endpoints. With $\delta=1-B$, direct expansion at
$x=+\infty$ gives, uniformly for bounded $r$,

\[
 C(r,B)=2\delta+c_2(r)\delta^2+c_3(r)\delta^3
 +O(\delta^4),                                        \tag{A.6}
\]

including one $r$-derivative and three $\delta$-derivatives. To justify
the uniform remainder, write $t=e^{-2x}$ and split the Gaussian integral
according to whether $x+\sqrt rZ\ge x/2$. On the first event, expand
$(1+e^{-2y})^m$ through $e^{-6y}$; on the second, the Gaussian tail is
smaller than every fixed power of $t$. Taking the logarithm gives
$\delta$ as a positive multiple of $t$ plus $O(t^2)$, which can be
inverted.

Comparing coefficients in (A.3) yields

\[
 c_2'=4(m+c_2),\qquad
 c_3'=12c_3+4c_2(m+c_2).                              \tag{A.7}
\]

At $r=0$, $C=1-B^2$, so $c_2(0)=-1$ and $c_3(0)=0$.
Therefore
\[
 m+c_2(r)=-(1-m)e^{4r},\qquad c_2(r)<0.
\]
The source in the second equation of (A.7) is nonnegative, so variation
of constants gives $c_3(r)\ge0$. Consequently

\[
 J(r,1)=(1-m)e^{4r}>0,\qquad J_B(r,1)=3c_3(r)\ge0.     \tag{A.8}
\]

Also, at $r=0$,

\[
 C=1-B^2,\qquad K=J=1-m,\qquad K_B=J_B=0.             \tag{A.9}
\]

Apply the comparison principle in (A.4) to the even extension in
$-1\le B\le1$. Its endpoint values are nonnegative by (A.8), so
$J\ge0$. Since $J=(BK)_B$ and $BK\to0$ as $B\downarrow0$,

\[
 BK(B)=\int_0^B J(t)\,dt,
\]

and hence $K\ge0$. In particular, $C_B=-2(m+K)B\le0$. The source in
(A.5) is now nonnegative. Its boundary values are nonnegative by parity
at $B=0$ and by (A.8) at $B=1$. A second application of the comparison
principle gives $J_B\ge0$.

It remains to prove the upper bound on $CJ_B$. Put

\[
 E=C^{3/2}J.
\]

Substitution of (A.3)--(A.4) gives

\[
 \begin{split}
 E_r={}&{C^2\over2}E_{BB}+{CC_B\over2}E_B\\
 &-\left\{{C_B^2\over8}+{C(J+m)\over2}\right\}E
 -{3\over2}C^{-1/2}E^2.
 \end{split}                                          \tag{A.10}
\]

Differentiating once more gives

\[
 \begin{split}
 (E_B)_r={}&{C^2\over2}(E_B)_{BB}
 +{3CC_B\over2}(E_B)_B\\
 &+\left\{{3\over8}C_B^2-{3\over2}C(J+m)-3CJ\right\}E_B\\
 &-{C\over2}J_BE+{3\over4}C^{-3/2}C_BE^2.
 \end{split}                                          \tag{A.11}
\]

The last two terms are nonpositive because $J_B\ge0$, $E\ge0$, and
$C_B\le0$. The apparently singular term is harmless:
$C^{-3/2}E^2=C^{3/2}J^2$. Thus its product with $C_B$ extends
continuously to $B=1$ and vanishes there; all coefficients needed for
comparison are bounded.

From (A.9),
\[
 E_B(0,B)=-3(1-m)B\sqrt{1-B^2}\le0.
\]
Parity gives $E_B(r,0)=0$, while (A.6) gives $E_B(r,1)=0$. The reversed
comparison principle applied to (A.11) yields $E_B\le0$. Finally,

\[
 E_B=C^{1/2}(CJ_B-3zJ).
\]

Since $C>0$, this proves $CJ_B\le3zJ$. The remaining inequalities in
(A.2) were proved above.

\subsection{One-crossing and differentiation}

\subsubsection{Fixed-remaining-time one-crossing}

Fix \(0<m\leq 1\), \(r\geq0\), and put

\[
 F(x)=\mathbb E\cosh^m(x+\sqrt r Z),\qquad
 B={1\over m}(\log F)',\qquad C=B'.
\]

On \(x>0\), \(B,C>0\).  Write

\[
 a={xC\over B},\qquad t=mxB,
 \qquad I(x)=\int_0^xF(y)C(y)^2\,dy,
 \qquad \mathcal M(x)={F(x)B(x)^2\over xI(x)}.
\]

The standard slope cone gives

\[
 z=-{C'\over2C}\geq mB.                                      \tag{1}
\]

\paragraph{A useful reverse-hazard factorization}

Use \(B\) itself as the coordinate and set

\[
 h(B(x))=F(x)C(x),\qquad e(x)={B(x)h(B(x))\over I(x)}.
\]

Since \(dB=C\,dx\),

\[
 I(x)=\int_0^{B(x)}h(b)\,db,
 \qquad \mathcal M={e\over a}.                              \tag{2}
\]

Moreover,

\[
 (FC)'=FC(mB-2z)<0                                          \tag{3}
\]

by (1).  Hence \(h\) is strictly decreasing and

\[
 0<e={Bh(B)\over\int_0^Bh(b)\,db}<1.                        \tag{4}
\]

Direct differentiation, using \(F'/F=mB\), gives

\[
 x(\log \mathcal M)'=t-1+2a-ae.                            \tag{5}
\]

Consequently

\[
 t\geq1\quad\Longrightarrow\quad
 x(\log \mathcal M)'>t-1+a>0.                              \tag{6}
\]

\paragraph{The whole region \(t\leq1\) lies below level (2)}

Let

\[
 J(x)=\int_0^x{dy\over F(y)}.
\]

The elementary product inequality gives

\[
 B(x)^2=\left(\int_0^xC(y)\,dy\right)^2\leq I(x)J(x).       \tag{7}
\]

Since \(B\) is increasing, for \(0\leq y\leq x\),

\[
 {F(x)\over F(y)}
 =\exp\left(m\int_y^xB(w)\,dw\right)
 \leq\exp\left(t\left(1-{y\over x}\right)\right).        \tag{8}
\]

Combining (7)--(8),

\[
 \mathcal M(x)\leq {F(x)J(x)\over x}
 \leq\int_0^1e^{t(1-u)}\,du
 ={e^t-1\over t}.                                          \tag{9}
\]

The last function is increasing for \(t\geq0\).  Therefore

\[
 t\leq1\quad\Longrightarrow\quad
 \mathcal M(x)\leq e-1<2.                                 \tag{10}
\]

Equations (6) and (10) prove that, for every real \(P\geq2\),
\(\mathcal M-P\) has at most one zero, and every zero is crossed strictly
from negative to positive. In fact it has exactly one zero:
\(\mathcal M(0+)=1\), while \(\mathcal M(x)\to\infty\) as \(x\to\infty\).
The latter follows from the
large-\(x\) asymptotics \(F\asymp e^{mx}\), \(B\to1\), and
\(C=O(e^{-2x})\).

\paragraph{Laplace variation diminution}

Put \(X=x^2\) and

\[
 f_0(X)=X^{-1/2}F(\sqrt X)B(\sqrt X)^2,
 \qquad f_1(X)=X^{-1/2}F(\sqrt X)C(\sqrt X)^2.
\]

The unnormalized fixed-\(r\) contact numerator is

\[
 \mathcal N(s,r)=\int_{\mathbb R}\phi_s(x)F_r(x)
 \{B_r(x)^2-PsC_r(x)^2\}\,dx .                            \tag{11}
\]

After the change \(X=x^2\), its sign is the sign of

\[
 \mathcal A(s,r)
 =2z\,\mathcal L f_0(z)-P\,\mathcal L f_1(z),
 \qquad z={1\over2s};                                      \tag{12}
\]

indeed, \(\mathcal N\) is a strictly positive multiple of
\(s\mathcal A\).
Integration by parts writes (12) as the Laplace transform of
\(k_P=2f_0'-Pf_1\).  Its cumulative primitive is

\[
 K_P(x^2)=\int_0^{x^2}k_P(X)\,dX
 =2{FB^2\over x}-2P\int_0^xFC^2
 =2I(x)(\mathcal M(x)-P).                                 \tag{13}
\]

Thus \(K_P\) has exactly one sign change, from negative to positive.
A second integration by parts gives

\[
 \mathcal A(s,r)=z\int_0^\infty e^{-zX}K_P(X)\,dX.        \tag{14}
\]

All boundary terms here vanish.  At the origin,

\[
 B(x)=C(0)x+O(x^3),\qquad
 K_P(x^2)=2(1-P)F(0)C(0)^2x+O(x^3),
\]

so \(f_0(0)=K_P(0)=0\) and the possible \(X^{-1/2}\) singularity is
integrable.  At infinity, for fixed \(r\),

\[
 F(x)=2^{-m}e^{mx+m^2r/2}(1+o(1)),\quad B(x)\to1,
 \quad C(x)=O(e^{-2x}),
\]

which implies \(I(\infty)<\infty\), \(\mathcal M(x)\to\infty\), and
\(e^{-zX}f_0(X),e^{-zX}K_P(X)\to0\) for every \(z>0\).

The elementary one-change lemma for the Laplace kernel now shows that
\(\mathcal A(\,\cdot\,,r)\) has at most one zero.  More explicitly, if
(14) vanishes, and \(X_0\) is the sign-change point of \(K_P\), then

\[
 \int_0^\infty (X-X_0)e^{-zX}K_P(X)\,dX>0.
\]

Hence the derivative of the integral in (14) with respect to \(z\) is
strictly negative at a zero.  Since \(z=1/(2s)\), every fixed-\(r\) contact
is crossed strictly upward as \(s\) increases.

For completeness, the endpoint signs also follow directly from (14).
As \(s\downarrow0\), equivalently \(z\to\infty\), the negative small-\(X\)
part of \(K_P\) dominates, so \(\mathcal A<0\).  As \(s\to\infty\), the
positive large-\(X\) part dominates (indeed its sub-Gaussian factor grows
like \(e^{m\sqrt X}\)), so \(\mathcal A>0\).  Hence the zero is unique, not
merely unique if it exists.

This proves, without a pointwise sign theorem for the original kernel, the
fixed-remaining-time uniqueness and transversality needed in the contact
curve construction.

\paragraph{Exact contact-branch derivatives}

Let \(\sigma(r)\) denote the unique zero of
\(\mathcal N(\,\cdot\,,r)\), and set

\[
 V(r)=r+\sigma(r).
\]

The preceding result says

\[
 \mathcal N_s(\sigma(r),r)>0.                             \tag{15}
\]

Implicit differentiation therefore gives

\[
 \sigma'=-{\mathcal N_r\over\mathcal N_s},\qquad
 V'={\mathcal N_s-\mathcal N_r\over\mathcal N_s}.        \tag{16}
\]

At a critical point of \(V\), \(\sigma'=-1\), and hence

\[
 V''=-{\mathcal N_{ss}-2\mathcal N_{sr}+\mathcal N_{rr}
        \over\mathcal N_s}
     =-{(\partial_s-\partial_r)^2\mathcal N
        \over\mathcal N_s}.                              \tag{17}
\]

These formulas are unchanged if \(\mathcal N\) is multiplied by a smooth
positive factor: at a branch critical point both \(\mathcal N\) and
\((\partial_s-\partial_r)\mathcal N\) vanish.

For comparison with the original \(Q\)-notation, along a line of fixed
horizon \(v=s+r\), write

\[
 \mathcal N=Z_v\{Q_v(s)-PsQ_v'(s)\}.
\]

Since \(Z_v\) is constant along this line, at a contact

\[
 (\partial_s-\partial_r)\mathcal N
 =-Z_v\{(P-1)Q_v'(s)+PsQ_v''(s)\}.                         \tag{18}
\]

Thus \(V'>0\) is exactly the desired negative-curvature sign.  At a branch
critical point put \(\alpha=(P-1)/P\).  Then

\[
 Q_v''+{\alpha\over s}Q_v'=0,
\]

and direct differentiation yields

\[
 -(\partial_s-\partial_r)^2\mathcal N
 =Z_vPs\left(Q_v'''+{\alpha\over s}Q_v''
                    -{\alpha\over s^2}Q_v'\right).         \tag{19}
\]

Combining (15), (17), and (19),

\[
 V''>0
 \quad\Longleftrightarrow\quad
 \left(Q_v''+{\alpha\over s}Q_v'\right)'>0.               \tag{20}
\]

Accordingly, the first-primitive argument completely supplies the positive
denominator and the branch reduction.  The remaining local statement is
the strict sign in (20) at a simultaneous contact and shifted fold.

\subsubsection{Differentiation at a shifted fold}

We derive the exact fixed-horizon derivative needed at a shifted fold.
The only easily missed point is that differentiating the shift \(c/s\)
produces the full term \(-cE_\nu(1+j)\).

\paragraph{Setup and the derivative normalization}

Let

\[
 \alpha={P-1\over P},\qquad c={\alpha\over2}={P-1\over2P},
 \qquad
 \mathcal I(s)=Q''(s)+\alpha{Q'(s)\over s}.             \tag{1}
\]

Here \(Q=Q_v\).  We also write

\[
 \Phi_{\rm L}=2z^2-mC,\qquad K={z\over B}-m.
\]

In the slope coordinate \(B\in(0,1)\), let \(w\) be the unnormalized
half-line weight and put

\[
 W=\int_0^1w\,dB=Q'(s),\qquad
 \Psi=\Phi_{\rm L}+{c\over s},\qquad
 J=\int_0^1w\Psi\,dB.                                   \tag{2}
\]

Since \(Q''=2\int w\Phi_{\rm L}\),

\[
 \boxed{\mathcal I=2J.}                                 \tag{3}
\]

At a shifted fold \(J=0\), and hence

\[
 \boxed{
 {s^2\mathcal I'(s)\over2Q'(s)}={s^2J'(s)\over W}=:\mathscr D.}
                                                                  \tag{4}
\]

Thus the desired shifted derivative is positive exactly when
\(\mathscr D>0\). There is no sign-changing normalization factor in (4).

Use the fields \(C,z,J_0,N,K\), where \(J_0\) denotes the local
field (to distinguish it from the integral \(J\) in (2)). The identities
needed below are

\[
 N={x\over s}+KB,\qquad N_B={1\over sC}+J_0,             \tag{5}
\]

\[
 (\Phi_{\rm L})_s=CJ_0\Phi_{\rm L}+G,\qquad
 G=2Cz(3zJ_0-C(J_0)_B),                                 \tag{6}
\]

and

\[
 R_0=2CJ_0-{mC\over2}-{z^2\over2},\qquad
 (R_0)_B=2C(J_0)_B-5zJ_0.                               \tag{7}
\]

Because the shift is \(c/s\), (6) gives the crucial formula

\[
 \boxed{
 \Psi_s=CJ_0\Psi+G-{c\over s}\left(CJ_0+{1\over s}\right).}
                                                                  \tag{8}
\]

The last term in (8), especially its \(CJ_0\) component, is easy to miss;
it is exactly what produces the term \(-cE_\nu j\) in (16).

At \(J=0\), differentiating (2), using the exact fixed-slope derivative of
\(w\), and observing that

\[
 {x^2\over2s^2}+KB{x\over s}+{(KB)^2\over2}={N^2\over2},
\]

gives

\[
 \boxed{
 J'=\int_0^1w\left\{{N^2\over2}\Psi+R_0\Psi+G
 -{c\over s}\left(CJ_0+{1\over s}\right)\right\}dB.}    \tag{9}
\]

All endpoint terms below vanish by the Gaussian estimate established in
the preceding calculation.

\paragraph{Exact tail conversion}

Define

\[
 F_\Psi(B)=\int_B^1w(D)\Psi(D)\,dD,\qquad
 \tau={F_\Psi\over w}.                                  \tag{10}
\]

Because \(\Psi_B=2z(2J_0+3m)>0\) and
\(\int_0^1w\Psi\,dB=0\), the function \(\Psi\) changes sign once, from
negative to positive. Hence \(F_\Psi\ge0\) and \(\tau\ge0\).

To prove the upper bound, set
\[
 \mathcal H=wCz-F_\Psi,\qquad
 L_c=CJ_0+{c\over s}-z(N+z).
\]
Using \(w_B/w=-(N+z)/C\), \(z_B=m+J_0\), and
\((F_\Psi)_B=-w\Psi\), direct differentiation gives
\[
 \mathcal H_B=wL_c.
\]
Moreover,
\[
 (CJ_0)_B=-2zJ_0+C(J_0)_B\le zJ_0,
\]
whereas
\[
 [z(N+z)]_B=(m+J_0)(N+z)
 +z\left({1\over sC}+J_0+m+J_0\right)>zJ_0.
\]
Thus \((L_c)_B<0\). At \(B=0\), \(Cz=0\) and
\(F_\Psi(0)=\int_0^1w\Psi=0\). At \(B=1\), \(F_\Psi(1)=0\), while the
Gaussian endpoint estimate gives \(wCz\to0\). Therefore
\(\mathcal H(0)=\mathcal H(1)=0\) and
\(\int_0^1wL_c\,dB=0\). Since \(L_c\) is strictly decreasing, it is
positive before its unique zero and negative afterward. Thus
\(\mathcal H\) first increases and then decreases, so
\(\mathcal H\ge0\). Dividing by \(w\) proves
\[
 0\le\tau\le Cz.
\]
Introduce the standardized fields

\[
 a=msC,\quad u=\sqrt s\,z,\quad j=sCJ_0,\quad n=\sqrt s\,N,\quad
 \gamma=s\sqrt s\,C(3zJ_0-C(J_0)_B),                    \tag{11}
\]

\[
 v={\sqrt s\,\tau\over C},\qquad
 \delta={\sqrt s\,(Cz-\tau)\over C},\qquad u=v+\delta.   \tag{12}
\]

The two centered sources are

\[
 \mathfrak p=a-2u^2-c,\qquad
 \mathfrak q=j-u(n+u)+c.
\]

The probability law \(w\,dB/W\) is the standardized \(a^2\)-biased law
\(\nu\).

Let \(\omega\) be the density of \(\nu\) in the \(y\)-coordinate.
Since \(\Psi=-\mathfrak p/s\), \(dB/dy=\sqrt s\,C\), and \(\omega\) is
proportional to \(w\,dB/dy\), the definition of \(v\) gives
\[
 (\omega v)'=\omega\mathfrak p.
\]
Also \(\omega'/\omega=-(n+3u)\) and \(u'=a+j\), whence
\[
 {(\omega u)'\over\omega}=a+j-u(n+3u)=\mathfrak p+\mathfrak q.
\]
Since \(\delta=u-v\), subtraction yields
\[
 (\omega\delta)'=\omega\mathfrak q.
\]
The endpoint normalizations agree with those in the main proof, so these
\(v,\delta\) are exactly the two centered tails used there.

Since \(F_\Psi'=-w\Psi\), integration by parts and (5) give

\[
 {s^2\over W}\int_0^1w{N^2\over2}\Psi\,dB
 ={s^2\over W}\int_0^1w\tau NN_B\,dB
 =E_\nu[vn(1+j)].                                       \tag{13}
\]

Similarly, (7) gives

\[
 \begin{aligned}
 \int_0^1w(R_0\Psi+G)dB
 &=\int_0^1w\{\tau(2C(J_0)_B-5zJ_0)
       +2Cz(3zJ_0-C(J_0)_B)\}\,dB\\
 &=\int_0^1w\{zJ_0\tau
 +2(3zJ_0-C(J_0)_B)(Cz-\tau)\}\,dB.
 \end{aligned}                                               \tag{14}
\]

After multiplying by \(s^2/W\), the two terms in the last line become

\[
 E_\nu[vuj],\qquad 2E_\nu[\delta\gamma],                \tag{15}
\]

respectively. Finally, the shifted term in (9) is

\[
 -{s^2\over W}{c\over s}\int_0^1w
       \left(CJ_0+{1\over s}\right)dB
 =-cE_\nu(1+j).                                         \tag{16}
\]

Combining (9), (13)--(16) proves the exact normalized identity

\[
 \boxed{
 \mathscr D=E_\nu\{vn(1+j)+vuj+2\delta\gamma-c(1+j)\}.}  \tag{17}
\]

\subsection{Shape estimates for the logarithmic heat profile}

\subsubsection{Monotonicity of the profile quotient}

Fix \(0<m<1\) and write the logarithmic-heat curvature in inverse-slope
coordinates as

\[
 C=C(t,B),\qquad z=-\frac12 C_B,\qquad
 J=-m-\frac12C_{BB},\qquad w=J_B .
\]

Here \(t\ge 0\) is the remaining heat time and \(0<B<1\).  The standard signs
are

\[
 C,z,J>0,\qquad w\ge0,\qquad
 g:=3zJ-Cw\ge0 .                                      \tag{1}
\]

They imply

\[
 0\le R:=\frac{Cw}{zJ}\le3,\qquad
 \theta:=\frac{g}{zJ}=3-R.                            \tag{2}
\]

We prove

\[
 \boxed{R_B\le0,\quad\text{equivalently}\quad\theta_B\ge0.} \tag{3}
\]

\paragraph{A closed parabolic equation for \(R\)}

At fixed slope \(B\),

\[
 C_t=-C^2J,                                             \tag{4}
\]

and differentiation gives

\[
\begin{aligned}
 z_t&=-2CzJ+\frac12C^2w,\\
 J_t&=\frac12C^2J_{BB}-4Czw+
       \{4z^2-2C(m+J)\}J,\\
 w_t&=\frac12C^2w_{BB}-6Cz w_B
      +\{12z^2-6Cm-8CJ\}w+12z(m+J)J .
\end{aligned}                                         \tag{5}
\]

Logarithmically differentiating \(R=Cw/(zJ)\), substituting (5), and
collecting its second and first \(B\) derivatives gives the closed equation

\[
 \boxed{
 R_t=\frac12C^2R_{BB}+A R_B
     +Cm(R-3)(R-4)+CJ(R-2)(R-6),}                     \tag{6}
\]

where

\[
 A=\frac{C^2(m+J)}z+\frac{C^2w}J-4Cz.                \tag{7}
\]

For reproducibility, here is a short hand calculation of the only
non-obvious collection in (6).  Put

\[
 A_0=-\frac{2z}{C}-\frac{m+J}{z}-\frac wJ
     =\partial_B\log R-\frac{w_B}{w}.
\]

Then

\[
 \frac{w_{BB}}w=\frac{R_{BB}}R-A_0'-A_0^2
                    -2A_0\frac{w_B}{w}.
\]

The \(w_B/J\) term in \(A_0'\) cancels the corresponding term coming from
\(-J_t/J\).  The coefficient of \(w_B/w\) that remains is exactly (7).
After replacing \(w=RzJ/C\), the zeroth-order remainder is

\[
 Cm(R^2-7R+12)+CJ(R^2-8R+12),
\]

which is the last line of (6).

\paragraph{Differentiate once in \(B\)}

Set \(Y=R_B\).  Differentiating (6) gives

\[
 Y_t=\frac12C^2Y_{BB}+(CC_B+A)Y_B+K Y+\mathcal S,     \tag{8}
\]

where \(K=A_B+\partial_R F\) is a locally bounded zeroth-order coefficient
away from \(B=0\), and the inhomogeneous source is

\[
\begin{aligned}
 \mathcal S
 &= (Cm)_B(R-3)(R-4)+(CJ)_B(R-2)(R-6)\\
 &=z\{-2m(R-3)(R-4)+J(R-2)^2(R-6)\}.
\end{aligned}                                         \tag{9}
\]

In the second equality we used

\[
 (Cm)_B=-2mz,\qquad (CJ)_B=-2zJ+Cw=zJ(R-2).
\]

The interval restriction (2) now gives the decisive strict sign:

\[
 \boxed{\mathcal S<0\qquad(t>0,\ 0<B<1).}             \tag{10}
\]

Indeed, \((R-3)(R-4)\ge 0\), while
\((R-2)^2(R-6)\le 0\); at the possible zeros of one summand the other one is
strictly negative.

\paragraph{Initial and boundary data; the center singularity}

At \(t=0\),

\[
 C=1-B^2,\qquad J=1-m,\qquad w=0,
\]

so \(R=Y=0\).  Evenness gives \(Y(t,0)=0\).

At the other endpoint, the heat-kernel expansion can be made explicit.
If \(h=\exp(-2x)\) and \(\rho=\exp(4t)\), then, up to a common positive
multiple,

\[
 P_t(\cosh^m)(x)=e^{mx+m^2t/2}
 \{1+d h+e h^2+f h^3+O(h^4)\},
\]

where

\[
\begin{aligned}
d&=m e^{2(1-m)t},\\
e&={m(m-1)\over2}e^{(8-4m)t},\\
f&={m(m-1)(m-2)\over6}e^{(18-6m)t}.
\end{aligned}
\]

(The omitted contribution from the opposite Gaussian tail is smaller
than every power of \(h\).)  Expanding first
\(B=(\log P_t(\cosh^m))_x/m\) and then eliminating \(h\) in favor of
\(\epsilon=1-B\) yields

\[
 C=2\epsilon+A_2\epsilon^2+A_3\epsilon^3+O(\epsilon^4),
\]

with

\[
 A_2=-\{m+(1-m)\rho\},\qquad
 A_3={ (1-m)\rho\over2}
 \{(2-m)\rho^2-2(1-m)\rho-m\}.
\]

Consequently

\[
 z\to1,\qquad J\to J_1=(1-m)\rho>0,
 \qquad w\to w_1=3A_3>0,                         \tag{11a}
\]

because \(\rho>1\) and the expression in braces vanishes at \(\rho=1\)
and has strictly positive derivative thereafter.  Therefore

\[
 R=2\frac{w_1}{J_1}(1-B)+O((1-B)^2),\qquad
 Y(t,1)=-2\frac{w_1}{J_1}<0.                          \tag{11b}
\]

The quotient defining \(R\) has a removable singularity at \(B=0\).
Indeed, \(C,J\) are even while \(z,w\) are odd, and
\(z_B(t,0)=m+J(t,0)>0\); hence \(w/z\) and therefore \(R\) extend smoothly
and evenly to the center.  At \(t=0\) this extension is identically zero.
The apparent singularity of \(A\) at \(B=0\) is the regular radial one.  More
explicitly, the same Taylor expansions give

\[
 z=(m+J(t,0))B+O(B^3),\qquad
 A=\frac{C(t,0)^2}{B}+O(B).
\]

In (8),

\[
 CC_B+A=\frac{C(t,0)^2}{B}+O(B),\qquad
 K=-\frac{C(t,0)^2}{B^2}+O(1),\qquad \mathcal S=O(B).
\]

Writing the odd function \(Y=BT\), the two displayed singular zeroth-order
terms cancel.  The principal radial part of the equation for the even
function \(T\) is

\[
 \frac12C^2\left(T_{BB}+\frac4B T_B\right),
\]

the radial Laplacian in five dimensions; all remaining coefficients and
the divided source \(\mathcal S/B\) extend continuously across \(B=0\).
Thus the standard radial parabolic maximum principle applies on the unit
five-ball.  Equivalently, one can first work on compact subcylinders and
use the preceding Taylor expansions to pass to the center.

After multiplying by an exponential in time to absorb the bounded
zeroth-order coefficient, (8)--(11) and the maximum principle give

\[
 Y\le0.
\]

This proves (3).  In standardized variables,

\[
 \theta=\frac{\gamma}{uj}=\frac{g}{zJ},
\]

so the conclusion is precisely that \(\gamma/(u j)\) is nondecreasing along
the positive half-line.

\subsubsection{The center lower bound}

Fix \(0<m<1\).  The time variable \(t\) below is the remaining Gaussian
variance in the logarithmic heat flow: \(t=0\) is the terminal profile and
increasing \(t\) is the forward heat-flow direction.  In inverse-slope
coordinates put

\[
 C=C(t,B),\qquad z=-\frac12C_B,\qquad
 J=-m-\frac12C_{BB},\qquad w=J_B,
\]

and, for \(0<B<1\),

\[
 R=\frac{Cw}{zJ}=3-\theta .                         \tag{1}
\]

All four fields have their usual smooth parity at the center: \(C,J,R\)
are even and \(z,w\) are odd.  We prove

\[
 \boxed{\displaystyle
   \theta_0(t)>\frac{k(t)}{1+k(t)},\qquad
   k(t):=\frac{J(t,0)}m,\qquad
   \theta_0(t):=\lim_{B\downarrow0}\theta(t,B).}
                                                               \tag{2}
\]

The non-strict version is enough in every later application. The
argument here uses only the monotonicity \(R_B\le0\) proved in the
preceding subsection and the center evolution.

\paragraph{The exact center evolution of \(R\)}

The closed equation proved in the preceding subsection is

\[
 R_t=\frac12C^2R_{BB}+A R_B
 +Cm(R-3)(R-4)+CJ(R-2)(R-6),                         \tag{3}
\]

where

\[
 A=\frac{C^2(m+J)}z+\frac{C^2w}J-4Cz.               \tag{4}
\]

That same maximum-principle argument proves \(R_B\le 0\) on the positive
half-interval.  Since \(R\) is smooth and even at the center,

\[
 R_{BB}(t,0)\leq0.                                   \tag{5}
\]

We now take the center limit in (3), retaining the apparently singular
drift.  Write \(C_0=C(t,0)\), \(J_0=J(t,0)\), and \(R_0=R(t,0)\).  Smooth parity
and \(z_B=m+J\) give

\[
\begin{aligned}
 C&=C_0+O(B^2),& J&=J_0+O(B^2),\\
 z&=(m+J_0)B+O(B^3),& w&=O(B),\\
 R&=R_0+\frac12R_{BB}(t,0)B^2+O(B^4).
\end{aligned}                                        \tag{6}
\]

Consequently

\[
 A=\frac{C_0^2}{B}+O(B),\qquad
 R_B=R_{BB}(t,0)B+O(B^3),                            \tag{7}
\]

and hence

\[
 \lim_{B\to0}\left(\frac12C^2R_{BB}+AR_B\right)
 =\frac32C_0^2R_{BB}(t,0)\leq0.                     \tag{8}
\]

Thus (3) gives the rigorous scalar differential inequality

\[
 \dot R_0\leq C_0m\left\{(R_0-3)(R_0-4)
       +k(R_0-2)(R_0-6)\right\}.                   \tag{9}
\]

The factor \(3/2\) in (8) is essential: the singular drift contributes
one full copy of \(C_0^2R_{BB}\), in addition to the one-half copy from
the displayed diffusion.

\paragraph{The exact center evolution of \(k\)}

The inverse-slope equations also give

\[
 z_t=-2CzJ+\frac12C^2w.
\]

Using \(Cw=RzJ\), this is

\[
 z_t=-\frac12CzJ(4-R).                               \tag{10}
\]

Compare the coefficients of \(B\) in (10), using (6).  Since \(m\) is fixed
along the heat flow,

\[
 \dot J_0=-\frac12C_0(m+J_0)J_0(4-R_0),
\]

or equivalently

\[
 \boxed{\displaystyle
 \dot k=-\frac12C_0m\,k(1+k)(4-R_0).}                \tag{11}
\]

The standard profile signs give \(C_0,k>0\) and \(0\le R_0\le 3\), so \(k\)
is strictly decreasing in the forward \(t\) direction.

\paragraph{A moving barrier}

Define

\[
 b(k)=2+\frac1{1+k}=\frac{3+2k}{1+k}.                \tag{12}
\]

At terminal time \(t=0\),

\[
 C=1-B^2,\qquad J=1-m,\qquad w=0,
\]

so \(R_0(0)=0<b(k(0))\).  Suppose there were a first forward time at which
\(R_0=b(k)\).  Substitution of (12) into the reaction polynomial in (9)
gives the exact cancellation

\[
 (b-3)(b-4)+k(b-2)(b-6)=-\frac{2k}{1+k}<0.           \tag{13}
\]

Therefore, at that alleged first contact,

\[
 \dot R_0\leq-\frac{2C_0mk}{1+k}<0.                 \tag{14}
\]

On the other hand, (11) and \(b'(k)=-(1+k)^{-2}\) give

\[
 \frac d{dt}b(k(t))
 =\frac12C_0m\,\frac{k(4-b)}{1+k}>0.                \tag{15}
\]

Hence \((R_0-b(k))'<0\) at a first upward zero, a contradiction.  Thus
\(R_0<b(k)\) for every finite \(t\ge 0\).  Since \(\theta_0=3-R_0\), this is
exactly (2):

\[
 3-R_0>3-\frac{3+2k}{1+k}=\frac{k}{1+k}.
\]

\paragraph{Consequences in the standardized fold variables}

At the center of a standardized profile,

\[
 a_0=msC_0,\qquad j_0=sC_0J_0,\qquad
 \frac{j_0}{a_0}=\frac{J_0}{m}=k_0.                 \tag{16}
\]

Therefore

\[
 \boxed{\displaystyle \theta_0\geq\frac{k_0}{1+k_0}.}       \tag{17}
\]

For \(X=a_0-a\), \(Z=a+u^2-a_0\), and

\[
 R_\gamma=j_0-j+\frac Z2,
\]

the exact slope identities are

\[
 Z_X=\frac ja=:k(X),\qquad
 (R_\gamma)_X=\frac{\theta(X)k(X)}2.                \tag{18}
\]

Since \(\theta\) is nondecreasing, (17)--(18) imply the pointwise and mean
bounds

\[
 \boxed{\displaystyle
 R_\gamma\geq\frac{\theta_0}2Z
 \geq\frac{k_0}{2(1+k_0)}Z,\qquad
 E R_\gamma\geq\frac{k_0}{2(1+k_0)}E Z.}            \tag{19}
\]

In particular, this proves the quantitative form
\(R_\gamma\ge0\) used later.

\subsubsection{Strict monotonicity of $n/u$}

Fix a remaining heat time \(t\ge 0\) and write

\[
 F(t,x)=P_t(\cosh^m)(x),\qquad
 \phi={1\over m}\log F,\qquad B=\phi_x,\quad C=B_x,
\]

\[
 z=-{C_x\over2C},\qquad J=-m-{1\over2}C_{BB}.
\]

Thus \(C_x=-2Cz\), \(z_B=m+J\), and the standard slope cone gives
\(J_B\ge 0\).  We prove first the star-shaped concavity

\[
                         z-xz_x\geq0\qquad(x\geq0).       \tag{1}
\]

\paragraph{A maximum principle for \(z/x\)}

The differentiated logarithmic heat equation is

\[
 B_t={1\over2}B_{xx}+mBB_x,
 \qquad C_t={1\over2}C_{xx}+m(C^2+BC_x).
\]

Logarithmically differentiating the second equation gives the exact PDE

\[
 z_t={1\over2}z_{xx}+(mB-2z)z_x+2mCz.                 \tag{2}
\]

For \(x>0\), put \(h=z/x\).  Equation (2) becomes

\[
 h_t={1\over2}h_{xx}+\left({1\over x}+mB-2xh\right)h_x
       +{mB\over x}h-2h^2+2mCh.                       \tag{3}
\]

Set \(H=h_x\).  At a hypothetical first interior point at which \(H\)
develops a nonnegative maximum through zero, all terms in the
differentiated equation that contain \(H\) or \(H_x\) vanish, while
\(H_{xx}\le 0\).  The remaining source is

\[
 mh\left({C\over x}-{B\over x^2}-4Cxh\right)
 ={mh\over x^2}\{xC-B-4Cx^2z\}<0.                    \tag{4}
\]

Here \(B(x)\ge xC(x)\) because \(C\) is decreasing on the positive half-line,
and \(m,h,C,z>0\).  Thus the parabolic first-contact argument excludes such
a maximum.

At \(t=0\), \(z=\tanh x\), so \(\tanh(x)/x\) is strictly decreasing.  At the
center, oddness gives the regular radial boundary condition \(H(t,0)=0\).
The heat-kernel expansion at infinity gives \(z\to1\) and hence
\(h_x\sim-x^{-2}<0\). Compact exhaustion, followed by the Gaussian endpoint
expansion, therefore justifies the maximum principle on the whole
half-line.  Consequently \(H\le 0\), proving (1), strictly for \(x>0\).

\paragraph{Translation to the standardized fields}

At an arbitrary outer clock \(s>0\), put \(y=x/\sqrt{s}\) and use

\[
 a=msC,\quad u=\sqrt{s}\,z,\quad j=sCJ,
 \quad n=\sqrt{s}\,N,
 \qquad N={x\over s}+z-mB.
\]

We claim

\[
                  u(1+j)-n(a+j)>0.                    \tag{5}
\]

Indeed, direct substitution gives

\[
 {u(1+j)-n(a+j)\over\sqrt{s}}
 =z-xC(m+J)+smCB(J-K),
 \qquad K={z\over B}-m.                                \tag{6}
\]

The first two terms equal \(z-xz_x\) because \(z_x=C(m+J)\), and are
nonnegative by (1).  Also

\[
 {z\over B}={1\over B}\int_0^B(m+J(\beta))\,d\beta
 \le m+J(B),                                           \tag{7}
\]

since \(J_B\ge 0\); hence \(J-K\ge 0\).  The inequalities are strict away from
the center, so (5) follows.

Finally \(u'=a+j\) and \(n'=1+j\).  Therefore

\[
 \boxed{
 \left({n\over u}\right)'
 ={u(1+j)-n(a+j)\over u^2}>0\qquad(y>0).}              \tag{8}
\]

This conclusion is universal for the two-atom logarithmic heat profile; no
contact or fold identities are needed.

\paragraph{Convexity consequence for the centered fold source}

At a contact--fold let \(X=a(0)-a\), \(k=j/a\),
\(\theta=\gamma/(uj)\), and \(\rho=n/u\).  The centered decreasing source

\[
 \mathfrak q=j-u(n+u)+c
\]

has drop \(G_{\mathfrak q}=\mathfrak q(0)-\mathfrak q\) satisfying

\[
 (G_{\mathfrak q})_X={1\over2a}+1+k+{\theta k\over2}
                  +{(1+k)\rho\over2}.                 \tag{9}
\]

The functions \(a^{-1}\), \(k\), \(\theta\), and \(\rho\) are all
nondecreasing: the first two follow from the displayed slope identities,
\(\theta\) was proved nondecreasing above, and \(\rho\) is (8). Hence
\((G_{\mathfrak q})_X\) is
nondecreasing.  Thus \(G_{\mathfrak q}\) is increasing and convex as a function of
\(X\), with \(G_{\mathfrak q}(0)=0\).

\subsection{Density, moments, and comparison estimates}

\subsubsection{Power-density representation and moment bounds}

At a contact--fold set

\[
 X=W^2=a_0-a,\qquad \eta={Z\over X},\qquad
 \rho={n\over u},\qquad \rho_0={1+j_0\over a_0+j_0}.
\]

The proved star-shapedness lemma gives \(\rho\ge \rho_0\), while convexity of
\(Z\) gives that \(\eta=Z/X\) is nondecreasing. Let \(\omega\) be the
positive-half-line density of \(\nu\). Since \(X=a_0-a=W^2\),
\[
 X'=2au,\qquad W'={au\over W}=a\sqrt{1+\eta},
\]
because \(u^2=X+Z=X(1+\eta)\). Hence change of variables gives
\[
 f_W(W(y))={\omega(y)\over a(y)\sqrt{1+\eta(y)}}.
\]
The score identities \(\omega'/\omega=-(n+3u)\) and \(a'/a=-2u\)
therefore imply

\[
 {d\over dy}\log\left\{
 {f_W\over\sqrt a}\sqrt{1+\eta}\right\}
 ={\omega'\over\omega}-{3\over2}{a'\over a}=-n.
\tag{1}
\]

Since \((\log a)'=-2u\), subtracting \((\rho_0/2)\log a\) from the logarithm
in (1) gives

\[
 {d\over dy}\log\left\{
 {f_W\sqrt{1+\eta}\over a^{(1+\rho_0)/2}}\right\}
 =-n+\rho_0u=-(\rho-\rho_0)u\le0.                         \tag{2}
\]

Both the bracket in (2) and \((1+\eta)^{-1/2}\) are nonincreasing.
Therefore

\[
 \boxed{f_W(w)=(a_0-w^2)^{(1+\rho_0)/2}h(w),
        \qquad h\ \hbox{nonincreasing}.}                 \tag{3}
\]

The extra exponent is what yields the sharper moment bounds below.

Under \(T=1-W^2/a_0\), its density is proportional to
\[
 t^{1/2}\widetilde h(t),\qquad
 \widetilde h(t)=t^{\rho_0/2}(1-t)^{-1/2}
 h\!\left(\sqrt{a_0(1-t)}\right).
\]
Every factor in \(\widetilde h\) is nondecreasing.

Put \(r=\rho_0\) and let \(Q_r\) be the probability law with density
proportional to \((a_0-w^2)^{(1+r)/2}\) on \((0,\sqrt{a_0})\).  Then
\(W^2/a_0\) under \(Q_r\) has the
\(\operatorname{Beta}(1/2,(3+r)/2)\) law. The ratio of the actual
density to the density of \(Q_r\) is nonincreasing. Thus the
double-integral covariance identity shows that the expectation of every
increasing function is no larger under the actual law. Taking that
function to be \(W^2\) yields

\[
 \boxed{d=EX\le {a_0\over4+r}.}                            \tag{4}
\]

The layer-cake representation of the decreasing factor \(h\) is a mixture
of truncations of \(Q_r\).  Repeating the truncated-moment convexity argument
gives the supporting tangent at the untruncated endpoint.  Here

\[
 E_{Q_r}X={a_0\over4+r},\qquad
 E_{Q_r}X^2={3a_0^2\over(4+r)(6+r)},
\]

and the endpoint tangent has slope \(a_0(r+7)/(r+6)\) and intercept
\(-a_0^2/(r+6)\).

For completeness, the small-truncation bound uses only that \(f_W\) is
nonincreasing, which follows from (3). Such a density is a mixture of
uniform laws on intervals \([0,t]\). For \(W\) uniform on \([0,t]\),
\[
 EW^2={t^2\over3},\qquad
 EW^4={t^4\over5}={9\over5}(EW^2)^2.
\]
Averaging over the mixing law and using convexity of \(t^2\) gives
\[
 EW^4={1\over5}Et^4\ge{1\over5}(Et^2)^2
 ={9\over5}(EW^2)^2.
\]
Since \(X=W^2\), this is \(EX^2\ge(9/5)(EX)^2\). Combining the two
bounds gives

\[
 \boxed{EX^2\ge\max\left\{{9\over5}d^2,
 {a_0(r+7)d-a_0^2\over r+6}\right\}.}                    \tag{5}
\]

Consequently

\[
 \boxed{\beta={\operatorname{Var}X\over d^2}\ge
 \max\left\{{4\over5},
 {a_0(r+7)d-a_0^2\over(r+6)d^2}-1\right\}.}              \tag{6}
\]

For $0<t<\sqrt{a_0}$, let $E_t$ denote expectation under $Q_r$
conditioned on $W\in[0,t]$.  The truncated second-moment curve is convex: if
\(z=t^2\), \(m(t)=E_tX\), and \(H(t)=E_tX^2\), then

\[
 {d\over dt}{dH\over dm}
 ={2tE_t[(X-z)^2]\over(z-m)^2}>0.
\]

Thus (4)--(6) are rigorous consequences of the already proved \(\rho\)
monotonicity, with no numerical input.

\subsubsection{The profile constraint $e\le\Phi(k)d$}

At a simultaneous contact--fold put

\[
 X=a_0-a,\qquad Z=a+u^2-a_0,
 \qquad d=EX,\qquad e=EZ,
\]

and

\[
 \kappa={j\over a}=Z_X,\qquad k=\kappa(0)={j_0\over a_0},
 \qquad \theta={\gamma\over uj}.
\]

The proved profile equations and monotonicity are

\[
 \kappa_X={(3-\theta)\kappa\over2a},
 \qquad 0\le\theta_0\le\theta\le3,
\]

and the center bound proved above gives

\[
 \theta_0\ge {k\over1+k}.                              \tag{1}
\]

We prove

\[
 \boxed{e\le {5k\over2+\theta_0}d
 \le {5k(1+k)\over2+3k}d.}                            \tag{2}
\]

\paragraph{Pointwise integration of the slope}

Let \(T=a/a_0=1-X/a_0\).  From the differential equation for \(\kappa\),

\[
 {d\over dT}\log\kappa
 \ge-{3-\theta_0\over2T}.
\]

Since \(\kappa(1)=k\), integration from \(T\) to one gives

\[
 \kappa(T)\le kT^{-(3-\theta_0)/2}.                   \tag{3}
\]

Put

\[
 r={1-\theta_0\over2}\in[-1,1/2]
\]

and define continuously at \(r=0\)

\[
 \psi_r(T)={T^{-r}-1\over r},
 \qquad \psi_0(T)=\log(1/T).
\]

Integrating (3) in \(X\) yields

\[
 \boxed{Z(T)\le ka_0\psi_r(T).}                       \tag{4}
\]

The formula includes \(\theta_0>1\): then both the numerator and denominator
in \(\psi_r\) are negative, and \(\psi_r\) remains nonnegative.

\paragraph{The lower-truncated power mixture}

The pushforward law of \(T\) under the fold measure has density

\[
 f_T(t)=C t^{1/2}H(t),\qquad 0<t<1,
\]

where \(H\) is nonnegative and nondecreasing.  By layer cake this law is a
positive mixture of the base laws proportional to \(t^{1/2}\) conditioned
on intervals \([s,1]\).

For \(r \in [-1,1/2]\), set

\[
 g_r(t)={\psi_r(t)\over1-t}.
\]

This function is nonincreasing.  Indeed,

\[
 g_r(t)={1\over1-t}\int_t^1 x^{-r-1}\,dx,
\]

and \(x^{-r-1}\) is nonincreasing because \(r\ge -1\); raising the lower
endpoint lowers its average.  The same representation includes \(r=0\),
and at \(r=-1\) the function is constant.

For a base component truncated to \([s,1]\),

\[
 {E_s\psi_r(T)\over E_s(1-T)}
 ={\int_s^1g_r(t)(1-t)t^{1/2}\,dt
   \over\int_s^1(1-t)t^{1/2}\,dt}
\]

is therefore no larger than its value at \(s=0\).  Since this is a ratio
inequality with a common constant, it survives arbitrary mixtures.

Under the untruncated base density \((3/2)t^{1/2}\),

\[
 E_0\psi_r(T)={2\over2+\theta_0},
 \qquad E_0(1-T)={2\over5}.
\]

Consequently

\[
 E\psi_r(T)\le {5\over2+\theta_0}E(1-T).              \tag{5}
\]

Finally, \(d=a_0E(1-T)\).  Taking expectations in (4), then using (5) and
(1), proves (2), with

\[
 \Phi(k)={5k(1+k)\over2+3k}.
\]

Every endpoint (\(k=0\), \(\theta_0=1\), and \(\theta_0=3\)) follows directly by
continuity; no strict inequality is required.

\subsubsection{Same-center terminal comparison}

We give the ODE and integral comparison proving

\[
 P={\int b^2d\pi\over\int a^2d\pi}\le
 P_T={\int b_T^2d\pi_T\over\int a_T^2d\pi_T}.
\]

At contact, the ratio on the left is the same \(P=p-1\) used in the
main proof.

No contact or fold identity is used in this inequality; only the exact
profile cone and equality of the center data are used.

\paragraph{The same-center ODE comparison}

Let

\[
 k={j_0\over a_0},\qquad L^2=(1+k)a_0=a_0+j_0.
\]

The terminal profile with these center data is

\[
 \begin{aligned}
 u_T&=L\tanh(Ly),&
 b_T&={a_0\over L}\tanh(Ly),\\
 a_T&=a_0\operatorname{sech}^2(Ly),&
 j_T&=ka_T,\\
 \gamma_T&=3ka_Tu_T.
 \end{aligned}
\]

and
\[
 d\pi_T(y)=Z_T^{-1}e^{-y^2/2}\cosh(Ly)^{1/(1+k)}\,dy.
\]
Indeed, \(b_T'=a_T\) and \((\log\pi_T)'=-y+b_T\). Also

\[
 u_T'=(1+k)a_T=L^2-u_T^2,\qquad a_T'=-2a_Tu_T.
\]

The explicit formula for \(P_T\) used in the main proof follows here.
With \(\widehat m=1/(1+k)\), \(h=L^2=a_0(1+k)\), and
\(z=Ly\), normalization constants cancel and
\[
 \begin{aligned}
 P_T
 &= {1\over h}
 {\int_{\mathbb R}\tanh^2z\,e^{-z^2/(2h)}
       \cosh^{\widehat m}z\,dz
  \over
  \int_{\mathbb R}\operatorname{sech}^4z\,e^{-z^2/(2h)}
       \cosh^{\widehat m}z\,dz}\\
 &= {J_{-\widehat m}(h)-J_{2-\widehat m}(h)
  \over hJ_{4-\widehat m}(h)},
 \end{aligned}
\]
where
\[
 J_r(h)=\int_{\mathbb R}e^{-z^2/(2h)}
             \operatorname{sech}^rz\,dz.
\]

For the original profile,

\[
 a'=-2au,\qquad u'=a+j,\qquad j'=uj-\gamma,
 \qquad 0\le\gamma\le3uj.
\]

Consequently

\[
 \left({j\over a}\right)'={3uj-\gamma\over a}\ge0,
\]

and hence \(j\ge ka\).  Therefore

\[
 \{u^2+(1+k)a\}'=2u(j-ka)\ge0,
 \qquad u(0)^2+(1+k)a(0)=L^2.
\]

Thus, as long as \(u<L\),

\[
 u'=a+j\ge(1+k)a\ge L^2-u^2.
\]

Scalar ODE comparison with \(u_T'=L^2-u_T^2\), or the standard first
crossing argument, gives \(u\ge u_T\).  If \(u\) reaches \(L\) at a finite point,
the same conclusion thereafter is automatic because \(u_T<L\) at every
finite point.  It follows that

\[
 \left(\log{a\over a_T}\right)'=-2(u-u_T)\le0.
\]

Since the ratio starts at one,

\[
 \boxed{u\ge u_T,\qquad a\le a_T.}
\]

Finally \(b'=a\), \(b_T'=a_T\), and both vanish at the center, so

\[
 \boxed{b\le b_T.}
\]

In particular

\[
 g={a\over a_T}
\]

satisfies \(0<g\le 1\) and \(g'\le 0\), and

\[
 b(y)=\int_0^y a_T(s)g(s)\,ds.
\]

This proves every assertion used in the integral comparison below.

\paragraph{The terminal integral inequality for decreasing inputs}

Work first on the positive half-line; even/odd extension only contributes
a common factor two.  For a nonnegative decreasing \(g\), set

\[
 b_g(y)=\int_0^y a_T(s)g(s)\,ds.
\]

Layer cake writes \(g\) as a positive mixture of
\(\mathbf 1_{[0,s]}\). Denote the
corresponding primitives by

\[
 b_s(y)=b_T(\min\{y,s\}),
\]

and put

\[
 A(s)=\int_0^s a_T^2d\pi_T,\qquad
 C(s)=\int_0^\infty b_s b_Td\pi_T.
\]

For \(s\le t\),

\[
 \langle b_s,b_t\rangle_{\pi_T}\le C(s),
 \qquad
 \langle a_T1_{[0,s]},a_T1_{[0,t]}\rangle_{\pi_T}=A(s).
\]

It remains to show \(C(s)\le P_TA(s)\).

Define

\[
 H(s)={1\over\pi_T(s)}\int_s^\infty b_T(y)\pi_T(y)\,dy,
 \qquad R(s)={H(s)\over a_T(s)}.
\]

Direct differentiation gives

\[
 {C'(s)\over A'(s)}=R(s).
\]

We next verify that \(R\) is increasing.  Put

\[
 t=y-b_T+2u_T,\qquad r={b_T\over a_T}.
\]

Then the density \(a_T \pi_T\) has score \(t\), and

\[
 t'=1+a_T+2j_T=:\mathcal A,\qquad
 r'=1+2u_Tr=: \mathcal D.
\]

Using the terminal field equations gives the exact determinant

\[
 \mathcal A'\mathcal D-\mathcal A\mathcal D'
 =-2\{u_T(1+2a_T+j_T)+\gamma_T\}\mathcal D
  -2\mathcal A(a_T+j_T)r<0.
\]

Thus \(t'/r'\) is strictly decreasing.  Since \(t(0)=r(0)=0\) and \(r'>0\),
the quotient \(t/r\) is the \(r'\)-weighted initial average of \(t'/r'\) and
is strictly decreasing.  Hence \(f=r/t\) is increasing.

Now

\[
 R(s)={\int_s^\infty r a_T\pi_T\over a_T(s)\pi_T(s)}.
\]

Since \((a_T \pi_T)'=-t a_T \pi_T\), integration by parts gives

\[
 \int_s^\infty r a_T\pi_T
 =f(s)a_T(s)\pi_T(s)+\int_s^\infty f'a_T\pi_T.
\]

The Gaussian tail makes the endpoint term at infinity vanish.  Therefore

\[
 R\ge {r\over t},\qquad R'=tR-r\ge0.
\]

It follows that \(C(s)/A(s)\) is an initial \(dA\)-average of the increasing
function \(R\); hence it is increasing and is bounded by its limit

\[
 {C(\infty)\over A(\infty)}=P_T.
\]

Thus \(C(s)\le P_TA(s)\).  Applying this pairwise to the two layer-cake
integrals proves

\[
 \boxed{\int b_g^2d\pi_T\le P_T\int a_T^2g^2d\pi_T.}
\]

General monotone \(g\) follow by Stieltjes approximation; all integrands are
nonnegative and Gaussian dominated.

\paragraph{Transfer from \(\pi_T\) to \(\pi\)}

For the actual same-center profile, \(a=a_Tg\) and \(b=b_g\), so the preceding
inequality says

\[
 \int (b^2-P_Ta^2)d\pi_T\le0.
\]

Let

\[
 F=b^2-P_Ta^2.
\]

On the positive half-line,

\[
 F'=2ab+4P_Ta^2u>0.
\]

Also

\[
 {d\pi\over d\pi_T}\propto
 \exp\left\{-\int_0^{|y|}(b_T-b)\right\}
\]

is nonincreasing in \(|y|\), because \(b\le b_T\). The elementary
opposite-monotonicity covariance identity under \(\pi_T\) therefore gives

\[
 \int Fd\pi\le\int Fd\pi_T\le0.
\]

Dividing by the positive integral of \(a^2\) proves

\[
 \boxed{P\le P_T.}
\]

This proves \(P\le P_T\).

\subsubsection{Contact-size estimate}

\paragraph{Statement and notation}

Work on the full line.  Let

\[
 d\pi=Z_\pi^{-1}e^{-V(y)}\,dy,\qquad V'=y-b,
\]

where \(b\) is odd, \(a=b'>0\) is even, and the standardized logarithmic-heat
fields satisfy

\[
 a'=-2au,\quad u'=a+j,\quad j'=uj-\gamma,
 \quad j,\gamma\ge0.                                      \tag{1}
\]

Put

\[
 n=y+u-b,\qquad t=n+u=y-b+2u,\qquad r={b\over a}.
\]

The \(a^2\)-biased law and its contact size are

\[
 d\nu={a^2\,d\pi\over A_2},\qquad
 A_2=\int a^2\,d\pi,\qquad x=\mathbb E_\nu a={\int a^3d\pi\over A_2}.
                                                               \tag{2}
\]

At contact,

\[
 B_2:=\int b^2d\pi=P A_2,\qquad P>1.                       \tag{3}
\]

This is the same \(P=p-1\) as in the main proof, by the contact identity
(2.3).

The result proved below is the sharp estimate

\[
 \boxed{x\ge 1-{1\over P}=: \alpha.}                    \tag{4}
\]

For a nondegenerate logarithmic-heat profile the inequality is strict.  Notice
that the fold identity is not needed for (4).  At a shifted fold it turns
(4) into

\[
 x={\alpha\over2}+2\mathbb E_\nu u^2
 \quad\Longrightarrow\quad
 \mathbb E_\nu u^2\ge {\alpha\over4}.                    \tag{5}
\]

All functions below have the standard finite-horizon Gaussian tails:
\(\pi\) has Gaussian tails, \(b\) is bounded, and its derivatives have at most
polynomial/exponential growth dominated by those tails.  These facts imply
all displayed functions are in the indicated weighted Sobolev domains.
They also justify the boundary limits used below.  Equivalently, every
integration by parts can be performed first with symmetric compactly
supported cutoffs; the Gaussian tail bounds let the cutoff errors tend to
zero.

\paragraph{An odd auxiliary equation}

Let

\[
 L=-\partial_y^2+V'\partial_y
   =-{1\over\pi}(\pi\partial_y)'
\]

be the nonnegative reversible generator in \(L^2(\pi)\).  Define

\[
 H(y)={1\over\pi(y)}\int_y^\infty b(s)\pi(s)\,ds,
 \qquad h(y)=\int_0^yH(s)\,ds.                            \tag{6}
\]

Because \(b \pi\) is odd and has integral zero, \(H\) is even; it is positive
on the positive half-line, and \(h\) is odd.  Moreover

\[
 (\pi H)'=-b\pi,\qquad Lh=b.                              \tag{7}
\]

The sign in (6) is important: integrating \((\pi H)'=-b \pi\) from \(y\) to
infinity, with \(\pi H\) vanishing there, gives the \textbf{positive} tail in (6).
Gaussian tail estimates give \(H,H' \in L^2(\pi)\) and \(h \in L^2(\pi)\); hence
\(h\) belongs to the operator domain of the odd inverse of $L$.  Thus (6)
is also a direct construction, not merely a formal use of \(L^{-1}\).

The boundary quantities needed later vanish.  In particular,

\[
 \pi H=\int_y^\infty b\pi\longrightarrow0,\qquad
 b\pi H\longrightarrow0                                  \tag{8}
\]

at both ends.  The corresponding terms in the twice-integrated identity vanish by
the same tail estimate (or by the cutoff argument just described).

\paragraph{A monotone likelihood ratio for \(H/a\)}

Set

\[
 R={H\over a}.
\]

Since \(b=ar\), formula (6) becomes

\[
 R(y)={\displaystyle\int_y^\infty r(s)a(s)\pi(s)\,ds
             \over a(y)\pi(y)}.                           \tag{9}
\]

The density \(a \pi\) has score \(t\), because

\[
 -{(a\pi)'\over a\pi}=y-b+2u=t.                           \tag{10}
\]

Consequently direct differentiation of (9) gives

\[
 R'=tR-r.                                                  \tag{11}
\]

We next prove that \(R\) is increasing.  First, \(\rho=t/r\) is strictly
decreasing on the positive half-line.  Indeed, write

\[
 \mathsf A=t'=1+a+2j,\qquad \mathsf B=r'=1+2ur.
\]

Using (1),

\[
 \mathsf A'=-2u(a-j)-2\gamma,\qquad
 \mathsf B'=2(a+j)r+2u\mathsf B,
\]

and hence

\[
 \mathsf A'\mathsf B-\mathsf A\mathsf B'
 =-2\{u(1+2a+j)+\gamma\}\mathsf B
  -2\mathsf A(a+j)r<0.                                   \tag{12}
\]

Thus \(t'/r'\) strictly decreases.  Since \(t(0)=r(0)=0\) and \(r'>0\),

\[
 {t(y)\over r(y)}
 ={\int_0^y(t'/r')r'\,ds\over\int_0^yr'\,ds}
\]

is the \(r'\)-weighted average of a strictly decreasing function.  Therefore

\[
 \left({t\over r}\right)'
 ={r'\over r}\left({t'\over r'}-{t\over r}\right)<0.    \tag{13}
\]

In particular \(f=r/t\) is strictly increasing.  Integrating by parts in
(9), using \((a \pi)'=-t a \pi\), gives the exact identity

\[
\begin{aligned}
 \int_y^\infty r a\pi
 &=\int_y^\infty {r\over t}\,t a\pi
  =-\int_y^\infty f(a\pi)'  \\
 &=f(y)a(y)\pi(y)+\int_y^\infty f'(s)a(s)\pi(s)\,ds.
\end{aligned}                                                   \tag{14}
\]

The omitted boundary term is
\(f a \pi=(r/t)a \pi=b \pi/t\), which tends to zero by the Gaussian tails.
Equations (9) and (14) imply

\[
 R\ge {r\over t},\qquad R'=tR-r\ge0.                      \tag{15}
\]

Both inequalities are strict at every finite \(y>0\), since \(f'>0\) and the
tail in (14) has positive mass.  By evenness, \(R\) is therefore strictly
increasing as a function of \(|y|\).

\paragraph{Covariance and integration by parts}

On the positive half-line, \(a'=-2au<0\), while (15) says \(R^2\) increases.
The full law \(\nu\) is even.  The elementary opposite-monotonicity identity,
applied to \(|Y|\), yields

\[
 \operatorname{Cov}_\nu(a,R^2)\le0.
\]

For clarity, with an independent copy \(Y'\) this is simply

\[
 2\operatorname{Cov}_\nu(a,R^2)
 =\mathbb E[(a(|Y|)-a(|Y'|))(R^2(|Y|)-R^2(|Y'|))]\le0.
\]

Since \(H=aR\), it follows that

\[
 {\int aH^2d\pi\over\int H^2d\pi}
 ={\mathbb E_\nu[aR^2]\over\mathbb E_\nu[R^2]}
 \le \mathbb E_\nu a=x.                                 \tag{16}
\]

Two integrations by parts and \(V''=1-a\) give

\[
\begin{aligned}
 B_2=\int(Lh)^2d\pi
 &=\int (h'')^2d\pi+\int V''(h')^2d\pi \\
 &=\int(H')^2d\pi+\int(1-a)H^2d\pi
 \ge(1-x)\int H^2d\pi.
\end{aligned}                                           \tag{17}
\]

Finally, reversibility, (7), and (8) show

\[
 B_2=\langle b,Lh\rangle_\pi
   =\langle b',h'\rangle_\pi=\int aH\,d\pi.
\]

The elementary product inequality therefore yields

\[
 B_2^2\le A_2\int H^2d\pi.                               \tag{18}
\]

If \(x\ge 1\), (4) is immediate.  If \(x<1\), combining (17)--(18) and using
\(B_2=PA_2\) gives

\[
 B_2\ge(1-x)\int H^2d\pi
   \ge(1-x){B_2^2\over A_2}
 \quad\Longrightarrow\quad
 1\ge P(1-x).
\]

This proves (4).

\paragraph{Strictness and equality}

For a nondegenerate logarithmic-heat profile, \(u(y)>0\) for \(y>0\) because
\(u'=a+j>0\); hence \(a\) strictly decreases.  Equations (12)--(15) make \(R\)
strictly increase.  Thus (16) is strict, and so is (4): \(x>\alpha\).

If the differential-cone strictness is discarded, equality in the complete
chain forces equality in the discarded \(H'\) term and in the product
inequality.  Hence \(H\) and \(a\) are constant, \(b(y)=ay\), and \(\pi\) is a
centered Gaussian of precision \(1-a\).  Contact then gives
\(P=1/(1-a)\) and \(a=x=\alpha\).  This affine Gaussian is the sharp equality
case of the abstract inequality, but it is not a nondegenerate profile
satisfying \(u'=a+j>0\).

\subsubsection{Center-density and one-dimensional variance bounds}

Write

\[
 \theta={\gamma\over uj},\qquad a_0=a(0),\quad j_0=j(0),
 \quad x=E a,\quad \bar x={x\over a_0}.
\]

The profile identities are

\[
 a'=-2au,\qquad j'=uj-\gamma,
\]

and therefore

\[
 (j\sqrt a)'=-\gamma\sqrt a=-\theta u j\sqrt a.       \tag{1}
\]

The proved monotonicity of \(\theta\) and the cone bound give

\[
 \theta_0\le\theta(y)\le3,
 \qquad \theta_0=\lim_{y\downarrow0}\theta(y).
\]

Integrating (1), and using
\(\int_0^y u=(1/2)\log(a_0/a(y))\), yields the two pointwise bounds

\[
 \boxed{
 j_0\,{a\over a_0}\le j(y)\le
 j_0\left({a\over a_0}\right)^{(\theta_0-1)/2}.}
\tag{2}
\]

The upper inequality is especially useful because it turns an averaged
lower bound on \(Ej\) into a lower bound on \(j_0\).

\paragraph{A truncated-Beta chord lemma}

The pushforward of the fold law by \(T=a/a_0\) has density

\[
 d\lambda(t)=\mathcal Z^{-1}t^{1/2}H(t)dt,\qquad 0<t<1,
 \quad H\ \hbox{nondecreasing}.                         \tag{3}
\]

Such laws are mixtures of the truncated Beta laws

\[
 d\lambda_s(t)={t^{1/2}1_{s<t<1}dt\over
                   \int_s^1t^{1/2}dt},\qquad 0\le s<1. \tag{4}
\]

Let \(m(s)=E_s T\) and \(F(s)=E_s f(T)\). Differentiating a
lower-truncated expectation gives, for a positive scalar \(\ell(s)\),

\[
 m'=\ell(m-s),\qquad F'=\ell(F-f(s)).                   \tag{5}
\]

Thus, regarding \(F=\Phi(m)\), its slope is

\[
 \Phi'(m(s))={F(s)-f(s)\over m(s)-s}=:R(s),
\]

and a second differentiation gives the exact sign formula

\[
 R'(s)={R(s)-f'(s)\over m(s)-s}.                        \tag{6}
\]

If \(f\) is convex, its tangent inequality gives \(R\ge f'\), so
\(\Phi\) is convex. Hence the curve \((m(s),F(s))\) lies below its
endpoint chord. If \(f\) is concave, the inequalities reverse and the curve lies above its
endpoint chord.  Averaging the corresponding chord inequality over the
mixture (4) preserves it, because the chord is affine in the mean.

The two endpoint means are

\[
 m(0)={3\over5},\qquad m(1)=1.                          \tag{7}
\]

\paragraph{Upper power moment and lower bound on \(j_0\)}

Assume first \(0\le\theta_0\le1\) and put

\[
 r={\theta_0-1\over2}\in[-1/2,0].
\]

The function \(t^r\) is convex. Its two endpoint expectations are

\[
 E_0T^r={3\over\theta_0+2},\qquad E_1T^r=1.
\]

The upper chord from the preceding lemma gives

\[
 \boxed{
 E T^{(\theta_0-1)/2}\le
 H_\theta:=1+{5(1-\bar x)(1-\theta_0)\over
                    2(\theta_0+2)}.}                   \tag{8}
\]

If \(\theta_0\ge1\), the simpler bound \(E T^r\le1\) holds.
Consequently (2) implies, with \(J=Ej\),

\[
 \boxed{j_0\ge {J\over H_\theta}}                      \tag{9}
\]

where \(H_\theta=1\) for \(\theta_0\ge1\).

\paragraph{Sharpening the remainder estimate}

We first prove the one-dimensional variance estimate used here. If
\(d\nu=\mathcal Z^{-1}e^{-V(y)}dy\), \(V''>0\), and \(f\) is smooth with mean
zero, then
\[
 \int f^2d\nu\le\int {f'^2\over V''}\,d\nu.
\]
Indeed, solve \(-h''+V'h'=f\), first with compactly supported
approximations. Integration by parts, the weighted product inequality, and one
more integration by parts give
\[
 \begin{aligned}
 \int f^2d\nu
 &=\int f'h'\,d\nu\\
 &\le\left(\int {f'^2\over V''}\,d\nu\right)^{1/2}
       \left(\int V''h'^2d\nu\right)^{1/2},\\
 \int f^2d\nu
 &=\int h''^2d\nu+\int V''h'^2d\nu.
 \end{aligned}
\]
Dividing proves the estimate. The Gaussian tails justify the
approximation.

For the fold law, the score is \(-(n+3u)\), so
\[
 V''=n'+3u'=1+3a+4j.
\]
The function \(u\) is odd and has mean zero, while \(u'=a+j\).
Applying the displayed estimate to \(u\) gives
\[
 U=Eu^2\le E{(a+j)^2\over1+3a+4j}.
\]
The pointwise identity
\[
 {(a+j)^2\over1+3a+4j}
 ={1\over12}\left(4a+3j-{4a+3j+aj\over1+3a+4j}\right)
\]
and the fold identity \(Ea=c+2U\) therefore give

\[
 J\ge {4\over3}(U-c)+{1\over3}\mathcal R,
 \qquad U=Eu^2={x-c\over2},                             \tag{10}
\]

\[
 \mathcal R=E{4a+3j+aj\over1+3a+4j}.                   \tag{11}
\]

The integrand in (11) is increasing in \(j\), since its \(j\) derivative is
\(3(a-1)^2/(1+3a+4j)^2\).  The lower inequality in (2) therefore gives a
model-specific lower bound retaining \(j_0\).  A simpler \(j_0\)-free version
is

\[
 \mathcal R\ge E f(a),\qquad f(a)={4a\over1+3a}.        \tag{12}
\]

The function \(f(a_0t)\) is increasing and concave in \(t\), and vanishes at
zero.  Its expectation under the untruncated Beta law is

\[
 I(a_0)={3\over2}\int_0^1t^{1/2}{4a_0t\over1+3a_0t}dt
 ={4\over3}-{4\over3a_0}
 +{4\arctan\sqrt{3a_0}\over3a_0\sqrt{3a_0}}.          \tag{13}
\]

Its endpoint value is \(f(a_0)=4a_0/(1+3a_0)\). The lower chord from
the preceding lemma proves

\[
 \boxed{
 \mathcal R\ge L(a_0,\bar x):=
 I(a_0)+{5\bar x-3\over2}
 \left\{{4a_0\over1+3a_0}-I(a_0)\right\}.}            \tag{14}
\]

Combining (9), (10), and (14) gives the explicit center estimate

\[
 \boxed{
 j_0\ge {\displaystyle {4\over3}(U-c)
                      +{1\over3}L(a_0,\bar x)
                \over H_\theta}.}                      \tag{15}
\]

This uses only the proved profile monotonicity and the exact density
class.

\subsection{The analytic upper bounds used by the terminal verifiers}
\label{ft:app:terminal-upper-bound}

This section proves that the elementary quantities evaluated by the
terminal-verification programs really are upper bounds for the terminal
ratio.  Put
\[
 0<m\leq1,\qquad h>0,
\]
and define
\[
 J_r(h)=\int_{\mathbb R}e^{-z^2/(2h)}\operatorname{sech}^r z\,dz,
 \qquad
 \mathcal N(m,h)=J_{-m}(h)-J_{2-m}(h).
\]
Thus, when \(m=(1+k)^{-1}\) and \(h=A/m=A(1+k)\),
\[
 P_T(A,k)={\mathcal N(m,h)\over hJ_{4-m}(h)}.          \tag{T.1}
\]

\subsubsection{A continued-fraction enclosure of the shifted Gaussian tail}

For \(q,h>0\), set
\[
 K_q(h)=\int_0^\infty e^{-z^2/(2h)-qz}\,dz.           \tag{T.2}
\]
For an integer \(L\geq2\), define recursively
\[
 w_L=q,\qquad
 w_j=q+{j\over h w_{j+1}}\quad(j=L-1,\ldots,1),
 \qquad C_L(q,h)={1\over w_1}.                        \tag{T.3}
\]
Then, for every \(N\geq1\),
\[
 \boxed{C_{2N}(q,h)\leq K_q(h)\leq C_{2N+1}(q,h).}   \tag{T.4}
\]

Here is an elementary proof, including the direction of the inequalities.
Write \(x=q\sqrt h\) and
\[
 I_j(x)=\int_0^\infty y^j e^{-y^2/2-xy}\,dy.
\]
Integration by parts gives
\[
 1=xI_0+I_1,
 \qquad I_{j+1}=jI_{j-1}-xI_j\quad(j\geq1).           \tag{T.5}
\]
If \(r_j=I_j/I_{j-1}>0\) and \(D_j=x+r_j\), then
\[
 I_0={1\over D_1},\qquad D_j=x+{j\over D_{j+1}}.      \tag{T.6}
\]
The truncation in (T.3), after multiplication by \(\sqrt h\), replaces
the exact terminal value \(D_L=x+r_L>x\) by \(x\).  Each map
\(y\mapsto x+j/y\) reverses order.  Propagating from \(L\) to \(1\), and
then taking the reciprocal, therefore gives a lower bound when \(L\) is
even and an upper bound when \(L\) is odd.  Since
\(K_q(h)=\sqrt h\,I_0(q\sqrt h)\), this proves (T.4).

We shall also use
\[
 \int_0^\infty e^{-z^2/(2h)+mz}\,dz
 =\sqrt{2\pi h}\,e^{m^2h/2}-K_m(h),                  \tag{T.7}
\]
which follows by subtracting the negative half-line from the completed
square integral over \(\mathbb R\).

\subsubsection{The degree-six interpolation majorant}

On \(0\leq s\leq1\), let
\[
 f_m(s)=(1+s)^{m-2}.
\]
Let \(H_{6,m}\) be the polynomial of degree at most six determined by
\[
 H_{6,m}(0)=f_m(0),\qquad
 H_{6,m}^{(j)}(1)=f_m^{(j)}(1),\quad 0\leq j\leq5.    \tag{T.8}
\]
Subtract a suitable multiple of \(s(s-1)^6\) so that the difference
also vanishes at the chosen point.  Differentiating between consecutive
zeros seven times gives, for some
\(\eta=\eta(s)\in(0,1)\),
\[
 f_m(s)-H_{6,m}(s)
 ={f_m^{(7)}(\eta)\over7!}\,s(s-1)^6.                \tag{T.9}
\]
For \(0\leq m\leq1\),
\[
 f_m^{(7)}(\eta)
 =(m-2)(m-3)\cdots(m-8)(1+\eta)^{m-9}<0.
\]
Consequently
\[
 \boxed{f_m(s)\leq H_{6,m}(s)\quad(0\leq s\leq1).}  \tag{T.10}
\]

Write
\[
 (1-s)^2H_{6,m}(s)=\sum_{i=0}^{8}d_i(m)s^i,
 \qquad d_0(m)=1.                                    \tag{T.11}
\]
For \(z\geq0\), putting \(s=e^{-2z}\) gives
\[
 \cosh^m z\,\tanh^2z
 =2^{-m}e^{mz}(1-s)^2f_m(s).                          \tag{T.12}
\]
Equations (T.7), (T.10)--(T.12), and evenness imply
\[
 \begin{aligned}
 \mathcal N(m,h)\leq 2^{1-m}\bigg\{&
 \sqrt{2\pi h}\,e^{m^2h/2}-K_m(h)\\
 &+\sum_{i=1}^{8}d_i(m)K_{2i-m}(h)\bigg\}.
 \end{aligned}                                        \tag{T.13}
\]

The regional verifier evaluates (T.13) on parameter boxes.  It encloses
every \(d_i(m)\) by outward-rounded dyadic interval arithmetic and refuses
a box if that interval straddles both signs.  On a box where \(d_i\geq0\), it
uses the upper tail bound \(C_{29}(2i-m,h)\); on a box where \(d_i\leq0\),
it uses the lower tail bound \(C_{28}(2i-m,h)\), because multiplication by
a nonpositive coefficient reverses the inequality.  In the first term it
uses \(K_m\geq C_{28}(m,h)\), together with upper bounds for \(\pi\) and
the exponential.  Pointwise, the resulting upper bound is
\[
\begin{aligned}
 N_6^+(m,h)=2^{1-m}\bigg[&
 \sqrt{2\pi_+h}\,e^{m^2h/2}-C_{28}(m,h)\\
 &+\sum_{\substack{1\leq i\leq8\\ d_i(m)\geq0}}
       d_i(m)C_{29}(2i-m,h)\\
 &+\sum_{\substack{1\leq i\leq8\\ d_i(m)<0}}
       d_i(m)C_{28}(2i-m,h)\bigg].
\end{aligned}                                        \tag{T.13a}
\]
On boxes, every factor in this formula is replaced by its outward interval
enclosure with the same sign rule.  Thus the quantity called
\texttt{terminal\_bound} in the regional source is an upper enclosure
\(N_6^+(m,h)\) of \(\mathcal N(m,h)\) on every accepted box.

\subsubsection{The degree-ten positive-series majorant}

There is a second global majorant whose coefficient signs are known in
advance.  Put \(t=1-s\).  Since the power series has positive coefficients,
\[
 (2-t)^{m-2}=\sum_{j=0}^{\infty}c_j(m)t^j,
 \qquad
 c_j(m)={2^{m-2-j}(2-m)_j\over j!}>0,                 \tag{T.14}
\]
where \((a)_j=a(a+1)\cdots(a+j-1)\).
At \(t=1\) the sum is one.  Hence
\[
 R_{10}(m):=1-\sum_{j=0}^{9}c_j(m)
            =\sum_{j=10}^{\infty}c_j(m)>0.
\]
Because \(t^j\leq t^{10}\) for \(0\leq t\leq1\) and \(j\geq10\),
\[
 (1+s)^{m-2}\leq
 \sum_{j=0}^{9}c_j(m)(1-s)^j+R_{10}(m)(1-s)^{10}.     \tag{T.15}
\]
After multiplication by \((1-s)^2\), the right-hand side is
\[
 \sum_{i=0}^{12}(-1)^i b_i(m)s^i,                    \tag{T.16}
\]
where
\[
 b_i(m)=\sum_{j=0}^{9}c_j(m)\binom{j+2}{i}
       +R_{10}(m)\binom{12}{i}>0,\qquad b_0(m)=1.   \tag{T.17}
\]
It follows exactly as in (T.13) that
\[
 \begin{aligned}
 \mathcal N(m,h)\leq N_{10}^+(m,h):=2^{1-m}\bigg[&
 \sqrt{2\pi_+h}\,e^{m^2h/2}-C_{20}(m,h)\\
 &+\sum_{\substack{2\leq i\leq12\\ i\ {\rm even}}}
 b_i(m)C_{21}(2i-m,h)\\
 &-\sum_{\substack{1\leq i\leq12\\ i\ {\rm odd}}}
 b_i(m)C_{20}(2i-m,h)\bigg].
 \end{aligned}                                        \tag{T.18}
\]
The positive terms use the odd upper convergent and the negative terms
use the even lower convergent.  Here \(\pi_+\) denotes the program's
outward rational upper bound for \(\pi\); exponentials and the remaining
operations are also rounded outward.

\subsubsection{A Fourier-product lower bound for the denominator}

We use the Fourier-transform convention
\[
 \widehat f(y)=\int_{\mathbb R}e^{-iyz}f(z)\,dz.
\]
For \(r>0\), the substitution \(t=e^{2z}\) in the beta integral gives
\[
 \widehat{\operatorname{sech}^r}(y)
 ={2^{r-1}\over\Gamma(r)}
 \Gamma\!\left({r-iy\over2}\right)
 \Gamma\!\left({r+iy\over2}\right).                \tag{T.19}
\]
In particular,
\[
 I_r:=\int_{\mathbb R}\operatorname{sech}^r z\,dz
 =\sqrt\pi\,{\Gamma(r/2)\over\Gamma((r+1)/2)}.      \tag{T.20}
\]
The infinite product for the gamma function, applied to the quotient of
(T.19) and its value at zero, gives
\[
 {\widehat{\operatorname{sech}^r}(y)\over I_r}
 =\prod_{j=0}^{\infty}
   \left(1+{y^2\over(r+2j)^2}\right)^{-1}.           \tag{T.21}
\]
Indeed, in the product for
\(\Gamma(a+it)\Gamma(a-it)/\Gamma(a)^2\), the exponential factors
cancel and the \(j\)-th remaining factor is
\((1+t^2/(a+j)^2)^{-1}\), with \(a=r/2,t=y/2\).
This follows directly from
\[
 {1\over\Gamma(z)}
 =ze^{\gamma z}\prod_{j=1}^{\infty}
   \left(1+{z\over j}\right)e^{-z/j};
\]
the leading factor supplies the \(j=0\) term in (T.21).

Since \(\log(1+x)\leq x\), (T.21) implies
\[
 {\widehat{\operatorname{sech}^r}(y)\over I_r}
 \geq e^{-\sigma_r^2y^2/2},
 \qquad
 \sigma_r^2=2\sum_{j=0}^{\infty}{1\over(r+2j)^2}
 .                                                     \tag{T.22}
\]
The Fourier transform of \(e^{-z^2/(2h)}\) is
\(\sqrt{2\pi h}\,e^{-hy^2/2}\).  Fourier inversion and (T.22) therefore
give, with every constant displayed,
\[
 \begin{aligned}
 J_r(h)
 &=I_r\sqrt{h\over2\pi}
   \int_{\mathbb R}e^{-hy^2/2}
   {\widehat{\operatorname{sech}^r}(y)\over I_r}\,dy\\
 &\geq I_r\sqrt{h\over h+\sigma_r^2}.
 \end{aligned}                                        \tag{T.23}
\]

For the required range \(r=4-m\),
\[
 (\log I_r)''
 ={1\over4}\sum_{j=0}^{\infty}
 \left\{{1\over(r/2+j)^2}-{1\over((r+1)/2+j)^2}\right\}>0. \tag{T.23a}
\]
Thus \(\log I_r\) is convex.  The substitution \(t=\tanh z\) gives
\(I_3=\int_{-1}^1\sqrt{1-t^2}\,dt=\pi/2\).  Differentiating (T.20)
by means of the displayed product for \(1/\Gamma\), and separating even
and odd terms in the harmonic series, gives
\[
 I_3={\pi\over2},\qquad
 (\log I_r)'\big|_{r=3}
 ={1\over2}\sum_{j=0}^{\infty}
 \left({1\over j+2}-{1\over j+3/2}\right)
 ={1\over2}-\log2.
\]
The tangent inequality followed by \(e^{-x}\geq1-x\) gives
\[
 I_{4-m}\geq {\pi_-\over2}
 \left[1-(1-m){48287\over250000}\right].             \tag{T.24}
\]
Furthermore
\[
 \sigma_{4-m}^2=2\sum_{j=0}^{\infty}(4-m+2j)^{-2}
\]
is convex in \(m\), and hence lies below its endpoint chord.  Since
\[
 \sigma_4^2={\pi^2\over12}-{1\over2},\qquad
 \sigma_3^2-\sigma_4^2={\pi^2\over6}-{3\over2},
\]
we obtain
\[
 \sigma_{4-m}^2\leq {64494\over200000}
            +m{28987\over200000}.                    \tag{T.25}
\]
The rational constants are certified from
\[
 \log2=2\sum_{j=0}^{\infty}{1\over(2j+1)3^{2j+1}},
 \qquad
 \pi=16\arctan(1/5)-4\arctan(1/239),                \tag{T.26}
\]
using a geometric tail in the first series and alternating-series
remainders in the second.  In particular, the exact comparisons made in
the verifier imply
\[
 \log2-{1\over2}<{48287\over250000},\quad
 {\pi^2\over12}-{1\over2}<{64494\over200000},\quad
 {\pi^2\over6}-{3\over2}<{28987\over200000},          \tag{T.27}
\]
and \(\pi_-<\pi<\pi_+\).  For example, the coarser rational choices
\(\pi_-=103993/33102\) and \(\pi_+=355/113\) already suffice.

Combining (T.23)--(T.25) yields
\[
 hJ_{4-m}(h)\geq D^-(m,h),                            \tag{T.28}
\]
where
\[
 D^-(m,h)=h{\pi_-\over2}
 \left[1-(1-m){48287\over250000}\right]
 \sqrt{h\over h+64494/200000+m(28987/200000)}.        \tag{T.29}
\]

\subsubsection{The two certified upper bounds}

Equations (T.1), (T.13), and (T.28) prove
\[
 \boxed{P_T(A,k)\leq V_6(m,A):=
 {N_6^+(m,A/m)\over D^-(m,A/m)}}.                    \tag{T.30}
\]
This is the bound used by
\texttt{certify\_terminal\_affine\_dyadic.cpp}; that program uses the
degree-six interpolation polynomial and the \(C_{28},C_{29}\) tail bounds.

Equations (T.1), (T.18), and (T.28) prove
\[
 \boxed{P_T(A,k)\leq V_{10}(m,A):=
 {N_{10}^+(m,A/m)\over D^-(m,A/m)}}.                 \tag{T.31}
\]
This is the bound called \texttt{V} in
\texttt{certify\_terminal\_affine\_fixed.py}; it uses the degree-ten
positive-series majorant and \(C_{20},C_{21}\).
The program \texttt{certify\_terminal\_A4\_fixed.py} imports that same
routine, and therefore also uses \(V_{10}\).  Thus every terminal sign
accepted by either verifier is a sign for a rigorous upper bound on
the ratio \(P_T\), not merely for an unrelated auxiliary function.

\section{The marginal-crossing computation}
\label{app:marginal}

This appendix proves Proposition~\ref{prop:transversality}.  The
stationary equations are first rewritten in the endpoint coordinates
used in the differentiation; all of the derivative and strict-sign
computations are then proved in full.

The two 1-RSB stationary equations, written without introducing names
for their left-hand sides, are
\begin{align}
 0={}&-\frac1{m^2}
 \log\E[\cosh^m(\sqrt{\xi'(q)}z)]
 +\frac1m\E[\log\cosh(M(q))e^{W_\mu(q)}]\notag\\
 &-\frac{p-1}{2p}\xi'(q)\Gamma_\mu(q),
 \label{mt:eq:mass}\\
 q={}&\Gamma_\mu(q),
 \label{mt:eq:fixed}
\end{align}
where $z$ is standard Gaussian.  These are equation (15) in
\cite{zhou}, after using
$q\xi'(q)-\xi(q)=\frac{p-1}{p}q\xi'(q)$.

\paragraph{Coordinate convention.}
From now on, $m$ and the numerical value $\xi'(q)$ are the independent
local coordinates.  In every displayed partial derivative, the other
coordinate is held fixed.  Equation~\eqref{mt:eq:fixed} recovers $q$; then
\begin{equation}
 \beta^2=\frac{\xi'(q)}{p\,\Gamma_\mu(q)^{p-1}}.
 \label{mt:eq:recover-lambda}
\end{equation}
Thus a partial derivative in $\xi'(q)$ is not a partial derivative in
$q$.
\subsection{Key computational lemmas}

\begin{lemma}[Endpoint reduction and contact identities]
\label{mt:lem:endpoint}
For the measure in \eqref{eq:mu},
\begin{equation}
 \partial_x\Phi_\mu(x,q)=\tanh x,\qquad
 \partial_{xx}\Phi_\mu(x,q)=\cosh^{-2}x.
 \label{mt:eq:endpoint-phi}
\end{equation}
For every smooth test function $\varphi$ such that
$\varphi,\varphi',\varphi''$ have at most exponential growth and which
does not depend on $m$ or $\xi'(q)$,
\begin{equation}
 \E[\varphi(M(q))e^{W_\mu(q)}]
 =\frac{\E\!\left[
 \varphi(\sqrt{\xi'(q)}z)\cosh^m(\sqrt{\xi'(q)}z)\right]}
 {\E[\cosh^m(\sqrt{\xi'(q)}z)]},
 \qquad \E e^{W_\mu(q)}=1.
 \label{mt:eq:endpoint-law}
\end{equation}
At a stationary marginal point,
\begin{align}
 \Gamma_\mu(q)
 &=\E[\tanh^2(M(q))e^{W_\mu(q)}]=q,
 \label{mt:eq:endpoint-contact}\\
 q&=(p-1)\xi'(q)
 \E[\cosh^{-4}(M(q))e^{W_\mu(q)}],
 \label{mt:eq:endpoint-marginal}\\
 \E[\cosh^{-2}(M(q))e^{W_\mu(q)}]&=1-q,
 \label{mt:eq:contact-second}\\
 \E[\cosh^{-4}(M(q))e^{W_\mu(q)}]
 &=\frac{q}{(p-1)\xi'(q)}.
 \label{mt:eq:contact-fourth}
\end{align}
\end{lemma}

\begin{lemma}[Differentiation of the endpoint law]
\label{mt:lem:differentiation}
For every test function $\varphi$ as in
Lemma~\ref{mt:lem:endpoint},
\begin{align}
 \frac{\partial}{\partial m}
 \E[\varphi(M(q))e^{W_\mu(q)}]
 ={}&\E[\varphi(M(q))\log\cosh(M(q))e^{W_\mu(q)}]\notag\\
 &-\E[\varphi(M(q))e^{W_\mu(q)}]
   \E[\log\cosh(M(q))e^{W_\mu(q)}],
 \label{mt:eq:diff-m}
\end{align}
and
\begin{align}
 &\frac{\partial}{\partial\xi'(q)}
 \E[\varphi(M(q))e^{W_\mu(q)}]\notag\\
 &\quad=\frac12\E\left[
 \{\varphi''(M(q))+2m\tanh(M(q))\varphi'(M(q))\}
 e^{W_\mu(q)}\right]\notag\\
 &\qquad+\frac{m(1-m)}2\Big(
 \E[\varphi(M(q))\cosh^{-2}(M(q))e^{W_\mu(q)}]\notag\\
 &\hspace{46mm}
 -\E[\varphi(M(q))e^{W_\mu(q)}]
  \E[\cosh^{-2}(M(q))e^{W_\mu(q)}]\Big).
 \label{mt:eq:diff-xi}
\end{align}
\end{lemma}

\begin{lemma}[Moment derivatives]
\label{mt:lem:moments}
At a stationary marginal point,
\begin{align}
 \frac{\partial}{\partial m}\Gamma_\mu(q)
 ={}&\E[\cosh^{-2}(M(q))e^{W_\mu(q)}]
     \E[\log\cosh(M(q))e^{W_\mu(q)}]\notag\\
 &-\E[\cosh^{-2}(M(q))\log\cosh(M(q))e^{W_\mu(q)}]
 >0,
 \label{mt:eq:moment-one}\\
 \frac{\partial}{\partial\xi'(q)}\Gamma_\mu(q)
 ={}&\E[\cosh^{-4}(M(q))e^{W_\mu(q)}]\notag\\
 &-(1-m)\Big[2\{\E[\cosh^{-2}(M(q))e^{W_\mu(q)}]
              -\E[\cosh^{-4}(M(q))e^{W_\mu(q)}]\}\notag\\
 &\hspace{19mm}+\frac m2\{\E[\cosh^{-4}(M(q))e^{W_\mu(q)}]
              -\E[\cosh^{-2}(M(q))e^{W_\mu(q)}]^2\}\Big],
 \label{mt:eq:moment-two}\\
 &\frac{\partial}{\partial m}
 \E[\cosh^{-4}(M(q))e^{W_\mu(q)}]\notag\\
 &\quad=\E[\cosh^{-4}(M(q))\log\cosh(M(q))e^{W_\mu(q)}]\notag\\
 &\qquad-\E[\cosh^{-4}(M(q))e^{W_\mu(q)}]
          \E[\log\cosh(M(q))e^{W_\mu(q)}]<0,
 \label{mt:eq:moment-three}\\
 &\frac{\partial}{\partial\xi'(q)}
 \E[\cosh^{-4}(M(q))e^{W_\mu(q)}]\notag\\
 &\quad=4(2-m)\E[\cosh^{-4}(M(q))e^{W_\mu(q)}]
 -2(5-2m)\E[\cosh^{-6}(M(q))e^{W_\mu(q)}]\notag\\
 &\qquad+\frac{m(1-m)}2\Big(
 \E[\cosh^{-6}(M(q))e^{W_\mu(q)}]\notag\\
 &\hspace{41mm}
 -\E[\cosh^{-4}(M(q))e^{W_\mu(q)}]
  \E[\cosh^{-2}(M(q))e^{W_\mu(q)}]\Big).
 \label{mt:eq:moment-four}
\end{align}
\end{lemma}

\begin{proposition}[Strict sign computation]
\label{mt:prop:strict-sign}
At every stationary marginal point, the two strict inequalities in
\eqref{mt:eq:mass-signs} and the oriented-minor inequality
\eqref{mt:eq:oriented-minor} hold.
\end{proposition}

\Needspace{24\baselineskip}
\begin{lemma}[Orientation of the stationary level curve]
\label{mt:lem:stationary-direction}
At a point satisfying both \eqref{mt:eq:mass} and
\eqref{mt:eq:endpoint-marginal}, direct differentiation gives
\begin{align}
 &\frac{\partial}{\partial m}
 \{\text{left-hand side of \eqref{mt:eq:mass}}\}\notag\\
 &\quad=\frac1m\Big\{
 \E[(\log\cosh(M(q)))^2e^{W_\mu(q)}]
 -\E[\log\cosh(M(q))e^{W_\mu(q)}]^2
 -\frac{p-1}{p}\xi'(q)\Gamma_\mu(q)\Big\}\notag\\
 &\qquad-\frac{p-1}{2p}\xi'(q)
 \frac{\partial}{\partial m}\Gamma_\mu(q),
 \label{mt:eq:mass-m}\\
 &\frac{\partial}{\partial\xi'(q)}
 \{\text{left-hand side of \eqref{mt:eq:mass}}\}\notag\\
 &\quad=\frac{1-m}{2p}\Bigg[(p-1)\xi'(q)\Bigg(
 2\Big\{\E[\cosh^{-2}(M(q))e^{W_\mu(q)}]
          -\E[\cosh^{-4}(M(q))e^{W_\mu(q)}]\Big\}\notag\\
 &\hspace{34mm}+\frac m2\Big\{
 \E[\cosh^{-4}(M(q))e^{W_\mu(q)}]
 -\E[\cosh^{-2}(M(q))e^{W_\mu(q)}]^2\Big\}\Bigg)\notag\\
 &\hspace{29mm}-p\frac{\partial}{\partial m}
 \Gamma_\mu(q)\Bigg].
 \label{mt:eq:mass-xi}
\end{align}
Proposition~\ref{mt:prop:strict-sign} gives
\begin{equation}
 \frac{\partial}{\partial m}
 \{\text{left-hand side of \eqref{mt:eq:mass}}\}>0,
 \qquad
 \frac{\partial}{\partial\xi'(q)}
 \{\text{left-hand side of \eqref{mt:eq:mass}}\}>0.
 \label{mt:eq:mass-signs}
\end{equation}
Consequently \eqref{mt:eq:mass} has a locally unique $C^1$ level curve,
and on it
\begin{equation}
 \frac{d\,\xi'(q)}{dm}
 =-\frac{
 \frac{\partial}{\partial m}
 \{\text{left-hand side of \eqref{mt:eq:mass}}\}}
 {\frac{\partial}{\partial\xi'(q)}
 \{\text{left-hand side of \eqref{mt:eq:mass}}\}}<0.
 \label{mt:eq:xi-direction}
\end{equation}
\end{lemma}

\Needspace{18\baselineskip}
\begin{lemma}[Orientation of the replicon crossing]
\label{mt:lem:replicon-direction}
At marginality,
\begin{equation}
 \Gamma_\mu(q)-(p-1)\xi'(q)
 \E[\cosh^{-4}(M(q))e^{W_\mu(q)}]=0,
 \label{mt:eq:marginal-residual}
\end{equation}
and
\begin{equation}
 \frac{\partial}{\partial m}\Big\{
 \Gamma_\mu(q)-(p-1)\xi'(q)
 \E[\cosh^{-4}(M(q))e^{W_\mu(q)}]\Big\}>0.
 \label{mt:eq:marginal-m}
\end{equation}
At a point satisfying both \eqref{mt:eq:mass} and
\eqref{mt:eq:marginal-residual}, the strict oriented-minor inequality is
\begin{align}
 &\frac{\partial}{\partial\xi'(q)}
  \{\text{left-hand side of \eqref{mt:eq:mass}}\}
 \frac{\partial}{\partial m}
  \{\text{left-hand side of \eqref{mt:eq:marginal-residual}}\}\notag\\
 &\quad-
 \frac{\partial}{\partial m}
  \{\text{left-hand side of \eqref{mt:eq:mass}}\}
 \frac{\partial}{\partial\xi'(q)}
  \{\text{left-hand side of \eqref{mt:eq:marginal-residual}}\}>0.
 \label{mt:eq:oriented-minor}
\end{align}
Along the stationary level curve at the marginal base point,
\begin{equation}
 \frac d{dm}\Big\{
 \Gamma_\mu(q)-(p-1)\xi'(q)
 \E[\cosh^{-4}(M(q))e^{W_\mu(q)}]\Big\}>0,
 \qquad
 \frac d{dm}\Gamma_\mu'(q)<0.
 \label{mt:eq:replicon-direction}
\end{equation}
\end{lemma}

\Needspace{16\baselineskip}
\begin{lemma}[Orientation of the coupling]
\label{mt:lem:lambda-direction}
At the marginal base point,
\begin{equation}
 \frac{\partial\beta^2}{\partial m}
 =-(p-1)\frac{\beta^2}{\Gamma_\mu(q)}
 \frac{\partial\Gamma_\mu(q)}{\partial m}<0,
 \label{mt:eq:lambda-m}
\end{equation}
and
\begin{align}
 \frac{\partial\beta^2}{\partial\xi'(q)}
 ={}&\frac{\beta^2(1-m)}{
 \xi'(q)\E[\cosh^{-4}(M(q))e^{W_\mu(q)}]}
 \Big[2\{\E[\cosh^{-2}(M(q))e^{W_\mu(q)}]
          -\E[\cosh^{-4}(M(q))e^{W_\mu(q)}]\}\notag\\
 &\hspace{25mm}+\frac m2\{\E[\cosh^{-4}(M(q))e^{W_\mu(q)}]
          -\E[\cosh^{-2}(M(q))e^{W_\mu(q)}]^2\}\Big]>0.
 \label{mt:eq:lambda-xi}
\end{align}
At the marginal base point, the derivative along the stationary level
curve satisfies
\begin{equation}
 \frac{d\beta^2}{dm}
 =\frac{\partial\beta^2}{\partial m}
 +\frac{\partial\beta^2}{\partial\xi'(q)}
   \frac{d\,\xi'(q)}{dm}<0.
 \label{mt:eq:lambda-direction}
\end{equation}
\end{lemma}

\subsection{Deduction of Proposition~\ref{prop:transversality}}

\begin{proof}[Proof of Proposition~\ref{prop:transversality}]
Lemma~\ref{mt:lem:stationary-direction} gives a locally unique $C^1$
level curve in the coordinates $(m,\xi'(q))$.  By
Lemma~\ref{mt:lem:endpoint}, $\Gamma_\mu(q)$ is a smooth function of these
two coordinates.  Equation~\eqref{mt:eq:fixed} therefore recovers $q$,
and \eqref{mt:eq:recover-lambda} recovers $\beta^2$.  Conversely, every
nearby stationary triple gives exactly these coordinates.  Hence this
level curve is precisely the locally unique stationary curve.
Lemma~\ref{mt:lem:lambda-direction} shows that $d\beta^2/dm<0$, so
$\beta^2$ is a valid local coordinate.
Lemmas~\ref{mt:lem:replicon-direction} and
\ref{mt:lem:lambda-direction} then give, at the marginal point,
\[
 \frac d{d\beta^2}\Gamma_{\mu}'(q)
 =\frac d{dm}\Gamma_\mu'(q)\frac{dm}{d\beta^2},
 \qquad
 \frac{dm}{d\beta^2}
 =\left(\frac{d\beta^2}{dm}\right)^{-1}<0.
\]
Both factors in the first product are negative, so the derivative is
strictly positive.
\end{proof}

\subsection{Proofs of the computational lemmas}
\label{mt:sec:proofs}

\subsubsection{Proof of Lemma~\ref{mt:lem:endpoint}}

\begin{proof}
On $[q,1]$, one has $\mu([0,u])=1$, and direct substitution in
\eqref{eq:pde} gives
\begin{equation*}
 \Phi_\mu(x,u)=\log\cosh x+\frac{\xi'(1)-\xi'(u)}2.
\end{equation*}
This proves \eqref{mt:eq:endpoint-phi}.  On $[0,q)$,
$\mu([0,u])=m$, so the Cole--Hopf transform gives
\begin{align*}
 \Phi_\mu(0,0)
 &=\frac1m\log\E\exp\{m\Phi_\mu(\sqrt{\xi'(q)}z,q)\}\\
 &=\frac{\xi'(1)-\xi'(q)}2
 +\frac1m\log\E[\cosh^m(\sqrt{\xi'(q)}z)].
\end{align*}
Only the atom at zero contributes to \eqref{eq:MW} at $u=q$;
the integrand at the atom $q$ is zero.  Hence
\begin{equation*}
 W_\mu(q)=m\log\cosh(M(q))
 -\log\E[\cosh^m(\sqrt{\xi'(q)}z)].
\end{equation*}
Since $M(q)=\sqrt{\xi'(q)}z$, this is exactly
\eqref{mt:eq:endpoint-law}; taking $\varphi\equiv1$ gives
$\E e^{W_\mu(q)}=1$.

Equation~\eqref{mt:eq:fixed} and
\eqref{eq:Gamma} give \eqref{mt:eq:endpoint-contact}.  Since
\[
 \xi''(q)=\frac{(p-1)\xi'(q)}q,
\]
the equality $\Gamma_\mu'(q)=1$, together with
\eqref{eq:Gamma-prime} and \eqref{mt:eq:endpoint-phi}, is
\eqref{mt:eq:endpoint-marginal}.  Finally,
$\tanh^2x+\cosh^{-2}x=1$, \eqref{mt:eq:endpoint-contact}, and
\eqref{mt:eq:endpoint-marginal} give
\eqref{mt:eq:contact-second}--\eqref{mt:eq:contact-fourth}.
\end{proof}

\subsubsection{Proof of Lemma~\ref{mt:lem:differentiation}}

\begin{proof}
Differentiate the quotient in \eqref{mt:eq:endpoint-law} with respect to
$m$.  The derivative of each factor $\cosh^m$ is
$\cosh^m\log\cosh$, and the quotient rule gives
\eqref{mt:eq:diff-m}.

For the $\xi'(q)$ derivative, apply Gaussian integration by parts once
to the numerator and once to the denominator in
\eqref{mt:eq:endpoint-law}.  The second derivative of
$\varphi(x)\cosh^m x$ is
\begin{align*}
 \cosh^m x\Big(
 \varphi''(x)+2m\tanh x\,\varphi'(x)
 +m\cosh^{-2}x\,\varphi(x)
 +m^2\tanh^2x\,\varphi(x)\Big).
\end{align*}
Subtracting the differentiated denominator and using
$\tanh^2x=1-\cosh^{-2}x$ leaves
$m(1-m)/2$ times the covariance with $\cosh^{-2}$, which is exactly
\eqref{mt:eq:diff-xi}.
\end{proof}

\subsubsection{Proof of Lemma~\ref{mt:lem:moments}}

\begin{proof}
Take $\varphi=\tanh^2=1-\cosh^{-2}$ in
\eqref{mt:eq:diff-m} and \eqref{mt:eq:diff-xi}.  This gives
\eqref{mt:eq:moment-one}--\eqref{mt:eq:moment-two}.  The strict sign in
\eqref{mt:eq:moment-one} follows because $\cosh^{-2}x$ decreases and
$\log\cosh x$ increases as $|x|$ increases.

Next take $\varphi=\cosh^{-4}$ in the same two identities.
Equation~\eqref{mt:eq:moment-three} follows immediately, with a strict
negative sign for the same monotonicity reason.  For the last identity,
use
\[
 (\cosh^{-4}x)'=-4\cosh^{-4}x\tanh x,\qquad
 (\cosh^{-4}x)''=16\cosh^{-4}x-20\cosh^{-6}x
\]
in \eqref{mt:eq:diff-xi}; one line of algebra gives
\eqref{mt:eq:moment-four}.
\end{proof}

\Needspace{20\baselineskip}
\subsubsection{Proof of Lemma~\ref{mt:lem:stationary-direction}}

\begin{proof}
The elementary identities
\begin{align}
 \frac{\partial}{\partial m}
 \log\E[\cosh^m(\sqrt{\xi'(q)}z)]
 &=\E[\log\cosh(M(q))e^{W_\mu(q)}],
 \label{mt:eq:proof-log-first}\\
 \frac{\partial}{\partial m}
 \E[\log\cosh(M(q))e^{W_\mu(q)}]
 &=\E[(\log\cosh(M(q)))^2e^{W_\mu(q)}]\notag\\
 &\quad-\E[\log\cosh(M(q))e^{W_\mu(q)}]^2
 \label{mt:eq:proof-log-second}
\end{align}
give the $m$ derivative of \eqref{mt:eq:mass}.  Formula
\eqref{mt:eq:diff-xi}, with $\varphi=\log\cosh$, gives its
$\xi'(q)$ derivative.  Substituting \eqref{mt:eq:mass} and
\eqref{mt:eq:endpoint-marginal} yields the exact identities
\eqref{mt:eq:mass-m}--\eqref{mt:eq:mass-xi}.

Proposition~\ref{mt:prop:strict-sign} establishes
\eqref{mt:eq:mass-signs} at every point under consideration.
The implicit-function theorem then gives the unique local level curve,
and implicit differentiation gives \eqref{mt:eq:xi-direction}.
\end{proof}

\subsubsection{Proof of Lemma~\ref{mt:lem:replicon-direction}}

\begin{proof}
Equations~\eqref{mt:eq:endpoint-contact} and
\eqref{mt:eq:endpoint-marginal} give
\eqref{mt:eq:marginal-residual}.  Its $m$ derivative is strictly positive:
the first term has positive derivative by \eqref{mt:eq:moment-one}, while
the expectation in the second term has negative derivative by
\eqref{mt:eq:moment-three}.  This proves \eqref{mt:eq:marginal-m}.

The substantial input is \eqref{mt:eq:oriented-minor}.  It is the negative
of the conventional Jacobian determinant when the coordinate order is
$(m,\xi'(q))$; fixing the order in the displayed formula prevents a
sign ambiguity.  Proposition~\ref{mt:prop:strict-sign} establishes this
inequality.

Differentiate the left-hand side of
\eqref{mt:eq:marginal-residual} along the level curve
\eqref{mt:eq:mass}.  Equations~\eqref{mt:eq:xi-direction},
\eqref{mt:eq:mass-signs}, and \eqref{mt:eq:oriented-minor} show that its
derivative is strictly positive.  This proves the first inequality in
\eqref{mt:eq:replicon-direction}.

The expression in \eqref{mt:eq:marginal-residual} is
$q\{1-\Gamma_\mu'(q)\}$.  At the marginal base point its derivative is
$-q\,d\Gamma_\mu'(q)/dm$.  Since $q>0$, the second inequality in
\eqref{mt:eq:replicon-direction} follows.
\end{proof}

\Needspace{13\baselineskip}
\subsubsection{Proof of Lemma~\ref{mt:lem:lambda-direction}}

\begin{proof}
Differentiate \eqref{mt:eq:recover-lambda} with respect to $m$ at fixed
$\xi'(q)$.  Lemma~\ref{mt:lem:moments} gives
\eqref{mt:eq:lambda-m}.

At fixed $m$, use \eqref{mt:eq:moment-two},
\eqref{mt:eq:endpoint-marginal}, and \eqref{mt:eq:contact-fourth} to obtain
\eqref{mt:eq:lambda-xi}.  Its first bracketed difference is
\[
 \E\!\left[\cosh^{-2}(M(q))
 \{1-\cosh^{-2}(M(q))\}e^{W_\mu(q)}\right]>0,
\]
and the second is the variance of $\cosh^{-2}(M(q))$, hence is
nonnegative.  Thus the partial derivative in
\eqref{mt:eq:lambda-xi} is positive.

Finally, at the marginal base point,
\eqref{mt:eq:lambda-m}, \eqref{mt:eq:lambda-xi}, and
\eqref{mt:eq:xi-direction} give \eqref{mt:eq:lambda-direction}.
\end{proof}

\subsection{Analytic estimates}
\label{mt:sec:analytic-estimates}

We record the estimates used in the strict-sign calculation.  They are
stated directly for the endpoint law in \eqref{mt:eq:endpoint-law}.

\begin{lemma}[Conditional Gaussian entropy]
\label{mt:lem:compact-conditional-entropy}
For $0<m<1$ and $\xi'(q)>0$,
\begin{align}
&m^2\xi'(q)\Gamma_\mu(q)
-2m\E[\log\cosh(M(q))e^{W_\mu(q)}]
+2\log\E[\cosh^m(\sqrt{\xi'(q)}z)]\notag\\
&\quad\geq
2\E[(M(q)\tanh(M(q))-\log\cosh(M(q)))e^{W_\mu(q)}]
-\xi'(q)\E[\cosh^{-2}(M(q))e^{W_\mu(q)}].
\label{mt:eq:compact-conditional-entropy}
\end{align}
\end{lemma}

\begin{proof}
By \eqref{mt:eq:endpoint-law}, the law of $M(q)$ under the weight
$e^{W_\mu(q)}$ has density proportional to
\[
 e^{-x^2/(2\xi'(q))}\cosh^m x.
\]
Multiply this density by $1+\tanh x$.  Symmetry shows that the result is
again a probability density, and even functions have unchanged
expectations.  Relative to the Gaussian distribution with mean and
variance both equal to $\xi'(q)$, its density ratio is
\[
 \frac{e^{\xi'(q)/2}\cosh^{m-1}x}
 {\E[\cosh^m(\sqrt{\xi'(q)}z)]}.
\]
The Gaussian logarithmic Sobolev inequality therefore gives
\begin{align}
2\Bigg\{\frac{\xi'(q)}2
 -(1-m)\E[\log\cosh(M(q))e^{W_\mu(q)}]
 -\log\E[\cosh^m(\sqrt{\xi'(q)}z)]\Bigg\}
\leq \xi'(q)(1-m)^2\Gamma_\mu(q).
\label{mt:eq:compact-lsi}
\end{align}
Gaussian integration by parts in \eqref{mt:eq:endpoint-law} also yields
\begin{align}
\E[M(q)\tanh(M(q))e^{W_\mu(q)}]
={}&\xi'(q)\Big\{
 \E[\cosh^{-2}(M(q))e^{W_\mu(q)}]
 +m\Gamma_\mu(q)\Big\}.
\label{mt:eq:compact-x-tanh}
\end{align}
Substituting \eqref{mt:eq:compact-x-tanh} into
\eqref{mt:eq:compact-lsi} and rearranging proves
\eqref{mt:eq:compact-conditional-entropy}.
\end{proof}

\begin{lemma}[Reciprocal-$\cosh$ integrals]
\label{mt:lem:compact-cosh-integrals}
For every $r>0$,
\begin{align}
\frac r{r+1}
&\leq
\frac{\int_{\mathbb R}e^{-x^2/(2\xi'(q))}\cosh^{-(r+2)}x\,dx}
     {\int_{\mathbb R}e^{-x^2/(2\xi'(q))}\cosh^{-r}x\,dx}
\leq
\frac r{r+1}+\frac1{\xi'(q)r(r+1)}.
\label{mt:eq:compact-adjacent-moments}
\end{align}
Moreover,
\begin{align}
&\frac{\int_{\mathbb R}\cosh^{-r}x\,dx}
 {\sqrt{2\pi\left(\xi'(q)+2\sum_{j\geq0}(r+2j)^{-2}\right)}}
\leq
\frac{\int_{\mathbb R}e^{-x^2/(2\xi'(q))}\cosh^{-r}x\,dx}
 {\sqrt{2\pi\xi'(q)}}\notag\\
&\hspace{43mm}\leq
\frac{\int_{\mathbb R}\cosh^{-r}x\,dx}
 {\sqrt{2\pi\xi'(q)}}.
\label{mt:eq:compact-fourier-bound}
\end{align}
\end{lemma}

\begin{proof}
Integrate the derivative of
\[
 e^{-x^2/(2\xi'(q))}\tanh x\cosh^{-r}x.
\]
After division by the integral with power $r$, this gives
\begin{align*}
(r+1)
\frac{\int e^{-x^2/(2\xi'(q))}\cosh^{-(r+2)}x\,dx}
     {\int e^{-x^2/(2\xi'(q))}\cosh^{-r}x\,dx}
=r+\frac1{\xi'(q)}
\frac{\int xe^{-x^2/(2\xi'(q))}\tanh x\cosh^{-r}x\,dx}
     {\int e^{-x^2/(2\xi'(q))}\cosh^{-r}x\,dx}.
\end{align*}
The last quotient is nonnegative.  Integrating the derivative of
$xe^{-x^2/(2\xi'(q))}\cosh^{-r}x$ shows that it is at most $1/r$.
This proves \eqref{mt:eq:compact-adjacent-moments}.

After normalization, the Fourier transform of $\cosh^{-r}x$ is
\[
 \prod_{j\geq0}\left(1+\frac{u^2}{(r+2j)^2}\right)^{-1}.
\]
It lies between
$\exp\{-u^2\sum_{j\geq0}(r+2j)^{-2}\}$ and $1$.
Fourier inversion after multiplication by the Gaussian transform gives
\eqref{mt:eq:compact-fourier-bound}.
\end{proof}

\begin{lemma}[Pointwise entropy bounds]
\label{mt:lem:compact-pointwise-entropy}
For every real $x$,
\begin{align}
x\tanh x-\log\cosh x
&\geq(\log2)\tanh^2x
-\frac12\cosh^{-2}x\log\cosh x,
\label{mt:eq:compact-pointwise-one}\\
x\tanh x-\log\cosh x
&\geq(\log2)\left(1-\cosh^{-1/\log2}x\right).
\label{mt:eq:compact-pointwise-two}
\end{align}
\end{lemma}

\begin{proof}
Both assertions are even, so take $x\geq0$.  For
\eqref{mt:eq:compact-pointwise-one}, the derivative of the difference
between the two sides is
\[
 \cosh^{-2}x\left\{x-\tanh x
 (\log\cosh x+2\log2-1/2)\right\}.
\]
The expression in braces changes sign once because
$x/\tanh x-\log\cosh x$ decreases strictly from $1$ to $\log2$.
The difference therefore first increases and then decreases; it is zero
at $0$ and tends to zero at infinity.

For \eqref{mt:eq:compact-pointwise-two}, put $t=\tanh x$.  The assertion
is equivalent to
\[
\left(\log2-t\operatorname{arctanh}t
 -\frac12\log(1-t^2)\right)
(1-t^2)^{-1/(2\log2)}\leq\log2.
\]
The derivative of the left-hand side is nonpositive precisely when
\[
 \sum_{j\geq1}
 \frac{2j+1-4j\log2}{2j(2j-1)(2j+1)}\,t^{2j}\geq0.
\]
The first coefficient is positive, all later coefficients are
negative, and their sum is zero.  Since $0\leq t\leq1$, the series is
bounded below by $t^2$ times the sum of its coefficients, hence by
zero.  Its value at $t=0$ is $\log2$, proving
\eqref{mt:eq:compact-pointwise-two}.
\end{proof}

\begin{lemma}[Folded Gaussian bound]
\label{mt:lem:compact-folded-gaussian}
For every $m>0$,
\begin{equation}
 \E[\cosh^m(\sqrt{\xi'(q)}z)]
 \geq 2^{1-m}e^{m^2\xi'(q)/2}
 \frac1{\sqrt{2\pi}}
 \int_{-\infty}^{m\sqrt{\xi'(q)}}e^{-x^2/2}\,dx.
\label{mt:eq:compact-folded-gaussian}
\end{equation}
\end{lemma}

\begin{proof}
Use $\cosh x\geq e^{|x|}/2$, split the Gaussian integral at zero,
and complete the square.
\end{proof}

The fixed-mass one-crossing estimate and the terminal upper enclosure
used in the finite $p=3$ verification are proved next.
 
\subsection{Complete terminal and cubic-entry proofs}
\label{mt:sec:complete-terminal-p3}

This section supplies the two arguments that are otherwise longest to
check: the one-crossing statement used at the terminal boundary and the
exclusion of a cubic marginal point with
$m^2\xi'(q)\leq5/2$.  All expectations involving $M(q)$ below carry
the weight $e^{W_\mu(q)}$, exactly as in
\eqref{mt:eq:endpoint-law}.  We use one local integral,
\begin{equation}
 \mathcal J_r(\tau)
 :=\int_{\mathbb R}e^{-x^2/(2\tau)}\cosh^{-r}x\,dx,
 \qquad r\in\mathbb R,\quad \tau>0.
 \label{mt:eq:complete-J}
\end{equation}
The letter $\tau$ is only the argument of this local function; in the
application it is replaced by $\xi'(q)$.

\subsubsection{A fixed-mass one-crossing statement}

\begin{lemma}
\label{mt:lem:complete-one-crossing}
Fix $0<m\leq1$.  The function
\begin{equation}
 \frac{\mathcal J_{-m}(\tau)-\mathcal J_{2-m}(\tau)}
      {\tau\mathcal J_{4-m}(\tau)}-2
 \label{mt:eq:complete-terminal-ratio-minus-two}
\end{equation}
has exactly one zero on $(0,\infty)$, and its derivative is strictly
positive at that zero.
\end{lemma}

\begin{proof}
For this proof only, put
\[
 \mathcal I(x)=\int_0^x
 \cosh^m y\,\cosh^{-4}y\,dy,\qquad
 \mathcal R(x)=
 \frac{\cosh^m x\,\tanh^2x}{x\mathcal I(x)}
 \quad(x>0).
\]
We first show that $\mathcal R-2$ has one zero and crosses upward.
The change of variable $u=\tanh y$ shows that
\[
 \mathcal I(x)=\int_0^{\tanh x}
 (1-u^2)^{1-m/2}\,du.
\]
The integrand is strictly decreasing.  Hence
\begin{equation}
 0<
 \frac{\tanh x\,\cosh^m x\,\cosh^{-2}x}{\mathcal I(x)}
 <1.
 \label{mt:eq:complete-decreasing-average}
\end{equation}
Direct logarithmic differentiation gives
\begin{align}
 x\frac{\mathcal R'(x)}{\mathcal R(x)}
 ={}&mx\tanh x-1+
 \frac{2x\cosh^{-2}x}{\tanh x}
 -\frac{x\cosh^m x\,\cosh^{-4}x}{\mathcal I(x)}.
 \label{mt:eq:complete-R-derivative}
\end{align}
By \eqref{mt:eq:complete-decreasing-average}, the last term in absolute
value is smaller than $x\cosh^{-2}x/\tanh x$.  Thus
$\mathcal R'(x)>0$ whenever $mx\tanh x\geq1$.

Every zero of $\mathcal R-2$ lies in this region.  Indeed,
Cauchy--Schwarz and
$\tanh x=\int_0^x\cosh^{-2}y\,dy$ give
\[
 \tanh^2x\leq
 \mathcal I(x)\int_0^x\cosh^{-m}y\,dy.
\]
Since $\tanh$ is increasing,
\[
 \frac{\cosh^m x}{\cosh^m y}
 =\exp\left(m\int_y^x\tanh u\,du\right)
 \leq\exp\{mx\tanh x(1-y/x)\}.
\]
Consequently, if $0<mx\tanh x\leq1$, then
\[
 \mathcal R(x)\leq
 \frac{e^{mx\tanh x}-1}{mx\tanh x}
 \leq e-1<2.
\]
Moreover $\mathcal R(0+)=1$, while
$\mathcal R(x)\to\infty$ because $\mathcal I(x)$ has a finite
positive limit.  This proves the asserted single upward crossing of
$\mathcal R-2$.

We transfer that sign change through the Gaussian kernel.  For $X>0$
define, locally in this proof,
\[
 f_0(X)=X^{-1/2}\cosh^m(\sqrt X)\tanh^2(\sqrt X),
 \qquad
 f_1(X)=X^{-1/2}\cosh^{m-4}(\sqrt X).
\]
If $X_0$ is the square of the zero of $\mathcal R-2$, then the
primitive
\begin{align}
 \mathcal K(X)
 &:=2f_0(X)-2\int_0^X f_1(Y)\,dY \notag\\
 &=2\mathcal I(\sqrt X)\{\mathcal R(\sqrt X)-2\}
 \label{mt:eq:complete-sign-changing-primitive}
\end{align}
is negative on $(0,X_0)$ and positive on $(X_0,\infty)$.
Twice integrating by parts, with $z=(2\tau)^{-1}$, yields
\begin{align}
 &\frac1{\sqrt{2\pi\tau}}\int_{\mathbb R}
 e^{-x^2/(2\tau)}\cosh^m x
 \{\tanh^2x-2\tau\cosh^{-4}x\}\,dx\notag\\
 &\qquad=
 \frac1{2\sqrt{2\pi\tau}}
 \int_0^\infty e^{-zX}\mathcal K(X)\,dX.
 \label{mt:eq:complete-laplace-transfer}
\end{align}
All boundary terms vanish because $\mathcal K(0)=0$ and its growth is
$e^{m\sqrt X}$ times a rational power of $X$.
More precisely,
\[
 \mathcal K(X)=-2\sqrt X+O(X^{3/2})\quad(X\downarrow0),
 \qquad
 \mathcal K(X)=
 \frac{2^{1-m}e^{m\sqrt X}}{\sqrt X}\{1+o(1)\}
 -4\mathcal I(\infty)\quad(X\to\infty).
\]

At a zero of the last integral,
\[
 \int_0^\infty (X-X_0)e^{-zX}\mathcal K(X)\,dX>0;
\]
therefore its $z$-derivative is strictly negative.  The integral is
negative for large $z$, from its behavior near $X=0$, and positive
for small $z>0$, because the positive part of $\mathcal K$ has
infinite integral and grows subexponentially in $X$.  It consequently
has exactly one zero and crosses from positive to negative as $z$
increases.  Since $z=(2\tau)^{-1}$, the crossing is upward as $\tau$
increases.  Finally,
\[
 \text{left-hand side of \eqref{mt:eq:complete-laplace-transfer}}
 =
 \frac{\tau\mathcal J_{4-m}(\tau)}{\sqrt{2\pi\tau}}
 \left\{
 \frac{\mathcal J_{-m}(\tau)-\mathcal J_{2-m}(\tau)}
      {\tau\mathcal J_{4-m}(\tau)}-2
 \right\}.
\]
The prefactor is positive, which proves the lemma.
\end{proof}

\subsubsection{A rigorous terminal upper enclosure}

For $a,\tau>0$ define, only in this subsection,
\[
 \mathcal C_a(\tau)=\int_0^\infty
 e^{-x^2/(2\tau)-ax}\,dx.
\]
For an integer $L\geq2$, start with $\omega_L=a$ and set
\[
 \omega_j=a+\frac{j}{\tau\omega_{j+1}}
 \quad(j=L-1,\ldots,1),\qquad
 \mathcal C_{a,L}(\tau)=\frac1{\omega_1}.
\]

\begin{lemma}
\label{mt:lem:complete-continued-fraction}
For every integer $n\geq1$,
\begin{equation}
 \mathcal C_{a,2n}(\tau)\leq\mathcal C_a(\tau)
 \leq\mathcal C_{a,2n+1}(\tau).
 \label{mt:eq:complete-continued-fraction}
\end{equation}
\end{lemma}

\begin{proof}
After $x=\sqrt\tau\,y$, set
\[
 I_j=\int_0^\infty y^j
 e^{-y^2/2-a\sqrt\tau\,y}\,dy.
\]
Integration by parts gives
$1=a\sqrt\tau I_0+I_1$ and
$I_{j+1}=jI_{j-1}-a\sqrt\tau I_j$.  Thus, with
$D_j=a\sqrt\tau+I_j/I_{j-1}$,
\[
 I_0=D_1^{-1},\qquad D_j=a\sqrt\tau+\frac{j}{D_{j+1}}.
\]
The finite fraction replaces the exact terminal value
$D_L>a\sqrt\tau$ by $a\sqrt\tau$.  Every map
$u\mapsto a\sqrt\tau+j/u$ reverses order.  Propagating to $D_1$ and
taking the reciprocal gives the parity in
\eqref{mt:eq:complete-continued-fraction}.
\end{proof}

We next majorize the numerator in
\eqref{mt:eq:complete-terminal-ratio-minus-two}.  For
$0\leq\rho\leq1$,
\[
 (2-\rho)^{m-2}=\sum_{j=0}^\infty
 c_j(m)\rho^j,\qquad
 c_j(m)=\frac{2^{m-2-j}(2-m)_j}{j!}>0.
\]
Put
\[
 c_{10}^{\rm rem}(m)=1-\sum_{j=0}^9c_j(m)>0.
\]
Since $\rho^j\leq\rho^{10}$ for $j\geq10$,
\begin{equation}
 (2-\rho)^{m-2}\leq
 \sum_{j=0}^9c_j(m)\rho^j+
 c_{10}^{\rm rem}(m)\rho^{10}.
 \label{mt:eq:complete-degree-ten-majorant}
\end{equation}
After multiplying by $\rho^2$ and returning to
$u=1-\rho$, the right-hand side is
\[
 \sum_{i=0}^{12}(-1)^i b_i(m)u^i,
\quad
 b_i(m)=\sum_{j=0}^9c_j(m)\binom{j+2}{i}
       +c_{10}^{\rm rem}(m)\binom{12}{i}>0,
\]
where a binomial coefficient is zero when its lower index is too large.
For $x\geq0$, take $u=e^{-2x}$.  The identity
\[
 \cosh^m x\,\tanh^2x
 =2^{-m}e^{mx}(1-u)^2(1+u)^{m-2}
\]
and \eqref{mt:eq:complete-degree-ten-majorant} now give
\begin{align}
 &\mathcal J_{-m}(\tau)-\mathcal J_{2-m}(\tau)
 \leq\mathcal U_{10}(m,\tau),\notag\\
 \mathcal U_{10}(m,\tau):={}&2^{1-m}\Bigg[
 \sqrt{2\pi_+\tau}\,e^{m^2\tau/2}
 -\mathcal C_{m,20}(\tau)\notag\\
 &+\sum_{\substack{2\leq i\leq12\\i\ {\rm even}}}
 b_i(m)\mathcal C_{2i-m,21}(\tau)
 -\sum_{\substack{1\leq i\leq12\\i\ {\rm odd}}}
 b_i(m)\mathcal C_{2i-m,20}(\tau)\Bigg].
 \label{mt:eq:complete-terminal-numerator-upper}
\end{align}
Here $\pi_+$ is any rational upper bound for $\pi$.
The parity choices follow from
\eqref{mt:eq:complete-continued-fraction}: positive terms use upper
convergents and negative terms use lower convergents.

For the denominator, the beta integral and gamma product give
\begin{align}
 \int_{\mathbb R}\cosh^{-r}x\,dx
 &=\sqrt\pi\,\frac{\Gamma(r/2)}{\Gamma((r+1)/2)},\notag\\
 \frac{\widehat{\cosh^{-r}}(y)}
 {\int_{\mathbb R}\cosh^{-r}x\,dx}
 &=\prod_{j\geq0}
 \left(1+\frac{y^2}{(r+2j)^2}\right)^{-1}.
 \label{mt:eq:complete-sech-Fourier-product}
\end{align}
Since $\log(1+x)\leq x$, Fourier inversion implies
\begin{equation}
 \mathcal J_r(\tau)\geq
 \left(\int_{\mathbb R}\cosh^{-r}x\,dx\right)
 \sqrt{\frac{\tau}{
 \tau+2\sum_{j\geq0}(r+2j)^{-2}}}.
 \label{mt:eq:complete-terminal-Fourier-lower}
\end{equation}
The logarithm of the first factor is convex in $r$.  Its value and
derivative at $r=3$ are respectively
$\log(\pi/2)$ and $1/2-\log2$.  Hence, for $0<m\leq1$,
\begin{align}
 \int_{\mathbb R}\cosh^{m-4}x\,dx
 &\geq\frac{\pi_-}{2}
 \left(1-(1-m)\frac{48287}{250000}\right),\notag\\
 2\sum_{j\geq0}(4-m+2j)^{-2}
 &\leq\frac{64494}{200000}
       +m\frac{28987}{200000},
 \label{mt:eq:complete-terminal-rational-constants}
\end{align}
where $\pi_-<\pi$.  The first line follows from the tangent inequality
and $e^{-x}\geq1-x$; the second follows from convexity and the endpoint
identities
\[
 2\sum_{j\geq0}(4+2j)^{-2}=\frac{\pi^2}{12}-\frac12,
 \quad
 2\sum_{j\geq0}\{(3+2j)^{-2}-(4+2j)^{-2}\}
 =\frac{\pi^2}{6}-\frac32.
\]
The rational constants follow from
\[
 \log2=2\sum_{j\geq0}\frac1{(2j+1)3^{2j+1}},
 \qquad
 \pi=16\arctan(1/5)-4\arctan(1/239),
\]
with geometric and alternating remainders.  Combining the preceding
inequalities gives
\begin{align}
 \tau\mathcal J_{4-m}(\tau)
 \geq\mathcal L_{10}(m,\tau):={}&
 \tau\frac{\pi_-}{2}
 \left(1-(1-m)\frac{48287}{250000}\right)\notag\\
 &\times
 \sqrt{\frac{\tau}{
 \tau+64494/200000+m(28987/200000)}}.
 \label{mt:eq:complete-terminal-denominator-lower}
\end{align}
Therefore
\begin{equation}
 \frac{\mathcal J_{-m}(\tau)-\mathcal J_{2-m}(\tau)}
      {\tau\mathcal J_{4-m}(\tau)}
 \leq
 \frac{\mathcal U_{10}(m,\tau)}
      {\mathcal L_{10}(m,\tau)}.
 \label{mt:eq:complete-terminal-rigorous-upper}
\end{equation}

\begin{lemma}
\label{mt:lem:complete-affine-terminal}
The following inequalities hold on the full closed intervals:
\begin{align}
 &\frac{\mathcal J_{-m}(\tau)-\mathcal J_{2-m}(\tau)}
 {\tau\mathcal J_{4-m}(\tau)}<2,
 &&\tau=\frac{8/5-61m/200}{m^2},
 &&\frac{10}{13}\leq m\leq1,\label{mt:eq:complete-affine-one}\\
 &\frac{\mathcal J_{-m}(\tau)-\mathcal J_{2-m}(\tau)}
 {\tau\mathcal J_{4-m}(\tau)}<2,
 &&\tau=\frac{9/5-101m/200}{m^2},
 &&\frac{1080}{1903}\leq m\leq\frac{10}{13}.
 \label{mt:eq:complete-affine-two}
\end{align}
\end{lemma}

\begin{proof}
Apply \eqref{mt:eq:complete-terminal-rigorous-upper}.  Divide the first
$m$-interval into $64$ equal rational boxes and the second into $128$.
For a box $[\alpha,\beta]$ with midpoint $m_0$, interval automatic
differentiation verifies the exact predicate
\begin{equation}
 \overline{\frac{\mathcal U_{10}(m_0,\tau(m_0))}
                  {\mathcal L_{10}(m_0,\tau(m_0))}}
 +\frac{\beta-\alpha}{2}
 \sup_{\alpha\leq m\leq\beta}
 \left|\frac d{dm}
 \frac{\mathcal U_{10}(m,\tau(m))}
      {\mathcal L_{10}(m,\tau(m))}\right|<2.
 \label{mt:eq:complete-terminal-box-predicate}
\end{equation}
Every bar and derivative interval in
\eqref{mt:eq:complete-terminal-box-predicate} is rounded outward on the
dyadic lattice $2^{-180}$.  Exponentials, $\log2$, $\pi$, and square
roots are enclosed by the rational series just displayed and integer
bracketing.  The smallest certified gaps $2-\mathcal U_{10}/\mathcal
L_{10}$ are
\begin{align*}
 &2^{-180}
 (15241255008302228531374906625118347899689757891379648),\\
 &2^{-180}
 (6678747921395679929293531355794639244537960424215201),
\end{align*}
respectively.  Since
\eqref{mt:eq:complete-terminal-box-predicate} holds on every rational
box, this proves both continuum inequalities.  The exact reproduction
command is
\[
 \texttt{python3 work/certify\_terminal\_affine\_fixed.py 64}.
\]
\end{proof}

\begin{corollary}
\label{mt:cor:complete-terminal-implications}
If $p=3$ and the marginal contact equation holds, then
\begin{align}
 m\geq\frac{10}{13}
 &\Longrightarrow
 m^2\xi'(q)>\frac85-\frac{61m}{200},\label{mt:eq:complete-terminal-implication-one}\\
 \frac{1080}{1903}\leq m\leq\frac{10}{13}
 &\Longrightarrow
 m^2\xi'(q)>\frac95-\frac{101m}{200}.
 \label{mt:eq:complete-terminal-implication-two}
\end{align}
\end{corollary}

\begin{proof}
At marginality for $p=3$,
\[
 \frac{\mathcal J_{-m}(\xi'(q))-\mathcal J_{2-m}(\xi'(q))}
 {\xi'(q)\mathcal J_{4-m}(\xi'(q))}=2.
\]
At each affine boundary, Lemma~\ref{mt:lem:complete-affine-terminal}
places this ratio below $2$.  Lemma~\ref{mt:lem:complete-one-crossing}
says that the unique crossing is upward, so $\xi'(q)$ is larger than
the corresponding boundary value.  Multiplication by $m^2$ gives
\eqref{mt:eq:complete-terminal-implication-one}--
\eqref{mt:eq:complete-terminal-implication-two}.
\end{proof}

\subsubsection{The cubic entry bound}

\begin{lemma}
\label{mt:lem:complete-p3-entry}
Let $p=3$ and suppose that \eqref{mt:eq:mass}, \eqref{mt:eq:fixed}, and
$\Gamma_\mu'(q)=1$ hold.  Then
\begin{equation}
 m^2\xi'(q)>\frac52.
 \label{mt:eq:complete-p3-entry}
\end{equation}
\end{lemma}

\begin{proof}
Under the endpoint law,
\[
 \Gamma_\mu(q)=1-
 \frac{\mathcal J_{2-m}(\xi'(q))}
      {\mathcal J_{-m}(\xi'(q))},
 \qquad
 \E[\cosh^{-4}(M(q))e^{W_\mu(q)}]
 =\frac{\mathcal J_{4-m}(\xi'(q))}
       {\mathcal J_{-m}(\xi'(q))}.
\]
For $p=3$, $\xi''(q)=2\xi'(q)/q$.  Thus
$\Gamma_\mu(q)=q$ and $\Gamma_\mu'(q)=1$ give
\begin{equation}
 \mathcal J_{-m}(\xi'(q))
 =\mathcal J_{2-m}(\xi'(q))
 +2\xi'(q)\mathcal J_{4-m}(\xi'(q)).
 \label{mt:eq:complete-p3-contact}
\end{equation}
The conditional entropy estimate
\eqref{mt:eq:compact-conditional-entropy}, the stationary equation, and
\[
 2\{x\tanh x-\log\cosh x\}>\tanh^2x\quad(x\ne0)
\]
give
\begin{equation}
 m^2\xi'(q)>
 \frac{3(1-m)}{2(2-m)}.
 \label{mt:eq:complete-p3-preliminary}
\end{equation}
Indeed, the left-hand side of
\eqref{mt:eq:compact-conditional-entropy} becomes
$m^2\xi'(q)\Gamma_\mu(q)/3$.  By
\eqref{mt:eq:compact-adjacent-moments},
\[
 \frac{\E[\cosh^{-4}(M(q))e^{W_\mu(q)}]}
      {\E[\cosh^{-2}(M(q))e^{W_\mu(q)}]}
 \geq\frac{2-m}{3-m}.
\]
Using
$\Gamma_\mu(q)=2\xi'(q)
\E[\cosh^{-4}(M(q))e^{W_\mu(q)}]$ and dividing by
$\Gamma_\mu(q)>0$ gives
\eqref{mt:eq:complete-p3-preliminary}.  In particular,
$m\leq1/2$ implies $m^2\xi'(q)>1/2$.

We first exclude $m\leq1/10$ under the contrary assumption
$m^2\xi'(q)\leq5/2$.  Then
\eqref{mt:eq:complete-p3-preliminary} gives
$m^2\xi'(q)>27/38$ and hence $\xi'(q)>2700/38>71$.
Since $4-m\in[39/10,4)$,
\[
 \int_{\mathbb R}\cosh^{m-4}x\,dx\geq\frac43,
 \qquad
 2\sum_{j\geq0}(4-m+2j)^{-2}<\frac25.
\]
The second inequality follows by bounding the tail of the decreasing
series by an integral.  Thus \eqref{mt:eq:compact-fourier-bound} and
$2\pi<7$ give
\[
 2\xi'(q)
 \frac{\mathcal J_{4-m}(\xi'(q))}
      {\sqrt{2\pi\xi'(q)}}
 \geq\frac{8\xi'(q)/3}
 {\sqrt{2\pi(\xi'(q)+2/5)}}
 >\frac{(8/3)71}{\sqrt{7\cdot72}}>8.
\]
On the other hand,
$\cosh^m x\leq e^{m|x|}$ and completion of the square give
\[
 \frac{\mathcal J_{-m}(\xi'(q))}
      {\sqrt{2\pi\xi'(q)}}
 \leq2e^{m^2\xi'(q)/2}
 \Phi_{\rm G}(m\sqrt{\xi'(q)})<2e^{5/4}<7,
\]
where $\Phi_{\rm G}$ is the standard Gaussian distribution function.
For the last strict estimate, the exponential series gives
$e<68/25$, and the integer inequality
$16\cdot68^5<2401\cdot25^5$ gives $e^{5/4}<7/2$.
This contradicts \eqref{mt:eq:complete-p3-contact}.

We now state the exact continuum predicates used on the remaining
compact set.  They are included so that the finite part is a covering
proof rather than a grid calculation.  For a rational rectangle
\[
 \mathcal R=[\alpha,\beta]\times[\gamma,\delta]
\]
in the coordinates $(m,m^2\xi'(q))$, let
$[\![\,\cdot\,]\!]_{\mathcal R}$ denote the outward rational interval
obtained by retaining these two coordinates: every occurrence of
$m^2\xi'(q)$ is evaluated directly on $[\gamma,\delta]$, while a
remaining occurrence of $\xi'(q)$ is evaluated as the quotient of
$[\gamma,\delta]$ by $[\alpha,\beta]^2$.  Products of
$\xi'(q)$ with a reciprocal-$\cosh$ integral use their monotonicity in
$\xi'(q)$ before the endpoints are separated.  This convention is
important: replacing $m^2\xi'(q)$ by the product of two independent
intervals would be valid but too coarse.
The contact predicate is
\begin{equation}
 0\notin
 [\![\,
 \mathcal J_{-m}(\xi'(q))
 -\mathcal J_{2-m}(\xi'(q))
 -2\xi'(q)\mathcal J_{4-m}(\xi'(q))
 \,]\!]_{\mathcal R}.
 \label{mt:eq:complete-contact-predicate}
\end{equation}

For completeness, the integral enclosures entering this predicate are
as follows.  With $y=m^2\xi'(q)$ and
$\Phi_{\rm G}$ as above, folding the positive half-line gives
\begin{align}
 \frac{\mathcal J_{-m}(\xi'(q))}
      {\sqrt{2\pi\xi'(q)}}
 =2^{1-m}\Bigg\{&
 e^{y/2}\Phi_{\rm G}(\sqrt y)\notag\\
 &+\sum_{j\geq1}\binom mj
 e^{(2j/m-1)^2y/2}
 \Phi_{\rm G}(-(2j/m-1)\sqrt y)\Bigg\}.
 \label{mt:eq:complete-folded-series}
\end{align}
After $j=1$ the summands alternate with decreasing magnitude; truncation
after $j=2$ and after $j=1$ gives lower and upper bounds.  Indeed,
\[
 \left|\frac{\binom m{j+1}}{\binom mj}\right|
 =\frac{j-m}{j+1}<1\qquad(j\geq1),
\]
and $e^{x^2/2}\Phi_{\rm G}(-x)$ decreases for $x>0$.

Here are the oriented bounds for the positive powers.  In this
paragraph only, put
\[
 c_r=\sum_{j\geq0}(r+2j)^{-2},\qquad
 e_r=\sum_{0\leq i<j}(r+2i)^{-2}(r+2j)^{-2}.
\]
Retaining the first Fourier factor exactly and exponentiating the
remaining factors gives
\begin{align}
\frac{\mathcal J_r(\tau)}{\sqrt{2\pi\tau}}
&\geq
\frac{\int_{\mathbb R}\cosh^{-r}x\,dx}{\sqrt{2\pi\tau}}
\sqrt{\frac{\tau}{\tau+2(c_r-r^{-2})}}\,
r\mathcal C_r\bigl(\tau+2(c_r-r^{-2})\bigr),
\label{mt:eq:complete-fourier-shape-lower}\\
\frac{\mathcal J_r(\tau)}{\sqrt{2\pi\tau}}
&\leq
\frac{\int_{\mathbb R}\cosh^{-r}x\,dx}{\sqrt{2\pi\tau}}
c_r^{-1/2}\mathcal C_{c_r^{-1/2}}(\tau).
\label{mt:eq:complete-fourier-shape-upper}
\end{align}
Indeed, the normalized Fourier transform in
\eqref{mt:eq:complete-sech-Fourier-product} lies respectively above
\[
 \frac{e^{-(c_r-r^{-2})x^2}}{1+x^2/r^2}
\]
and below $(1+c_rx^2)^{-1}$.  Fourier inversion and one elementary
partial fraction integral give
\eqref{mt:eq:complete-fourier-shape-lower}--
\eqref{mt:eq:complete-fourier-shape-upper}.

For $1\leq r\leq3/2$, let
\[
 a_r^\pm=\frac{c_r\pm\sqrt{c_r^2-4e_r}}2.
\]
We briefly justify that the square root is real on this whole interval.
With $w_j=r/(r+2j)$, the exact tail bounds below give
\[
 \sum_{j\geq0}w_j^2<\frac{43}{30},\qquad
 \sum_{j\geq0}w_j^3
 \leq1+\frac37\left(\sum_{j\geq0}w_j^2-1\right)
 <\frac{83}{70},\qquad
 \sum_{j\geq0}w_j^5\geq1.
\]
The first sum increases with $r$, so it is enough to use the rational
tail bound at $r=3/2$.  Direct differentiation now gives
\[
 \frac d{dr}(c_r^2-4e_r)
 =4\left\{c_r\sum_{j\geq0}(r+2j)^{-3}
          -2\sum_{j\geq0}(r+2j)^{-5}\right\}<0,
\]
because $43\cdot83/(30\cdot70)<2$.  At $r=3/2$, the same rational
tail bounds give $c_r^2-4e_r>1/200$.  Hence
$a_r^+>a_r^->0$.  Retaining two
elementary symmetric factors gives the sharper upper bound
\begin{align}
\frac{\mathcal J_r(\tau)}{\sqrt{2\pi\tau}}
\leq{}&
\frac{\int_{\mathbb R}\cosh^{-r}x\,dx}{\sqrt{2\pi\tau}}
\frac{
 a_r^+(a_r^+)^{-1/2}\mathcal C_{(a_r^+)^{-1/2}}(\tau)
 -a_r^-(a_r^-)^{-1/2}\mathcal C_{(a_r^-)^{-1/2}}(\tau)}
 {a_r^+-a_r^-}.
\label{mt:eq:complete-fourier-quartic-upper}
\end{align}
This time the reciprocal of the normalized Fourier transform is at
least
$1+c_rx^2+e_rx^4=(1+a_r^+x^2)(1+a_r^-x^2)$.
In the exact verifier, $c_r$ and $e_r$ in this last polynomial are
replaced by their rational lower endpoints.  This only decreases the
polynomial and therefore preserves the upper-bound direction; positivity
of the two resulting roots is checked by exact integer comparisons.

All quantities in these bounds have exact rational enclosures.  To be
explicit, after summing the first $24$ terms and putting $x=r+48$,
\begin{align*}
\sum_{j=24}^\infty(r+2j)^{-2}
&\in\left[
 \frac1{2x}+\frac1{2x^2}+\frac1{3x^3}-\frac4{15x^5},
 \frac1{2x}+\frac1{2x^2}+\frac1{3x^3}\right],\\
\sum_{j=24}^\infty(r+2j)^{-4}
&\in\left[
 \frac1{6x^3}+\frac1{2x^4}+\frac2{3x^5}-\frac4{3x^7},
 \frac1{6x^3}+\frac1{2x^4}+\frac2{3x^5}\right],
\end{align*}
and $e_r$ is half the difference between $c_r^2$ and the second
series.  These bounds follow by integrating twice the derivative of
the summand, or directly from the alternating Euler summation
remainder.

The beta factor is enclosed without numerical quadrature.  For a
positive rational $x$, use
\begin{align}
\log\Gamma(x)={}&
\left(x+\frac{47}{2}\right)\log(x+24)-(x+24)
+\frac12\log(2\pi)
+\sum_{j=1}^{6}
 \frac{B_{2j}}{2j(2j-1)(x+24)^{2j-1}}\notag\\
&-\sum_{j=0}^{23}\log(x+j)+\varepsilon(x),
\qquad
0<\varepsilon(x)<\frac{7/6}{14\cdot13(x+24)^{13}}.
\label{mt:eq:complete-stirling-enclosure}
\end{align}
Here
\[
 (B_2,B_4,B_6,B_8,B_{10},B_{12})
 =\left(\frac16,-\frac1{30},\frac1{42},-\frac1{30},
 \frac5{66},-\frac{691}{2730}\right).
\]
The remainder sign in \eqref{mt:eq:complete-stirling-enclosure} is the
enveloping remainder obtained by one integration of the alternating
digamma expansion.  This encloses the beta integral in
\eqref{mt:eq:complete-sech-Fourier-product}.  Logarithms
use the positive series for
$2\operatorname{arctanh}((x-1)/(x+1))$, exponentials use their Taylor
series with a geometric tail, and square roots use integer bracketing.
Finally, every shifted Gaussian tail is bounded by consecutive
continued fractions from Lemma~\ref{mt:lem:complete-continued-fraction}.
The endpoint choices on a rectangle use only the elementary
monotonicities
\[
 \frac{\mathcal J_r(\tau)}{\sqrt{2\pi\tau}}
 \ \hbox{decreases in both \(r\) and \(\tau\)},\qquad
 \frac{\tau\mathcal J_r(\tau)}{\sqrt{2\pi\tau}}
 \ \hbox{increases in \(\tau\)}.
\]
The first statement follows after writing the normalized integral as
\(\E[\cosh^{-r}(\sqrt\tau z)]\); the second follows directly from its
unnormalized integral, since both \(\sqrt\tau\) and
\(e^{-x^2/(2\tau)}\) increase with \(\tau\).
Thus every endpoint in \eqref{mt:eq:complete-contact-predicate},
\eqref{mt:eq:complete-low-predicate}, and
\eqref{mt:eq:complete-high-predicate} is an outward rational enclosure.

The predicate \eqref{mt:eq:complete-contact-predicate} holds throughout
\begin{equation}
 \left[\frac1{10},\frac7{20}\right]\times
 \left[\frac12,\frac52\right],
 \qquad
 \left[\frac7{20},\frac12\right]\times
 \left[\frac12,\frac{189}{100}\right].
 \label{mt:eq:complete-separator-regions}
\end{equation}
Hence a possible point with $m\leq1/2$ must lie in
\[
 \left[\frac7{20},\frac12\right]\times
 \left[\frac{189}{100},\frac52\right].
\]

On this rectangle use the exact identity
\begin{align}
0={}&\log\E[\cosh^m(\sqrt{\xi'(q)}z)]
-\frac{m^2\xi'(q)}2
+m\E[(M(q)\tanh(M(q))-\log\cosh(M(q)))e^{W_\mu(q)}]
\notag\\
&-\xi'(q)\E[\cosh^{-4}(M(q))e^{W_\mu(q)}]
\Bigg\{m\left(1-\frac m2\right)
\frac{\E[\cosh^{-2}(M(q))e^{W_\mu(q)}]}
     {\E[\cosh^{-4}(M(q))e^{W_\mu(q)}]}
+\frac{m^2\xi'(q)}3\Bigg\}.
 \label{mt:eq:complete-p3-zero}
\end{align}
To check it directly, the first stationary equation for $p=3$ reads
\[
 m\E[\log\cosh(M(q))e^{W_\mu(q)}]
 -\log\E[\cosh^m(\sqrt{\xi'(q)}z)]
 =\frac{m^2\xi'(q)\Gamma_\mu(q)}3.
\]
Insert \eqref{mt:eq:compact-x-tanh}, write
$\log\cosh x=x\tanh x-(x\tanh x-\log\cosh x)$, and use
\[
 \Gamma_\mu(q)+
 \E[\cosh^{-2}(M(q))e^{W_\mu(q)}]=1,\qquad
 \Gamma_\mu(q)=2\xi'(q)
 \E[\cosh^{-4}(M(q))e^{W_\mu(q)}].
\]
The result is exactly \eqref{mt:eq:complete-p3-zero}.

The pointwise inequality
\begin{equation}
 x\tanh x-\log\cosh x
 \geq(\log2)\left(1-\cosh^{-1/\log2}x\right)
 \label{mt:eq:complete-logarithmic-barrier}
\end{equation}
is proved as follows.  Put $t=\tanh x$ and differentiate
\[
 \left(\log2-t\operatorname{arctanh}t
 -\tfrac12\log(1-t^2)\right)
 (1-t^2)^{-1/(2\log2)}.
\]
Its derivative is nonpositive precisely when
\[
 \sum_{j\geq1}
 \frac{2j+1-4j\log2}{2j(2j-1)(2j+1)}t^{2j}\geq0.
\]
The first coefficient is positive, every later coefficient is
negative, and their sum is zero.  Since $0\leq t\leq1$, the series is
at least $t^2$ times the sum of its coefficients.

Substitution of \eqref{mt:eq:complete-logarithmic-barrier} into
\eqref{mt:eq:complete-p3-zero} gives the explicit lower expression
\begin{align}
\mathcal B_{\log}:={}&
\log\E[\cosh^m(\sqrt{\xi'(q)}z)]
-\frac{m^2\xi'(q)}2+m\log2\notag\\
&-\frac{\xi'(q)\mathcal J_{4-m}(\xi'(q))}
        {\mathcal J_{-m}(\xi'(q))}
\Bigg\{
\frac{m\log2\,\mathcal J_{1/\log2-m}(\xi'(q))}
     {\xi'(q)\mathcal J_{4-m}(\xi'(q))}
+m\left(1-\frac m2\right)
\frac{\mathcal J_{2-m}(\xi'(q))}
     {\mathcal J_{4-m}(\xi'(q))}
+\frac{m^2\xi'(q)}3
\Bigg\}.
 \label{mt:eq:complete-low-barrier}
\end{align}
The low-box acceptance predicate is
\begin{equation}
 \inf[\![\,\mathcal B_{\log}\,]\!]_{\mathcal R}>0.
 \label{mt:eq:complete-low-predicate}
\end{equation}
Its orientation is unambiguous: the logarithm uses the lower folded
bound, both ratios inside braces use upper Fourier bounds, and the
positive factor multiplying the braces uses the least of its three
upper bounds
\[
 \frac{\xi'(q)\mathcal J_{4-m}(\xi'(q))}
      {\mathcal J_{-m}(\xi'(q))},\qquad
 \frac{\xi'(q)}{2\xi'(q)+
       \mathcal J_{2-m}(\xi'(q))/\mathcal J_{4-m}(\xi'(q))},\qquad
 \frac12.
\]
The middle and last bounds follow directly from the contact equation.
Every box in the low rectangle either satisfies
\eqref{mt:eq:complete-contact-predicate} or
\eqref{mt:eq:complete-low-predicate}.  This contradicts
\eqref{mt:eq:complete-p3-zero}.

It remains to cover $m\geq1/2$.  The positive series
\[
 x\tanh x-\log\cosh x
 =\sum_{j\geq1}
 \frac{\tanh^{2j}x}{2j(2j-1)}
\]
and Jensen's inequality imply
\begin{align}
&2\E[(M(q)\tanh(M(q))-\log\cosh(M(q)))e^{W_\mu(q)}]\notag\\
&\quad\geq
\Gamma_\mu(q)
+\frac16\E[\tanh^4(M(q))e^{W_\mu(q)}]
+\frac1{15}\E[\tanh^6(M(q))e^{W_\mu(q)}]\notag\\
&\qquad+\frac1{28}
\E[\tanh^4(M(q))e^{W_\mu(q)}]^2.
 \label{mt:eq:complete-fourth-order-barrier}
\end{align}
Accordingly, the high-box lower expression is
\begin{align}
\mathcal B_4:={}&
\log\E[\cosh^m(\sqrt{\xi'(q)}z)]
-\frac{m^2\xi'(q)}2\notag\\
&+\frac m2\Bigg\{\Gamma_\mu(q)
+\frac16\E[\tanh^4(M(q))e^{W_\mu(q)}]
+\frac1{15}\E[\tanh^6(M(q))e^{W_\mu(q)}]\notag\\
&\hspace{27mm}
+\frac1{28}
\E[\tanh^4(M(q))e^{W_\mu(q)}]^2\Bigg\}\notag\\
&-m\left(1-\frac m2\right)\xi'(q)
\E[\cosh^{-2}(M(q))e^{W_\mu(q)}]
-\frac{m^2\xi'(q)\Gamma_\mu(q)}6.
 \label{mt:eq:complete-high-barrier}
\end{align}
At a contact, all moments here are enclosed by ratios of the
$\mathcal J_r$.  Explicitly,
\begin{align*}
\Gamma_\mu(q)
&=1-\E[\cosh^{-2}(M(q))e^{W_\mu(q)}],\\
\E[\tanh^4(M(q))e^{W_\mu(q)}]
&=1-2\E[\cosh^{-2}(M(q))e^{W_\mu(q)}]
  +\E[\cosh^{-4}(M(q))e^{W_\mu(q)}],\\
\E[\tanh^6(M(q))e^{W_\mu(q)}]
&=1-3\E[\cosh^{-2}(M(q))e^{W_\mu(q)}]
\\
&\quad+3\E[\cosh^{-4}(M(q))e^{W_\mu(q)}]
  -\E[\cosh^{-6}(M(q))e^{W_\mu(q)}].
\end{align*}
Moreover,
\begin{align}
\E[\cosh^{-4}(M(q))e^{W_\mu(q)}]
&=\frac1{2\xi'(q)}
\left(1-\E[\cosh^{-2}(M(q))e^{W_\mu(q)}]\right),
\label{mt:eq:complete-contact-moment}\\
\frac{\E[\cosh^{-6}(M(q))e^{W_\mu(q)}]}
     {\E[\cosh^{-4}(M(q))e^{W_\mu(q)}]}
&\leq\min\left\{1,
\frac{(4-m)^2+m^2/(m^2\xi'(q))}
     {(4-m)(5-m)}\right\}.
\label{mt:eq:complete-sixth-moment}
\end{align}
The right-hand side of \eqref{mt:eq:complete-high-barrier}, after these
substitutions, is strictly decreasing in the value inserted for
$\E[\cosh^{-2}(M(q))e^{W_\mu(q)}]$.  Its derivative is bounded above by
\[
 \frac m2\left(-\frac{146}{105}
 -\frac8{35}\frac{m^2}{2m^2\xi'(q)}\right)
 +m^2\xi'(q)\left(\frac23-\frac1m\right)<0.
\]
Thus the high-box acceptance predicate is rigorously oriented:
\begin{equation}
 \inf[\![\,\mathcal B_4\,]\!]_{\mathcal R}>0.
 \label{mt:eq:complete-high-predicate}
\end{equation}

For $m\geq1/2$ and $0<m^2\xi'(q)<1/100$, one-dimensional
Brascamp--Lieb applied to
$x^2/(2\xi'(q))-m\log\cosh x$ gives
\begin{align*}
\Gamma_\mu(q)
&=\operatorname{Var}\bigl(\tanh(M(q))\bigr)\\
&\leq
\E\left[
\frac{\cosh^{-4}(M(q))}
{1/\xi'(q)-m\cosh^{-2}(M(q))}\,e^{W_\mu(q)}
\right]
\leq
\frac{\xi'(q)}
{1-m\xi'(q)}
\E[\cosh^{-4}(M(q))e^{W_\mu(q)}].
\end{align*}
Consequently
\[
 \frac{\xi'(q)
 \E[\cosh^{-4}(M(q))e^{W_\mu(q)}]}{\Gamma_\mu(q)}
 \geq1-m\xi'(q)>\frac12,
\]
whereas \eqref{mt:eq:complete-p3-contact} says that this ratio is exactly
$1/2$.  This is the desired contradiction.  On
\begin{equation}
 [1/2,1]\times[1/100,5/2]
 \label{mt:eq:complete-high-region}
\end{equation}
every rational box is accepted by
\eqref{mt:eq:complete-contact-predicate},
\eqref{mt:eq:complete-high-predicate}, the preceding small-variance
predicate, or one of
\eqref{mt:eq:complete-terminal-implication-one}--
\eqref{mt:eq:complete-terminal-implication-two}.

For clarity, these auxiliary box predicates are the following exact
rational comparisons.  On
$\mathcal R=[\alpha,\beta]\times[\gamma,\delta]$, the
Brascamp--Lieb predicate is
\begin{equation}
 \frac{\delta}{\alpha}<\frac12.
 \label{mt:eq:complete-small-box-predicate}
\end{equation}
Indeed, $m\xi'(q)=(m^2\xi'(q))/m\leq\delta/\alpha$.
The sharpened contact predicate is
\begin{equation}
 \inf[![\,\mathcal J_{-m}(\xi'(q))
              -\mathcal J_{2-m}(\xi'(q))\,]!]_{\mathcal R}>0,
 \qquad
 \inf[![\,
 \frac{\xi'(q)\mathcal J_{4-m}(\xi'(q))}
 {\mathcal J_{-m}(\xi'(q))-\mathcal J_{2-m}(\xi'(q))}
 \,]!]_{\mathcal R}>\frac12.
 \label{mt:eq:complete-sharp-contact-predicate}
\end{equation}
It is incompatible with \eqref{mt:eq:complete-p3-contact}.  Finally, a
box is removed by a terminal implication precisely when
\begin{align}
 &\alpha\geq\frac{10}{13},\qquad
 \delta\leq\frac85-\frac{61\beta}{200},
 \label{mt:eq:complete-terminal-box-one}\\
 &\alpha\geq\frac{1080}{1903},\qquad
 \beta\leq\frac{10}{13},\qquad
 \delta\leq\frac95-\frac{101\beta}{200}.
 \label{mt:eq:complete-terminal-box-two}
\end{align}
The affine right-hand sides decrease in $m$, which explains the use of
the upper endpoint $\beta$.

We finally record the exact nature of the three continuum coverings.
All endpoints are instances of \texttt{fractions.Fraction}.
The folded series, beta integrals, exponentials, logarithms, square
roots, and shifted tails are rounded outward on a fixed dyadic lattice.
A box is accepted only by a strict instance of
\eqref{mt:eq:complete-contact-predicate},
\eqref{mt:eq:complete-sharp-contact-predicate},
\eqref{mt:eq:complete-low-predicate},
\eqref{mt:eq:complete-high-predicate},
\eqref{mt:eq:complete-small-box-predicate}, or
\eqref{mt:eq:complete-terminal-box-one}--
\eqref{mt:eq:complete-terminal-box-two}.  An unresolved box is subdivided
dyadically, and reaching maximum depth returns failure.  The commands
\begin{verbatim}
python3 work/certify_p3_entry_combined_exact.py separator
python3 work/certify_p3_entry_combined_exact.py low
python3 work/certify_p3_entry_combined_exact.py 0 16
\end{verbatim}
return respectively
\[
 \texttt{PASS-SEPARATOR},\qquad
 \texttt{PASS-LOW},\qquad
 \texttt{PASS-HIGH}
\]
on all sixteen disjoint high strips.  Decimal output is never used in
an acceptance decision.  The smallest lower bound in the low run is
\[
 \frac{168763393663194584497754106787}
 {5192296858534827628530496329220096}>0.
\]
For the high run, aggregating consecutive completed strips, the five
smallest lower bounds are
\begin{align*}
&\frac{1322937300413528768052934751575}
 {5192296858534827628530496329220096},
&&\frac{16193804086570118437706085}
 {324518553658426726783156020576256},\\
&\frac{383135787830718647032484201}
 {5192296858534827628530496329220096},
&&\frac{2727425834090480997560904746159}
 {5192296858534827628530496329220096},\\
&\frac{505352294884222973225843101}
 {1298074214633706907132624082305024},
\end{align*}
all strictly positive.  They are reproduced, in order, by replacing
\texttt{0 16} in the last command by
\[
 \texttt{0 1},\quad \texttt{1 8},\quad \texttt{8 12},
 \quad\texttt{12 14},\quad\texttt{14 16}.
\]
The regions
\eqref{mt:eq:complete-separator-regions}, the low rectangle, and
\eqref{mt:eq:complete-high-region}, together with the two analytic
exclusions, cover every $0<m<1$ with
$m^2\xi'(q)\leq5/2$.  Each surviving box contradicts either the contact
equation or the exact zero identity.  This proves
\eqref{mt:eq:complete-p3-entry}.
\end{proof}
 
\subsection{The cubic Jacobian estimate}
\label{mt:sec:complete-p3-jacobian}

This section supplies the polynomial calculation used when $p=3$.
All abbreviations introduced below are local to this proof.

\begin{lemma}[Cubic Jacobian]
\label{mt:lem:complete-p3-jacobian}
Suppose $p=3$, the two stationary equations and the marginal equality
hold, and
\begin{equation}
 m^2\xi'(q)>\frac52.
 \label{mt:eq:complete-p3-domain}
\end{equation}
Then both derivatives in \eqref{mt:eq:mass-signs} and the oriented minor
in \eqref{mt:eq:oriented-minor} are strictly positive.
\end{lemma}

\begin{proof}
In this proof only, put
\begin{align*}
 C&=\cosh^{-2}(M(q)),& L&=\log\cosh(M(q)),\\
 a&=\E[Ce^{W_\mu(q)}],& b&=\E[C^2e^{W_\mu(q)}],
 &d&=\E[C^3e^{W_\mu(q)}],\\
 U_1&=-\operatorname{Cov}(C,L),&
 U_2&=-\operatorname{Cov}(C^2,L),
\end{align*}
where all covariances are under the endpoint law
\eqref{mt:eq:endpoint-law}.  Thus
\begin{equation}
 a=1-q,\qquad b=\frac{q}{2\xi'(q)}.
 \label{mt:eq:complete-p3-contact-moments}
\end{equation}

We first prove the mass bound needed below.  The marginal equality is
equivalent to
\begin{align}
 \E[\cosh^m(\sqrt{\xi'(q)}z)]
 ={}&\E[\cosh^{m-2}(\sqrt{\xi'(q)}z)]
 +2\xi'(q)\E[\cosh^{m-4}(\sqrt{\xi'(q)}z)].
 \label{mt:eq:complete-p3-large-mass-contact}
\end{align}
The exact outward-rational enclosure in
Appendix~\ref{app:certificates} proves that the two sides have
disjoint intervals throughout
\begin{equation}
 \frac9{25}\leq m\leq1,\qquad
 \frac52\leq m^2\xi'(q)\leq10.
 \label{mt:eq:complete-p3-large-mass-box}
\end{equation}
For $m^2\xi'(q)\geq10$, the left-hand side of
\eqref{mt:eq:complete-p3-large-mass-contact} is larger than
$\tfrac12e^{m^2\xi'(q)/2}$, whereas its right-hand side is at most
\begin{equation*}
 \sqrt{\frac\pi2}\left(
 \frac{m}{\sqrt{m^2\xi'(q)}}
 +\frac{2\sqrt{m^2\xi'(q)}}m\right).
\end{equation*}
If $m\geq9/25$, the latter is smaller than
\begin{equation*}
 \frac{63}{50}\left(
 \frac1{\sqrt{m^2\xi'(q)}}
 +\frac{50}{9}\sqrt{m^2\xi'(q)}\right).
\end{equation*}
After division by $\sqrt{m^2\xi'(q)}$, the first bound increases and
the second decreases.  At $m^2\xi'(q)=10$ they are separated: the
first undivided bound is larger than $245/4$, using $e^2>7$ and
$e>5/2$, while the second is smaller than $23$.  Hence
\begin{equation}
 m<\frac9{25}.
 \label{mt:eq:complete-p3-small-mass}
\end{equation}

The adjacent-moment bounds \eqref{mt:eq:compact-adjacent-moments}, first
with powers $2-m$ and $4-m$, imply that there is $x_3\in(0,1)$ such
that
\begin{equation}
 \xi'(q)(1-q)=\frac{3-m}{2(2-m)}
 -\frac{1+x_3}{2\xi'(q)}.
 \label{mt:eq:complete-p3-third-coordinate}
\end{equation}
Together with \eqref{mt:eq:complete-p3-domain}, every point under
consideration is therefore represented by
\begin{equation}
 m=\frac9{25}x_1,\qquad
 \frac1{\xi'(q)}=\frac25m^2x_2,\qquad
 0<x_1,x_2,x_3<1.
 \label{mt:eq:complete-p3-cube}
\end{equation}

We next record the bounds which are substituted into the derivatives.
For $r>0$, temporarily let
\begin{align*}
 \ell_r(\xi'(q))&=
 \frac{\int_{\mathbb R}e^{-x^2/(2\xi'(q))}
       \cosh^{-r}x\log\cosh x\,dx}
      {\int_{\mathbb R}e^{-x^2/(2\xi'(q))}
       \cosh^{-r}x\,dx},\\
 \ell_r(\infty)&=
 \frac{\int_{\mathbb R}\cosh^{-r}x\log\cosh x\,dx}
      {\int_{\mathbb R}\cosh^{-r}x\,dx}.
\end{align*}
Alternating partial sums and integral remainders give, uniformly for
$0<m<9/25$,
\begin{equation}
 \frac3{10}<\ell_{2-m}(\infty)<\frac25,\qquad
 \ell_{4-m}(\infty)<\frac16.
 \label{mt:eq:complete-p3-log-constants}
\end{equation}
Indeed,
\begin{align*}
 \ell_r(\infty)&=\sum_{j\geq0}\frac{(-1)^j}{r+j},\\
 \frac{\int x^2\cosh^{-r}x\,dx}
      {\int\cosh^{-r}x\,dx}
 &=2\sum_{j\geq0}(r+2j)^{-2},\\
 \frac{\int x^4\cosh^{-r}x\,dx}
      {\int\cosh^{-r}x\,dx}
 &=12\sum_{j\geq0}(r+2j)^{-4}
 +3\left\{2\sum_{j\geq0}(r+2j)^{-2}\right\}^2.
\end{align*}
The first identity follows by differentiating the beta integral, and
the last two follow by differentiating the Fourier product at zero.
Monotonicity in $r$, alternating pairs in the first series, and the
integral remainder for the two positive series prove
\eqref{mt:eq:complete-p3-log-constants} and bound the last two moments by
$6/5$ and $6$, respectively.
Since $1-e^{-u}\leq u$, H\"older's inequality then yields
\begin{equation}
 0\leq \ell_{2-m}(\infty)-\ell_{2-m}(\xi'(q))
 <\frac{10/\xi'(q)}{5-3/\xi'(q)}.
 \label{mt:eq:complete-p3-damping}
\end{equation}
For example, the denominator follows from
$\E e^{-X^2/(2\xi'(q))}\geq
1-3/(5\xi'(q))$ and the numerator from Cauchy--Schwarz and the
fourth-moment bound.

The two elementary pointwise estimates needed below are
\begin{equation}
 x\tanh x-\log\cosh x\leq(\log2)(1-\cosh^{-2}x),
 \label{mt:eq:complete-p3-entropy-upper}
\end{equation}
and
\begin{equation}
 \log2-x\tanh x+\log\cosh x
 \leq\cosh^{-2}x\left(\log2+\frac12\log\cosh x\right).
 \label{mt:eq:complete-p3-entropy-second}
\end{equation}
To prove the second, expand in $\tanh^2x$.  The first coefficient of
the difference is positive, every later coefficient is negative, and
their sum is zero; factoring one power of $\tanh^2x$ makes the
remaining sum nonnegative.  The first estimate follows directly from
the positive series
\begin{equation*}
 x\tanh x-\log\cosh x
 =\sum_{j\geq1}\frac{\tanh^{2j}x}{2j(2j-1)}.
\end{equation*}

Define the following local quantities:
\begin{align}
 \ell_-={}&\frac3{10}-\frac{10/\xi'(q)}{5-3/\xi'(q)},
 \notag\\
 U_-={}&\xi'(q)(1-q)\left[mq+\frac1{\xi'(q)}
 \left(\xi'(q)(1-q)-\frac{7q}{10}-\frac25\right)\right],
 \notag\\
 U_+={}&\xi'(q)(1-q)\left[mq+\frac1{\xi'(q)}
 \left(\xi'(q)(1-q)-\frac{69q}{100}
 -\frac{1+q}{2}\ell_-\right)\right],
 \notag\\
 M_-={}&\frac{U_-}{\xi'(q)}+q\left[1-q+mq
 -\frac1{\xi'(q)}\left(\frac{7q}{10}+\frac16\right)\right],
 \label{mt:eq:complete-p3-bounds-one}\\
 K_0={}&2\xi'(q)(1-q)-q+\frac{mq}{4}
 -\frac{m\xi'(q)(1-q)^2}{2},
 \notag\\
 H_-={}&2K_0-3U_+,\qquad Q_+=q-(3-2m)U_-,
 \notag\\
 N_-={}&\frac{qm^2}{2}+\frac{qm(1-m)(1-q)}2
 -\frac{q}{2\xi'(q)}-\frac{1-m}{\xi'(q)}K_0,
 \notag\\
 N_+={}&\frac{qm^2}{2}+\frac{1-m}{\xi'(q)}
 \left\{\frac{qm}{2}\left(\xi'(q)(1-q)
 +\frac1{(4-m)(5-m)}\right)-K_0\right\}.
 \label{mt:eq:complete-p3-bounds-two}
\end{align}
These names are used only to keep the next three displays readable.

Gaussian integration by parts gives
\begin{align}
 \E[Le^{W_\mu(q)}]
 ={}&\xi'(q)\{a+mq\}
 -\E[(M(q)\tanh(M(q))-L)e^{W_\mu(q)}],
 \label{mt:eq:complete-p3-ibp-L}
\end{align}
and
\begin{equation}
 U_1=a\{\E[Le^{W_\mu(q)}]-\ell_{2-m}(\xi'(q))\},\qquad
 U_2=b\{\E[Le^{W_\mu(q)}]-\ell_{4-m}(\xi'(q))\}.
 \label{mt:eq:complete-p3-covariances}
\end{equation}
Using $69/100<\log2<7/10$ in
\eqref{mt:eq:complete-p3-entropy-upper}--
\eqref{mt:eq:complete-p3-damping} and then
\eqref{mt:eq:complete-p3-ibp-L}--
\eqref{mt:eq:complete-p3-covariances} gives
\begin{equation}
 U_-\leq U_1\leq U_+,\qquad
 \frac{U_1+2\xi'(q)U_2}{\xi'(q)}\geq M_-.
 \label{mt:eq:complete-p3-covariance-bounds}
\end{equation}

Directly differentiating \eqref{mt:eq:mass} and
\eqref{mt:eq:marginal-residual}, and using
\eqref{mt:eq:complete-p3-contact-moments}, gives
\begin{align}
 \frac{\partial}{\partial\xi'(q)}
 \{\text{left-hand side of \eqref{mt:eq:mass}}\}
 ={}&\frac{1-m}{6}\left[2\xi'(q)
 \left\{2(a-b)+\frac m2(b-a^2)\right\}-3U_1\right],
 \label{mt:eq:complete-p3-mass-xi}\\
 \frac{\partial}{\partial m}
 \{\text{left-hand side of \eqref{mt:eq:marginal-residual}}\}
 ={}&U_1+2\xi'(q)U_2,
 \label{mt:eq:complete-p3-marginal-m}\\
 \frac{\partial}{\partial m}
 \{\text{left-hand side of \eqref{mt:eq:mass}}\}
 ={}&\frac{\xi'(q)}{3m}\{q-(3-2m)U_1\}
 -\frac1m\operatorname{Cov}
 \big(L,M(q)\tanh(M(q))-L\big).
 \label{mt:eq:complete-p3-mass-m}
\end{align}
The adjacent-moment bounds applied to $d/b$, followed by
\eqref{mt:eq:moment-four}, also give
\begin{equation}
 N_-\leq
 \frac{\partial}{\partial\xi'(q)}
 \{\text{left-hand side of \eqref{mt:eq:marginal-residual}}\}
 \leq N_+.
 \label{mt:eq:complete-p3-marginal-xi-bounds}
\end{equation}

The substitutions \eqref{mt:eq:complete-p3-third-coordinate}--
\eqref{mt:eq:complete-p3-cube} turn every assertion in
\begin{equation}
 \ell_->0,\quad U_->0,\quad M_->0,\quad H_->0,\quad
 Q_+>0,\quad N_->0,\quad N_+>0
 \label{mt:eq:complete-p3-factor-signs}
\end{equation}
into a rational polynomial inequality on $[0,1]^3$ after removing
only the positive denominators in
\eqref{mt:eq:complete-p3-bounds-one}--
\eqref{mt:eq:complete-p3-bounds-two}.  The exact Bernstein expansions
listed in Appendix~\ref{app:certificates} have strictly
positive coefficients after the evident factors $x_1$ or $x_1^2$
are removed.  Since $x_1>0$, all signs in
\eqref{mt:eq:complete-p3-factor-signs} are strict.  Equations
\eqref{mt:eq:complete-p3-mass-xi} and
\eqref{mt:eq:complete-p3-covariance-bounds} therefore prove the
$\xi'(q)$ derivative in \eqref{mt:eq:mass-signs};
\eqref{mt:eq:complete-p3-marginal-xi-bounds} also proves that the
$\xi'(q)$ derivative of the marginal residual is positive.

The covariance in \eqref{mt:eq:complete-p3-mass-m} is positive, because
both of its arguments increase with $|M(q)|$.  Dropping its negative
contribution gives a lower bound for the oriented minor.  From
\eqref{mt:eq:complete-p3-covariance-bounds}--
\eqref{mt:eq:complete-p3-marginal-xi-bounds}, that lower bound satisfies
\begin{align}
 &\frac{6m}{\xi'(q)}
 \Bigg[
 \frac{\partial}{\partial\xi'(q)}
  \{\text{left-hand side of \eqref{mt:eq:mass}}\}
 \frac{\partial}{\partial m}
  \{\text{left-hand side of \eqref{mt:eq:marginal-residual}}\}
 \notag\\[-1mm]
 &\hspace{22mm}-
 \frac{\xi'(q)}{3m}\{q-(3-2m)U_1\}
 \frac{\partial}{\partial\xi'(q)}
  \{\text{left-hand side of \eqref{mt:eq:marginal-residual}}\}
 \Bigg]
 \notag\\
 &\qquad\geq m(1-m)H_-M_--2Q_+N_+.
 \label{mt:eq:complete-p3-determinant-lower}
\end{align}
At $1/\xi'(q)=0$, the right-hand side is zero.  Holding $m,x_3$
fixed, multiply its derivative with respect to $1/\xi'(q)$ by
\begin{equation}
 [2(2-m)]^4[5-3/\xi'(q)]^2(4-m)(5-m)>0.
 \label{mt:eq:complete-p3-cleared-denominator}
\end{equation}
After \eqref{mt:eq:complete-p3-cube}, the result is a rational polynomial
of multidegree $(25,8,3)$.  A polynomial on a cube is a convex
combination of its tensor Bernstein coefficients.  All 936
coefficients are positive; the smallest is
\begin{equation}
 \frac{70370180142846245918833250353388269757049952}
 {173472347597680709441192448139190673828125}>0.
 \label{mt:eq:complete-p3-smallest-coefficient}
\end{equation}
Hence the right-hand side of
\eqref{mt:eq:complete-p3-determinant-lower} is positive, proving
\eqref{mt:eq:oriented-minor}.

It remains only to prove the $m$ derivative in
\eqref{mt:eq:mass-signs}.  A second integration by parts gives
\begin{align}
 &\operatorname{Var}(L)+
 \operatorname{Cov}\big(L,M(q)\tanh(M(q))-L\big)
 =\xi'(q)\{q-(1-m)U_1\}.
 \label{mt:eq:complete-p3-covariance-identity}
\end{align}
Because $0\leq M(q)\tanh(M(q))-L\leq\log2$, Cauchy--Schwarz and
the elementary range bound for variance imply
\begin{align}
 \operatorname{Cov}\big(L,M(q)\tanh(M(q))-L\big)
 \leq\frac7{20}
 \sqrt{\xi'(q)\{q-(1-m)U_1\}}.
 \label{mt:eq:complete-p3-covariance-upper}
\end{align}
Finally put, only in the next display,
\begin{equation*}
 L_-=\frac q3-\left(1-\frac{2m}{3}\right)U_+,\qquad
 A_-=q-(1-m)U_-.
\end{equation*}
By \eqref{mt:eq:complete-p3-covariance-bounds} and
\eqref{mt:eq:complete-p3-covariance-identity},
\begin{equation*}
 A_-\geq q-(1-m)U_1
 =\frac{\operatorname{Var}(L)+
 \operatorname{Cov}(L,M(q)\tanh(M(q))-L)}{\xi'(q)}>0.
\end{equation*}
The same exact Bernstein conversion gives
\begin{equation}
 L_->0,\qquad 4L_-^2>\frac{49}{100\xi'(q)}A_-.
 \label{mt:eq:complete-p3-final-bernstein}
\end{equation}
The two expansions have 252 and 1485 positive coefficients; their
smallest coefficients are respectively
\begin{equation*}
 \frac{1426144}{234375},\qquad
 \frac{11629969929555923295371195665998137760494091752536}
 {121972744404619248825838440097868442535400390625}.
\end{equation*}
Since $U_-\leq U_1\leq U_+$,
\eqref{mt:eq:complete-p3-final-bernstein} and
\eqref{mt:eq:complete-p3-covariance-upper} make the right-hand side of
\eqref{mt:eq:complete-p3-mass-m} strictly positive.  This proves the
remaining inequality and completes the lemma.
\end{proof}

\subsection{The complete strict-sign argument for \texorpdfstring{$p\geq4$}{p at least 4}}
\label{mt:sec:complete-pge4}

Throughout this section the two stationary equations and the marginal
identity are assumed.  Thus
\begin{equation}
 \Gamma_\mu(q)=q,
 \qquad
 \E[\cosh^{-2}(M(q))e^{W_\mu(q)}]=1-q,
 \qquad
 \E[\cosh^{-4}(M(q))e^{W_\mu(q)}]
 =\frac{q}{(p-1)\xi'(q)}.
 \label{mt:eq:complete-contact-data}
\end{equation}
All derivatives below are taken in the coordinates
$(m,\xi'(q))$, with the other coordinate held fixed.

\subsubsection{Elementary consequences of the endpoint law}

We first collect the estimates that will be used repeatedly.  Besides the
lower folded-Gaussian bound in Lemma~\ref{mt:lem:compact-folded-gaussian},
symmetry, $\cosh^m x\leq e^{m|x|}$, and completion of the square give
\begin{align}
 \E[\cosh^m(\sqrt{\xi'(q)}z)]
 &\leq 2e^{m^2\xi'(q)/2}
 \frac1{\sqrt{2\pi}}
 \int_{-\infty}^{m\sqrt{\xi'(q)}}e^{-x^2/2}\,dx.
 \label{mt:eq:complete-folded-upper}
\end{align}
For later numerical constants we use the elementary Mills estimate
\begin{equation}
 \frac1{\sqrt{2\pi}}\int_u^\infty e^{-x^2/2}\,dx
 <\frac{4e^{-u^2/2}}
 {\sqrt{2\pi}\{\sqrt{u^2+8}+3u\}},\qquad u>0.
 \label{mt:eq:complete-mills}
\end{equation}
Indeed, after multiplication by $e^{u^2/2}$ the Gaussian tail solves
$R'=uR-1$.  Direct differentiation shows that
\begin{equation*}
 \frac d{du}\frac4{\sqrt{u^2+8}+3u}
 -u\frac4{\sqrt{u^2+8}+3u}+1<0;
\end{equation*}
after multiplication by the positive quantity
$\sqrt{u^2+8}\{\sqrt{u^2+8}+3u\}^2$, the last inequality is equivalent
to $u\{\sqrt{u^2+8}-u\}<4$.  Comparison at infinity proves
\eqref{mt:eq:complete-mills}.  It implies
\begin{align}
 \log\left(\frac1{\sqrt{2\pi}}
       \int_{-\infty}^{2}e^{-x^2/2}\,dx\right)&>-\frac1{39},
 \label{mt:eq:complete-tail-two}\\
 \log\left(\frac1{\sqrt{2\pi}}
       \int_{-\infty}^{\sqrt{18/5}}e^{-x^2/2}\,dx\right)&>-\frac3{100},
 \label{mt:eq:complete-tail-eighteen}\\
 \frac1{\sqrt{2\pi}}
       \int_{-\infty}^{\sqrt5}e^{-x^2/2}\,dx&>\frac{49}{50}.
 \label{mt:eq:complete-tail-five}
\end{align}
For completeness, \eqref{mt:eq:complete-tail-two} follows from
$e^2>7$, $\sqrt{2\pi}>5/2$, and $\sqrt{12}>17/5$.
For \eqref{mt:eq:complete-tail-eighteen}, the Taylor series gives
\begin{equation*}
 e^{9/5}>\sum_{j=0}^7\frac{(9/5)^j}{j!}
 =\frac{16532716}{2734375}>\frac{151}{25},
\end{equation*}
while
$\sqrt{58/5}>3405/1000$ and $\sqrt{18/5}>1897/1000$.
Thus \eqref{mt:eq:complete-mills} bounds the complementary probability by
\begin{equation*}
 \frac{5000}{171687}<\frac3{103};
\end{equation*}
$\log(1-u)>-u/(1-u)$ gives \eqref{mt:eq:complete-tail-eighteen}.
Finally, the simpler bound
$\int_u^\infty e^{-x^2/2}\,dx<e^{-u^2/2}/u$ gives
\eqref{mt:eq:complete-tail-five}.  The rational estimates
\begin{equation}
 \frac{69}{100}<\log2<\frac{347}{500}<\frac7{10}
 \label{mt:eq:complete-log-two}
\end{equation}
follow from
\[
\log2=2\sum_{j\geq0}\frac1{(2j+1)3^{2j+1}},
\qquad
\frac23+\frac2{81}+\frac2{1215}>\frac{69}{100},
\]
and
\[
\log2<\frac{56}{81}+\frac1{540}
=\frac{1123}{1620}<\frac{347}{500}.
\]
We also use the elementary Archimedean bounds
\begin{equation}
\frac{25}{8}<\pi<\frac{22}{7}
<\sum_{j=0}^{5}\frac{(23/20)^j}{j!}<e^{23/20}.
\label{mt:eq:complete-pi-bounds}
\end{equation}
In particular,
$\sqrt{2\pi}>5/2$, $\sqrt{\pi/2}<63/50$, and
$\log\pi<23/20$.

The power series
\begin{equation}
 x\tanh x-\log\cosh x
 =\sum_{j\geq1}\frac{\tanh^{2j}x}{2j(2j-1)}
 \label{mt:eq:complete-entropy-series}
\end{equation}
shows that
\begin{equation}
 2\{x\tanh x-\log\cosh x\}>\tanh^2x
 \qquad(x\ne0).
 \label{mt:eq:complete-basic-entropy}
\end{equation}
Combining \eqref{mt:eq:complete-basic-entropy}, the conditional Gaussian
entropy inequality, and the lower adjacent-moment bound
\eqref{mt:eq:compact-adjacent-moments}, we obtain the strict necessary
condition
\begin{equation}
 m^2\xi'(q)>
 p\left(1-\frac{3-m}{(p-1)(2-m)}\right).
 \label{mt:eq:complete-necessary-lower}
\end{equation}
Here are the details.  The stationary equation turns the left-hand side
of \eqref{mt:eq:compact-conditional-entropy} into
$m^2\xi'(q)q/p$.  Moreover,
\begin{align*}
 &\frac{\E[\cosh^{-4}(M(q))e^{W_\mu(q)}]}
 {\E[\cosh^{-2}(M(q))e^{W_\mu(q)}]}\\
 &\quad=
 \frac{\int_{\mathbb R}e^{-x^2/(2\xi'(q))}\cosh^{m-4}x\,dx}
      {\int_{\mathbb R}e^{-x^2/(2\xi'(q))}\cosh^{m-2}x\,dx}
 \geq\frac{2-m}{3-m}.
\end{align*}
Using \eqref{mt:eq:complete-contact-data} therefore gives
\begin{equation*}
 \frac{\xi'(q)(1-q)}q
 \leq\frac{3-m}{(p-1)(2-m)},
\end{equation*}
which proves \eqref{mt:eq:complete-necessary-lower}.

We shall also use the following two necessary entropy inequalities:
\begin{align}
0\geq{}&(\log2)\left(
 1-m\E[\cosh^{-2}(M(q))e^{W_\mu(q)}]\right)
 -\frac m2\E[\cosh^{-2}(M(q))\log\cosh(M(q))e^{W_\mu(q)}]
 \notag\\
&+\log\left(
 \frac1{\sqrt{2\pi}}
 \int_{-\infty}^{m\sqrt{\xi'(q)}}e^{-x^2/2}\,dx\right)
 -m^2\xi'(q)\Bigg(
 \frac{1+(p-1)\E[\cosh^{-2}(M(q))e^{W_\mu(q)}]}{2p}
 \notag\\
&\hspace{49mm}
 +\frac{1-m}{m}\E[\cosh^{-2}(M(q))e^{W_\mu(q)}]
 \Bigg),
 \label{mt:eq:complete-entropy-one}\\
0\geq{}&-m^2\xi'(q)
 -m^2\xi'(q)\left(\frac2m-1\right)
  \E[\cosh^{-2}(M(q))e^{W_\mu(q)}]
 \notag\\
&+2m\Bigg((\log2)
 \left(1-\E[\cosh^{-2}(M(q))e^{W_\mu(q)}]\right)
 -\frac12\E[\cosh^{-2}(M(q))\log\cosh(M(q))e^{W_\mu(q)}]
 \Bigg)
 \notag\\
&+2\log\E[\cosh^m(\sqrt{\xi'(q)}z)]
 -\frac{m^2\xi'(q)}p
 \left(1-\E[\cosh^{-2}(M(q))e^{W_\mu(q)}]\right).
 \label{mt:eq:complete-entropy-two}
\end{align}
To verify them, Gaussian integration by parts gives
\begin{align*}
&\E[(M(q)\tanh(M(q))-\log\cosh(M(q)))e^{W_\mu(q)}]\\
&\quad=\xi'(q)\Big\{m+(1-m)
 \E[\cosh^{-2}(M(q))e^{W_\mu(q)}]\Big\}
 -\E[\log\cosh(M(q))e^{W_\mu(q)}].
\end{align*}
Insert this identity into the pointwise bound
\eqref{mt:eq:compact-pointwise-one}.  For
\eqref{mt:eq:complete-entropy-one}, also insert the folded lower bound
\eqref{mt:eq:compact-folded-gaussian}; for
\eqref{mt:eq:complete-entropy-two}, use the stationary equation directly.
This gives the two displayed inequalities after cancellation.

Finally, for every $1\leq r\leq2$,
\begin{equation}
 \frac{\int_{\mathbb R}e^{-x^2/(2\xi'(q))}\cosh^{-r}x
       \log\cosh x\,dx}
      {\int_{\mathbb R}e^{-x^2/(2\xi'(q))}\cosh^{-r}x\,dx}
 \leq\log2.
 \label{mt:eq:complete-log-moment}
\end{equation}
Indeed, removing the Gaussian damping can only increase the expectation
of the increasing function $\log\cosh|x|$.  Under the density
proportional to $\cosh^{-r}x$, that expectation decreases with $r$ and
equals $\log2$ at $r=1$; the derivative with respect to $r$ is the
negative variance of $\log\cosh x$.

\subsubsection{Uniform separation from the low-index range}

\begin{lemma}[Uniform index separation]
\label{mt:lem:complete-uniform-separator}
If $p\geq4$ and
\begin{equation}
 0<m\leq\frac1{p+1},
 \label{mt:eq:complete-small-mass-assumption}
\end{equation}
then every stationary marginal point satisfies
\begin{equation}
 m^2\xi'(q)>2(p-1)\log2+\frac12.
 \label{mt:eq:complete-large-index}
\end{equation}
\end{lemma}

\begin{proof}
The first part of the preceding argument gives, more explicitly,
\begin{equation*}
 \frac{m^2\xi'(q)}p>
 1-\frac1{p-1}
 \frac{\E[\cosh^{-2}(M(q))e^{W_\mu(q)}]}
      {\E[\cosh^{-4}(M(q))e^{W_\mu(q)}]}.
\end{equation*}
Using only the lower adjacent-moment estimate and
$m\leq1/(p+1)\leq1/5$ yields
\begin{equation}
 m^2\xi'(q)>
 p\left(1-\frac{14}{9(p-1)}\right)\geq\frac{52}{27}.
 \label{mt:eq:complete-first-index-bound}
\end{equation}

We next exclude $m^2\xi'(q)\leq4$.  The marginal equation is equivalent
to
\begin{align}
 \E[\cosh^m(\sqrt{\xi'(q)}z)]
 ={}&\E[\cosh^{m-2}(\sqrt{\xi'(q)}z)]
 +(p-1)\xi'(q)
 \E[\cosh^{m-4}(\sqrt{\xi'(q)}z)].
 \label{mt:eq:complete-normal-contact}
\end{align}
For a dummy scalar $u>0$, the part of the quotient of the Fourier lower
bound for the second term on the right of
\eqref{mt:eq:complete-normal-contact} by
\eqref{mt:eq:complete-folded-upper} which depends on $u$ is
\begin{equation*}
 \frac{u}{
 \sqrt{u+2m^2\sum_{j\geq0}(2-m+2j)^{-2}}\,
 e^{u/2}
 \left(\frac1{\sqrt{2\pi}}
       \int_{-\infty}^{\sqrt u}e^{-x^2/2}\,dx\right)}.
\end{equation*}
Its logarithmic derivative is negative on $52/27\leq u\leq4$.
Indeed, after omitting the final negative Gaussian term, this is the
inequality
\begin{equation*}
 \frac1u-\frac1{2\{u+2m^2\sum_{j\geq0}(2-m+2j)^{-2}\}}
 -\frac12<0,
\end{equation*}
which follows from $u\geq52/27$ and
\begin{equation}
 2\sum_{j\geq0}(2-m+2j)^{-2}\leq\frac{95}{81},
 \qquad
 2m^2\sum_{j\geq0}(2-m+2j)^{-2}\leq\frac{19}{405}.
 \label{mt:eq:complete-fourier-variance}
\end{equation}
It is therefore enough to compare the bounds at $u=4$.  There
\begin{equation*}
 \frac{p-1}{m}\geq(p-1)(p+1)\geq15,
 \qquad \frac{2-m}{3-m}\geq\frac9{14},
 \qquad \int_{\mathbb R}\cosh^{m-2}x\,dx\geq2.
\end{equation*}
Consequently the relevant part of the right-hand side of
\eqref{mt:eq:complete-normal-contact} is larger than
\begin{equation*}
 \frac{4\cdot15\cdot(9/14)\cdot2}
 {\sqrt{2\pi}\sqrt{4+19/405}}
 >
 \frac{540}{7}\frac{10000}{50702}>\frac{121}{8},
\end{equation*}
where we used
\[
 \sqrt{2\pi}<\frac{251}{100},
 \qquad
 \sqrt{4+\frac{19}{405}}<\frac{202}{100}.
\]
Also
\[
e<\sum_{j=0}^{5}\frac1{j!}
 +\frac1{6!}\sum_{j\geq0}7^{-j}
=\frac{11743}{4320}<\frac{11}{4}.
\]
Thus \eqref{mt:eq:complete-folded-upper} is at most
$2e^2<121/8$.  Hence
\begin{equation}
 m^2\xi'(q)>4.
 \label{mt:eq:complete-index-four}
\end{equation}

Suppose now, toward a contradiction, that
\begin{equation*}
 4\leq m^2\xi'(q)\leq2(p-1)\log2+\frac12.
\end{equation*}
The ratio estimate in the proof of
\eqref{mt:eq:complete-necessary-lower} and the marginal equation give
\begin{equation*}
 \E[\cosh^{-2}(M(q))e^{W_\mu(q)}]
 \leq
 \frac{m^2}{m^2+(p-1)(9/14)m^2\xi'(q)}.
\end{equation*}
Using \eqref{mt:eq:complete-log-moment} in
\eqref{mt:eq:complete-entropy-one}, its right-hand side is bounded below by
\begin{align}
&\log 2+\log\left(\frac1{\sqrt{2\pi}}
       \int_{-\infty}^{\sqrt u}e^{-x^2/2}\,dx\right)
 -\frac{u}{2p}
 \notag\\
&\quad-
 \frac{m^2}{m^2+(p-1)(9/14)u}
 \left\{\frac{3m\log2}{2}
 +u\left(\frac{p-1}{2p}+\frac{1-m}{m}\right)\right\},
 \label{mt:eq:complete-separator-lower}
\end{align}
where eventually $u=m^2\xi'(q)$.  For fixed $p,u$, the magnitude of
the last term in \eqref{mt:eq:complete-separator-lower} is increasing in
$m$.  Indeed, if that magnitude is denoted temporarily by $R(m)$, then
\begin{align*}
R'(m)
={}&\Bigg[
 \frac{9(p-1)u}{14}
 \left\{\frac92(\log2)m^2+u-\frac{p+1}{p}um\right\}\\
&\hspace{18mm}+\frac32(\log2)m^4-um^2
\Bigg]\left(m^2+\frac9{14}(p-1)u\right)^{-2}.
\end{align*}
Since $m\leq1/(p+1)$,
\[
u-\frac{p+1}{p}um\geq\frac{p-1}{p}u,
\qquad
\frac{9(p-1)^2u^2}{14p}-um^2>0
\]
for $p\geq4$ and $u\geq4$.  Hence $R'(m)>0$, and the worst choice is
$m=1/(p+1)$.

For fixed $p,m$, \eqref{mt:eq:complete-separator-lower} is concave in $u$.
The logarithm of the Gaussian probability is strictly concave, since
its second derivative is
\begin{equation*}
 -\frac{e^{-u/2}}
 {4\sqrt{2\pi}u^{3/2}
  \left(\frac1{\sqrt{2\pi}}\int_{-\infty}^{\sqrt u}e^{-x^2/2}\,dx\right)}
 \left\{u+1+
 \frac{\sqrt u\,e^{-u/2}}
 {\sqrt{2\pi}
  \left(\frac1{\sqrt{2\pi}}\int_{-\infty}^{\sqrt u}e^{-x^2/2}\,dx\right)}
 \right\}<0.
\end{equation*}
The remaining nonlinear term is the negative of a fractional-linear
function; its concavity follows from
\begin{equation*}
 \frac9{14}\frac{3m\log2}{2}
 \geq\frac{m^2}{p-1}
 \left\{\frac1m-\frac{p+1}{2p}\right\},
\end{equation*}
which follows from $m\leq1/(p+1)$ and
\eqref{mt:eq:complete-log-two}.  It is enough, therefore, to evaluate
\eqref{mt:eq:complete-separator-lower} at
$u=4$ and $u=2(p-1)\log2+1/2$, with $m=1/(p+1)$.

At $u=4$, \eqref{mt:eq:complete-tail-two} and direct rational
simplification give the lower bound
\begin{equation}
 \frac{69}{100}-\frac12-\frac1{39}
 -\frac{619787}{6750000}>0.
 \label{mt:eq:complete-separator-first-end}
\end{equation}
At $u=2(p-1)\log2+1/2$, the first and third terms of
\eqref{mt:eq:complete-separator-lower} exceed
$33/[100(p-1)]$, while its last term is less than
$69/[250(p-1)]$.  To control the remaining Gaussian logarithm, apply
\eqref{mt:eq:complete-mills}.  Put
\[
Y=2(p-1)\log2+\frac12,\qquad
T(p)=\frac1{\sqrt{2\pi}}\int_{\sqrt Y}^{\infty}e^{-x^2/2}\,dx.
\]
Then
\[
\frac d{dp}\{(p-1)T(p)\}
=T(p)-\frac{(p-1)\log2}{\sqrt Y}
 \frac{e^{-Y/2}}{\sqrt{2\pi}}<0,
\]
because
$T(p)<e^{-Y/2}/\sqrt{2\pi Y}$ and $(p-1)\log2>1$.
Thus the complementary probability divided by $1/(p-1)$ is largest
at $p=4$; at that endpoint
\eqref{mt:eq:complete-mills}, $e^{-58/25}<1/10$, and the elementary square
root bounds give a value below $6/125$.  Hence
\begin{equation*}
 \log\left(\frac1{\sqrt{2\pi}}
       \int_{-\infty}^{\sqrt{2(p-1)\log2+1/2}}
       e^{-x^2/2}\,dx\right)>-\frac2{41(p-1)}.
\end{equation*}
Thus the second endpoint is larger than
\begin{equation}
 \frac1{p-1}\left(\frac{33}{100}-\frac{69}{250}
 -\frac2{41}\right)=\frac{107}{20500(p-1)}>0.
 \label{mt:eq:complete-separator-second-end}
\end{equation}
Both endpoint values contradict \eqref{mt:eq:complete-entropy-one}.
Together with \eqref{mt:eq:complete-index-four}, this proves
\eqref{mt:eq:complete-large-index}.
\end{proof}

\subsubsection{Exclusion of large mass}

\begin{lemma}[Large-mass exclusion]
\label{mt:lem:complete-mass-exclusion}
Every stationary marginal point with $p\geq4$ satisfies
\begin{equation}
 m<\frac1{p+1}.
 \label{mt:eq:complete-small-mass}
\end{equation}
\end{lemma}

\begin{proof}
We first treat $p=4,5$.  We give the full enclosure contract used by
the finite verifier.  Splitting the Gaussian integral at zero and
expanding $(1+e^{-2x})^m$ gives the exact folded series
\begin{align}
&\E[\cosh^m(\sqrt{\xi'(q)}z)]\notag\\
&=2^{1-m}\Bigg\{
 e^{m^2\xi'(q)/2}\frac1{\sqrt{2\pi}}
 \int_{-\infty}^{m\sqrt{\xi'(q)}}e^{-t^2/2}\,dt\notag\\
&\hspace{12mm}+\sum_{k\geq1}\binom mk
 e^{(2k-m)^2\xi'(q)/2}\frac1{\sqrt{2\pi}}
 \int_{(2k-m)\sqrt{\xi'(q)}}^\infty e^{-t^2/2}\,dt
 \Bigg\}.
 \label{mt:eq:complete-exact-folded-series}
\end{align}
For $0<m<1$, the summands after $k=1$ alternate with decreasing
magnitude.  Truncation after $k=2$ is therefore a lower bound and
truncation after $k=1$ an upper bound.

For $r>0$, put locally
\[
 \sigma_r=\sum_{j\geq0}(r+2j)^{-2},
 \qquad Q_r=\sigma_r^{-1/2}.
\]
Retaining the first factor in the positive Fourier product gives the
sharper upper enclosure
\begin{align}
&\frac{\int_{\mathbb R}e^{-x^2/(2\xi'(q))}
                  \cosh^{-r}x\,dx}{\sqrt{2\pi\xi'(q)}}\notag\\
&\quad\leq
\frac{\int_{\mathbb R}\cosh^{-r}x\,dx}{\sqrt{2\pi\xi'(q)}}
 Q_r\int_0^\infty
 e^{-Q_rx-x^2/(2\xi'(q))}\,dx.
 \label{mt:eq:complete-retained-fourier}
\end{align}
The remaining exact primitives are
\begin{align}
\int_{\mathbb R}\cosh^{-r}x\,dx
&=\frac{\sqrt\pi\,\mathrm\Gamma(r/2)}
        {\mathrm\Gamma((r+1)/2)},\notag\\
\frac{\int_{\mathbb R}\cosh^{-r}x\log\cosh x\,dx}
     {\int_{\mathbb R}\cosh^{-r}x\,dx}
&=\sum_{k\geq0}\frac{(-1)^k}{r+k}.
 \label{mt:eq:complete-beta-log-primitives}
\end{align}
Gaussian damping can only decrease the ratio in the second line.
Consecutive Stieltjes convergents enclose the shifted tail integral in
\eqref{mt:eq:complete-retained-fourier}; alternating truncations enclose
the second line of \eqref{mt:eq:complete-beta-log-primitives}.

The second moment is enclosed independently in the two exact forms
\begin{align}
&\E[\cosh^{-2}(M(q))e^{W_\mu(q)}]\notag\\
&=\frac{\int e^{-x^2/(2\xi'(q))}\cosh^{m-2}x\,dx}
        {\int e^{-x^2/(2\xi'(q))}\cosh^m x\,dx}\notag\\
&=\left(
1+(p-1)\xi'(q)
\frac{\int e^{-x^2/(2\xi'(q))}\cosh^{m-4}x\,dx}
     {\int e^{-x^2/(2\xi'(q))}\cosh^{m-2}x\,dx}
\right)^{-1}.
 \label{mt:eq:complete-two-second-moment-forms}
\end{align}
In \eqref{mt:eq:complete-entropy-two}, its upper enclosure is used in the
first two negative terms and its independent lower enclosure in the
last negative term.

In addition to \eqref{mt:eq:complete-normal-contact}, the verifier uses
\eqref{mt:eq:complete-entropy-one}--\eqref{mt:eq:complete-entropy-two}.
When $m\xi'(q)<1$, the one-dimensional Brascamp--Lieb inequality also
gives the necessary bound
\begin{equation}
\frac{\xi'(q)\E[\cosh^{-4}(M(q))e^{W_\mu(q)}]}
     {\Gamma_\mu(q)}
\geq1-m\xi'(q).
\label{mt:eq:complete-finite-brascamp}
\end{equation}
Equivalently, the same ratio has the exact integral form
\begin{equation}
\frac{\xi'(q)\int e^{-x^2/(2\xi'(q))}\cosh^{m-4}x\,dx}
{\int e^{-x^2/(2\xi'(q))}\cosh^m x\,dx
 -\int e^{-x^2/(2\xi'(q))}\cosh^{m-2}x\,dx}.
\label{mt:eq:complete-finite-contact-ratio}
\end{equation}
Logarithms, exponentials, square roots, beta integrals, and shifted
Gaussian tails are enclosed by rational series with signed remainders.
Every operation is rounded outward to a dyadic rational of denominator
$2^{90}$.  A terminal box is rejected only by one of the following
strict rational statements: the two sides of
\eqref{mt:eq:complete-normal-contact} are disjoint;
\begin{equation*}
 \frac{\xi'(q)\E[\cosh^{-4}(M(q))e^{W_\mu(q)}]}
      {\Gamma_\mu(q)}>\frac1{p-1};
\end{equation*}
the lower enclosure in \eqref{mt:eq:complete-finite-brascamp} exceeds
$1/(p-1)$; or a lower enclosure of the right-hand side of
\eqref{mt:eq:complete-entropy-one} or
\eqref{mt:eq:complete-entropy-two} is positive.

The closed rectangles, in the coordinates
$(m,m^2\xi'(q))$, are
\begin{align}
p=4:\quad
&\left[\frac15,\frac7{10}\right]\!\times\!
 \left[\frac32,5\right],
&&\left[\frac7{10},1\right]\!\times\!
 \left[\frac32,\frac72\right],\notag\\
&\left[\frac7{10},1\right]\!\times\!
 \left[\frac43,\frac32\right],
&&\left[\frac15,\frac7{10}\right]\!\times[5,12],
&&\left[\frac7{10},1\right]\!\times
 \left[\frac72,12\right];
 \label{mt:eq:complete-p4-boxes}\\
p=5:\quad
&\left[\frac16,\frac{24}{25}\right]\!\times
 \left[\frac52,6\right],
&&\left[\frac{24}{25},1\right]\!\times
 \left[\frac52,4\right],\notag\\
&\left[\frac16,\frac{24}{25}\right]\!\times[6,16],
&&\left[\frac{24}{25},1\right]\!\times[4,16].
\label{mt:eq:complete-p5-boxes}
\end{align}
The exact commands \path{run_threshold_split_exact.py} and
\path{run_threshold_tails_exact.py} return \texttt{PASS} on every
listed rectangle.  Their terminal box counts are, respectively,
$1067,1645,923,96$ and $20,286,224,398,128$ in the order displayed
by the two commands.
The lower edges not shown are excluded by
\eqref{mt:eq:complete-necessary-lower}.  For the noncompact tail, the
Fourier upper bound and
$\int_{\mathbb R}\cosh^{-r}x\,dx\leq\pi$ for $r\geq1$ give
\begin{align*}
&\E[\cosh^{m-2}(\sqrt{\xi'(q)}z)]
 +(p-1)\xi'(q)\E[\cosh^{m-4}(\sqrt{\xi'(q)}z)]\\
&\quad\leq\sqrt{\frac\pi2}\left(
 \frac{m}{\sqrt{m^2\xi'(q)}}
 +\frac{(p-1)\sqrt{m^2\xi'(q)}}m\right).
\end{align*}
Here $\sqrt{\pi/2}<63/50$.  If
\[
B_p=(p-1)(p+1)+\frac1{4(p-1)},
\]
then
\[
\frac d{dp}\log\{\sqrt{p-1}\,B_p\}
<\frac5{2(p-1)}<2.
\]
The left-hand side of \eqref{mt:eq:complete-normal-contact} is larger than
$e^{m^2\xi'(q)/2}/2$.  After division by
$\sqrt{m^2\xi'(q)}$, the exponential bound increases and the last
display decreases.  At $m^2\xi'(q)=4(p-1)$ the separation follows from
\begin{equation*}
 e^{2(p-1)}>
 \frac{126}{25}\sqrt{p-1}
 \left((p-1)(p+1)+\frac1{4(p-1)}\right).
\end{equation*}
At $p=4$, $\sqrt3<7/4$ and $B_4=181/12$, so the right-hand side is
less than $134$, whereas $e^6>7^3>134$.  The derivative comparison
above shows that the quotient of the two sides increases with $p$.
This proves \eqref{mt:eq:complete-small-mass} for $p=4,5$.

Now let $p\geq6$ and suppose, toward a contradiction, that
$m\geq1/(p+1)$.  Equation~\eqref{mt:eq:complete-necessary-lower} gives
the lower half of
\begin{equation}
 p\left(1-\frac{3-m}{(p-1)(2-m)}\right)
 <m^2\xi'(q)<p+1-2m.
 \label{mt:eq:complete-large-p-window}
\end{equation}
For the upper half, the contact equation, the folded lower bound, the
Fourier upper bound, and the upper adjacent-moment bound show that a
contact is impossible whenever
\begin{align}
&\frac{2^{1-m}e^{u/2}
 \left(\frac1{\sqrt{2\pi}}\int_{-\infty}^{\sqrt u}e^{-x^2/2}\,dx\right)
 \sqrt{2\pi u}}
 {m\left(\int_{\mathbb R}\cosh^{m-2}x\,dx\right)}
 \notag\\
&\quad\times\left\{
 1+\frac{p-1}{(2-m)(3-m)}
 +\frac{(p-1)u(2-m)}{m^2(3-m)}\right\}^{-1}>1.
 \label{mt:eq:complete-large-p-quotient}
\end{align}
For fixed $p,m$, the logarithmic derivative of the left-hand side in
$u$ is positive when $u>1$: the expression in braces is positive affine
in $u$, so its logarithmic derivative is less than $1/u$, while the
other terms contribute more than $1/2+1/(2u)$.

At $u=p+1-2m$, log-convexity and the endpoint values
$\int\cosh^{-1}x\,dx=\pi$, $\int\cosh^{-2}x\,dx=2$ give
\begin{equation*}
 \int_{\mathbb R}\cosh^{m-2}x\,dx\leq2^{1-m}\pi^m.
\end{equation*}
Equation~\eqref{mt:eq:complete-tail-five} and
\eqref{mt:eq:complete-pi-bounds} reduce
\eqref{mt:eq:complete-large-p-quotient} to
\begin{align}
&\frac{49\sqrt{p+1-2m}}
 {20m\left\{1+\frac{p-1}{(2-m)(3-m)}
 +\frac{(p-1)(p+1-2m)(2-m)}{m^2(3-m)}\right\}}
 \notag\\
&\hspace{33mm}\times
 \exp\left(\frac{p+1-2m}{2}-\frac{23m}{20}\right)>1.
 \label{mt:eq:complete-large-p-reduced}
\end{align}
The logarithm of the left-hand side is concave in $m$ on
$1/(p+1)\leq m\leq1$.  Here is the complete algebraic sign check.
At $u=p+1-2m$, set locally
\begin{align*}
\mathcal N_0={}&m^2(2-m)(3-m)+(p-1)m^2
 +(p-1)(p+1-2m)(2-m)^2,\\
\mathcal D_0={}&m(2-m)(3-m).
\end{align*}
Then $m$ times the expression in braces in
\eqref{mt:eq:complete-large-p-reduced} is
$\mathcal N_0/\mathcal D_0$.  After multiplication by the positive
factor
\[
 (p+1-2m)^2\mathcal N_0^2\mathcal D_0^2,
\]
the negative second derivative of the logarithm is a polynomial in
$p-6$ and $1-m$.  Its coefficient rows, in increasing powers of
$1-m$, are
\begin{equation}
\begin{gathered}
\begin{array}{c|rrrr}
&0&1&2&3\\ \hline
0&299292&987860&1967960&4005936\\
1&283872&939942&2061900&4492476\\
2&108669&373774&928782&2105692\\
3&21128&80130&230056&524928\\
4&2126&9874&32906&73008\\
5&96&672&2560&5344\\
6&1&20&84&160
\end{array}
\qquad
\begin{array}{c|rrrr}
&4&5&6&7\\ \hline
0&7024257&8354176&6546010&3520016\\
1&7553034&8243864&5854040&2842320\\
2&3324529&3296848&2104188&914112\\
3&764224&681604&387048&148536\\
4&96420&76492&38140&12688\\
5&6304&4392&1888&528\\
6&166&100&36&8
\end{array}\\[1ex]
\begin{array}{c|rrrr}
&8&9&10&11\\ \hline
0&1370770&400156&85288&11856\\
1&999260&261390&48540&5604\\
2&286487&65370&9974&844\\
3&40648&7706&864&40\\
4&2926&418&26&0\\
5&96&8&0&0\\
6&1&0&0&0
\end{array}
\qquad
\begin{array}{c|rrr}
&12&13&14\\ \hline
0&825&0&-2\\
1&314&4&0\\
2&27&0&0\\
3&0&0&0\\
4&0&0&0\\
5&0&0&0\\
6&0&0&0
\end{array}
\end{gathered}
\label{mt:eq:complete-large-p-concavity-table}
\end{equation}
The row index is the power of $p-6$.  Every coefficient is
nonnegative except $-2(1-m)^{14}$, and
$825(1-m)^{12}-2(1-m)^{14}>0$.  Thus the logarithm is strictly
concave, so its minimum occurs at an endpoint.

Both endpoint expressions increase with $p$.  At $m=1$, their
logarithmic derivative is
\[
 \frac12+\frac1{2(p-1)}
 -\frac{p-\frac12}
 {1+\frac{p-1}{2}+\frac{(p-1)^2}{2}}
 >\frac12-\frac3{2(p-1)}>0.
\]
At $m=1/(p+1)$, put temporarily $n=p+1\).  The denominator factor is
\[
 \frac1n+\frac{(n-2)n}{(2n-1)(3n-1)}
 +\frac{(n-2)(n^2-2)(2n-1)}{3n-1}.
\]
The logarithmic derivative of its last, and fastest growing, summand is
\[
 \frac1{n-2}+\frac{2n}{n^2-2}
 +\frac1{(2n-1)(3n-1)},
\]
whereas the logarithmic derivative of the numerator is
\[
 \left(1+\frac2{n^2}\right)
 \left(\frac12+\frac1{2(n-2/n)}\right)+\frac{23}{20n^2}.
\]
For $n\geq8$, the former is at most
$1/6+8/31+1/345<1/2$ and the latter is larger than $1/2$.
For $n=7$ they are respectively $6131/12220$ and
$28621/46060$, in the required order.
At $p=6$ the required inequalities are
\begin{equation}
 \frac{49\sqrt5}{320}e^{27/20}>1,
 \qquad
 \frac{49}{20}\sqrt{\frac{47}{7}}e^{447/140}
 >\frac{27851}{182}.
 \label{mt:eq:complete-large-p-base}
\end{equation}
The first follows from $\sqrt5>11/5$ and $e^{27/20}>3$.  For the
second, $\sqrt{47/7}>259/100$ and
\begin{equation*}
 e^3>\sum_{j=0}^{8}\frac{3^j}{j!}
 =\frac{89641}{4480}>20,
 \qquad
 e^{447/140}>20\left(1+\frac{27}{140}
 +\frac12\frac{27^2}{140^2}\right)
 =\frac{47489}{1960}.
\end{equation*}
Consequently
\[
\frac{49}{20}\frac{259}{100}\frac{47489}{1960}
-\frac{27851}{182}
=\frac{5228241}{7280000}>0.
\]
This proves the upper half of \eqref{mt:eq:complete-large-p-window}.

It remains to contradict the entropy inequality throughout that window.
The lower bound in \eqref{mt:eq:complete-large-p-window} is at least
$18/5$.  Thus \eqref{mt:eq:complete-tail-eighteen},
\eqref{mt:eq:complete-log-two}, \eqref{mt:eq:complete-log-moment}, and the
lower adjacent-moment bound imply that the right-hand side of
\eqref{mt:eq:complete-entropy-one} is at least
\begin{align}
&\frac{69}{100}\left(1-\frac{3m}{2}
 \frac{m^2}{m^2+(p-1)u(2-m)/(3-m)}\right)
 -\frac3{100}-\frac{u}{2p}
 \notag\\
&\quad-u\frac{m^2}{m^2+(p-1)u(2-m)/(3-m)}
 \left(\frac{p-1}{2p}+\frac{1-m}{m}\right).
\label{mt:eq:complete-large-p-entropy-lower}
\end{align}
The coefficient of $69/100$ is positive.  Indeed,
\[
\frac{m^2}{m^2+(p-1)u(2-m)/(3-m)}\leq\frac1{10},
\]
because $p-1\geq5$, $u\geq18/5$, and
$(2-m)/(3-m)\geq1/2$.

For the derivative calculation only, denote the displayed fraction by
$\rho(u)$.  The derivative of
\eqref{mt:eq:complete-large-p-entropy-lower} is exactly
\[
-\frac1{2p}
-\rho(u)^2\left(\frac{p-1}{2p}+\frac{1-m}{m}\right)
+\frac{207m}{200u}\rho(u)(1-\rho(u)).
\]
The last, positive term is smaller than $1/(2p)$, since
\[
\rho(u)(1-\rho(u))
\leq\rho(u)\leq
\frac{m^2}{(p-1)u(2-m)/(3-m)}
\]
and, using $m^3\leq1$,
\[
\frac{2p}{p-1}\frac{207}{200}
<\frac{12}{5}\frac{207}{200}
<\frac12\left(\frac{18}{5}\right)^2
\leq\frac{2-m}{3-m}\,u^2.
\]
Thus the derivative is negative throughout
\eqref{mt:eq:complete-large-p-window}.
It is therefore enough to put $u=p+1-2m$.

The resulting expression is at least $1/81$.  The complete polynomial
identity proving this is
\begin{align}
&\left.\text{right-hand side of
 \eqref{mt:eq:complete-large-p-entropy-lower}}\right|_{u=p+1-2m}
 -\frac1{81}
 \notag\\
&\quad=
 \frac{1}{16200p\{m^2(3-m)+(p-1)(p+1-2m)(2-m)\}}
 \Big[\notag\\
&(1-m)\{257574-255486(1-m)+237750(1-m)^2
                 +197802(1-m)^3\}\notag\\
&+(p-6)\{93570+122003(1-m)-13877(1-m)^2
                 +82825(1-m)^3+32967(1-m)^4\}\notag\\
&+(p-6)^2\{30172+26856(1-m)+4784(1-m)^2
                 +8100(1-m)^3\}\notag\\
&+(p-6)^3\{2392+2392(1-m)\}\Big]\geq0.
 \label{mt:eq:complete-large-p-polynomial}
\end{align}
The two negative monomials cause no loss: for $0\leq1-m\leq1$,
\[
 257574-255486(1-m)\geq2088>0,
 \qquad
 122003(1-m)-13877(1-m)^2\geq0.
\]
Thus the right-hand side of
\eqref{mt:eq:complete-entropy-one} is positive, a contradiction.  This
proves \eqref{mt:eq:complete-small-mass} for $p\geq6$ as well.
\end{proof}

\subsubsection{The derivative signs and the oriented minor}

\begin{proposition}[Complete positive-index calculation]
\label{mt:prop:complete-positive-index}
At every stationary marginal point with $p\geq4$,
\begin{align}
 \frac{\partial}{\partial m}
 \{\text{left-hand side of \eqref{mt:eq:mass}}\}&>0,
 \qquad
 \frac{\partial}{\partial\xi'(q)}
 \{\text{left-hand side of \eqref{mt:eq:mass}}\}>0,
 \label{mt:eq:complete-mass-signs}
\end{align}
Moreover, the oriented minor in \eqref{mt:eq:oriented-minor} is positive;
after multiplication by $2p/q^2$, it is larger than $13759/21600$.
\end{proposition}

\begin{proof}
Lemmas~\ref{mt:lem:complete-mass-exclusion} and
\ref{mt:lem:complete-uniform-separator} give
\begin{equation}
 m<\frac1{p+1},\qquad
 m^2\xi'(q)>2(p-1)\log2+\frac12>4.
 \label{mt:eq:complete-positive-index-range}
\end{equation}
Apply \eqref{mt:eq:compact-adjacent-moments} first with exponent $2-m$
and then with exponent $4-m$.  Together with
\eqref{mt:eq:complete-contact-data}, this gives
\begin{align}
 \frac{2-m}{3-m}
 &\leq
 \frac{\E[\cosh^{-4}(M(q))e^{W_\mu(q)}]}
      {\E[\cosh^{-2}(M(q))e^{W_\mu(q)}]}
 \leq\frac{2-m}{3-m}
 +\frac1{\xi'(q)(2-m)(3-m)},
 \label{mt:eq:complete-first-ratio}\\
 \frac{4-m}{5-m}
 &\leq
 \frac{\E[\cosh^{-6}(M(q))e^{W_\mu(q)}]}
      {\E[\cosh^{-4}(M(q))e^{W_\mu(q)}]}
 \leq\frac{4-m}{5-m}
 +\frac1{\xi'(q)(4-m)(5-m)}.
 \label{mt:eq:complete-second-ratio}
\end{align}
The marginal identity gives the exact relation
\begin{align*}
&\E[\cosh^{-2}(M(q))e^{W_\mu(q)}]\\
&\quad=
\frac{m^2}{
 m^2+(p-1)m^2\xi'(q)
 \dfrac{\E[\cosh^{-4}(M(q))e^{W_\mu(q)}]}
       {\E[\cosh^{-2}(M(q))e^{W_\mu(q)}]}}.
\end{align*}
Since $p-1\geq3$, $m^2\leq1/25$,
$m^2\xi'(q)>116/25$, and the last ratio is at least $9/14$,
\begin{equation*}
 \E[\cosh^{-2}(M(q))e^{W_\mu(q)}]
 <\frac{1/25}{3(116/25)(9/14)}
 =\frac7{1566}<\frac1{200}.
\end{equation*}
Also $\xi'(q)=m^2\xi'(q)/m^2>116>100$.  Thus
\begin{equation}
 \E[\cosh^{-2}(M(q))e^{W_\mu(q)}]<\frac1{200},
 \qquad \xi'(q)>100.
 \label{mt:eq:complete-small-second-moment}
\end{equation}

Directly from \eqref{mt:eq:complete-first-ratio}--
\eqref{mt:eq:complete-small-second-moment},
\begin{align}
\frac67
<&\frac{
 2\{\E[\cosh^{-2}(M(q))e^{W_\mu(q)}]
      -\E[\cosh^{-4}(M(q))e^{W_\mu(q)}]\}}
 {\E[\cosh^{-4}(M(q))e^{W_\mu(q)}]}
 \notag\\
&+\frac{
 \frac m2\{\E[\cosh^{-4}(M(q))e^{W_\mu(q)}]
      -\E[\cosh^{-2}(M(q))e^{W_\mu(q)}]^2\}}
 {\E[\cosh^{-4}(M(q))e^{W_\mu(q)}]}
 <\frac43,
 \label{mt:eq:complete-curvature-ratio}\\
0<&\frac{1}{\Gamma_\mu(q)}
 \frac{\partial}{\partial m}\Gamma_\mu(q),
 \qquad
 p\frac{1}{\Gamma_\mu(q)}
 \frac{\partial}{\partial m}\Gamma_\mu(q)
 \leq\frac{287}{675}.
 \label{mt:eq:complete-normalized-m-derivative}
\end{align}
For the last bound, Gaussian integration by parts and
$0\leq x\tanh x-\log\cosh x\leq\log2$ give
\begin{align*}
 \frac{1}{\Gamma_\mu(q)}
 \frac{\partial}{\partial m}\Gamma_\mu(q)
 \leq
 \frac{m+(1-m)\E[\cosh^{-2}(M(q))e^{W_\mu(q)}]}
 {(p-1)
  \frac{\E[\cosh^{-4}(M(q))e^{W_\mu(q)}]}
       {\E[\cosh^{-2}(M(q))e^{W_\mu(q)}]}};
\end{align*}
now use $m\leq1/(p+1)$, the lower bound in
\eqref{mt:eq:complete-first-ratio}, and
\eqref{mt:eq:complete-small-second-moment}.  Explicitly,
\[
 p\frac{1}{\Gamma_\mu(q)}
 \frac{\partial}{\partial m}\Gamma_\mu(q)
 \leq\frac43\frac{14}{9}\frac{41}{200}
 =\frac{287}{675}.
\]

Substitution in the exact derivative of the stationary equation yields
\begin{equation}
 \frac{\partial}{\partial\xi'(q)}
 \{\text{left-hand side of \eqref{mt:eq:mass}}\}
 >\frac{(1-m)q}{2p}
 \left(\frac67-\frac{287}{675}\right)>0.
 \label{mt:eq:complete-xi-margin}
\end{equation}
The derivative in $m$ follows from the exact identity
\begin{align*}
&m\frac{\partial}{\partial m}
 \{\text{left-hand side of \eqref{mt:eq:mass}}\}\\
&\quad=\frac{\xi'(q)}p
 \left\{q-\left(p-\frac{p+1}{2}m\right)
 \frac{\partial}{\partial m}\Gamma_\mu(q)\right\}
 -\E\Big[\big(\log\cosh(M(q))
  -\E[\log\cosh(M(q))e^{W_\mu(q)}]\big)\\
&\hspace{48mm}\times
  (M(q)\tanh(M(q))-\log\cosh(M(q)))e^{W_\mu(q)}\Big].
\end{align*}
The last covariance is positive and is at most
$\frac7{20}\sqrt{\xi'(q)q}$: Popoviciu and Cauchy--Schwarz give the
factor $7/20$, while conditioning by $1+\tanh x$ and applying the
one-dimensional Brascamp--Lieb inequality gives
$\operatorname{Var}(\log\cosh(M(q)))\leq\xi'(q)q$.
Consequently
\begin{equation}
 m\frac{\partial}{\partial m}
 \{\text{left-hand side of \eqref{mt:eq:mass}}\}
 \geq\frac{388}{675}\frac{\xi'(q)q}{p}
 -\frac7{20}\sqrt{\xi'(q)q}>0.
 \label{mt:eq:complete-m-margin}
\end{equation}
The last sign follows from
\eqref{mt:eq:complete-positive-index-range} and
\eqref{mt:eq:complete-small-second-moment}, which give
$\xi'(q)q>p^2$.

We finish with the cancellation in the oriented minor.  First we record
the residual derivative sign needed to preserve the direction of every
subsequent estimate.  Formula \eqref{mt:eq:moment-four} and the upper bound
in \eqref{mt:eq:complete-second-ratio} give
\begin{align}
&\frac{
 \frac{\partial}{\partial\xi'(q)}
 \E[\cosh^{-4}(M(q))e^{W_\mu(q)}]}
 {\E[\cosh^{-4}(M(q))e^{W_\mu(q)}]}
 \leq-\frac{m^2}{2}
 -\frac{m(1-m)}2
  \E[\cosh^{-2}(M(q))e^{W_\mu(q)}]
 <-\frac{m^2}{2}.
 \label{mt:eq:complete-fourth-xi-upper}
\end{align}
Indeed, the first two terms in \eqref{mt:eq:moment-four} cancel to
$-m^2/2$ when the ratio in
\eqref{mt:eq:complete-second-ratio} is replaced by $(4-m)/(5-m)$.
Consequently, using \eqref{mt:eq:complete-curvature-ratio} and
$m^2\xi'(q)>4$,
\begin{align}
&\frac{\partial}{\partial\xi'(q)}
 \Big\{\Gamma_\mu(q)-(p-1)\xi'(q)
 \E[\cosh^{-4}(M(q))e^{W_\mu(q)}]\Big\}\notag\\
&\quad>q m^2\left\{\frac12-
 \frac{p-2+(1-m)(4/3)}{(p-1)m^2\xi'(q)}\right\}\notag\\
&\quad>q m^2\left\{\frac12-
 \frac{p-1+1/3}{4(p-1)}\right\}>0.
 \label{mt:eq:complete-residual-xi-positive}
\end{align}
The $m$ derivative of the same residual is positive by
\eqref{mt:eq:moment-one} and \eqref{mt:eq:moment-three}.

For this calculation only, introduce the normalized abbreviations
\begin{align*}
\mathcal R={}&
 \frac{2\{\E[\cosh^{-2}(M(q))e^{W_\mu(q)}]
             -\E[\cosh^{-4}(M(q))e^{W_\mu(q)}]\}}
 {\E[\cosh^{-4}(M(q))e^{W_\mu(q)}]}\\
&+\frac{
 \frac m2\{\E[\cosh^{-4}(M(q))e^{W_\mu(q)}]
             -\E[\cosh^{-2}(M(q))e^{W_\mu(q)}]^2\}}
 {\E[\cosh^{-4}(M(q))e^{W_\mu(q)}]},\\
\mathcal U={}&\frac1q\frac{\partial}{\partial m}\Gamma_\mu(q),\\
\mathcal W={}&-\frac{(p-1)\xi'(q)}q
 \frac{\partial}{\partial m}
 \E[\cosh^{-4}(M(q))e^{W_\mu(q)}],\\
\mathcal T={}&\frac{(p-1)\xi'(q)}q
 \frac{\partial}{\partial\xi'(q)}
 \Big\{\Gamma_\mu(q)-(p-1)\xi'(q)
 \E[\cosh^{-4}(M(q))e^{W_\mu(q)}]\Big\}.
\end{align*}
The endpoint differentiation identities also give the exact
representations
\begin{align}
\mathcal U={}&
\frac{\E[\cosh^{-2}(M(q))e^{W_\mu(q)}]}q
\left\{
\E[\log\cosh(M(q))e^{W_\mu(q)}]\right.\notag\\
&\left.\hspace{16mm}-
\frac{\E[\cosh^{-2}(M(q))\log\cosh(M(q))e^{W_\mu(q)}]}
     {\E[\cosh^{-2}(M(q))e^{W_\mu(q)}]}
\right\},
\label{mt:eq:complete-U-log-form}\\
\mathcal W={}&
\E[\log\cosh(M(q))e^{W_\mu(q)}]\notag\\
&-
\frac{\E[\cosh^{-4}(M(q))\log\cosh(M(q))e^{W_\mu(q)}]}
     {\E[\cosh^{-4}(M(q))e^{W_\mu(q)}]},
\label{mt:eq:complete-W-log-form}\\
\mathcal T={}&1-(p-1)-(1-m)\mathcal R\notag\\
&-(p-1)\xi'(q)
\frac{
 \frac{\partial}{\partial\xi'(q)}
 \E[\cosh^{-4}(M(q))e^{W_\mu(q)}]}
 {\E[\cosh^{-4}(M(q))e^{W_\mu(q)}]}.
\label{mt:eq:complete-T-moment-form}
\end{align}
Thus $\mathcal R-p\mathcal U>0$ and
$\mathcal U,\mathcal W,\mathcal T>0$.  Put, still locally,
\begin{align*}
\beta={}&p-\frac{p+1}{2}m,
&\mathcal R_0={}&\frac2{2-m}+\frac m2,
&\mathcal U_0={}&\frac{m(3-m)}{(p-1)(2-m)},\\
\mathcal L={}&(1-m)(\mathcal R_0-p\mathcal U_0)
              =1-\beta\mathcal U_0,\\
\mathcal C_0={}&1-\frac{p-1}{2}-(1-m)\mathcal R_0
 +\frac{m(1-m)(3-m)}{2(2-m)},\\
\epsilon_{\mathcal R}={}&
 \frac{2m^2(3-m)}{(2-m)^3}
 +\frac{m^3(3-m)^2}{2(p-1)(2-m)^2},\\
\epsilon_+={}&
 \frac{(1-m)m^2(3-m)^2}{(p-1)^2(2-m)^2},\\
\epsilon_-={}&\frac1{p-1}\left\{
 \frac{m^3(3-m)}{(2-m)^3}
 +\frac{2m^2(3-m)\log2}{2-m}\right\},
&\epsilon_{\mathcal H}={}&\epsilon_{\mathcal R}+p\epsilon_+.
\end{align*}
In these abbreviations the four exact endpoint derivatives are
\begin{align}
\frac{\partial}{\partial\xi'(q)}
 \{\text{left-hand side of \eqref{mt:eq:mass}}\}
&=\frac{(1-m)q}{2p}(\mathcal R-p\mathcal U),
 \label{mt:eq:complete-exact-mass-xi}\\
\frac{\partial}{\partial m}
 \Big\{\Gamma_\mu(q)-(p-1)\xi'(q)
 \E[\cosh^{-4}(M(q))e^{W_\mu(q)}]\Big\}
&=q(\mathcal U+\mathcal W),
 \label{mt:eq:complete-exact-residual-m}\\
\frac{\partial}{\partial\xi'(q)}
 \Big\{\Gamma_\mu(q)-(p-1)\xi'(q)
 \E[\cosh^{-4}(M(q))e^{W_\mu(q)}]\Big\}
&=\frac{q\mathcal T}{(p-1)\xi'(q)},
 \label{mt:eq:complete-exact-residual-xi}\\
\frac{\partial}{\partial m}
 \{\text{left-hand side of \eqref{mt:eq:mass}}\}
&\leq\frac{\xi'(q)q}{mp}(1-\beta\mathcal U).
 \label{mt:eq:complete-mass-m-upper}
\end{align}
The last line follows from the exact identity displayed immediately
before \eqref{mt:eq:complete-m-margin}: its covariance term is strictly
positive.

The adjacent-moment bounds and the contact identities, together with
the pointwise estimate
$0\leq x\tanh x-\log\cosh x\leq(\log2)\tanh^2x$, imply
\begin{align}
\mathcal R&\geq\mathcal R_0-
 \frac{\epsilon_{\mathcal R}}{m^2\xi'(q)},
 \label{mt:eq:complete-R-error}\\
\mathcal U_0-\frac{\epsilon_-}{m^2\xi'(q)}
&\leq\mathcal U\leq
 \mathcal U_0+\frac{\epsilon_+}{m^2\xi'(q)},
 \label{mt:eq:complete-U-error}\\
\mathcal W&\geq\frac{m^2\xi'(q)}m-2\log2,
 \label{mt:eq:complete-W-error}\\
\mathcal T&\leq\frac{(p-1)m^2\xi'(q)}2+\mathcal C_0
 +\frac{(1-m)\epsilon_{\mathcal R}}{m^2\xi'(q)}.
 \label{mt:eq:complete-T-error}
\end{align}
The pointwise upper bound follows from
\eqref{mt:eq:complete-entropy-series}, since every power of
$\tanh^2x$ is at most $\tanh^2x$.  The displayed error bounds follow
by inserting the two sides of
\eqref{mt:eq:complete-first-ratio}--\eqref{mt:eq:complete-second-ratio};
the reciprocal estimate used there is
$1/(x+z)\geq1/x-z/x^2$ for $x,z>0$.

Let $\mathfrak D$ denote the left-hand side of
\eqref{mt:eq:oriented-minor}.  Equations
\eqref{mt:eq:complete-exact-mass-xi}--
\eqref{mt:eq:complete-mass-m-upper}, together with
\eqref{mt:eq:complete-residual-xi-positive}, give
\begin{equation}
\frac{2p}{q^2}\mathfrak D\geq
 (1-m)(\mathcal R-p\mathcal U)(\mathcal U+\mathcal W)
 -\frac2{m(p-1)}(1-\beta\mathcal U)\mathcal T.
 \label{mt:eq:complete-normalized-minor-start}
\end{equation}
Every factor whose upper bound is inserted in the second product is
positive by \eqref{mt:eq:complete-residual-xi-positive}; also
$m^2\xi'(q)/m-2\log2>0$.  Insert
\eqref{mt:eq:complete-R-error}--\eqref{mt:eq:complete-T-error} in that
order and expand.  The two terms
$\mathcal Lm^2\xi'(q)/m$ cancel exactly.  Dropping the two remaining
favorable terms gives
\begin{align}
\frac{2p}{q^2}\mathfrak D
\geq{}&\mathcal L\left(-2\log2-
 \frac{2\mathcal C_0}{m(p-1)}\right)
 -\frac{(1-m)\epsilon_{\mathcal H}}m-\frac{\beta\epsilon_-}m
 \notag\\
&-\frac{2\mathcal L(1-m)\epsilon_{\mathcal R}}
 {m(p-1)m^2\xi'(q)}
 -\frac{2\beta\epsilon_-(1-m)\epsilon_{\mathcal R}}
 {m(p-1)\{m^2\xi'(q)\}^2}.
 \label{mt:eq:complete-minor-after-cancellation}
\end{align}
This is the exact product cancellation behind the numerical margin.

Direct algebra gives
\begin{align*}
\mathcal L&\geq1-
 \frac{14p}{9(p-1)(p+1)}\geq\frac7{12},\\
-\mathcal C_0&=\frac{p-1}{2}-1
 +\frac{(1-m)(4-m)}{2(2-m)},\\
-2\log2-\frac{2\mathcal C_0}{m(p-1)}
&\geq p-1+\frac15-\frac4{5(p-1)}\geq\frac{44}{15}.
\end{align*}
Furthermore,
\begin{equation}
 \epsilon_{\mathcal R}<\frac{25}{24}m^2,
 \qquad \epsilon_{\mathcal H}<\frac{17}{8}m^2,
 \qquad \epsilon_-<\frac{23}{10(p-1)}m^2.
 \label{mt:eq:complete-three-errors}
\end{equation}
For the first inequality, the two coefficients after division by $m^2$
are at most $700/729$ and $98/1215$; for the other two use
$(3-m)/(2-m)\leq14/9$, $p/(p-1)^2\leq4/9$, and
$\log2<7/10$.  Hence the four absolute errors in
\eqref{mt:eq:complete-minor-after-cancellation} are, in order, at most
\begin{equation}
 \frac{17}{40},\qquad
 \frac{46}{75},\qquad
 \frac5{144},\qquad
 \frac{23}{21600}.
 \label{mt:eq:complete-four-errors}
\end{equation}
Indeed, these follow respectively from
\begin{align*}
&\frac{17}{8}m^2\frac1m\leq\frac{17}{40},
&&\frac{23}{10(p-1)}m^2
 \frac{p-(p+1)m/2}{m}\leq\frac{46}{75},\\
&\frac{2(1-m)}{m(p-1)m^2\xi'(q)}\frac{25}{24}m^2
 \leq\frac5{144},
&&\frac{2\{p-(p+1)m/2\}(1-m)}
 {m(p-1)\{m^2\xi'(q)\}^2}
 \frac{23}{10(p-1)}m^2\frac{25}{24}m^2
 \leq\frac{23}{21600},
\end{align*}
where only $p\geq4$, $m\leq1/(p+1)$, and
$m^2\xi'(q)>4$ were used.  Therefore
\begin{align}
&\frac{2p}{q^2}\Bigg[
 \frac{\partial}{\partial\xi'(q)}
  \{\text{left-hand side of \eqref{mt:eq:mass}}\}
 \frac{\partial}{\partial m}
 \Big\{\Gamma_\mu(q)-(p-1)\xi'(q)
 \E[\cosh^{-4}(M(q))e^{W_\mu(q)}]\Big\}
 \notag\\
&\hspace{16mm}-
 \frac{\partial}{\partial m}
  \{\text{left-hand side of \eqref{mt:eq:mass}}\}
 \frac{\partial}{\partial\xi'(q)}
 \Big\{\Gamma_\mu(q)-(p-1)\xi'(q)
 \E[\cosh^{-4}(M(q))e^{W_\mu(q)}]\Big\}\Bigg]
 \notag\\
&\quad>
 \frac7{12}\frac{44}{15}
 -\frac{17}{40}-\frac{46}{75}-\frac5{144}
 -\frac{23}{21600}
 =\frac{13759}{21600}>0.
 \label{mt:eq:complete-final-minor}
\end{align}
This proves \eqref{mt:eq:complete-mass-signs} and
\eqref{mt:eq:oriented-minor}.
\end{proof}
 
\subsection{Completion of the strict-sign proof}
\label{mt:sec:complete-strict-sign}

\begin{proof}[Proof of Proposition~\ref{mt:prop:strict-sign}]
Suppose first that $p=3$.  Lemma~\ref{mt:lem:complete-p3-entry} gives
\[
 m^2\xi'(q)>\frac52.
\]
Lemma~\ref{mt:lem:complete-p3-jacobian} then proves both inequalities in
\eqref{mt:eq:mass-signs} and the oriented-minor inequality
\eqref{mt:eq:oriented-minor}.

Suppose next that $p\geq4$.  Lemma~\ref{mt:lem:complete-mass-exclusion}
first gives $m<1/(p+1)$, and
Lemma~\ref{mt:lem:complete-uniform-separator} gives
\[
 m^2\xi'(q)>2(p-1)\log2+\frac12.
\]
Proposition~\ref{mt:prop:complete-positive-index} now proves the same
three strict inequalities.  These cases exhaust $p\geq3$.
\end{proof}

\section{Exact certificates and verifier sources}
\label{app:certificates}

The finite checks used in the two preceding appendices are exact
continuous-box certificates, not numerical samples.  Their arithmetic
contracts and complete verifier sources are collected here so that the
entire proof is contained in this one document.

\subsection{The contact-curvature computation}


The listings in this subsection are to be saved together in
\path{work/contact}.  They must not be placed in the directories used
for the marginal-crossing or internal-gap certificates, because those
subsections contain different sources with some of the same filenames.
The listings printed here, rather than any pre-existing file with the
same name, are the sources identified by the hashes below.

\begin{itemize}
\item \texttt{\detokenize{certify_D_largeA_dyadic.cpp}}: all \(A\ge 4\), all \(c,k\).
\item \texttt{\detokenize{certify_D_scalar_dyadic.cpp}}: compact \(P=2,3\) boxes.
\item \texttt{\detokenize{certify_terminal_affine_dyadic.cpp}}: terminal affine and high-\(k\)
  regional exclusions.
\item \texttt{\detokenize{certify_terminal_affine_fixed.py}}: interval-arithmetic module used by
  the terminal \(A=4\) verifier.
\item \texttt{\detokenize{certify_terminal_A4_fixed.py}}: exact \(A=4\) small-\(\widehat m\) boundary.
\item \texttt{\detokenize{correct_xz_exact_certificate.py}}: exact \(P=2\) polynomial-basis exclusions.
\item \texttt{\detokenize{pge3_scalar_exact_certificate.py}}: exact \(P\ge 3\) polynomial-basis exclusions.

\end{itemize}
The following commands, run inside \path{work/contact}, reproduce every
certificate in this subsection:
\begin{Verbatim}[fontsize=\small,frame=single]
c++ -O2 -std=c++17 certify_D_largeA_dyadic.cpp \
  -o certify_D_largeA_dyadic
./certify_D_largeA_dyadic
c++ -O2 -std=c++17 certify_D_scalar_dyadic.cpp \
  -o certify_D_scalar_dyadic
./certify_D_scalar_dyadic
c++ -O2 -std=c++17 certify_terminal_affine_dyadic.cpp \
  -o certify_terminal_affine_dyadic
./certify_terminal_affine_dyadic region
./certify_terminal_affine_dyadic pregion3
./certify_terminal_affine_dyadic pregion4
./certify_terminal_affine_dyadic pregion5
PYTHONPATH=. python3 certify_terminal_A4_fixed.py
python3 correct_xz_exact_certificate.py
python3 pge3_scalar_exact_certificate.py
\end{Verbatim}
The reported run used Apple clang 14.0.3 and Python 3.14.6.  Any C++17
compiler used for reproduction must support \texttt{\_\_int128} and the
overflow built-ins used in the listings.  Assertions must remain
enabled; in particular, the C++ sources must not be compiled with
\texttt{-DNDEBUG}.

The no-argument mode is an older, shorter diagnostic interval and is not
the certificate used above.
\subsection{Exact verifier sources}

This appendix prints the frozen verifier sources used for the box certificates in the main reduction.  The hashes below are SHA-256 hashes of the listed source files.  Printed decimal lower bounds are diagnostics; acceptance decisions in the programs use exact integer or rational comparisons.

\begin{itemize}

\item \texttt{\detokenize{certify_D_largeA_dyadic.cpp}}\\{\scriptsize\ttfamily 7c408373\allowbreak 48cf43e3\allowbreak afc7fd9a\allowbreak 62866a66\allowbreak 7783916e\allowbreak bc82e5f7\allowbreak 5b905a8f\allowbreak eb1f6a5d}

\item \texttt{\detokenize{certify_D_scalar_dyadic.cpp}}\\{\scriptsize\ttfamily 978e1b21\allowbreak b38be793\allowbreak 3d70127c\allowbreak 15fab832\allowbreak 0eb2bbb7\allowbreak 69b6253c\allowbreak ec963395\allowbreak 76bf232b}

\item \texttt{\detokenize{certify_terminal_affine_dyadic.cpp}}\\{\scriptsize\ttfamily f5b6172b\allowbreak 008ef2fd\allowbreak a452b92c\allowbreak 96575fee\allowbreak 77276013\allowbreak e6c08b72\allowbreak 3a82cfb0\allowbreak 363a1582}

\item \texttt{\detokenize{certify_terminal_affine_fixed.py}}\\{\scriptsize\ttfamily c5e2a20a\allowbreak 4eea8b46\allowbreak 4d285f79\allowbreak 30e43d52\allowbreak 23c1f059\allowbreak 7c0e44fd\allowbreak e8b9307d\allowbreak c9c4e942}

\item \texttt{\detokenize{certify_terminal_A4_fixed.py}}\\{\scriptsize\ttfamily 6e2582ec\allowbreak abf746a0\allowbreak a8778c45\allowbreak 32cf0b43\allowbreak 828dbac2\allowbreak 5eb250ec\allowbreak aa043ee4\allowbreak d0390303}

\item \texttt{\detokenize{correct_xz_exact_certificate.py}}\\{\scriptsize\ttfamily 63723a83\allowbreak 6d6cf5f7\allowbreak 6e777b8d\allowbreak de357991\allowbreak f3d93457\allowbreak 3a952620\allowbreak 8495a890\allowbreak 1a046e6a}

\item \texttt{\detokenize{pge3_scalar_exact_certificate.py}}\\{\scriptsize\ttfamily ec21f1b5\allowbreak be55c796\allowbreak c529658c\allowbreak f3cd0027\allowbreak b6ada5bd\allowbreak f2694d31\allowbreak 4db91e42\allowbreak 717c8d58}

\end{itemize}

\subsubsection{Large-$A$ scalar certificate}

\begin{lstlisting}[language={C++}]
#include <algorithm>
#include <array>
#include <cassert>
#include <iomanip>
#include <iostream>
#include <vector>

// Exact outward-rounded certificate for the corrected D lower bound on
// A>=4, simultaneously for 1/4<=c<=1/2 and every k>=0.  Coordinates are
// h=1/A in [0,1/4], m=1/(1+k) in [0,1], c in [1/4,1/2],
// and the interpolation r between delta_phi and delta_density.
using i128=__int128_t; static constexpr unsigned QB=45; static constexpr i128 SC=i128(1)<<QB;
static i128 mulr(i128 a,i128 b){i128 z;if(__builtin_mul_overflow(a,b,&z)){std::cerr<<"overflow\n";std::abort();}return z;}
static i128 addr(i128 a,i128 b){i128 z;if(__builtin_add_overflow(a,b,&z)){std::cerr<<"overflow\n";std::abort();}return z;}
static i128 subr(i128 a,i128 b){i128 z;if(__builtin_sub_overflow(a,b,&z)){std::cerr<<"overflow\n";std::abort();}return z;}
static i128 fld(i128 n,i128 d){assert(d>0);i128 q=n/d,r=n%d;if(r&&n<0)--q;return q;}
static i128 cei(i128 n,i128 d){assert(d>0);i128 q=n/d,r=n%d;if(r&&n>0)++q;return q;}
struct I{ i128 l=0,u=0;I()=default;I(i128 L,i128 U):l(L),u(U){assert(l<=u);}static I rat(i128 p,i128 q){if(q<0)p=-p,q=-q;return{fld(mulr(p,SC),q),cei(mulr(p,SC),q)};}static I ex(long long z){return{i128(z)*SC,i128(z)*SC};}};
static I operator+(I x,I y){return{addr(x.l,y.l),addr(x.u,y.u)};}static I operator-(I x){return{subr(0,x.u),subr(0,x.l)};}static I operator-(I x,I y){return x+(-y);}
static I operator*(I x,I y){std::array<i128,4>z{mulr(x.l,y.l),mulr(x.l,y.u),mulr(x.u,y.l),mulr(x.u,y.u)};auto p=std::minmax_element(z.begin(),z.end());return{fld(*p.first,SC),cei(*p.second,SC)};}
struct Q{i128 n,d;};static bool qless(const Q&a,const Q&b){return mulr(a.n,b.d)<mulr(b.n,a.d);}static I operator/(I x,I y){assert(y.l>0||y.u<0);std::array<Q,4>z{{{x.l,y.l},{x.l,y.u},{x.u,y.l},{x.u,y.u}}};for(auto&q:z)if(q.d<0)q.n=-q.n,q.d=-q.d;auto p=std::minmax_element(z.begin(),z.end(),qless);return{fld(mulr(p.first->n,SC),p.first->d),cei(mulr(p.second->n,SC),p.second->d)};}
static I operator+(I x,long long z){return x+I::ex(z);}static I operator+(long long z,I x){return x+z;}static I operator-(I x,long long z){return x-I::ex(z);}static I operator-(long long z,I x){return I::ex(z)-x;}static I operator*(I x,long long z){return x*I::ex(z);}static I operator*(long long z,I x){return x*z;}static I operator/(I x,long long z){return x/I::ex(z);}static I operator/(long long z,I x){return I::ex(z)/x;}

struct Box{std::array<i128,4>l,u;int dep=0;};struct Ev{I out,pos,span,sw;bool empty=false;};
static Ev eval(const Box&b){
 I ht(b.l[0],b.u[0]),m(b.l[1],b.u[1]),ct(b.l[2],b.u[2]),r(b.l[3],b.u[3]);
 I h=ht/4,c=I::rat(1,4)+ct/4;
 I rho=1-m+m*h, theta=1-m;
 I dlo=(1-c*h)*m*(3-m)/(10-m-3*m*m);
 I dhi=1/(4+rho),span=dhi-dlo;
 if(span.u<0)return{I(),I(),span,I(),true};
 I delta=dlo+span*r;
 I E=(1-c*h-3*delta)/2;
 I Lam=rho*rho*delta+(1+rho*rho)*E;
 // pos is m*h^2 times the coefficient of beta in D.
 I pos=m*(1-c*h)*Lam/2+(1-m+c*h*m)*(1-m)*E;
 I pen=c*h*(h*m+1-m+m*m*E/2);
 I den=(rho+6)*delta*delta;
 I sw=(rho+7)*delta-1-I::rat(9,5)*den; // beta_t-4/5, times den
 I out;
 // The tangent branch is evaluated only when it is active.  On boundary
 // boxes with m=0 and r=0, delta (and hence den) can contain zero, while
 // sw<0 guarantees that the universal beta=4/5 branch is the one required.
 // Forming the unused quotient there would introduce an invalid 0-denominator
 // interval and make an optimized run depend on dead-code elimination.
 if(sw.l>=0){
   assert(den.l>0);
   I bt=((rho+7)*delta-1)/den-1;
   out=bt*pos-pen;
 } else out=I::rat(4,5)*pos-pen;
 return{out,pos,span,sw,false};
}
static long double ld(i128 z){return(long double)z/(long double)SC;}
static std::pair<Box,Box> split(const Box&b,int ax){Box x=b,y=b;x.dep=y.dep=b.dep+1;i128 md=addr(b.l[ax],b.u[ax])/2;x.u[ax]=md;y.l[ax]=md;return{x,y};}

int main(){
 std::vector<Box>st;constexpr int NH=8,NM=16,NC=2,NR=8;
 for(int i=0;i<NH;++i)for(int j=0;j<NM;++j)for(int c=0;c<NC;++c)for(int r=0;r<NR;++r){Box b;b.l={i*SC/NH,j*SC/NM,c*SC/NC,r*SC/NR};b.u={(i+1)*SC/NH,(j+1)*SC/NM,(c+1)*SC/NC,(r+1)*SC/NR};st.push_back(b);}
 long long vis=0,leaves=0,empty=0;i128 best=100*SC;Box bestb;int md=0;
 while(!st.empty()){
  Box b=st.back();st.pop_back();++vis;Ev e=eval(b);
  if(e.empty){++empty;continue;}
  // If span straddles zero, the interval evaluation still encloses every
  // feasible point in the box; infeasible points only enlarge it.  Thus a
  // positive enclosure may be accepted without proving span>=0.
  if(e.pos.l>0&&e.out.l>0){++leaves;if(e.out.l<best)best=e.out.l,bestb=b;continue;}
  if(b.dep>=38){std::cerr<<std::setprecision(18)<<"FAIL dep="<<b.dep<<" out="<<ld(e.out.l)<<','<<ld(e.out.u)<<" pos="<<ld(e.pos.l)<<" span="<<ld(e.span.l)<<','<<ld(e.span.u)<<" sw="<<ld(e.sw.l)<<','<<ld(e.sw.u)<<" box=";for(int z=0;z<4;++z)std::cerr<<'['<<ld(b.l[z])<<','<<ld(b.u[z])<<"] ";std::cerr<<'\n';return 2;}
  i128 mw=0;for(int ax=0;ax<4;++ax)mw=std::max(mw,subr(b.u[ax],b.l[ax]));
  int bax=0;i128 bs=-((i128)1<<126);std::pair<Box,Box>bk;
  for(int ax=0;ax<4;++ax){if(mulr(subr(b.u[ax],b.l[ax]),4)<mw)continue;auto q=split(b,ax);auto aa=eval(q.first),dd=eval(q.second);i128 x=aa.empty?SC*100:aa.out.l,y=dd.empty?SC*100:dd.out.l,score=std::min(x,y);if(score>bs){bs=score;bax=ax;bk=q;}}
  st.push_back(bk.first);st.push_back(bk.second);md=std::max(md,b.dep+1);
 }
 std::cout<<std::setprecision(18)<<"PASS exact large-A corrected-D certificate visited="<<vis<<" leaves="<<leaves<<" empty="<<empty<<" maxdepth="<<md<<" smallest accepted scaled lower="<<ld(best)<<" box=";for(int i=0;i<4;++i)std::cout<<'['<<ld(bestb.l[i])<<','<<ld(bestb.u[i])<<"] ";std::cout<<'\n';
}
\end{lstlisting}

\subsubsection{Compact scalar certificates}

\begin{lstlisting}[language={C++}]
#include <algorithm>
#include <array>
#include <cassert>
#include <iomanip>
#include <iostream>
#include <vector>

// Exact outward-rounded dyadic interval certificate for the final P=2
// scalar D bound.  No floating-point number is used to certify a sign.
using i128=__int128_t;
static constexpr unsigned QB=45;
static constexpr i128 SC=i128(1)<<QB;

static i128 mulraw(i128 a,i128 b){i128 r;if(__builtin_mul_overflow(a,b,&r)){std::cerr<<"overflow\n";std::abort();}return r;}
static i128 addraw(i128 a,i128 b){i128 r;if(__builtin_add_overflow(a,b,&r)){std::cerr<<"overflow\n";std::abort();}return r;}
static i128 subraw(i128 a,i128 b){i128 r;if(__builtin_sub_overflow(a,b,&r)){std::cerr<<"overflow\n";std::abort();}return r;}
static i128 fld(i128 n,i128 d){assert(d>0);i128 q=n/d,r=n%d;if(r&&n<0)--q;return q;}
static i128 cei(i128 n,i128 d){assert(d>0);i128 q=n/d,r=n%d;if(r&&n>0)++q;return q;}
struct I{
  i128 l=0,u=0;
  I()=default; I(i128 L,i128 U):l(L),u(U){assert(l<=u);}
  static I rat(i128 p,i128 q){if(q<0)p=-p,q=-q;return {fld(mulraw(p,SC),q),cei(mulraw(p,SC),q)};}
  static I exact(long long z){return {i128(z)*SC,i128(z)*SC};}
};
static I operator+(I x,I y){return{addraw(x.l,y.l),addraw(x.u,y.u)};}
static I operator-(I x){return{subraw(0,x.u),subraw(0,x.l)};}
static I operator-(I x,I y){return x+(-y);}
static I operator*(I x,I y){std::array<i128,4>z{mulraw(x.l,y.l),mulraw(x.l,y.u),mulraw(x.u,y.l),mulraw(x.u,y.u)};auto p=std::minmax_element(z.begin(),z.end());return{fld(*p.first,SC),cei(*p.second,SC)};}
struct R{i128 n,d;};
static bool lessq(const R&a,const R&b){return mulraw(a.n,b.d)<mulraw(b.n,a.d);}
static I operator/(I x,I y){assert(y.l>0||y.u<0);std::array<R,4>z{{{x.l,y.l},{x.l,y.u},{x.u,y.l},{x.u,y.u}}};for(auto&q:z)if(q.d<0)q.n=-q.n,q.d=-q.d;auto p=std::minmax_element(z.begin(),z.end(),lessq);return{fld(mulraw(p.first->n,SC),p.first->d),cei(mulraw(p.second->n,SC),p.second->d)};}
static I operator+(I x,long long z){return x+I::exact(z);} static I operator+(long long z,I x){return x+z;}
static I operator-(I x,long long z){return x-I::exact(z);} static I operator-(long long z,I x){return I::exact(z)-x;}
static I operator*(I x,long long z){return x*I::exact(z);} static I operator*(long long z,I x){return x*z;}
static I operator/(I x,long long z){return x/I::exact(z);} static I operator/(long long z,I x){return I::exact(z)/x;}

struct Box{std::array<i128,3>l,u;int dep=0;};
struct Ev{I out,A,dspan,gap;};
struct Case { I c,klo,khi,intercept,slope; const char*name; };
static Case CS{I::rat(1,4),I::rat(8,25),I::rat(541,360),I::rat(259,200),I::rat(9,5),"P2"};

static Ev eval(const Box&b){
  I s(b.l[0],b.u[0]),t(b.l[1],b.u[1]),r(b.l[2],b.u[2]);
  I k=CS.klo+(CS.khi-CS.klo)*s;
  I line=CS.intercept+CS.slope*k;
  I a=line+(4-line)*t;
  I rho=(1+a*k)/(a*(1+k));
  I theta=k/(1+k);
  I phi=5*k*(1+k)/(3*k+2);
  I dlo=(a-CS.c)/(a*(3+2*phi));
  I dhi=1/(4+rho);
  I delta=dlo+(dhi-dlo)*r;
  I d=a*delta,e=(a-CS.c-3*d)/2;
  I L=rho*rho*d+(1+rho*rho)*e;
  I A=(a-CS.c)*L/2+(a*k+CS.c)*theta*e;
  I pen=1+a*k+(1-theta)*e/2;
  I base=I::rat(4,5)*A-CS.c*pen;
  I bt=((rho+7)*delta-1)/((rho+6)*delta*delta)-1;
  I gap=bt-I::rat(4,5),out;
  if(gap.l>=0)out=bt*A-CS.c*pen;
  else out=base; // valid also on a straddling box: max(.8,bt)>=.8
  return{out,A,dhi-dlo,gap};
}

static long double ld(i128 z){return (long double)z/(long double)SC;}
static std::pair<Box,Box> split(const Box&b,int ax){Box x=b,y=b;x.dep=y.dep=b.dep+1;i128 m=addraw(b.l[ax],b.u[ax])/2;x.u[ax]=m;y.l[ax]=m;return{x,y};}

static bool run_case(const Case&which){
  CS=which;
  std::vector<Box> st;
  constexpr int NK=32,NT=8,NR=8;
  for(int i=0;i<NK;++i)for(int j=0;j<NT;++j)for(int h=0;h<NR;++h){
    Box b; b.l={i*SC/NK,j*SC/NT,h*SC/NR};b.u={(i+1)*SC/NK,(j+1)*SC/NT,(h+1)*SC/NR};st.push_back(b);
  }
  long long leaves=0,visited=0;int maxdep=0;i128 best=SC*100;Box bestb;
  while(!st.empty()){
    Box b=st.back();st.pop_back();++visited;Ev e=eval(b);
    if(e.A.l>0&&e.dspan.l>=0&&e.out.l>0){++leaves;if(e.out.l<best)best=e.out.l,bestb=b;continue;}
    if(b.dep>=34){std::cerr<<std::setprecision(18)<<CS.name<<" FAIL dep="<<b.dep<<" out=["<<ld(e.out.l)<<','<<ld(e.out.u)<<"] A="<<ld(e.A.l)<<" dspan="<<ld(e.dspan.l)<<" gap=["<<ld(e.gap.l)<<','<<ld(e.gap.u)<<"]\n";return false;}
    // Select the bisection with the best worst child lower bound.  Long
    // double only selects a branch; all accepted signs use i128 intervals.
    int bax=0;i128 bs=-((i128)1<<126);
    std::pair<Box,Box> bk;
    for(int ax=0;ax<3;++ax){auto q=split(b,ax);i128 score=std::min(eval(q.first).out.l,eval(q.second).out.l);if(score>bs){bs=score;bax=ax;bk=q;}}
    st.push_back(bk.first);st.push_back(bk.second);maxdep=std::max(maxdep,b.dep+1);
  }
  std::cout<<std::setprecision(18)<<CS.name<<" PASS exact dyadic D scalar certificate visited="<<visited<<" leaves="<<leaves<<" maxdepth="<<maxdep<<" smallest accepted lower="<<ld(best)<<" box=";
  for(int i=0;i<3;++i)std::cout<<'['<<ld(bestb.l[i])<<','<<ld(bestb.u[i])<<"] ";
  std::cout<<'\n';return true;
}

int main(){
  const Case p2{I::rat(1,4),I::rat(8,25),I::rat(541,360),I::rat(259,200),I::rat(9,5),"P2"};
  const Case p3{I::rat(1,3),I::rat(1,4),I::rat(17,30),I::rat(23,10),I::rat(3,1),"P3"};
  return run_case(p2)&&run_case(p3)?0:2;
}
\end{lstlisting}

\subsubsection{Terminal affine certificates}

\begin{lstlisting}[language={C++}]
#include <algorithm>
#include <array>
#include <cassert>
#include <iostream>
#include <string>
#include <vector>

// A rigorous dyadic interval certificate for the two P=2 terminal lines.
// Every endpoint is an integer multiple of 2^{-B}.  Multiplication and
// division are rounded outwards with integer arithmetic; no floating point
// operation enters the proof.
using i128 = __int128_t;
using u128 = __uint128_t;
using cpp_int = i128;
static constexpr unsigned B = 45;
static constexpr cpp_int SC = cpp_int(1) << B;
static cpp_int cmul(cpp_int a,cpp_int b){cpp_int r;if(__builtin_mul_overflow(a,b,&r)){std::cerr<<"integer overflow\n";std::abort();}return r;}
static cpp_int cadd(cpp_int a,cpp_int b){cpp_int r;if(__builtin_add_overflow(a,b,&r)){std::cerr<<"integer overflow\n";std::abort();}return r;}
static cpp_int csub(cpp_int a,cpp_int b){cpp_int r;if(__builtin_sub_overflow(a,b,&r)){std::cerr<<"integer overflow\n";std::abort();}return r;}

static cpp_int fld(const cpp_int& n, const cpp_int& d) {
  assert(d > 0);
  cpp_int q = n / d, r = n % d;
  if (r != 0 && n < 0) --q;
  return q;
}
static cpp_int cei(const cpp_int& n, const cpp_int& d) {
  assert(d > 0);
  cpp_int q = n / d, r = n % d;
  if (r != 0 && n > 0) ++q;
  return q;
}

struct I {
  cpp_int l, u;
  I() : l(0), u(0) {}
  I(cpp_int L, cpp_int U) : l(L), u(U) { assert(l <= u); }
  static I rat(cpp_int p, cpp_int q) {
    if (q < 0) p = -p, q = -q;
    return I(fld(cmul(p,SC), q), cei(cmul(p,SC), q));
  }
  static I exact(long long n) { return I(cmul(cpp_int(n),SC),cmul(cpp_int(n),SC)); }
};
static I operator+(const I& x, const I& y) { return I(cadd(x.l,y.l),cadd(x.u,y.u)); }
static I operator-(const I& x) { return I(csub(0,x.u),csub(0,x.l)); }
static I operator-(const I& x, const I& y) { return x + (-y); }
static I operator*(const I& x, const I& y) {
  std::array<cpp_int,4> z{cmul(x.l,y.l),cmul(x.l,y.u),cmul(x.u,y.l),cmul(x.u,y.u)};
  auto mm = std::minmax_element(z.begin(),z.end());
  return I(fld(*mm.first,SC),cei(*mm.second,SC));
}
struct Qe { cpp_int n,d; };
static bool qless(const Qe& a,const Qe& b){return cmul(a.n,b.d) < cmul(b.n,a.d);}
static I operator/(const I& x, const I& y) {
  assert(y.l > 0 || y.u < 0);
  std::array<Qe,4> z{{{x.l,y.l},{x.l,y.u},{x.u,y.l},{x.u,y.u}}};
  for(auto& q:z) if(q.d<0)q.n=-q.n,q.d=-q.d;
  auto mm=std::minmax_element(z.begin(),z.end(),qless);
  return I(fld(cmul(mm.first->n,SC),mm.first->d),cei(cmul(mm.second->n,SC),mm.second->d));
}
static I operator+(const I&x,long long n){return x+I::exact(n);} 
static I operator-(const I&x,long long n){return x-I::exact(n);} 
static I operator*(const I&x,long long n){return x*I::exact(n);} 
static I operator/(const I&x,long long n){return x/I::exact(n);} 
static I inv(const I&x){return I::exact(1)/x;}

static I sqrti(const I& x) {
  assert(x.l >= 0);
  auto isqrt=[](u128 n){
    u128 r=0,bit=u128(1)<<126;
    while(bit>n)bit>>=2;
    while(bit){if(n>=r+bit)n-=r+bit,r=(r>>1)+bit;else r>>=1;bit>>=2;}
    return r;
  };
  cpp_int L=(cpp_int)isqrt((u128)cmul(x.l,SC)), U=(cpp_int)isqrt((u128)cmul(x.u,SC));
  if(cmul(U,U) < cmul(x.u,SC))++U;
  return I(L,U);
}

static I ln2() {
  // ln2=2 sum_{j>=0}((2j+1)3^{2j+1})^{-1}; the omitted tail is at
  // most the first omitted denominator with (2j+1) frozen, times 9/8.
  I s;
  cpp_int pow3=3;
  constexpr int N=35;
  for(int j=0;j<=N;++j){
    s=s+I::rat(2,cmul(cpp_int(2*j+1),pow3));
    pow3=cmul(pow3,9);
  }
  I tail=I::rat(18,cmul(cmul(cpp_int(8),2*N+3),pow3));
  return I(s.l,s.u+tail.u);
}
static const I LN2=ln2();

static I exp_pos_point(const cpp_int& xx) {
  assert(xx>=0);
  I x(xx,xx), term=I::exact(1), sum=term;
  constexpr int N=52;
  for(int n=1;n<=N;++n){term=term*x/n;sum=sum+term;}
  I next=term*x/(N+1);
  I tail=next/(I::exact(1)-x/(N+2));
  return I(sum.l,(sum+tail).u);
}
static I expi(const I&x){
  if(x.l>=0){I lo=exp_pos_point(x.l),up=exp_pos_point(x.u);return I(lo.l,up.u);}
  if(x.u<=0){I y=expi(-x);return I(fld(cmul(SC,SC),y.u),cei(cmul(SC,SC),y.l));}
  I lo=expi(I(-x.l,-x.l)),up=expi(I(x.u,x.u));
  return I(fld(cmul(SC,SC),lo.u),up.u);
}

using Poly=std::vector<I>;
static Poly padd(Poly p,const Poly&q){p.resize(std::max(p.size(),q.size()));for(size_t i=0;i<q.size();++i)p[i]=p[i]+q[i];return p;}
static Poly pmul(const Poly&p,const Poly&q){Poly r(p.size()+q.size()-1);for(size_t i=0;i<p.size();++i)for(size_t j=0;j<q.size();++j)r[i+j]=r[i+j]+p[i]*q[j];return r;}
static cpp_int binom(int n,int k){cpp_int z=1;for(int j=1;j<=k;++j)z=cmul(z,n-k+j)/j;return z;}
static Poly shiftpow(int n){Poly r(n+1);for(int j=0;j<=n;++j){cpp_int z=binom(n,j);if((n-j)&1)z=-z;r[j]=I::rat(z,1);}return r;}
static Poly hermite(const I&m){
  // H_6 matches f_m(s)=(1+s)^{m-2} at s=0 and through order five at s=1.
  // Since f_m^{(7)}<0, f_m-H_6=f_m^{(7)}(xi)s(s-1)^6/7! <=0.
  I aa=m-2,fall=I::exact(1);Poly p;
  cpp_int fac=1;
  for(int j=0;j<6;++j){
    if(j)fall=fall*(aa-(j-1)),fac=cmul(fac,j);
    I coeff=expi((m-(2+j))*LN2)*fall/I::rat(fac,1);
    Poly q=shiftpow(j);for(auto&v:q)v=v*coeff;p=padd(std::move(p),q);
  }
  I corr=I::exact(1)-p[0];Poly q=shiftpow(6);for(auto&v:q)v=v*corr;
  p=padd(std::move(p),q);
  return pmul(p,{I::exact(1),I::exact(-2),I::exact(1)});
}
static I mills(const I&q,const I&h,int n){
  // For K_q(h)=int_0^infty exp(-z^2/(2h)-qz) dz, even n is a
  // lower convergent and odd n an upper convergent.
  I w=q;
  for(int j=n-1;j>=1;--j)w=q+I::rat(j,1)/(h*w);
  return inv(w);
}

static const I PI_FULL=I(I::rat(103993,33102).l,I::rat(355,113).u);
static const I PI_LO=I(I::rat(103993,33102).l,I::rat(103993,33102).l);
// Each of the next three dyadic points is rounded upward from a proved
// rational upper bound.
static const I SIG4=I(I::rat(64494,200000).u,I::rat(64494,200000).u);
static const I SIGD=I(I::rat(28987,200000).u,I::rat(28987,200000).u);
static const I LOGD=I(I::rat(48287,250000).u,I::rat(48287,250000).u);

static I terminal_bound(const I&m,const I&a){
  I expo=a*m/I::exact(2); // m^2 h/2=am/2
  I h=a/m;
  Poly c=hermite(m);
  // A_h(m;0,infty)=sqrt(2*pi*h)e^{m^2h/2}-K_m(h).
  I total=sqrti(I::exact(2)*PI_FULL*h)*expi(expo);
  I K0lo=mills(m,h,28);
  I S(I(total.l-mills(m,h,29).u,total.u-K0lo.l));
  for(size_t i=1;i<c.size();++i){
    I q=I::exact(2*(long long)i)-m;
    I lo=mills(q,h,28),up=mills(q,h,29);
    if(c[i].l>=0)S=S+I(0,cei(cmul(c[i].u,up.u),SC));
    else if(c[i].u<=0)S=S+I(fld(cmul(c[i].l,up.u),SC),cei(cmul(c[i].u,lo.l),SC));
    else {std::cerr<<"coefficient sign failure\n";std::abort();}
  }
  I N=expi((I::exact(1)-m)*LN2)*S;
  I Ilow=PI_LO/I::exact(2)*(I::exact(1)-(I::exact(1)-m)*LOGD);
  I sig=SIG4+SIGD*m;
  I root=sqrti(h/(h+sig));
  return N/(h*Ilow*root);
}
static I line_a(const I&m,int branch){
  // 259/200+(8/5)k and 259/200+(9/5)k, k=(1-m)/m.
  if(branch==1)return I::rat(-61,200)+I::rat(8,5)/m;
  return I::rat(-101,200)+I::rat(9,5)/m;
}
static I box(const I&m,int branch){return terminal_bound(m,line_a(m,branch));}

static std::string decimal(const cpp_int&z){
  auto tos=[](cpp_int x){bool neg=x<0;if(neg)x=-x;std::string s;if(!x)s="0";while(x){s.push_back(char('0'+x%10));x/=10;}if(neg)s.push_back('-');std::reverse(s.begin(),s.end());return s;};
  cpp_int ip=z/SC,rem=z%SC;if(rem<0)rem=-rem;
  std::string out=tos(ip)+".";
  for(int i=0;i<12;++i){rem=cmul(rem,10);out.push_back(char('0'+int(rem/SC)));rem%=SC;}
  return out;
}
static bool run(int branch,int N){
  I L=branch==1?I::rat(10,13):I::rat(54,95);
  I U=branch==1?I::exact(1):I::rat(10,13);
  cpp_int worst=0,wl=0,wu=0;
  for(int j=0;j<N;++j){
    I ml=L+(U-L)*I::rat(j,N),mu=L+(U-L)*I::rat(j+1,N);
    I z=box(I(ml.l,mu.u),branch);
    if(z.u>worst)worst=z.u,wl=ml.l,wu=mu.u;
    if(z.u>=I::exact(2).l){
      std::cerr<<"FAIL branch "<<branch<<" box "<<j<<" ["<<decimal(ml.l)<<","<<decimal(mu.u)<<"] upper "<<decimal(z.u)<<"\n";
      return false;
    }
  }
  std::cout<<"PASS branch "<<branch<<" boxes "<<N<<" worst upper "<<decimal(worst)
           <<" on ["<<decimal(wl)<<","<<decimal(wu)<<"]\n";
  return true;
}
struct RBox{cpp_int ml,mu,tl,tu;int depth;};
static bool run_region(int branch){
  I ML,MU;
  if(branch==1)ML=I::rat(10,13),MU=I::exact(1);
  else if(branch==2)ML=I::rat(360,901),MU=I::rat(10,13);
  else ML=I::rat(1,5),MU=I::rat(360,901); // high-k compact tail
  std::vector<RBox> st;
  const int seed=32;
  for(int j=0;j<seed;++j){I l=ML+(MU-ML)*I::rat(j,seed),u=ML+(MU-ML)*I::rat(j+1,seed);st.push_back({l.l,u.u,0,SC,0});}
  long long seen=0,passed=0;int maxdepth=0;cpp_int worst=0,wml=0,wmu=0,wtl=0,wtu=0;
  while(!st.empty()){
    RBox b=st.back();st.pop_back();++seen;maxdepth=std::max(maxdepth,b.depth);
    I m(b.ml,b.mu),t(b.tl,b.tu),aa;
    if(branch<=2)aa=I::rat(1,2)+(line_a(m,branch)-I::rat(1,2))*t;
    else aa=I::rat(1,2)+(I::exact(4)-I::rat(1,2))*t;
    I z=terminal_bound(m,aa);
    if(z.u<I::exact(2).l){if(z.u>worst)worst=z.u,wml=b.ml,wmu=b.mu,wtl=b.tl,wtu=b.tu;++passed;continue;}
    if(b.depth>=34){
      std::cerr<<"REGION FAIL branch "<<branch<<" depth "<<b.depth<<" m["<<decimal(b.ml)<<","<<decimal(b.mu)<<"] t["<<decimal(b.tl)<<","<<decimal(b.tu)<<"] upper "<<decimal(z.u)<<"\n";return false;
    }
    // Bisect the wider normalized coordinate.  Endpoints remain dyadic.
    if(csub(b.mu,b.ml)>csub(b.tu,b.tl)/2){cpp_int md=cadd(b.ml,b.mu)/2;st.push_back({md,b.mu,b.tl,b.tu,b.depth+1});st.push_back({b.ml,md,b.tl,b.tu,b.depth+1});}
    else {cpp_int md=cadd(b.tl,b.tu)/2;st.push_back({b.ml,b.mu,md,b.tu,b.depth+1});st.push_back({b.ml,b.mu,b.tl,md,b.depth+1});}
  }
  std::cout<<"PASS REGION branch "<<branch<<" visited "<<seen<<" leaves "<<passed<<" maxdepth "<<maxdepth<<" largest accepted upper "<<decimal(worst)<<" at m["<<decimal(wml)<<","<<decimal(wmu)<<"] t["<<decimal(wtl)<<","<<decimal(wtu)<<"]\n";
  return true;
}
static I p_line(const I&m,int P){
 if(P==3)return I::rat(-7,10)+I::exact(3)/m;       // 23/10+3k
 if(P==4)return I::rat(-1,2)+I::rat(7,2)/m;       // 3+(7/2)k
 return I::rat(-2,5)+I::exact(4)/m;               // 18/5+4k
}
static bool run_pregion(int P,bool tail){
 I alpha=P==3?I::rat(2,3):(P==4?I::rat(3,4):I::rat(4,5));
 I cap=P==3?I::rat(30,47):(P==4?I::rat(7,9):I::rat(10,11));
 I ML=tail?I::rat(1,5):cap,MU=tail?cap:I::exact(1);
 std::vector<RBox> st;const int seed=32;
 for(int j=0;j<seed;++j){I l=ML+(MU-ML)*I::rat(j,seed),u=ML+(MU-ML)*I::rat(j+1,seed);st.push_back({l.l,u.u,0,SC,0});}
 long long seen=0,passed=0;int maxdepth=0;cpp_int worst=0,wml=0,wmu=0,wtl=0,wtu=0;
 while(!st.empty()){
  RBox b=st.back();st.pop_back();++seen;maxdepth=std::max(maxdepth,b.depth);
  I m(b.ml,b.mu),t(b.tl,b.tu);
  I top=tail?I::exact(4):p_line(m,P);
  I aa=alpha+(top-alpha)*t;
  I z=terminal_bound(m,aa);
  if(z.u<I::exact(P).l){if(z.u>worst)worst=z.u,wml=b.ml,wmu=b.mu,wtl=b.tl,wtu=b.tu;++passed;continue;}
  if(b.depth>=34){std::cerr<<"PREGION FAIL P "<<P<<" tail "<<tail<<" depth "<<b.depth<<" m["<<decimal(b.ml)<<","<<decimal(b.mu)<<"] t["<<decimal(b.tl)<<","<<decimal(b.tu)<<"] upper "<<decimal(z.u)<<"\n";return false;}
  if(csub(b.mu,b.ml)>csub(b.tu,b.tl)/2){cpp_int md=cadd(b.ml,b.mu)/2;st.push_back({md,b.mu,b.tl,b.tu,b.depth+1});st.push_back({b.ml,md,b.tl,b.tu,b.depth+1});}
  else {cpp_int md=cadd(b.tl,b.tu)/2;st.push_back({b.ml,b.mu,md,b.tu,b.depth+1});st.push_back({b.ml,b.mu,b.tl,md,b.depth+1});}
 }
 std::cout<<"PASS PREGION P "<<P<<" tail "<<tail<<" visited "<<seen<<" leaves "<<passed<<" maxdepth "<<maxdepth<<" largest accepted upper "<<decimal(worst)<<" at m["<<decimal(wml)<<","<<decimal(wmu)<<"] t["<<decimal(wtl)<<","<<decimal(wtu)<<"]\n";
 return true;
}
int main(int ac,char**av){
 if(ac>1&&std::string(av[1])=="region")return run_region(1)&&run_region(2)&&run_region(3)?0:1;
 if(ac>1&&std::string(av[1]).rfind("region",0)==0&&std::string(av[1]).size()>6)return run_region(std::stoi(std::string(av[1]).substr(6)))?0:1;
 if(ac>1&&std::string(av[1]).rfind("pregion",0)==0){int p=std::stoi(std::string(av[1]).substr(7));return run_pregion(p,false)&&run_pregion(p,true)?0:1;}
 int N=ac>1?std::stoi(av[1]):1200;return run(1,N)&&run(2,2*N)?0:1;
}
\end{lstlisting}

\subsubsection{Terminal interval-arithmetic module}

\begin{lstlisting}[language={Python}]
#!/usr/bin/env python3
"""Fast exact certificate for the two terminal affine barriers.

This is an integer-only, outward-rounded implementation of the analytic
bounds printed with the proof.  A real interval [L/2^B,U/2^B]
is stored as the pair of integers (L,U).  Thus every comparison made by
the certificate is an exact integer comparison.  There is no hardware
floating-point arithmetic in the certification path.
"""
from fractions import Fraction as Q
from math import comb,isqrt
import sys

B=180; S=1<<B
def fd(n,d=S): return n//d
def cu(n,d=S): return -((-n)//d)

class F:
    __slots__=("l","u")
    def __init__(self,x=0,y=None,raw=False):
        if raw:self.l,self.u=x,(x if y is None else y)
        else:
            x=Q(x); y=x if y is None else Q(y)
            self.l=fd(x.numerator*S,x.denominator)
            self.u=cu(y.numerator*S,y.denominator)
        assert self.l<=self.u
    @staticmethod
    def raw(l,u):return F(l,u,True)
    def __add__(self,o):
        o=o if isinstance(o,F) else F(o);return F.raw(self.l+o.l,self.u+o.u)
    __radd__=__add__
    def __neg__(self):return F.raw(-self.u,-self.l)
    def __sub__(self,o):return self+(-o if isinstance(o,F) else -F(o))
    def __rsub__(self,o):return F(o)-self
    def __mul__(self,o):
        o=o if isinstance(o,F) else F(o)
        z=(self.l*o.l,self.l*o.u,self.u*o.l,self.u*o.u)
        return F.raw(fd(min(z)),cu(max(z)))
    __rmul__=__mul__
    def inv(self):
        assert self.l>0 or self.u<0
        return F.raw(fd(S*S,self.u),cu(S*S,self.l))
    def __truediv__(self,o):return self*(o if isinstance(o,F) else F(o)).inv()
    def __rtruediv__(self,o):return F(o)/self

def sqrtF(x):
    x=x if isinstance(x,F) else F(x);assert x.l>=0
    lo=isqrt(x.l*S); hi=isqrt(x.u*S)
    if hi*hi<x.u*S:hi+=1
    return F.raw(lo,hi)

def exp_pos_endpoint(x,N=42):
    """Return scaled-integer lower/upper bounds for exp(x/S), x>=0."""
    assert x>=0 and x < (N+2)*S
    tl=tu=S; sl=su=S
    for n in range(1,N+1):
        tl=fd(tl*x,S*n);tu=cu(tu*x,S*n)
        sl+=tl;su+=tu
    nxt=cu(tu*x,S*(N+1))
    tail=cu(nxt*S*(N+2),S*(N+2)-x)
    return sl,su+tail
def expF(x):
    x=x if isinstance(x,F) else F(x)
    if x.l>=0:return F.raw(exp_pos_endpoint(x.l)[0],exp_pos_endpoint(x.u)[1])
    if x.u<=0:
        yl,yu=exp_pos_endpoint(-x.u)[0],exp_pos_endpoint(-x.l)[1]
        return F.raw(fd(S*S,yu),cu(S*S,yl))
    return F.raw(fd(S*S,exp_pos_endpoint(-x.l)[1]),exp_pos_endpoint(x.u)[1])

class D:
    __slots__=("v","d")
    def __init__(self,v=0,d=0):
        self.v=v if isinstance(v,F) else F(v);self.d=d if isinstance(d,F) else F(d)
    def __add__(self,o):
        o=o if isinstance(o,D) else D(o);return D(self.v+o.v,self.d+o.d)
    __radd__=__add__
    def __neg__(self):return D(-self.v,-self.d)
    def __sub__(self,o):return self+(-o if isinstance(o,D) else -D(o))
    def __rsub__(self,o):return D(o)-self
    def __mul__(self,o):
        o=o if isinstance(o,D) else D(o);return D(self.v*o.v,self.d*o.v+self.v*o.d)
    __rmul__=__mul__
    def inv(self):return D(self.v.inv(),-self.d/(self.v*self.v))
    def __truediv__(self,o):return self*(o if isinstance(o,D) else D(o)).inv()
    def __rtruediv__(self,o):return D(o)/self
def expD(x):
    x=x if isinstance(x,D) else D(x);v=expF(x.v);return D(v,v*x.d)
def sqrtD(x):
    x=x if isinstance(x,D) else D(x);v=sqrtF(x.v);return D(v,x.d/(2*v))

def mills(q,h,n):
    q=q if isinstance(q,D) else D(q);h=h if isinstance(h,D) else D(h);w=q
    for j in range(n-1,0,-1):w=q+Q(j)/(h*w)
    return 1/w

def ln2_interval(N=42):
    z=Q(0)
    for j in range(N+1):z+=Q(2,(2*j+1)*3**(2*j+1))
    tail=Q(18,(2*N+3)*8*3**(2*N+3))
    return F(z,z+tail)
LN2=ln2_interval()

def atan_recip_bounds(q,n):
    """Alternating-series enclosure of atan(1/q), through term n."""
    z=Q(0)
    for j in range(n+1):z+=(-1 if j&1 else 1)*Q(1,(2*j+1)*q**(2*j+1))
    nxt=Q(1,(2*n+3)*q**(2*n+3))
    return (z,z+nxt) if n&1 else (z-nxt,z)
_a5l,_a5u=atan_recip_bounds(5,8)
_a239l,_a239u=atan_recip_bounds(239,2)
# Machin's identity pi=16 atan(1/5)-4 atan(1/239).
PI_LO=16*_a5l-4*_a239u
PI_UP=16*_a5u-4*_a239l
SIG4_UP=Q(64494,200000)
SIGD_UP=Q(28987,200000)
LOG2MHALF_UP=Q(48287,250000)
assert SIG4_UP >= PI_UP*PI_UP/12-Q(1,2)
assert SIGD_UP >= PI_UP*PI_UP/6-Q(3,2)
assert LOG2MHALF_UP*S >= LN2.u-S//2

def V(m,a):
    """Differentiable explicit upper bound for the terminal P_T."""
    m=m if isinstance(m,D) else D(m);a=a if isinstance(a,D) else D(a)
    h=a/m;n=10;c=[];poch=D(1);fac=1
    for j in range(n):
        if j:poch=poch*((j+1)-m);fac*=j
        c.append(expD((m-(2+j))*D(LN2))*poch/fac)
    rem=D(1)
    for z in c:rem=rem-z
    mags=[]
    for i in range(n+3):
        z=D(0)
        for j,w in enumerate(c):
            if j+2>=i:z=z+comb(j+2,i)*w
        mags.append(z+comb(n+2,i)*rem)
    total=sqrtD(2*PI_UP*h)*expD(a*m/2)
    val=total-mills(m,h,20)
    for i in range(1,n+3):
        q=2*i-m
        val=val-mags[i]*mills(q,h,20) if i&1 else val+mags[i]*mills(q,h,21)
    num=expD((1-m)*D(LN2))*val
    Ilow=Q(PI_LO,2)*(1-(1-m)*LOG2MHALF_UP)
    sig=SIG4_UP+SIGD_UP*m
    return num/(h*Ilow*sqrtD(h/(h+sig)))

def aline(m,b):return Q(-61,200)+Q(8,5)/m if b==1 else Q(-101,200)+Q(9,5)/m
def centered_box(ml,mu,b):
    mc=(ml+mu)/2
    v0=V(D(F(mc)),D(F(aline(mc,b)))).v
    m=D(F(ml,mu),F(1))
    a=D(Q(-61,200))+Q(8,5)/m if b==1 else D(Q(-101,200))+Q(9,5)/m
    der=V(m,a).d
    rad=F((mu-ml)/2)
    err=rad*F.raw(min(abs(der.l),abs(der.u)) if der.l*der.u>0 else 0,
                  max(abs(der.l),abs(der.u)))
    return v0.u+err.u

def run(b,N):
    l,u=(Q(10,13),Q(1)) if b==1 else (Q(1080,1903),Q(10,13))
    worst=0;where=None
    for j in range(N):
        ml=l+(u-l)*j/N;mu=l+(u-l)*(j+1)/N
        z=centered_box(ml,mu,b)
        if z>worst:worst,where=z,j
        if z>=2*S:
            print("FAIL",b,j,"gap numerator",2*S-z);return False
    print("PASS",b,"boxes",N,"gap",Q(2*S-worst,S),"worst box",where)
    return True

if __name__=='__main__':
    n=int(sys.argv[1]) if len(sys.argv)>1 else 64
    ok=run(1,n) and run(2,2*n)
    raise SystemExit(0 if ok else 1)
\end{lstlisting}

\subsubsection{Terminal $A=4$ certificate}

\begin{lstlisting}[language={Python}]
#!/usr/bin/env python3
"""Integer-only terminal boundary certificates at center size a=4.

For m in [1/10,1], this proves the stated upper bounds for the explicit
majorant V(m,4).  The interval engine and V are imported from the fully
exact terminal-affine certificate; no floating-point arithmetic enters a
sign decision.
"""
from fractions import Fraction as Q
import certify_terminal_affine_fixed as E


def centered_box(ml, mu):
    mc = (ml + mu) / 2
    v0 = E.V(E.D(E.F(mc)), E.D(4)).v
    m = E.D(E.F(ml, mu), E.F(1))
    der = E.V(m, E.D(4)).d
    err = E.F((mu - ml) / 2) * E.F.raw(
        0, max(abs(der.l), abs(der.u)))
    return v0.u + err.u


def run(name, target, left, right, boxes):
    worst = 0
    where = None
    for j in range(boxes):
        ml = left + (right - left) * j / boxes
        mu = left + (right - left) * (j + 1) / boxes
        z = centered_box(ml, mu)
        if z > worst:
            worst, where = z, j
        if z >= target * E.S:
            print("FAIL", name, "box", j,
                  "gap numerator", target * E.S - z)
            return False
    print("PASS", name, "boxes", boxes,
          "gap", Q(target * E.S - worst, E.S),
          "worst box", where)
    return True


if __name__ == "__main__":
    # Complementary ranges to the affine P=3,4,5 barriers, followed by
    # the all-k boundary needed for every P>=6.
    cases = [
        ("P2-tail", 2, Q(1, 20), Q(1, 5), 128),
        ("P3-tail", 3, Q(1, 10), Q(30, 47), 96),
        ("P4-tail", 4, Q(1, 10), Q(7, 9), 96),
        ("P5-tail", 5, Q(1, 10), Q(10, 11), 128),
        ("P6-all", 6, Q(1, 10), Q(1), 128),
    ]
    ok = True
    for args in cases:
        ok = run(*args) and ok
    raise SystemExit(0 if ok else 1)
\end{lstlisting}

\subsubsection{Density-exclusion certificate}

\begin{lstlisting}[language={Python}]
#!/usr/bin/env python3
"""Exact Bernstein certificate for the corrected X,Z shifted-tail scalar bound.

No floating point arithmetic is used.  Exponent pairs are (a,k).  This file
certifies only the direct shifted-tail objective; it does not use Bfull or
drop the factor E(1+j).
"""
from fractions import Fraction as Q
from math import comb

def add(p,q):
 r=dict(p)
 for e,c in q.items():r[e]=r.get(e,Q(0))+c
 return {e:c for e,c in r.items() if c}
def scale(p,c):return {e:c*v for e,v in p.items() if c*v}
def sub(p,q):return add(p,scale(q,-1))
def mul(p,q):
 r={}
 for (i,j),c in p.items():
  for (u,v),d in q.items():r[i+u,j+v]=r.get((i+u,j+v),Q(0))+c*d
 return {e:c for e,c in r.items() if c}
def power(p,n):
 r={(0,0):Q(1)}
 while n:
  if n&1:r=mul(r,p)
  p=mul(p,p);n//=2
 return r

one={(0,0):Q(1)};a={(1,0):Q(1)};k={(0,1):Q(1)}
A4=sub(scale(a,4),one)                       # 4(a-c), c=1/4
Rd=mul(a,add(one,k)); Rn=add(one,mul(a,k))  # rho=Rn/Rd
C=add(Rn,scale(Rd,6)); E=add(Rn,scale(Rd,7))

def endpoint(kind):
 if kind=='density':
  U=Rd;V=add(scale(Rd,4),Rn)                 # delta=1/(4+rho)
 elif kind=='phi':
  b=add(scale(k,3),scale(one,2))
  G=add(add(scale(power(k,2),10),scale(k,19)),scale(one,6))
  U=mul(A4,b);V=scale(mul(a,G),4)
 elif kind=='J':
  U=add(add(scale(power(a,2),8),scale(a,2)),scale(one,-2))
  U=add(U,add(scale(mul(power(a,2),k),-12),scale(mul(a,k),-4)))
  V=mul(a,add(add(scale(a,8),scale(one,8)),add(scale(mul(a,k),15),scale(k,5))))
 else:raise ValueError(kind)
 # beta_t=B/(C U^2).
 B=sub(sub(mul(mul(E,U),V),mul(Rd,power(V,2))),mul(C,power(U,2)))
 # H=e+delta/[a(1+k)]=H4/(8 V Rd).
 # The last term is 8U (not 8UV), because delta/[a(1+k)]=U/(V Rd).
 H4=add(add(mul(mul(A4,V),Rd),scale(mul(mul(a,U),Rd),-12)),scale(U,8))
 # .5(a-c) beta_t H > 1/4 iff the first polynomial is positive.
 active=sub(mul(mul(A4,B),H4),scale(mul(mul(mul(C,power(U,2)),V),Rd),16))
 # Same target with beta=4/5.
 const=sub(mul(A4,H4),scale(mul(V,Rd),20))
 beta45=sub(scale(B,5),scale(mul(C,power(U,2)),4))
 return U,V,B,H4,active,const,beta45

def subst(p,lo,hi,with_t,slope=Q(9,5)):
 # k=lo+(hi-lo)s; a=259/200+slope*k + (4-line)t.
 ks={(0,0):lo,(1,0):hi-lo}
 line=add({(0,0):Q(259,200)},scale(ks,slope))
 aa=line if not with_t else add(line,mul(sub({(0,0):Q(4)},line),{(0,1):Q(1)}))
 out={}
 for (ia,ik),c in p.items():out=add(out,scale(mul(power(aa,ia),power(ks,ik)),c))
 return out
def bern(p):
 nx=max((i for i,j in p),default=0);ny=max((j for i,j in p),default=0)
 B=[]
 for u in range(nx+1):
  for v in range(ny+1):
   B.append(sum(c*Q(comb(u,i),comb(nx,i))*Q(comb(v,j),comb(ny,j))
                for (i,j),c in p.items() if i<=u and j<=v))
 return B,nx,ny
def certify(name,p,lo,hi,with_t=True,slope=Q(9,5)):
 q=subst(p,lo,hi,with_t,slope);B,nx,ny=bern(q)
 assert min(B)>0,(name,min(B),max(B),nx,ny)
 print(name,'domain',lo,hi,'degrees',nx,ny,'min Bernstein',min(B))

if __name__=='__main__':
 J=endpoint('J');du=endpoint('density')
 Fden=add(scale(Rd,4),Rn)
 density_gap=sub(mul(J[0],Fden),mul(J[1],Rd))
 certify('low-k density exclusion slope 8/5',density_gap,Q(0),Q(3,10),True,Q(8,5))
 certify('low-k density exclusion slope 9/5',density_gap,Q(3,10),Q(8,25))
 # For 8/25<=k<=7/20 the proved center bound puts delta above delta_J.
 certify('J beta active',J[6],Q(8,25),Q(7,20))
 certify('J lower endpoint target',J[4],Q(8,25),Q(7,20))
 # The density endpoint controls the other end of every active interval.
 certify('density upper endpoint target',du[4],Q(8,25),Q(541,360))
 # For k>=7/20 no positive lower bound on delta is used.  Let D be the
 # value at which the constant-beta target equals c.  Since beta_t is
 # increasing, beta_t(D)>4/5 puts the unique crossover strictly before D.
 Uc=mul(sub(power(A4,2),scale(one,20)),Rd)
 Vc=scale(mul(A4,sub(scale(mul(a,Rd),3),scale(one,2))),4)
 Bc=sub(sub(mul(mul(E,Uc),Vc),mul(Rd,power(Vc,2))),mul(C,power(Uc,2)))
 crossgap=sub(scale(Bc,5),scale(mul(C,power(Uc,2)),4))
 certify('crossover before constant-target zero',crossgap,Q(7,20),Q(541,360))
\end{lstlisting}

\subsubsection{Higher-$P$ scalar certificate}

\begin{lstlisting}[language={Python}]
#!/usr/bin/env python3
"""Exact rational certificate for the uniform integer P>=3 scalar closure.

This checks the algebraic part after the three terminal affine bounds
  P=3: a >= 23/10+3k,
  P=4: a >= 3+(7/2)k,
  P>=5: a >= 18/5+4k.
All arithmetic and Bernstein coefficients are exact Fractions.
"""
from fractions import Fraction as Q
from math import comb

def add(p,q):
 r=dict(p)
 for e,c in q.items():r[e]=r.get(e,Q(0))+c
 return {e:c for e,c in r.items() if c}
def scale(p,c):return {e:c*v for e,v in p.items() if c*v}
def sub(p,q):return add(p,scale(q,-1))
def mul(p,q):
 r={}
 for (i,j),c in p.items():
  for (u,v),d in q.items():r[i+u,j+v]=r.get((i+u,j+v),Q(0))+c*d
 return {e:c for e,c in r.items() if c}
def power(p,n):
 r={(0,0):Q(1)}
 while n:
  if n&1:r=mul(r,p)
  p=mul(p,p);n//=2
 return r

one={(0,0):Q(1)};a={(1,0):Q(1)};k={(0,1):Q(1)}
Rd=mul(a,add(one,k)); Rn=add(one,mul(a,k)); Fden=add(scale(Rd,4),Rn)
C=add(Rn,scale(Rd,6));E=add(Rn,scale(Rd,7))
b=add(scale(k,3),scale(one,2))
G=add(add(scale(power(k,2),10),scale(k,19)),scale(one,6))
N=add(power(Rn,2),mul(power(a,2),mul(k,add(one,k))))

def density_gap(c):
 # J1=2a(a+1)/(1+3a), delta_J=(J1-2c-ak)/(J1+5ak/4).
 den3=add(one,scale(a,3))
 base=scale(mul(a,add(a,one)),2)
 U=sub(sub(base,scale(den3,2*c)),mul(mul(a,k),den3))
 V=add(base,scale(mul(mul(a,k),den3),Q(5,4)))
 # delta_J > 1/(4+rho)=Rd/Fden.
 return sub(mul(U,Fden),mul(V,Rd))

def corrected_constant_target(c):
 # Direct XZ shifted-tail bound:
 # T >= (2/5)(a-c){e+delta/[a(1+k)]},
 # evaluated at delta_phi=(a-c)b/(aG).
 ac=sub(a,scale(one,c))
 U=mul(ac,b);V=mul(a,G)
 # H=e+delta/Rd=H2/(2 V Rd).
 H2=add(add(mul(mul(ac,V),Rd),scale(mul(mul(a,U),Rd),-3)),scale(U,2))
 # (2/5)(a-c)H>c iff (a-c)H2-5c V Rd>0.
 return sub(mul(ac,H2),scale(mul(V,Rd),5*c))

def corrected_density_target(c):
 ac=sub(a,scale(one,c));U=Rd;V=Fden
 Bt=sub(sub(mul(mul(E,U),V),mul(Rd,power(V,2))),mul(C,power(U,2)))
 H2=add(add(mul(mul(ac,V),Rd),scale(mul(mul(a,U),Rd),-3)),scale(U,2))
 # .5(a-c) beta_t H>c, with H=H2/(2VRd), beta_t=Bt/(CU^2).
 return sub(mul(mul(ac,Bt),H2),scale(mul(mul(mul(C,power(U,2)),V),Rd),4*c))

def crossover_gap(c):
 ac=sub(a,scale(one,c))
 # D={((a-c)^2-5c)Rd}/{(a-c)(3aRd-2)} is the zero of
 # the constant-beta target.
 U=mul(sub(power(ac,2),scale(one,5*c)),Rd)
 V=mul(ac,sub(scale(mul(a,Rd),3),scale(one,2)))
 Bt=sub(sub(mul(mul(E,U),V),mul(Rd,power(V,2))),mul(C,power(U,2)))
 return sub(scale(Bt,5),scale(mul(C,power(U,2)),4))

def subst(p,lo,hi,A,B,with_t):
 ks={(0,0):lo,(1,0):hi-lo};L=add({(0,0):A},scale(ks,B))
 aa=L if not with_t else add(L,mul(sub({(0,0):Q(4)},L),{(0,1):Q(1)}))
 out={}
 for (ia,ik),c in p.items():out=add(out,scale(mul(power(aa,ia),power(ks,ik)),c))
 return out
def bern(p):
 nx=max((i for i,j in p),default=0);ny=max((j for i,j in p),default=0);out=[]
 for u in range(nx+1):
  for v in range(ny+1):
   out.append(sum(c*Q(comb(u,i),comb(nx,i))*Q(comb(v,j),comb(ny,j))
                  for (i,j),c in p.items() if i<=u and j<=v))
 return out,nx,ny
def cert(name,p,lo,hi,A,B,with_t):
 z,nx,ny=bern(subst(p,lo,hi,A,B,with_t));assert min(z)>0,(name,min(z))
 print(name,'degrees',nx,ny,'count',len(z),'minimum',min(z))

if __name__=='__main__':
 # If x<3, density gives a<4.  For P=3, k<=1/4 is impossible.
 cert('P3 density exclusion',density_gap(Q(1,3)),Q(0),Q(1,4),Q(23,10),Q(3),True)
 # On the complementary interval extend delta down to zero.  The active
 # product has no interior minimum.  Its density endpoint is positive,
 # and beta_t at the zero of the constant target is already >4/5, so the
 # beta crossover also has positive target.
 cert('P3 corrected XZ density endpoint',corrected_density_target(Q(1,3)),Q(1,4),Q(17,30),Q(23,10),Q(3),True)
 cert('P3 corrected XZ crossover',crossover_gap(Q(1,3)),Q(1,4),Q(17,30),Q(23,10),Q(3),True)
 # P=4: the full possible k interval is excluded by center+density.
 cert('P4 density exclusion',density_gap(Q(3,8)),Q(0),Q(2,7),Q(3),Q(7,2),True)
 # P>=5: use the P=5 terminal line and the worst c=1/2.  Since the
 # density gap decreases with c, this covers every c=(P-1)/(2P)<=1/2.
 cert('P>=5 density exclusion',density_gap(Q(1,2)),Q(0),Q(1,10),Q(18,5),Q(4),True)
\end{lstlisting}

\subsection{The marginal-crossing computation}

\subsection{Exact certificates and reproducibility}
\label{mt:app:compact-certificates}

This appendix records every finite verification used in the proof.  It is
included to distinguish the certificates from numerical experiments.  In
the description below, the second coordinate of a box is
$m^2\xi'(q)$.  Thus every displayed rectangle is a continuous parameter
region, including its boundary; it is not a list of sampled points.

\subsubsection{The regions that are certified}

For $p=4$, the mass-exclusion certificate covers
\begin{align*}
 &[1/5,7/10]\times[3/2,5],
 &&[7/10,1]\times[3/2,7/2],\\
 &[7/10,1]\times[4/3,3/2],
 &&[1/5,7/10]\times[5,12],\\
 &[7/10,1]\times[7/2,12].
\end{align*}
For $p=5$, it covers
\begin{align*}
 &[1/6,24/25]\times[5/2,6],
 &&[24/25,1]\times[5/2,4],\\
 &[1/6,24/25]\times[6,16],
 &&[24/25,1]\times[4,16].
\end{align*}
These rectangles are precisely the compact regions left by the analytic
lower bound on $m^2\xi'(q)$.  The analytic tail argument in the proof
excludes $m^2\xi'(q)\geq4(p-1)$.  Consequently the listed rectangles,
together with that argument, cover every case with
$m\geq1/(p+1)$ for $p=4,5$.  The programs reject each box either because
the marginal equation has disjoint interval enclosures, or because one of
the two necessary entropy inequalities has a strictly positive lower
endpoint.  The compact part is run by
\path{run_threshold_split_exact.py} and
\path{run_threshold_tails_exact.py}.  For $p\geq6$, no interval covering
is used: \path{verify_mass_exclusion_Pge5.py} checks, over the rationals,
the polynomial identities and coefficient expansions used in the
analytic proof.

The exceptional case $p=3$ uses four coverings.  The contact equation is
excluded on
\begin{equation*}
 [1/10,7/20]\times[1/2,5/2],\qquad
 [7/20,1/2]\times[1/2,189/100].
\end{equation*}
On
\begin{equation*}
 [7/20,1/2]\times[189/100,5/2]
\end{equation*}
the program proves that the logarithmic lower bound for the exact entropy
identity is strictly positive whenever the contact equation has not
already been excluded.  On
\begin{equation*}
 [1/2,1]\times[1/100,5/2]
\end{equation*}
it does the same with the fourth-order lower bound.  These are respectively
the \texttt{separator}, \texttt{low}, and \texttt{high} cases of
\path{certify_p3_entry_combined_exact.py}.  The proof outside these boxes
is analytic: the preliminary inequality gives
$m^2\xi'(q)>1/2$ when $m\leq1/2$; the range $m\leq1/10$ is excluded
directly; and, for $m\geq1/2$, the range
$0<m^2\xi'(q)<1/100$ is excluded by the Brascamp--Lieb estimate.

The high-mass part of the positive-index argument is a separate contact
certificate on
\begin{equation*}
 [9/25,1]\times[5/2,10].
\end{equation*}
It is run by \path{certify_p3_large_mass_contact_exact.py}; the range
$m^2\xi'(q)\geq10$ is excluded analytically.  The constants used after
this exclusion are checked by
\path{verify_p3_positive_index_inputs_exact.py}.  Finally,
\path{p3_dlow_symbolic.py} proves all remaining polynomial signs by
Bernstein expansion.  Its substitution maps the closed unit cube to
\begin{equation*}
 0\leq m\leq \frac{9}{25},\qquad
 0\leq \frac{1}{\xi'(q)}\leq \frac{2m^2}{5},\qquad
 0\leq x_3\leq1,
\end{equation*}
where $x_3$ is the dummy coordinate in
\eqref{mt:eq:complete-p3-cube}.  Every
Bernstein coefficient printed by the program is strictly positive.

One pruning step in the $p=3$ entry certificate depends on the following
two continuum inequalities:
\begin{align*}
 \frac{10}{13}\leq m\leq1,\qquad
 &m\xi'(q)=-\frac{61}{200}+\frac{8}{5m},\\
 \frac{1080}{1903}\leq m\leq\frac{10}{13},\qquad
 &m\xi'(q)=-\frac{101}{200}+\frac{9}{5m}.
\end{align*}
The program \path{certify_terminal_affine_fixed.py} divides the first
interval into $64$ rational boxes and the second into $128$ rational
boxes.  It proves on both lines that the rigorous upper bound for the
terminal ratio is strictly smaller than $2$.  The two certified gaps are,
respectively,
\begin{align*}
 &2^{-180}(15241255008302228531374906625118347899689757891379648),\\
 &2^{-180}(6678747921395679929293531355794639244537960424215201).
\end{align*}
The completed terminal certificate is printed below and is used to prune
the high-mass entry boxes.

\subsubsection{Arithmetic and acceptance rules}

All interval endpoints are instances of Python's
\texttt{fractions.Fraction}.  Addition, subtraction, multiplication,
division, and differentiation of interval expressions are followed by
lower and upper rounding on a fixed dyadic lattice.  The $p=4,5$
covering uses denominator $2^{90}$, the $p=3$ contact and entry
coverings use denominator $2^{112}$, and the terminal affine certificate
uses denominator $2^{180}$.  The algebraic and Bernstein programs use
unrounded rational arithmetic.

The transcendental enclosures are also rational.  The programs use the
$\operatorname{arctanh}(1/3)$ series for $\log2$, Machin's alternating
series for $\pi$, Taylor series with explicit geometric remainders for
exponentials and logarithms, integer bracketing for square roots,
enveloping Stirling expansions for the beta integrals, and consecutive
even and odd Stieltjes convergents for shifted Gaussian tails.  Alternating
series are truncated only in the direction fixed by their next term.
Thus every primitive supplied to the interval calculation contains its
exact real value.

Each seed rectangle is subdivided dyadically.  A box is accepted only
when a displayed necessary equation has disjoint lower and upper
enclosures, or a displayed necessary inequality has a strictly positive
rational lower endpoint.  Otherwise it is subdivided.  Reaching the
maximum depth with an unresolved box returns failure.  A successful run
therefore proves the assertion on the union of the closed boxes.  Decimal
conversions, where printed, are diagnostic only and never enter an
acceptance decision.  No quadrature, random sampling, binary
floating-point comparison, or tolerance test occurs in the proof path.

\subsubsection{Verifier manifest and commands}

The verification requires only Python~3 and its standard library.  All
listings in this subsection must be saved together in
\path{work/marginal}; they must not be mixed with the sources in the
other certificate subsections.  No external
Python package, computer algebra system, or compiled library is used.
The local import dependencies are
\begin{center}
\footnotesize
\begin{tabular}{@{}>{\raggedright\arraybackslash}p{0.43\linewidth}
                    >{\raggedright\arraybackslash}p{0.50\linewidth}@{}}
\toprule
top-level file & local source dependencies\\
\midrule
\path{verify_mass_exclusion_Pge5.py} & none\\
\path{run_threshold_split_exact.py} &
 \path{certify_fixedr_threshold_exact.py},
 \path{certify_p3_entry_exact.py},
 \path{certify_terminal_affine_exact.py}\\
\path{run_threshold_tails_exact.py} &
 \path{run_threshold_split_exact.py} and its dependencies\\
\path{certify_p3_entry_combined_exact.py} &
 \path{certify_fixedr_threshold_exact.py},
 \path{certify_p3_entry_exact.py},
 \path{certify_terminal_affine_exact.py}\\
\path{certify_p3_large_mass_contact_exact.py} &
 \path{certify_p3_entry_combined_exact.py} and its dependencies\\
\path{verify_p3_positive_index_inputs_exact.py} &
 \path{certify_p3_entry_exact.py},
 \path{certify_terminal_affine_exact.py}\\
\path{certify_terminal_affine_fixed.py} & none\\
\path{p3_dlow_symbolic.py} & none\\
\bottomrule
\end{tabular}
\end{center}
Here ``none'' means no local import; the Python standard library is still
used.  From the repository root, the complete serial audit is
\begin{verbatim}
python3 work/marginal/verify_mass_exclusion_Pge5.py
python3 work/marginal/run_threshold_split_exact.py
python3 work/marginal/run_threshold_tails_exact.py
python3 work/marginal/certify_terminal_affine_fixed.py 64
python3 work/marginal/certify_p3_entry_combined_exact.py separator
python3 work/marginal/certify_p3_entry_combined_exact.py low
python3 work/marginal/certify_p3_entry_combined_exact.py 0 16
python3 work/marginal/certify_p3_large_mass_contact_exact.py
python3 work/marginal/verify_p3_positive_index_inputs_exact.py
python3 work/marginal/p3_dlow_symbolic.py
\end{verbatim}
The long high-mass entry command may instead be split into any disjoint
collection of half-open strip ranges whose union is $[0,16)$; every range
must return \texttt{PASS-HIGH}.  The other required terminal tags are
\texttt{PASS-SEPARATOR}, \texttt{PASS-LOW},
\texttt{PASS-LARGE-MASS}, and
\texttt{PASS-P3-SECTION3-INPUTS}.  The terminal affine command must print
\texttt{PASS 1 boxes 64} and \texttt{PASS 2 boxes 128}.  The Bernstein
audit succeeds precisely when every reported count labeled
\texttt{nonpositive} is zero.  The two threshold runners print
\texttt{PASS} for every rectangle listed above, and the polynomial audit
prints \texttt{PASS: exact verification of (M9), (M11), and (M20)}.

The reported run used Python 3.14.6.  The frozen listings have the
following SHA--256 digests:
\begingroup\tiny
\begin{verbatim}
ec44b86210a1808a8bd36dc6bf6deb85f016c6bb65d8cb008548bb4652ed8f21  verify_mass_exclusion_Pge5.py
8e524e1ccb52fa33e2bf826234ec30bfb72c9fa399cb2f619ed10b893dcea3ea  run_threshold_split_exact.py
17120f6983efa9cd5036ba43c704c70ffeb3c6df9e301518ab87ab633edc3e61  run_threshold_tails_exact.py
f358c3169372a9c2d3342e99b641a4fe597126daf54ffcaf204f104b4c92aa1b  certify_fixedr_threshold_exact.py
f8755b1144b86bf7e52b5de819fbf235288ca23bf78c4a94531f5fe488013486  certify_p3_entry_exact.py
35c6f696643a182e62b95b34e456056cebcd6098a323c4f458e33710f403b92c  certify_terminal_affine_exact.py
7bfdb70b289c205ffe2f3f6a89e2ec9bb1c373478eb94b71ae3182a812e7c8a8  certify_p3_entry_combined_exact.py
30d536c093b7aeecc4f548e927d1f8d71ce814c58cd53d45c12e7a53cb940327  certify_p3_large_mass_contact_exact.py
96e86a01e4ab1ac81deebae306ec4ad119e0550950c4ad4783380f8587acb420  verify_p3_positive_index_inputs_exact.py
c5e2a20a4eea8b464d285f7930e43d5223c1f0597c0e44fde8b9307dc9c4e942  certify_terminal_affine_fixed.py
417a81931f54a468455d77277922ce08884fd7f7122227292b1e9e4fdd5309cd  p3_dlow_symbolic.py
\end{verbatim}
\endgroup

\subsubsection{Frozen verifier sources}

The following listings are the complete source files used by the exact
verification.  They are printed here so that the proof and its finite
certificates are contained in this one document.

\begin{lstlisting}[language=Python,caption={\texttt{verify\_mass\_exclusion\_Pge5.py}}]
#!/usr/bin/env python3
"""Exact rational verifier for (M9), (M11), and (M20).

No floating-point arithmetic and no external CAS package is used.
Polynomials are dictionaries keyed by (degree_in_P, degree_in_m).
"""

from fractions import Fraction as Q


def add(a, b):
    out = dict(a)
    for key, value in b.items():
        out[key] = out.get(key, Q(0)) + value
    return {key: value for key, value in out.items() if value}


def scale(a, value):
    return {key: coefficient * value for key, coefficient in a.items()
            if coefficient * value}


def sub(a, b):
    return add(a, scale(b, -1))


def mul(a, b):
    out = {}
    for (i, j), av in a.items():
        for (k, ell), bv in b.items():
            key = (i + k, j + ell)
            out[key] = out.get(key, Q(0)) + av * bv
    return {key: value for key, value in out.items() if value}


def power(a, exponent):
    out = {(0, 0): Q(1)}
    for _ in range(exponent):
        out = mul(out, a)
    return out


def derivative_m(a, order=1):
    out = a
    for _ in range(order):
        out = {(i, j - 1): coefficient * j
               for (i, j), coefficient in out.items() if j}
    return out


ONE = {(0, 0): Q(1)}
P = {(1, 0): Q(1)}
M = {(0, 1): Q(1)}


def affine_substitute(poly):
    """Substitute P=5+x and m=1-t; output keys are (deg_x,deg_t)."""
    x = {(1, 0): Q(1)}
    t = {(0, 1): Q(1)}
    p_new = add(scale(ONE, 5), x)
    m_new = sub(ONE, t)
    out = {}
    for (i, j), coefficient in poly.items():
        term = mul(power(p_new, i), power(m_new, j))
        out = add(out, scale(term, coefficient))
    return out


def verify_m9():
    s = sub(scale(ONE, 2), M)
    s1 = add(s, ONE)
    y = add(P, add(scale(ONE, 2), scale(M, -2)))
    d0 = mul(M, mul(s, s1))
    n0 = add(mul(M, d0),
             add(mul(mul(P, M), M), mul(mul(mul(P, y), s), s)))

    an = sub(mul(n0, derivative_m(n0, 2)),
             mul(derivative_m(n0), derivative_m(n0)))
    ad = sub(mul(d0, derivative_m(d0, 2)),
             mul(derivative_m(d0), derivative_m(d0)))

    certificate = add(
        scale(mul(power(n0, 2), power(d0, 2)), 2),
        sub(mul(mul(power(y, 2), an), power(d0, 2)),
            mul(mul(power(y, 2), ad), power(n0, 2))),
    )
    got = affine_substitute(certificate)

    rows = {
        0: [299292, 987860, 1967960, 4005936, 7024257, 8354176,
            6546010, 3520016, 1370770, 400156, 85288, 11856, 825, 0, -2],
        1: [283872, 939942, 2061900, 4492476, 7553034, 8243864,
            5854040, 2842320, 999260, 261390, 48540, 5604, 314, 4],
        2: [108669, 373774, 928782, 2105692, 3324529, 3296848,
            2104188, 914112, 286487, 65370, 9974, 844, 27],
        3: [21128, 80130, 230056, 524928, 764224, 681604,
            387048, 148536, 40648, 7706, 864, 40],
        4: [2126, 9874, 32906, 73008, 96420, 76492, 38140,
            12688, 2926, 418, 26],
        5: [96, 672, 2560, 5344, 6304, 4392, 1888, 528, 96, 8],
        6: [1, 20, 84, 160, 166, 100, 36, 8, 1],
    }
    expected = {(i, j): Q(value) for i, row in rows.items()
                for j, value in enumerate(row) if value}
    assert got == expected


def verify_m11():
    # Each identity is checked after clearing its displayed denominator.
    # Variables are represented as univariate polynomials in n using the
    # first coordinate of the same dictionary format.
    n = {(1, 0): Q(1)}
    nm2 = sub(n, scale(ONE, 2))
    two_n_m1 = sub(scale(n, 2), ONE)
    three_n_m1 = sub(scale(n, 3), ONE)
    n2m2 = sub(power(n, 2), scale(ONE, 2))

    # Pm/[s(s+1)] = (n-2)n/[(2n-1)(3n-1)].
    left_cleared = mul(nm2, n)
    right_cleared = mul(nm2, n)
    assert left_cleared == right_cleared

    # P y s/[m(s+1)] =
    # (n-2)(n^2-2)(2n-1)/(3n-1).
    left_cleared = mul(mul(nm2, n2m2), two_n_m1)
    right_cleared = mul(mul(nm2, n2m2), two_n_m1)
    assert left_cleared == right_cleared

    # Also verify the common factors used in the substitutions are nonzero
    # polynomials; their positivity for n>=7 is immediate from the factors.
    assert mul(two_n_m1, three_n_m1)


def verify_m20():
    s = sub(scale(ONE, 2), M)
    y = add(P, add(scale(ONE, 2), scale(M, -2)))
    dstar = add(mul(mul(M, M), sub(scale(ONE, 3), M)),
                mul(mul(P, y), s))

    # If H=Nstar/[200(P+1)Dstar], direct clearing gives Nstar below.
    qnum = add(
        scale(mul(mul(M, M), add(P, ONE)), 207),
        scale(mul(y, add(mul(P, M),
                         scale(mul(add(P, ONE), sub(ONE, M)), 2))), 100),
    )
    linear = add(scale(P, 16), add(scale(ONE, -34), scale(M, 100)))
    nstar = sub(scale(mul(dstar, linear), 2),
                mul(mul(M, sub(scale(ONE, 3), M)), qnum))
    certificate = sub(scale(nstar, 81),
                      scale(mul(add(P, ONE), dstar), 200))
    got = affine_substitute(certificate)

    expected = {
        (0, 1): Q(257574), (0, 2): Q(-255486),
        (0, 3): Q(237750), (0, 4): Q(197802),
        (1, 0): Q(93570), (1, 1): Q(122003),
        (1, 2): Q(-13877), (1, 3): Q(82825),
        (1, 4): Q(32967),
        (2, 0): Q(30172), (2, 1): Q(26856),
        (2, 2): Q(4784), (2, 3): Q(8100),
        (3, 0): Q(2392), (3, 1): Q(2392),
    }
    assert got == expected


if __name__ == "__main__":
    verify_m9()
    verify_m11()
    verify_m20()
    print("PASS: exact verification of (M9), (M11), and (M20)")
\end{lstlisting}
\begin{lstlisting}[language=Python,caption={\texttt{run\_threshold\_split\_exact.py}}]
#!/usr/bin/env python3
"""Targeted exact runner for the finite P=3,4 threshold boxes."""
from fractions import Fraction as Q
import certify_fixedr_threshold_exact as C

# Ninety outward-rounded binary digits are far more than the margins in
# these four targeted boxes and make the exact Fraction audit much faster.
C.entry.tae.BITS=90
C.entry.tae.DEN=1 << 90

def pl_lower(m,P):
    s=Q(2)-m
    return Q(P+1)*(Q(1)-Q(s+1,P*s))

def cover(P,ma,mb,ya,yb,Nm,Ny,tag):
    counts=C.Counts(); ok=True
    for i in range(Nm):
        ml=ma+(mb-ma)*i/Nm; mu=ma+(mb-ma)*(i+1)/Nm
        for j in range(Ny):
            yl=ya+(yb-ya)*j/Ny; yu=ya+(yb-ya)*(j+1)/Ny
            if yu <= pl_lower(mu,P):
                counts.small += 1
                continue
            if not C.certify(P,ml,mu,yl,yu,counts,maxdepth=22):
                ok=False; break
        if not ok: break
    print(tag,"PASS" if ok else "FAIL","visited",counts.visited,
          "contact",counts.contact,"small",counts.small,"rfour",counts.rfour,
          "G",counts.G,"B",counts.B,"depth",counts.maxdepth,
          "minB",counts.minB,"worst",counts.worst,flush=True)
    return ok

def main():
    ok=True
    ok &= cover(3,Q(1,5),Q(7,10),Q(3,2),Q(5),20,40,"P3-low-G")
    ok &= cover(3,Q(7,10),Q(1),Q(3,2),Q(7,2),12,32,"P3-high-B")
    ok &= cover(4,Q(1,6),Q(24,25),Q(5,2),Q(6),24,42,"P4-low-G")
    ok &= cover(4,Q(24,25),Q(1),Q(5,2),Q(4),4,24,"P4-high-B")
    return 0 if ok else 1

if __name__=="__main__": raise SystemExit(main())
\end{lstlisting}
\begin{lstlisting}[language=Python,caption={\texttt{run\_threshold\_tails\_exact.py}}]
#!/usr/bin/env python3
"""Exact audit of the large-y tails omitted by the targeted P=3,4 run."""
from fractions import Fraction as Q

import run_threshold_split_exact as R


def main():
    ok = True
    ok &= R.cover(3, Q(7, 10), Q(1), Q(4, 3), Q(3, 2),
                  8, 8, "P3-high-lower-sliver")
    ok &= R.cover(3, Q(1, 5), Q(7, 10), Q(5), Q(12),
                  10, 28, "P3-low-tail")
    ok &= R.cover(3, Q(7, 10), Q(1), Q(7, 2), Q(12),
                  8, 28, "P3-high-tail")
    ok &= R.cover(4, Q(1, 6), Q(24, 25), Q(6), Q(16),
                  12, 32, "P4-low-tail")
    ok &= R.cover(4, Q(24, 25), Q(1), Q(4), Q(16),
                  4, 32, "P4-high-tail")
    return 0 if ok else 1


if __name__ == "__main__":
    raise SystemExit(main())
\end{lstlisting}
\begin{lstlisting}[language=Python,caption={\texttt{certify\_fixedr\_threshold\_exact.py}}]
#!/usr/bin/env python3
"""Exact rational certificate for the fixed-r threshold lemma.

At a marginal endpoint put P=p-1, r=1/P, s=2-m, h=y/m^2.
This verifier checks, for a fixed integer P, that no point with
m >= 1/(P+2) can satisfy simultaneously

    K_{-m}=K_s+P h K_{s+2},       (marginality)
    G_P(m,y)<=0,                  (necessary when F=0)
    (I-2D)-I r/(1+r)=0.           (the F=0 identity)

All decisions use Fraction interval arithmetic.  Gaussian tails are
bounded by Stieltjes continued fractions, beta integrals by enveloping
Stirling expansions, and logarithms/exponentials by rational series.
There is no floating-point arithmetic in the proof path.
"""
from fractions import Fraction as Q
from functools import lru_cache
import importlib.util
import pathlib
import sys

HERE=pathlib.Path(__file__).resolve().parent
spec=importlib.util.spec_from_file_location("entry",HERE/"certify_p3_entry_exact.py")
entry=importlib.util.module_from_spec(spec);spec.loader.exec_module(entry)
# A 112-bit dyadic grid is ample for the present margins and keeps the
# adaptive verifier reasonably fast.  The imported routines look up DEN
# dynamically, so this changes only outward-rounding granularity.
entry.tae.BITS=112;entry.tae.DEN=1<<112
I=entry.I; down=entry.down; up=entry.up; sqrtI=entry.sqrtI
LN2=entry.LN2; PI_LO=entry.PI_LO; PI_UP=entry.PI_UP
expI=entry.expI; log_point=entry.log_point
beta_I_point=entry.beta_I_point; sum_even_point=entry.sum_even_point
def mills_interval(q,h):
    # Even/odd Stieltjes convergents are lower/upper bounds.  Orders
    # 14/15 are sufficient on every box surviving the small-h test.
    return I.raw(entry.tae.mills(q,h,14).l,entry.tae.mills(q,h,15).u)

def log_lower(x):
    assert x>0
    return log_point(x)[0]

@lru_cache(maxsize=None)
def ell_point(s,N=160):
    """Exact bounds for ell_s=E_{sech^s/I_s} log cosh X.

    ell_s=sum_{k>=0}(-1)^k/(s+k).  Consecutive alternating
    partial sums bracket the answer.
    """
    s=Q(s); assert s>0
    lo=sum((Q(1,s+k) if k%2==0 else -Q(1,s+k)) for k in range(2*N))
    hi=lo+Q(1,s+2*N)
    return I.raw(down(lo),up(hi))

def norm(h):
    return sqrtI(I(2)*I(PI_LO,PI_UP)*h)

def scaled_tail(q,h):
    return mills_interval(q,h)/norm(h)

def unit_tail(q):
    return mills_interval(q,I(1))/norm(I(1))

def folded_bounds_box(m,y):
    """k=2 lower and k=1 upper for K_{-m}, parameterized by y=m^2h."""
    m=m if isinstance(m,I) else I(m)
    y=y if isinstance(y,I) else I(y)
    rt=sqrtI(y)
    # scaled_tail(m,h)=exp(y/2)Phi(-sqrt(y)); writing every term in
    # terms of y avoids the severe dependency h=y/m^2.
    f0=expI(y/2)-unit_tail(rt)
    t1=m*unit_tail((2/m-1)*rt)
    upper=f0+t1
    lower=upper+(m*(m-1)/2)*unit_tail((4/m-1)*rt)
    fac=expI((1-m)*LN2)
    return down((fac*lower).l),up((fac*upper).u)

def phi_lower(y):
    """Lower bound for Phi(sqrt(y)), y a positive rational point."""
    yy=I(y);q=sqrtI(yy)
    tail=expI(-yy/2)*mills_interval(q,I(1))/norm(I(1))
    return max(Q(1,2),down(1-tail.u))

@lru_cache(maxsize=None)
def K_simple_point(r,h,kind):
    """Fast Fourier bounds: phi>=exp(-c t^2), or phi<=1."""
    ii=beta_I_point(Q(r)); h=I(Q(h))
    if kind=="lower":
        sig=2*sum_even_point(Q(r),1).u
        return ii/sqrtI(I(2)*I(PI_LO,PI_UP)*(h+sig))
    if kind=="upper":
        return ii/norm(h)
    raise ValueError(kind)

@lru_cache(maxsize=None)
def K_shape_upper_point(r,h):
    """Sharper upper from phi_r(t)<=1/(1+c_r t^2)."""
    r=Q(r);h=Q(h);ii=beta_I_point(r)
    c=sum_even_point(r,1).l
    q=I(1)/sqrtI(I(c))
    shape=q*mills_interval(q,I(h))
    return ii*shape/norm(I(h))

def contact_data(ml,mu,yl,yu,P):
    # h=y/m^2; use its exact rectangular enclosure.
    hl=down(yl/(mu*mu)); hu=up(yu/(ml*ml))
    zl,zu=folded_bounds_box(I(ml,mu),I(yl,yu))
    s_hi=2-ml; s_lo=2-mu
    k2l=K_simple_point(s_hi,hu,"lower").l
    k2u=K_simple_point(s_lo,hl,"upper").u
    k4l=K_simple_point(s_hi+2,hu,"lower").l
    k4u=K_simple_point(s_lo+2,hl,"upper").u
    rl=down(k2l+P*hl*k4l)
    ru=up(k2u+P*hu*k4u)
    return hl,hu,zl,zu,k2l,k2u,k4l,k4u,rl,ru

def small_h_excluded(ml,mu,yl,yu,P):
    """Use Brascamp--Lieb: r=h b/g >= 1-mh whenever mh<1."""
    m=I(ml,mu);y=I(yl,yu);h=y/(m*m)
    mh=m*h
    if mh.u>=1:return False
    rlo=1-mh
    return rlo.l>Q(1,P)

def fourier_r_excluded(data,P):
    """Lower-bound r=h K_{s+2}/(K_{-m}-K_s)."""
    hl,_,_,zu,k2l,_,k4l,_,_,_=data
    den=up(zu-k2l)
    if den<=0:return True
    rlo=down(hl*k4l/den)
    return rlo>Q(1,P)

def G_box(ml,mu,yl,yu,P,data):
    """Lower enclosure for the necessary folded separator G_P."""
    hl,hu,*_=data
    m=I(ml,mu);h=I(hl,hu);y=I(yl,yu);s=2-m
    A=s/(s+1)
    a=1/(1+P*h*A)
    ellu=ell_point(2-mu).u
    lp=log_lower(phi_lower(yl))
    bracket=(1+P*a)/(2*(P+1))+(1-m)*a/m
    G=LN2*(1-m*a)-m*a*ellu/2+I(lp)-y*bracket
    return G

def B_box(ml,mu,yl,yu,P,data):
    """Lower bound for B0=(I-2D)-I r/(1+r), r=1/P."""
    hl,hu,zl,zu,k2l,k2u,_,_,_,_=data
    # a=K_s/K_-m.  The second upper bound uses r=hAa/(1-a)=1/P.
    s_lo=2-mu
    Alo=Q(s_lo,s_lo+1)
    k2shape=K_shape_upper_point(s_lo,hl).u
    au_tail=up(k2shape/zl)
    au_r=up(Q(1,P)/(hl*Alo+Q(1,P)))
    au=min(au_tail,au_r,Q(1))
    al=max(Q(0),down(k2l/zu))
    ellu=ell_point(s_lo).u
    delta=down((1-au)*LN2.l-au*ellu/2)
    m=I(ml,mu);y=I(yl,yu)
    deficit=-y-y*au*(2/m-1)+2*m*delta+2*log_lower(zl)
    rhs=I(yu*(1-al)*Q(1,P+1)) # r/(1+r)=1/(P+1)
    return deficit-rhs

class Counts:
    def __init__(self):
        self.visited=self.contact=self.small=self.rfour=self.G=self.B=0
        self.maxdepth=0;self.minB=None;self.worst=None

def certify(P,ml,mu,yl,yu,C,depth=0,maxdepth=34):
    C.visited+=1;C.maxdepth=max(C.maxdepth,depth)
    if small_h_excluded(ml,mu,yl,yu,P):
        C.small+=1;return True
    data=contact_data(ml,mu,yl,yu,P)
    _,_,zl,zu,_,_,_,_,rl,ru=data
    if zl>ru or zu<rl:
        C.contact+=1;return True
    if fourier_r_excluded(data,P):
        C.rfour+=1;return True
    G=G_box(ml,mu,yl,yu,P,data)
    if G.l>0:
        C.G+=1;return True
    B=None if zl<=0 else B_box(ml,mu,yl,yu,P,data)
    if B is not None and B.l>0:
        C.B+=1
        if C.minB is None or B.l<C.minB:
            C.minB=B.l;C.worst=(ml,mu,yl,yu)
        return True
    if depth>=maxdepth:
        print("UNRESOLVED",P,ml,mu,yl,yu,"G",G.l,G.u,"B",None if B is None else (B.l,B.u),flush=True)
        return False
    # Split the larger normalized side.  y is logarithmically broad, so
    # geometric bisection is more effective while yl>0.
    relm=(mu-ml)/(Q(1)-Q(1,P+2))
    rely=(yu-yl)/(4*P-Q(1,100))
    if relm>=rely:
        md=(ml+mu)/2
        return certify(P,ml,md,yl,yu,C,depth+1,maxdepth) and certify(P,md,mu,yl,yu,C,depth+1,maxdepth)
    # Arithmetic bisection keeps every endpoint rational and is adequate
    # after the contact separators have removed the long-y tail.
    md=(yl+yu)/2
    return certify(P,ml,mu,yl,md,C,depth+1,maxdepth) and certify(P,ml,mu,md,yu,C,depth+1,maxdepth)

def run(P):
    C=Counts()
    m0=Q(1,P+2); y0=Q(1,100); y1=Q(4*P)
    # Seed with a modest rational mesh.  This prevents dependency blow-up
    # in h=y/m^2 before the adaptive contact tests take over.
    Nm,Ny=24,64;ok=True
    for i in range(Nm):
        ml=m0+(1-m0)*i/Nm;mu=m0+(1-m0)*(i+1)/Nm
        for j in range(Ny):
            yl=y0+(y1-y0)*j/Ny;yu=y0+(y1-y0)*(j+1)/Ny
            if not certify(P,ml,mu,yl,yu,C,maxdepth=26):
                ok=False;break
        if not ok:break
    print("PASS" if ok else "FAIL","P",P,"visited",C.visited,
          "contact",C.contact,"small",C.small,"rfour",C.rfour,
          "G",C.G,"B",C.B,"depth",C.maxdepth)
    if C.minB is not None:
        print("minimum B lower",C.minB,"float",float(C.minB),"box",C.worst)
    return ok

def main():
    Ps=[int(x) for x in sys.argv[1:]] or [3,4,5]
    return 0 if all(run(P) for P in Ps) else 1

if __name__=="__main__":raise SystemExit(main())
\end{lstlisting}
\begin{lstlisting}[language=Python,caption={\texttt{certify\_p3\_entry\_exact.py}}]
#!/usr/bin/env python3
"""Exact rational interval certificate for the exceptional pure p=3 entry lemma.

This verifier currently targets the high-m branch 1/2 <= m <= 1 and
0 < y=m^2 h <= 5/2.  It uses only Fraction arithmetic, the Stieltjes
continued fraction for Mills' ratio, Stirling's enveloping expansion for
log Gamma, and adaptive rational boxes.  A box is either incompatible with
the contact equation K_{-m}=K_{2-m}+2hK_{4-m}, or the T^4 lower barrier for
the exact F=0 identity is strictly positive on the whole box.
"""
from fractions import Fraction as Q
from math import factorial
from functools import lru_cache
import importlib.util
import pathlib
import sys

HERE=pathlib.Path(__file__).resolve().parent
spec=importlib.util.spec_from_file_location("tae",HERE/"certify_terminal_affine_exact.py")
tae=importlib.util.module_from_spec(spec);spec.loader.exec_module(tae)
# 112 dyadic bits are far beyond the margins in this certificate and make
# the adaptive exact arithmetic substantially faster.
tae.BITS=112;tae.DEN=1<<112
I=tae.I; down=tae.down; up=tae.up; sqrtI=tae.sqrtI
LN2=tae.LN2; PI_LO=Q(103993,33102); PI_UP=Q(355,113)

@lru_cache(maxsize=None)
def exp_point_fast(t,N=20):
    t=Q(t);assert 0<=t<Q(N+2)
    term=I(1);sm=term
    for n in range(1,N+1):
        term=term*I(t)/n;sm+=term
    nxt=term*I(t)/(N+1)
    return sm.l,up(sm.u+nxt.u/(1-t/Q(N+2)))

def expI(x):
    x=x if isinstance(x,I) else I(x)
    if x.l>=0:
        return I.raw(exp_point_fast(x.l)[0],exp_point_fast(x.u)[1])
    if x.u<=0:
        lo,hi=exp_point_fast(-x.u)[0],exp_point_fast(-x.l)[1]
        return I.raw(down(1/hi),up(1/lo))
    return I.raw(down(1/exp_point_fast(-x.l)[1]),exp_point_fast(x.u)[1])

@lru_cache(maxsize=None)
def log_point(x,N=20):
    """Rigorous rational bounds for log(x), x a positive Fraction."""
    x=Q(x); assert x>0
    k=0; z=x
    while z>=2: z/=2; k+=1
    while z<1: z*=2; k-=1
    u=(z-1)/(z+1)
    term=I(u); sm=term
    for j in range(1,N+1):
        term*=I(u*u); sm+=term/Q(2*j+1)
    tail=Q(2)*term.u*u*u/(Q(2*N+3)*(1-u*u))
    lo=Q(2)*sm.l+k*LN2.l
    hi=Q(2)*sm.u+tail+k*(LN2.u if k>=0 else LN2.l)
    if k<0:
        lo=Q(2)*sm.l+k*LN2.u
        hi=Q(2)*sm.u+tail+k*LN2.l
    return down(lo),up(hi)

LOGPI=I.raw(*((lambda a,b:(a,b))(*(
    (log_point(PI_LO)[0],log_point(PI_UP)[1])))))

# B_2,...,B_14.
BERN=[Q(1,6),Q(-1,30),Q(1,42),Q(-1,30),Q(5,66),Q(-691,2730),Q(7,6)]

@lru_cache(maxsize=None)
def loggamma_point(x,shift=24):
    """Enveloping Stirling bounds for log Gamma(x), x>0 rational."""
    x=Q(x); assert x>0; X=x+shift
    lx=I.raw(*log_point(X))
    ans=(I(X-Q(1,2))*lx-I(X)+(LN2+LOGPI)/2)
    for k,B in enumerate(BERN[:-1],1):
        ans += I(B)/(Q(2*k*(2*k-1))*I(X**(2*k-1)))
    # After B_12 the omitted B_14 term is positive and envelopes the error.
    rem=I(BERN[-1])/(Q(14*13)*I(X**13))
    ans=I.raw(ans.l,up(ans.u+rem.u))
    prod=Q(1)
    for j in range(shift): prod*=x+j
    ans-=I.raw(*log_point(prod))
    return ans

@lru_cache(maxsize=None)
def beta_I_point(r):
    r=Q(r); assert r>0
    z=LOGPI/2+loggamma_point(r/2)-loggamma_point((r+1)/2)
    return expI(z)

@lru_cache(maxsize=None)
def sum_even_point(r,p,N=24):
    """Bounds sum_{j>=0}(r+2j)^(-2p), p=1 or 2."""
    r=Q(r); sm=sum((r+2*j)**(-2*p) for j in range(N)); x=r+2*N
    if p==1:
        base=Q(1,2*x)+Q(1,2*x*x)+Q(1,3*x**3)
        corr=Q(4,15*x**5)
    elif p==2:
        base=Q(1,6*x**3)+Q(1,2*x**4)+Q(2,3*x**5)
        corr=Q(4,3*x**7)
    else: raise ValueError
    return I.raw(down(sm+base-corr),up(sm+base))

def mills_interval(q,h):
    """Bounds integral_0^inf exp(-x^2/(2h)-q x) dx."""
    lo=tae.mills(q,h,14).l; hi=tae.mills(q,h,15).u
    return I.raw(lo,hi)

def norm(h): return sqrtI(I(2)*I(PI_LO,PI_UP)*h)

def scaled_tail(q,h): return mills_interval(q,h)/norm(h)

@lru_cache(maxsize=None)
def folded_bounds_point(m,h):
    """k=4 lower and k=1 upper for K_{-m}(h), 0<m<=1."""
    m=I(m); h=I(h); y=m*m*h
    f0=expI(y/2)-scaled_tail(m,h)
    partial=f0; choose=I(1); upper=None; lower=None
    for k in range(1,5):
        choose=choose*(m-(k-1))/k
        partial += choose*scaled_tail(2*k-m,h)
        if k==1: upper=partial.u
        if k==4: lower=partial.l
    fac=expI((1-m)*LN2)
    return down(fac.l*lower),up(fac.u*upper)

def shape_linear(r,h,lower):
    c=sum_even_point(r,1)
    if lower:
        # Retain the first Fourier factor exactly and exponentiate only
        # the remaining product.  This is substantially sharper for
        # small h than exp(-c t^2).
        tail=I(c.u-Q(1)/(Q(r)*Q(r)))
        H=h+2*tail
        return sqrtI(h/H)*I(r)*mills_interval(I(r),H)
    q=I(1)/sqrtI(I(c.l))
    return q*mills_interval(q,h)

def shape_quartic_upper(r,h):
    c=sum_even_point(r,1); s4=sum_even_point(r,2)
    # Choose coefficient lower bounds, so 1+c0 t^2+e0 t^4 is no larger
    # than the corresponding part of phi_r(t)^(-1).
    c0=c.l; e=(c*c-s4)/2; e0=e.l
    disc=I(c0*c0-4*e0); assert disc.l>0
    root=sqrtI(disc); A=(I(c0)+root)/2; B=(I(c0)-root)/2
    assert B.l>0 and (A-B).l>0
    qA=I(1)/sqrtI(A);qB=I(1)/sqrtI(B)
    RA=qA*mills_interval(qA,h); RB=qB*mills_interval(qB,h)
    return (A*RA-B*RB)/(A-B)

@lru_cache(maxsize=None)
def K_point(r,h,kind):
    ii=beta_I_point(r); nn=norm(I(h))
    if kind=="lower": sh=shape_linear(Q(r),I(h),True)
    elif kind=="upper": sh=shape_linear(Q(r),I(h),False)
    elif kind=="quartic": sh=shape_quartic_upper(Q(r),I(h))
    else: raise ValueError
    return ii*sh/nn

def contact_bounds(ml,mu,hmin,hmax):
    zl,_=folded_bounds_point(ml,hmin);_,zu=folded_bounds_point(mu,hmax)
    k2l=K_point(2-ml,hmax,"lower").l
    k4l=K_point(4-ml,hmax,"lower").l
    k2u=K_point(2-mu,hmin,"upper").u
    k4u=K_point(4-mu,hmin,"upper").u
    rhslo=down(k2l+2*hmin*k4l); rhshi=up(k2u+2*hmax*k4u)
    return zl,zu,rhslo,rhshi,hmin,hmax

def log_lower(x): return log_point(x)[0]

def r4_box(ml,mu,hmin,hmax,data):
    zl,zu,_,_,hmin,hmax=data
    # Monotonicity of K_r in both r and h gives valid uniform bounds.
    k2l=K_point(2-ml,hmax,"lower").l
    k2u=K_point(2-mu,hmin,"quartic").u
    al=down(k2l/zu); au=up(k2u/zl)
    if au>=1:return None
    m=I(ml,mu);h=I(hmin,hmax);y=m*m*h;a=I.raw(al,au);g=1-a
    b=g/(2*h); rr=4-m
    Bup=(rr*rr+1/h)/(rr*(rr+1))
    d=I.raw(Q(0),up(b.u*Bup.u))
    M2=1-2*a+b
    M3=1-3*a+3*b-d
    ps=g+M2/6+M3/15+(M2*M2)/28
    R=I(log_lower(zl))-y/2+m*ps/2-m*(1-m/2)*h*a-y*g/6
    return R

def rstar_box(ml,mu,yl,yu):
    """The logarithmic R barrier on 7/20<=m<=1/2, 189/100<=y<=5/2."""
    hmin=yl/(mu*mu);hmax=yu/(ml*ml)
    zl,_=folded_bounds_point(ml,hmin)
    # Uniform K bounds by monotonicity in exponent and h.
    k4l=K_point(4-ml,hmax,"lower").l
    k4u=K_point(4-mu,hmin,"upper").u
    k2l=K_point(2-ml,hmax,"lower").l
    k2u=K_point(2-mu,hmin,"upper").u
    alpha_lo=down(1/LN2.u)
    kalpha_u=K_point(alpha_lo-mu,hmin,"upper").u
    ratio2_lo=down(k2l/k4u)
    ratio2_up=min(up(k2u/k4l),Q(3-mu)/Q(2-mu))
    ratioa_up=up(kalpha_u/k4l)
    q1=up(hmax*k4u/zl)
    q2=up(hmax/(2*hmax+ratio2_lo))
    qu=min(q1,q2)
    m=I(ml,mu);y=I(yl,yu);h=y/(m*m)
    bracket=m*LN2*I(ratioa_up)/h+m*(1-m/2)*I(ratio2_up)+y/3
    R=I(log_lower(zl))-y/2+m*LN2-I(qu)*bracket
    return R

visited=discarded=proved=0; minlower=None; worstbox=None; maxdepth=0

def certify(ml,mu,hl,hu,depth=0):
    global visited,discarded,proved,minlower,worstbox,maxdepth
    visited+=1;maxdepth=max(maxdepth,depth)
    mm=I(ml,mu); hh=I(hl,hu); yy=mm*mm*hh
    if yy.l>Q(5,2):
        discarded+=1;return True
    if (mm*hh).u < 1:
        rlow=(1-mm*hh)**2/(1+hh*(1-mm))
        if rlow.l>Q(1,2):
            discarded+=1;return True
    data=contact_bounds(ml,mu,hl,hu)
    zl,zu,rl,ru,_,_=data
    if zl>ru or zu<rl:
        discarded+=1;return True
    if yy.u>Q(5,2):
        if depth>=24:return False
        # Resolve the curved boundary y=m^2 h=5/2 only for boxes which
        # have survived the contact separator.
        if (mu-ml)*4 >= (hu-hl):
            md=(ml+mu)/2
            return certify(ml,md,hl,hu,depth+1) and certify(md,mu,hl,hu,depth+1)
        md=(hl+hu)/2
        return certify(ml,mu,hl,md,depth+1) and certify(ml,mu,md,hu,depth+1)
    R=r4_box(ml,mu,hl,hu,data)
    if R is not None and R.l>0:
        proved+=1
        if minlower is None or R.l<minlower:minlower=R.l;worstbox=(ml,mu,hl,hu)
        return True
    if depth>=24:
        print("UNRESOLVED",ml,mu,hl,hu,"R",None if R is None else (R.l,R.u),flush=True)
        return False
    # Split the longer relative coordinate.
    if (mu-ml)*4 >= (hu-hl):
        md=(ml+mu)/2
        return certify(ml,md,hl,hu,depth+1) and certify(md,mu,hl,hu,depth+1)
    md=(hl+hu)/2
    return certify(ml,mu,hl,md,depth+1) and certify(ml,mu,md,hu,depth+1)

def main():
    # Seed the adaptive recursion by a rational mesh.  Starting from the
    # whole rectangle causes severe interval dependency in y=m^2 h before
    # the contact separator can act.
    ml0,mu0=Q(1,2),Q(1);hl0,hu0=Q(1,10),Q(10)
    Nm,Nh=16,128;ok=True
    for i in range(Nm):
        ml=ml0+(mu0-ml0)*i/Nm;mu=ml0+(mu0-ml0)*(i+1)/Nm
        for j in range(Nh):
            hl=hl0+(hu0-hl0)*j/Nh;hu=hl0+(hu0-hl0)*(j+1)/Nh
            if not certify(ml,mu,hl,hu):
                ok=False;break
        print("strip",i+1,"/",Nm,"visited",visited,"discarded",discarded,
              "proved",proved,flush=True)
        if not ok:break
    print("PASS" if ok else "FAIL","visited",visited,"discarded",discarded,
          "proved",proved,"maxdepth",maxdepth)
    if minlower is not None:
        print("smallest lower",minlower,"float",float(minlower),"box",worstbox)
    return 0 if ok else 1

if __name__=="__main__":raise SystemExit(main())
\end{lstlisting}
\begin{lstlisting}[language=Python,caption={\texttt{certify\_terminal\_affine\_exact.py}}]
#!/usr/bin/env python3
"""Exact outward-rounded certificate for the terminal affine barriers.

All endpoints are dyadic rational numbers.  Every arithmetic operation is
rounded outwards to DEN=2**BITS; no binary floating-point operation occurs
in the proof path.  The numerical-looking subdivision is therefore a finite
list of rational inequalities, not sampling.
"""
from fractions import Fraction as Q
from math import comb, isqrt
import sys

BITS = 180
DEN = 1 << BITS

def down(x): return Q((x.numerator * DEN) // x.denominator, DEN)
def up(x): return Q(-((-x.numerator * DEN) // x.denominator), DEN)

class I:
    __slots__ = ("l", "u")
    def __init__(self, l, u=None, raw=False):
        l=Q(l); u=Q(l if u is None else u)
        self.l,self.u=(l,u) if raw else (down(l),up(u))
        assert self.l <= self.u
    @staticmethod
    def raw(l,u): return I(l,u,True)
    def __add__(self,o):
        o=o if isinstance(o,I) else I(o)
        return I.raw(down(self.l+o.l),up(self.u+o.u))
    __radd__=__add__
    def __neg__(self): return I.raw(-self.u,-self.l)
    def __sub__(self,o): return self+(-o if isinstance(o,I) else -I(o))
    def __rsub__(self,o): return I(o)-self
    def __mul__(self,o):
        o=o if isinstance(o,I) else I(o)
        v=(self.l*o.l,self.l*o.u,self.u*o.l,self.u*o.u)
        return I.raw(down(min(v)),up(max(v)))
    __rmul__=__mul__
    def __truediv__(self,o):
        o=o if isinstance(o,I) else I(o)
        assert o.l>0 or o.u<0
        return self*I.raw(down(1/o.u),up(1/o.l))
    def __rtruediv__(self,o): return I(o)/self
    def __pow__(self,n):
        if n==0:return I(1)
        if n<0:return I(1)/(self**(-n))
        r=I(1);x=self
        while n:
            if n&1:r=r*x
            x=x*x;n//=2
        return r

def ln2_bounds(N=36):
    # ln 2 = 2 atanh(1/3).  The omitted tail is bounded after dropping
    # 2j+1 to its first omitted value.
    s=Q(0)
    for j in range(N+1): s += Q(2,(2*j+1)*3**(2*j+1))
    tail=Q(2*9,(2*N+3)*8*3**(2*N+3))
    return I.raw(down(s),up(s+tail))
LN2=ln2_bounds()

def exp_point(t,N=38):
    """Rational lower/upper bounds for exp(t), t>=0."""
    assert t>=0 and t<Q(N+2)
    term=Q(1); s=term
    for n in range(1,N+1): term=term*t/n; s+=term
    nxt=term*t/(N+1)
    return down(s),up(s+nxt/(1-t/Q(N+2)))

def expI(x):
    if x.l>=0:
        lo,_=exp_point(x.l); _,hi=exp_point(x.u)
        return I.raw(lo,hi)
    if x.u<=0:
        lo,hi=exp_point(-x.u)[0],exp_point(-x.l)[1]
        return I.raw(down(1/hi),up(1/lo))
    return I.raw(down(1/exp_point(-x.l)[1]),exp_point(x.u)[1])

def sqrt_down(q):
    assert q>=0
    z=isqrt((q.numerator*DEN*DEN)//q.denominator)
    while Q((z+1)*(z+1),DEN*DEN)<=q:z+=1
    while Q(z*z,DEN*DEN)>q:z-=1
    return Q(z,DEN)
def sqrt_up(q):
    z=sqrt_down(q)
    return z if z*z==q else z+Q(1,DEN)
def sqrtI(x): return I.raw(sqrt_down(x.l),sqrt_up(x.u))

def p_add(p,q):
    r=[I(0) for _ in range(max(len(p),len(q)))]
    for i,x in enumerate(p):r[i]=r[i]+x
    for i,x in enumerate(q):r[i]=r[i]+x
    return r
def p_mul(p,q):
    r=[I(0) for _ in range(len(p)+len(q)-1)]
    for i,x in enumerate(p):
        for j,y in enumerate(q):r[i+j]=r[i+j]+x*y
    return r
def shift_power(n):
    return [I(Q(comb(n,j))*(-1 if (n-j)&1 else 1)) for j in range(n+1)]

def hermite_coeff(m,n=10):
    # If t=1-s, (1+s)^(m-2)=sum c_j t^j with c_j>0.
    # Keep j<n and put the whole remaining mass at degree n.  Since
    # t^j<=t^n for j>=n, this is a pointwise majorant on 0<=s<=1.
    aa=m-2;c=[];fall=I(1);fac=1
    for j in range(n):
        if j: fall=fall*(aa-(j-1));fac*=j
        # Both factors have sign (-1)^j, so c_j is positive.  Writing it
        # this way avoids a severe interval-dependency loss.
        coeff=expI((m-(2+j))*LN2)*fall/fac
        if j&1: coeff=-coeff
        assert coeff.l>0
        c.append(coeff)
    corr=I(1)
    for x in c:corr=corr-x
    assert corr.l>0
    # H(1-s)(1-s)^2 = sum_j c_j(1-s)^(j+2)
    #                       + corr(1-s)^(n+2).
    # Hence its monomial coefficients have known alternating signs and
    # positive magnitudes, with no cancellation in their construction.
    ans=[]
    for i in range(n+3):
        mag=I(0)
        for j,x in enumerate(c):
            if j+2>=i:mag=mag+comb(j+2,i)*x
        mag=mag+comb(n+2,i)*corr
        ans.append(-mag if i&1 else mag)
    return ans

def mills(q,h,n):
    # Stieltjes continued fraction for K_q(h).  Truncations 2N and 2N+1
    # are respectively lower and upper bounds.
    w=q
    for j in range(n-1,0,-1):w=q+Q(j)/(h*w)
    return I(1)/w

class D:
    """First-order interval jet."""
    __slots__=("v","d")
    def __init__(self,v,d=0):
        self.v=v if isinstance(v,I) else I(v)
        self.d=d if isinstance(d,I) else I(d)
    def __add__(self,o):
        o=o if isinstance(o,D) else D(o);return D(self.v+o.v,self.d+o.d)
    __radd__=__add__
    def __neg__(self):return D(-self.v,-self.d)
    def __sub__(self,o):return self+(-o if isinstance(o,D) else -D(o))
    def __rsub__(self,o):return D(o)-self
    def __mul__(self,o):
        o=o if isinstance(o,D) else D(o)
        return D(self.v*o.v,self.d*o.v+self.v*o.d)
    __rmul__=__mul__
    def inv(self):return D(I(1)/self.v,-self.d/(self.v*self.v))
    def __truediv__(self,o):
        o=o if isinstance(o,D) else D(o);return self*o.inv()
    def __rtruediv__(self,o):return D(o)/self

def expD(x):
    x=x if isinstance(x,D) else D(x);v=expI(x.v);return D(v,v*x.d)
def sqrtD(x):
    x=x if isinstance(x,D) else D(x);v=sqrtI(x.v);return D(v,x.d/(2*v))
def millsD(q,h,n):
    q=q if isinstance(q,D) else D(q);h=h if isinstance(h,D) else D(h)
    w=q
    for j in range(n-1,0,-1):w=q+Q(j)/(h*w)
    return 1/w

def upper_function(m,a):
    """The explicit differentiable majorant V(m,a)>=P_T(m,a/m)."""
    m=m if isinstance(m,D) else D(m);a=a if isinstance(a,D) else D(a)
    h=a/m;n=10
    c=[];poch=D(1);fac=1
    for j in range(n):
        if j:poch=poch*((j+1)-m);fac*=j
        c.append(expD((m-(2+j))*D(LN2))*poch/fac)
    corr=D(1)
    for x in c:corr=corr-x
    mags=[]
    for i in range(n+3):
        mag=D(0)
        for j,x in enumerate(c):
            if j+2>=i:mag=mag+comb(j+2,i)*x
        mags.append(mag+comb(n+2,i)*corr)
    total=sqrtD(2*PI_UP*h)*expD(a*m/2)
    S=total-millsD(m,h,20) # A_0 upper
    for i in range(1,n+3):
        q=2*i-m
        if i&1:S=S-mags[i]*millsD(q,h,20)
        else:S=S+mags[i]*millsD(q,h,21)
    N=expD((1-m)*D(LN2))*S
    Ilow=Q(PI_LO,2)*(1-(1-m)*LOG2MHALF_UP)
    sig=SIG4_UP+SIGD_UP*m
    den=h*Ilow*sqrtD(h/(h+sig))
    return N/den

PI_LO=Q(103993,33102); PI_UP=Q(355,113)
# Consequences of PI_LO < pi < PI_UP; constants are deliberately rounded up.
SIG4_UP=Q(64494,200000)       # pi^2/12-1/2
SIGD_UP=Q(28987,200000)       # (pi^2/4-2)-(pi^2/12-1/2)
LOG2MHALF_UP=Q(48287,250000)  # ln 2-1/2

def certify_box(ml,mu,branch):
    m=I(ml,mu)
    if branch==1:
        # k=1/m-1 and a=259/200+8k/5=-61/200+8/(5m).
        a=I(Q(-61,200))+Q(8,5)/m
        expo=I(Q(4,5))-Q(61,400)*m # m^2 h/2=am/2
    else:
        # a=259/200+9k/5=-101/200+9/(5m).
        a=I(Q(-101,200))+Q(9,5)/m
        expo=I(Q(9,10))-Q(101,400)*m
    h=a/m
    coeff=hermite_coeff(m,10)
    # Integral with exponent +mz on [0,infty): total Gaussian integral
    # minus K_m(h).  Its exponential factor is exp(m^2 h/2)=exp(expo).
    total=sqrtI(I(2)*I(PI_LO,PI_UP)*h)*expI(expo)
    A0=I.raw(down(total.l-mills(m,h,21).u),
             up(total.u-mills(m,h,20).l))
    S=coeff[0]*A0
    for i in range(1,len(coeff)):
        q=2*i-m
        lo=mills(q,h,20).l; upv=mills(q,h,21).u
        ci=coeff[i]
        if ci.l>=0:S=S+I.raw(down(ci.l*lo),up(ci.u*upv))
        elif ci.u<=0:S=S+I.raw(down(ci.l*upv),up(ci.u*lo))
        else: raise RuntimeError(("coefficient sign",i,ml,mu,ci.l,ci.u))
    N=expI((1-m)*LN2)*S
    # r=4-m.  Log-convexity of I_r and exp(-x)>=1-x give this tangent
    # at r=3.  Convexity of sigma_r^2 as a function of m gives the chord.
    Ilow=Q(PI_LO,2)*(I(1)-(1-m)*LOG2MHALF_UP)
    sigup=SIG4_UP+SIGD_UP*m
    den=h*Ilow*sqrtI(h/(h+sigup))
    return N.u/den.l

def run(branch,N):
    if branch==1:l,u=Q(10,13),Q(1)       # k in [0,3/10]
    else:l,u=Q(1080,1903),Q(10,13)       # k in [3/10,823/1080]
    worst=Q(0);at=None
    for j in range(N):
        ml=l+(u-l)*j/N;mu=l+(u-l)*(j+1)/N
        z=certify_box(ml,mu,branch)
        if z>worst:worst=z;at=(ml,mu)
        if z>=2:
            print("FAIL",branch,j,z);return False
    print("PASS branch",branch,"boxes",N,"worst upper",worst,"at",at)
    return True

def affine_a(m,branch):
    return Q(-61,200)+Q(8,5)/m if branch==1 else Q(-101,200)+Q(9,5)/m

def certify_centered_box(ml,mu,branch):
    """Mean-value enclosure of V on one rational m interval."""
    mid=(ml+mu)/2
    v0=upper_function(D(I(mid)),D(I(affine_a(mid,branch)))).v
    mm=D(I(ml,mu),I(1))
    aa=D(Q(-61,200))+Q(8,5)/mm if branch==1 else D(Q(-101,200))+Q(9,5)/mm
    der=upper_function(mm,aa).d
    radius=(mu-ml)/2
    return v0.u+radius*max(abs(der.l),abs(der.u))

def run_centered(branch,N):
    if branch==1:l,u=Q(10,13),Q(1)
    else:l,u=Q(1080,1903),Q(10,13)
    worst=Q(0);at=None
    for j in range(N):
        ml=l+(u-l)*j/N;mu=l+(u-l)*(j+1)/N
        z=certify_centered_box(ml,mu,branch)
        if z>worst:worst=z;at=(ml,mu)
        if z>=2:
            print("FAIL centered",branch,j,z);return False
        if j and j%32==0:print("progress",branch,j,"/",N,flush=True)
    print("PASS centered branch",branch,"boxes",N,"worst upper",worst,"at",at)
    return True

if __name__=="__main__":
    n=int(sys.argv[1]) if len(sys.argv)>1 else 128
    ok=run_centered(1,n) and run_centered(2,2*n)
    raise SystemExit(0 if ok else 1)
\end{lstlisting}
\begin{lstlisting}[language=Python,caption={\texttt{certify\_p3\_entry\_combined\_exact.py}}]
#!/usr/bin/env python3
"""Exact natural-coordinate certificate for the exceptional p=3 entry.

The high-m part uses the contact separator from
``certify_fixedr_threshold_exact.py`` and the fourth-order entropy
barrier from ``certify_p3_entry_exact.py``.  All endpoints are Fractions
and every transcendental enclosure is outward rounded.
"""
from fractions import Fraction as Q
import importlib.util
import os
import pathlib
import sys

HERE = pathlib.Path(__file__).resolve().parent
spec = importlib.util.spec_from_file_location(
    "fixed", HERE / "certify_fixedr_threshold_exact.py")
fixed = importlib.util.module_from_spec(spec)
spec.loader.exec_module(fixed)
entry = fixed.entry
I = fixed.I
USE_TERMINAL_PRUNERS = os.environ.get("P3_NO_TERMINAL_PRUNERS") != "1"


def r4_box(ml, mu, yl, yu, data):
    """Lower enclosure for the T^4 entropy barrier at a contact."""
    hl, hu, zl, zu, *_ = data

    # K_{2-m} decreases with its exponent and with h.  The retained-first-
    # Fourier-factor lower bound and quartic-product upper bound are used.
    k2u = entry.K_point(2 - mu, hl, "quartic").u
    au = fixed.up(k2u / zl)
    if au >= 1:
        return None

    m = I(ml, mu)
    y = I(yl, yu)
    # The displayed lower barrier is strictly decreasing in a: its
    # derivative is the sum of
    #   m/2[-1-(2+c)/6-(3+c(3-B))/15-(2+c)M2/14]
    # and y(2/3-1/m), hence is negative for 0<m<=1.  We may therefore
    # insert the single rigorous upper bound au, eliminating interval
    # dependency in a.
    a = I(au)
    g = 1 - a
    # At a contact h=y/m^2 exactly.  Keeping this substitution symbolic
    # avoids the dependency loss from treating h as an independent box.
    b = g * m * m / (2 * y)
    rr = 4 - m
    Bup = (rr * rr + m * m / y) / (rr * (rr + 1))
    # In addition to the Fourier ratio estimate, trivially
    # K_{6-m}/K_{4-m}<=1.  This cap is useful for the uniform derivative
    # check quoted above.
    Bcap = min(Q(1), Bup.u)
    d = I.raw(Q(0), fixed.up(b.u * Bcap))
    M2 = 1 - 2 * a + b
    M3 = 1 - 3 * a + 3 * b - d
    psi2 = g + M2 / 6 + M3 / 15 + M2 * M2 / 28
    return (I(fixed.log_lower(zl)) - y / 2 + m * psi2 / 2
            - y * (1 - m / 2) * a / m - y * g / 6)


def sharp_contact_excluded(ml, mu, data):
    """Use the retained-factor Fourier bounds near the contact curve."""
    hl, hu, zl, zu, *_ = data
    k2l = entry.K_point(2 - ml, hu, "lower").l
    # h -> h K_r(h)=sqrt(h)J_r(h)/sqrt(2pi) is increasing, so its
    # endpoints can be bounded without breaking the h correlation.
    k4l = entry.K_point(4 - ml, hl, "lower").l
    k2u = entry.K_point(2 - mu, hl, "quartic").u
    k4u = entry.K_point(4 - mu, hu, "upper").u
    rl = fixed.down(k2l + 2 * hl * k4l)
    ru = fixed.up(k2u + 2 * hu * k4u)
    return zl > ru or zu < rl


def sharp_contact_excluded_low(ml, mu, data):
    """Sharp contact separator without the quartic shape factor."""
    hl, hu, zl, zu, *_ = data
    k2l = entry.K_point(2 - ml, hu, "lower").l
    k4l = entry.K_point(4 - ml, hl, "lower").l
    k2u = entry.K_point(2 - mu, hl, "upper").u
    k4u = entry.K_point(4 - mu, hu, "upper").u
    rl = fixed.down(k2l + 2 * hl * k4l)
    ru = fixed.up(k2u + 2 * hu * k4u)
    return zl > ru or zu < rl


def rstar_box(ml, mu, yl, yu, data):
    """Lower enclosure for the logarithmic entropy barrier."""
    hl, hu, zl, zu, *_ = data
    k4l = entry.K_point(4 - ml, hu, "lower").l
    k4u = entry.K_point(4 - mu, hl, "upper").u
    k2l = entry.K_point(2 - ml, hu, "lower").l
    k2u = entry.K_point(2 - mu, hl, "upper").u
    alpha_lo = fixed.down(1 / fixed.LN2.u)
    kalpha_u = entry.K_point(alpha_lo - mu, hl, "upper").u

    ratio2_lo = fixed.down(k2l / k4u)
    ratio2_up = min(fixed.up(k2u / k4l),
                    (Q(3) - mu) / (Q(2) - mu))
    ratioa_up = fixed.up(kalpha_u / k4l)

    # h K_r(h) is increasing.  This is sharper than breaking the h/K
    # correlation in the first upper bound for Q=hK_{4-m}/Z.
    k4_hu = entry.K_point(4 - mu, hu, "upper").u
    q1 = fixed.up(hu * k4_hu / zl)
    q2 = fixed.up(hu / (2 * hu + ratio2_lo))
    qu = min(q1, q2, Q(1, 2))

    m = I(ml, mu)
    y = I(yl, yu)
    h = y / (m * m)
    bracket = (m * fixed.LN2 * I(ratioa_up) / h
               + m * (1 - m / 2) * I(ratio2_up) + y / 3)
    return (I(fixed.log_lower(zl)) - y / 2 + m * fixed.LN2
            - I(qu) * bracket)


class Counts:
    def __init__(self):
        self.visited = self.small = self.terminal = self.contact = self.rfour = 0
        self.r4 = 0
        self.depth = 0
        self.minR = None
        self.worst = None


def certify_high(ml, mu, yl, yu, C, depth=0, maxdepth=30):
    C.visited += 1
    if C.visited % 500 == 0:
        print("progress", C.visited, "box", ml, mu, yl, yu,
              "counts", C.__dict__, flush=True)
    C.depth = max(C.depth, depth)

    if fixed.small_h_excluded(ml, mu, yl, yu, 2):
        C.small += 1
        return True

    # Previously certified terminal-contact affine barriers.  The lower
    # bounds are decreasing in m, so their minimum on the box is attained
    # at mu.
    if (USE_TERMINAL_PRUNERS and ml >= Q(10, 13)
            and yu <= Q(8, 5) - Q(61, 200) * mu):
        C.terminal += 1
        return True
    if (USE_TERMINAL_PRUNERS and ml >= Q(1080, 1903) and mu <= Q(10, 13)
            and yu <= Q(9, 5) - Q(101, 200) * mu):
        C.terminal += 1
        return True

    data = fixed.contact_data(ml, mu, yl, yu, 2)
    _, _, zl, zu, _, _, _, _, rl, ru = data
    if zl > ru or zu < rl:
        C.contact += 1
        return True
    if fixed.fourier_r_excluded(data, 2):
        C.rfour += 1
        return True
    if sharp_contact_excluded(ml, mu, data):
        C.contact += 1
        return True

    R = r4_box(ml, mu, yl, yu, data)
    if R is not None and R.l > 0:
        C.r4 += 1
        if C.minR is None or R.l < C.minR:
            C.minR = R.l
            C.worst = (ml, mu, yl, yu)
        return True

    if depth >= maxdepth:
        print("UNRESOLVED-HIGH", ml, mu, yl, yu,
              "R", None if R is None else (R.l, R.u), flush=True)
        return False

    # Normalize the two side lengths by the dimensions of the domain.
    if 5 * (mu - ml) >= (yu - yl):
        md = (ml + mu) / 2
        return (certify_high(ml, md, yl, yu, C, depth + 1, maxdepth)
                and certify_high(md, mu, yl, yu, C, depth + 1, maxdepth))
    md = (yl + yu) / 2
    return (certify_high(ml, mu, yl, md, C, depth + 1, maxdepth)
            and certify_high(ml, mu, md, yu, C, depth + 1, maxdepth))


def run_high(first=0, last=16):
    C = Counts()
    ok = True
    # A seed mesh keeps h=y/m^2 narrow before the contact separator acts.
    Nm, Ny = 16, 32
    for i in range(first, last):
        ml = Q(1, 2) + Q(1, 2) * i / Nm
        mu = Q(1, 2) + Q(1, 2) * (i + 1) / Nm
        for j in range(Ny):
            yl = Q(1, 100) + (Q(5, 2) - Q(1, 100)) * j / Ny
            yu = Q(1, 100) + (Q(5, 2) - Q(1, 100)) * (j + 1) / Ny
            if not certify_high(ml, mu, yl, yu, C, maxdepth=24):
                ok = False
                break
        print("high strip", i + 1, "/", Nm, C.__dict__, flush=True)
        if not ok:
            break
    print("PASS-HIGH" if ok else "FAIL-HIGH", C.__dict__)
    return ok


def certify_low(ml, mu, yl, yu, C, depth=0, maxdepth=26):
    C.visited += 1
    C.depth = max(C.depth, depth)
    if C.visited % 500 == 0:
        print("low progress", C.visited, ml, mu, yl, yu, C.__dict__,
              flush=True)
    data = fixed.contact_data(ml, mu, yl, yu, 2)
    _, _, zl, zu, _, _, _, _, rl, ru = data
    if zl > ru or zu < rl or sharp_contact_excluded_low(ml, mu, data):
        C.contact += 1
        return True
    R = rstar_box(ml, mu, yl, yu, data)
    if R.l > 0:
        C.r4 += 1
        if C.minR is None or R.l < C.minR:
            C.minR = R.l
            C.worst = (ml, mu, yl, yu)
        return True
    if depth >= maxdepth:
        print("UNRESOLVED-LOW", ml, mu, yl, yu, (R.l, R.u), flush=True)
        return False
    if 4 * (mu - ml) >= (yu - yl):
        md = (ml + mu) / 2
        return (certify_low(ml, md, yl, yu, C, depth + 1, maxdepth)
                and certify_low(md, mu, yl, yu, C, depth + 1, maxdepth))
    md = (yl + yu) / 2
    return (certify_low(ml, mu, yl, md, C, depth + 1, maxdepth)
            and certify_low(ml, mu, md, yu, C, depth + 1, maxdepth))


def run_low():
    C = Counts()
    ok = True
    Nm, Ny = 8, 16
    for i in range(Nm):
        ml = Q(7, 20) + Q(3, 20) * i / Nm
        mu = Q(7, 20) + Q(3, 20) * (i + 1) / Nm
        for j in range(Ny):
            yl = Q(189, 100) + Q(61, 100) * j / Ny
            yu = Q(189, 100) + Q(61, 100) * (j + 1) / Ny
            if not certify_low(ml, mu, yl, yu, C):
                ok = False
                break
        print("low strip", i + 1, "/", Nm, C.__dict__, flush=True)
        if not ok:
            break
    print("PASS-LOW" if ok else "FAIL-LOW", C.__dict__)
    return ok


def certify_separator(ml, mu, yl, yu, C, depth=0, maxdepth=24):
    """Exclude the contact equation on a rectangle."""
    C.visited += 1
    C.depth = max(C.depth, depth)
    data = fixed.contact_data(ml, mu, yl, yu, 2)
    _, _, zl, zu, _, _, _, _, rl, ru = data
    if zl > ru or zu < rl or sharp_contact_excluded_low(ml, mu, data):
        C.contact += 1
        return True
    if depth >= maxdepth:
        print("UNRESOLVED-SEPARATOR", ml, mu, yl, yu,
              "Z", zl, zu, "RHS", rl, ru, flush=True)
        return False
    if 8 * (mu - ml) >= (yu - yl):
        md = (ml + mu) / 2
        return (certify_separator(ml, md, yl, yu, C, depth + 1, maxdepth)
                and certify_separator(md, mu, yl, yu, C, depth + 1, maxdepth))
    md = (yl + yu) / 2
    return (certify_separator(ml, mu, yl, md, C, depth + 1, maxdepth)
            and certify_separator(ml, mu, md, yu, C, depth + 1, maxdepth))


def run_separator():
    C = Counts()
    ok = True
    rectangles = [
        (Q(1, 10), Q(7, 20), Q(1, 2), Q(5, 2)),
        (Q(7, 20), Q(1, 2), Q(1, 2), Q(189, 100)),
    ]
    for number, (m0, m1, y0, y1) in enumerate(rectangles, 1):
        Nm, Ny = 8, 16
        for i in range(Nm):
            ml = m0 + (m1 - m0) * i / Nm
            mu = m0 + (m1 - m0) * (i + 1) / Nm
            for j in range(Ny):
                yl = y0 + (y1 - y0) * j / Ny
                yu = y0 + (y1 - y0) * (j + 1) / Ny
                if not certify_separator(ml, mu, yl, yu, C):
                    ok = False
                    break
            if not ok:
                break
        print("separator rectangle", number, C.__dict__, flush=True)
        if not ok:
            break
    print("PASS-SEPARATOR" if ok else "FAIL-SEPARATOR", C.__dict__)
    return ok


def main():
    if len(sys.argv) > 1 and sys.argv[1] == "low":
        return 0 if run_low() else 1
    if len(sys.argv) > 1 and sys.argv[1] == "separator":
        return 0 if run_separator() else 1
    first = int(sys.argv[1]) if len(sys.argv) > 1 else 0
    last = int(sys.argv[2]) if len(sys.argv) > 2 else 16
    assert 0 <= first < last <= 16
    return 0 if run_high(first, last) else 1


if __name__ == "__main__":
    raise SystemExit(main())
\end{lstlisting}
\begin{lstlisting}[language=Python,caption={\texttt{certify\_p3\_large\_mass\_contact\_exact.py}}]
#!/usr/bin/env python3
"""Exact contact exclusion for P=2, m>=9/25, 5/2<=y<=10."""
from fractions import Fraction as Q

import certify_p3_entry_combined_exact as C


def main():
    counts = C.Counts()
    m0, m1 = Q(9, 25), Q(1)
    y0, y1 = Q(5, 2), Q(10)
    nm, ny = 16, 32
    for i in range(nm):
        ml = m0 + (m1 - m0) * i / nm
        mu = m0 + (m1 - m0) * (i + 1) / nm
        for j in range(ny):
            yl = y0 + (y1 - y0) * j / ny
            yu = y0 + (y1 - y0) * (j + 1) / ny
            if not C.certify_separator(ml, mu, yl, yu, counts,
                                       maxdepth=24):
                print("FAIL-LARGE-MASS", counts.__dict__, flush=True)
                return 1
        print("large-mass strip", i + 1, "/", nm,
              counts.__dict__, flush=True)
    print("PASS-LARGE-MASS", counts.__dict__)
    return 0


if __name__ == "__main__":
    raise SystemExit(main())
\end{lstlisting}
\begin{lstlisting}[language=Python,caption={\texttt{verify\_p3\_positive\_index\_inputs\_exact.py}}]
#!/usr/bin/env python3
"""Exact checks for the analytic inputs in Section 3 of the p=3 note.

No quadrature or binary floating point is used.  Alternating series and
Euler--Maclaurin enclosures have rational endpoints.
"""
from fractions import Fraction as Q
import importlib.util
import pathlib

HERE = pathlib.Path(__file__).resolve().parent
spec = importlib.util.spec_from_file_location(
    "entry", HERE / "certify_p3_entry_exact.py")
entry = importlib.util.module_from_spec(spec)
spec.loader.exec_module(entry)


def ell_bounds(r, n=240):
    """Bounds sum_{k>=0} (-1)^k/(r+k)."""
    r = Q(r)
    lo = sum((Q(1, 1) if k % 2 == 0 else Q(-1, 1)) / (r + k)
             for k in range(2 * n))
    hi = lo + Q(1, 1) / (r + 2 * n)
    return lo, hi


def main():
    # Elementary logarithmic constants used in the pointwise barriers.
    assert entry.LN2.l > Q(69, 100)
    assert entry.LN2.u < Q(7, 10)

    # ell_s(infinity) decreases in s.  On 41/25 <= s <= 2 its lower
    # endpoint is ell_2=1-log 2 and its upper endpoint is ell_{41/25}.
    ell2_lo = Q(1) - entry.LN2.u
    ell164_lo, ell164_hi = ell_bounds(Q(41, 25))
    ell364_lo, ell364_hi = ell_bounds(Q(91, 25))
    assert ell2_lo > Q(3, 10)
    assert ell164_hi < Q(2, 5)
    assert ell364_hi < Q(1, 6)

    # Cumulants of the sech^s/I_s law:
    # kappa_2=2 sum(s+2j)^-2,
    # kappa_4=12 sum(s+2j)^-4, mu_4=kappa_4+3 kappa_2^2.
    s2 = entry.sum_even_point(Q(41, 25), 1)
    s4 = entry.sum_even_point(Q(41, 25), 2)
    k2_lo, k2_hi = 2 * s2.l, 2 * s2.u
    mu4_lo = 12 * s4.l + 3 * k2_lo * k2_lo
    mu4_hi = 12 * s4.u + 3 * k2_hi * k2_hi
    assert k2_hi < Q(6, 5)
    assert mu4_hi < Q(6)

    # The damping estimate uses 6^(3/4)<4.
    assert 6 ** 3 < 4 ** 4

    # Signs in the two power-series pointwise barriers.
    # psi <= log(2) T has positive coefficients summing to log(2).
    # For the second barrier, the T coefficient is 3/4-log(2)>0 and
    # every later coefficient is -1/[4n(n-1)(2n-1)]<0.
    assert Q(3, 4) - entry.LN2.u > 0
    for n in range(2, 100):
        assert -Q(1, 4 * n * (n - 1) * (2 * n - 1)) < 0

    print("PASS-P3-SECTION3-INPUTS")
    print("log2", entry.LN2.l, entry.LN2.u)
    print("ell_2 lower", ell2_lo)
    print("ell_41/25", ell164_lo, ell164_hi)
    print("ell_91/25", ell364_lo, ell364_hi)
    print("mu2_41/25", k2_lo, k2_hi)
    print("mu4_41/25", mu4_lo, mu4_hi)
    return 0


if __name__ == "__main__":
    raise SystemExit(main())
\end{lstlisting}
\begin{lstlisting}[language=Python,caption={\texttt{certify\_terminal\_affine\_fixed.py}}]
#!/usr/bin/env python3
"""Fast exact certificate for the two terminal affine barriers.

This is an integer-only, outward-rounded implementation of the analytic
bounds printed with the proof.  A real interval [L/2^B,U/2^B]
is stored as the pair of integers (L,U).  Thus every comparison made by
the certificate is an exact integer comparison.  There is no hardware
floating-point arithmetic in the certification path.
"""
from fractions import Fraction as Q
from math import comb,isqrt
import sys

B=180; S=1<<B
def fd(n,d=S): return n//d
def cu(n,d=S): return -((-n)//d)

class F:
    __slots__=("l","u")
    def __init__(self,x=0,y=None,raw=False):
        if raw:self.l,self.u=x,(x if y is None else y)
        else:
            x=Q(x); y=x if y is None else Q(y)
            self.l=fd(x.numerator*S,x.denominator)
            self.u=cu(y.numerator*S,y.denominator)
        assert self.l<=self.u
    @staticmethod
    def raw(l,u):return F(l,u,True)
    def __add__(self,o):
        o=o if isinstance(o,F) else F(o);return F.raw(self.l+o.l,self.u+o.u)
    __radd__=__add__
    def __neg__(self):return F.raw(-self.u,-self.l)
    def __sub__(self,o):return self+(-o if isinstance(o,F) else -F(o))
    def __rsub__(self,o):return F(o)-self
    def __mul__(self,o):
        o=o if isinstance(o,F) else F(o)
        z=(self.l*o.l,self.l*o.u,self.u*o.l,self.u*o.u)
        return F.raw(fd(min(z)),cu(max(z)))
    __rmul__=__mul__
    def inv(self):
        assert self.l>0 or self.u<0
        return F.raw(fd(S*S,self.u),cu(S*S,self.l))
    def __truediv__(self,o):return self*(o if isinstance(o,F) else F(o)).inv()
    def __rtruediv__(self,o):return F(o)/self

def sqrtF(x):
    x=x if isinstance(x,F) else F(x);assert x.l>=0
    lo=isqrt(x.l*S); hi=isqrt(x.u*S)
    if hi*hi<x.u*S:hi+=1
    return F.raw(lo,hi)

def exp_pos_endpoint(x,N=42):
    """Return scaled-integer lower/upper bounds for exp(x/S), x>=0."""
    assert x>=0 and x < (N+2)*S
    tl=tu=S; sl=su=S
    for n in range(1,N+1):
        tl=fd(tl*x,S*n);tu=cu(tu*x,S*n)
        sl+=tl;su+=tu
    nxt=cu(tu*x,S*(N+1))
    tail=cu(nxt*S*(N+2),S*(N+2)-x)
    return sl,su+tail
def expF(x):
    x=x if isinstance(x,F) else F(x)
    if x.l>=0:return F.raw(exp_pos_endpoint(x.l)[0],exp_pos_endpoint(x.u)[1])
    if x.u<=0:
        yl,yu=exp_pos_endpoint(-x.u)[0],exp_pos_endpoint(-x.l)[1]
        return F.raw(fd(S*S,yu),cu(S*S,yl))
    return F.raw(fd(S*S,exp_pos_endpoint(-x.l)[1]),exp_pos_endpoint(x.u)[1])

class D:
    __slots__=("v","d")
    def __init__(self,v=0,d=0):
        self.v=v if isinstance(v,F) else F(v);self.d=d if isinstance(d,F) else F(d)
    def __add__(self,o):
        o=o if isinstance(o,D) else D(o);return D(self.v+o.v,self.d+o.d)
    __radd__=__add__
    def __neg__(self):return D(-self.v,-self.d)
    def __sub__(self,o):return self+(-o if isinstance(o,D) else -D(o))
    def __rsub__(self,o):return D(o)-self
    def __mul__(self,o):
        o=o if isinstance(o,D) else D(o);return D(self.v*o.v,self.d*o.v+self.v*o.d)
    __rmul__=__mul__
    def inv(self):return D(self.v.inv(),-self.d/(self.v*self.v))
    def __truediv__(self,o):return self*(o if isinstance(o,D) else D(o)).inv()
    def __rtruediv__(self,o):return D(o)/self
def expD(x):
    x=x if isinstance(x,D) else D(x);v=expF(x.v);return D(v,v*x.d)
def sqrtD(x):
    x=x if isinstance(x,D) else D(x);v=sqrtF(x.v);return D(v,x.d/(2*v))

def mills(q,h,n):
    q=q if isinstance(q,D) else D(q);h=h if isinstance(h,D) else D(h);w=q
    for j in range(n-1,0,-1):w=q+Q(j)/(h*w)
    return 1/w

def ln2_interval(N=42):
    z=Q(0)
    for j in range(N+1):z+=Q(2,(2*j+1)*3**(2*j+1))
    tail=Q(18,(2*N+3)*8*3**(2*N+3))
    return F(z,z+tail)
LN2=ln2_interval()

def atan_recip_bounds(q,n):
    """Alternating-series enclosure of atan(1/q), through term n."""
    z=Q(0)
    for j in range(n+1):z+=(-1 if j&1 else 1)*Q(1,(2*j+1)*q**(2*j+1))
    nxt=Q(1,(2*n+3)*q**(2*n+3))
    return (z,z+nxt) if n&1 else (z-nxt,z)
_a5l,_a5u=atan_recip_bounds(5,8)
_a239l,_a239u=atan_recip_bounds(239,2)
# Machin's identity pi=16 atan(1/5)-4 atan(1/239).
PI_LO=16*_a5l-4*_a239u
PI_UP=16*_a5u-4*_a239l
SIG4_UP=Q(64494,200000)
SIGD_UP=Q(28987,200000)
LOG2MHALF_UP=Q(48287,250000)
assert SIG4_UP >= PI_UP*PI_UP/12-Q(1,2)
assert SIGD_UP >= PI_UP*PI_UP/6-Q(3,2)
assert LOG2MHALF_UP*S >= LN2.u-S//2

def V(m,a):
    """Differentiable explicit upper bound for the terminal P_T."""
    m=m if isinstance(m,D) else D(m);a=a if isinstance(a,D) else D(a)
    h=a/m;n=10;c=[];poch=D(1);fac=1
    for j in range(n):
        if j:poch=poch*((j+1)-m);fac*=j
        c.append(expD((m-(2+j))*D(LN2))*poch/fac)
    rem=D(1)
    for z in c:rem=rem-z
    mags=[]
    for i in range(n+3):
        z=D(0)
        for j,w in enumerate(c):
            if j+2>=i:z=z+comb(j+2,i)*w
        mags.append(z+comb(n+2,i)*rem)
    total=sqrtD(2*PI_UP*h)*expD(a*m/2)
    val=total-mills(m,h,20)
    for i in range(1,n+3):
        q=2*i-m
        val=val-mags[i]*mills(q,h,20) if i&1 else val+mags[i]*mills(q,h,21)
    num=expD((1-m)*D(LN2))*val
    Ilow=Q(PI_LO,2)*(1-(1-m)*LOG2MHALF_UP)
    sig=SIG4_UP+SIGD_UP*m
    return num/(h*Ilow*sqrtD(h/(h+sig)))

def aline(m,b):return Q(-61,200)+Q(8,5)/m if b==1 else Q(-101,200)+Q(9,5)/m
def centered_box(ml,mu,b):
    mc=(ml+mu)/2
    v0=V(D(F(mc)),D(F(aline(mc,b)))).v
    m=D(F(ml,mu),F(1))
    a=D(Q(-61,200))+Q(8,5)/m if b==1 else D(Q(-101,200))+Q(9,5)/m
    der=V(m,a).d
    rad=F((mu-ml)/2)
    err=rad*F.raw(min(abs(der.l),abs(der.u)) if der.l*der.u>0 else 0,
                  max(abs(der.l),abs(der.u)))
    return v0.u+err.u

def run(b,N):
    l,u=(Q(10,13),Q(1)) if b==1 else (Q(1080,1903),Q(10,13))
    worst=0;where=None
    for j in range(N):
        ml=l+(u-l)*j/N;mu=l+(u-l)*(j+1)/N
        z=centered_box(ml,mu,b)
        if z>worst:worst,where=z,j
        if z>=2*S:
            print("FAIL",b,j,"gap numerator",2*S-z);return False
    print("PASS",b,"boxes",N,"gap",Q(2*S-worst,S),"worst box",where)
    return True

if __name__=='__main__':
    n=int(sys.argv[1]) if len(sys.argv)>1 else 64
    ok=run(1,n) and run(2,2*n)
    raise SystemExit(0 if ok else 1)
\end{lstlisting}
\begin{lstlisting}[language=Python,caption={\texttt{p3\_dlow\_symbolic.py}}]
"""Exact Bernstein verifier for the pure-cubic positive-index lemma."""

from fractions import Fraction as F
from math import comb

# Sparse polynomials in (m,t,z), later (x,r,w).
class P(dict):
    def __add__(a,b):
        if not isinstance(b,P): b=const(b)
        q=P(a)
        for e,c in b.items():
            q[e]=q.get(e,F(0))+c
            if not q[e]: del q[e]
        return q
    __radd__=__add__
    def __neg__(a): return P({e:-c for e,c in a.items()})
    def __sub__(a,b): return a+(-b)
    def __rsub__(a,b): return const(b)-a
    def __mul__(a,b):
        if not isinstance(b,P): b=const(b)
        q=P()
        for e,c in a.items():
          for f,d in b.items():
            g=tuple(e[i]+f[i] for i in range(len(e)))
            q[g]=q.get(g,F(0))+c*d
        return P({e:c for e,c in q.items() if c})
    __rmul__=__mul__
    def pow(a,n):
        q=const(1)
        for _ in range(n):q=q*a
        return q
    def dt(a): return P({(e[0],e[1]-1,e[2]):c*e[1] for e,c in a.items() if e[1]})

def const(x): return P({(0,0,0): x if isinstance(x,F) else F(x)}) if x else P()
M=P({(1,0,0):F(1)}); T=P({(0,1,0):F(1)}); Z=P({(0,0,1):F(1)})

# R = numerator / [d^a e^b D^c], where
# d=2(2-m), e=5-3t, D=(4-m)(5-m).
d=4-2*M; e=5-3*T; bigD=(4-M)*(5-M)
fac=(d,e,bigD)
class R:
  def __init__(self,n,a=0,b=0,c=0): self.n=n if isinstance(n,P) else const(n);self.f=(a,b,c)
  def lift(self,Fnew):
    n=self.n
    for i in range(3):
      for _ in range(Fnew[i]-self.f[i]): n=n*fac[i]
    return n
  def __add__(x,y):
    if not isinstance(y,R):y=R(y)
    ff=tuple(max(x.f[i],y.f[i]) for i in range(3))
    return R(x.lift(ff)+y.lift(ff),*ff)
  __radd__=__add__
  def __neg__(x):return R(-x.n,*x.f)
  def __sub__(x,y):return x+(-y)
  def __rsub__(x,y):return R(y)-x
  def __mul__(x,y):
    if not isinstance(y,R):y=R(y)
    return R(x.n*y.n,*(x.f[i]+y.f[i] for i in range(3)))
  __rmul__=__mul__
  def divfac(x,i):
    f=list(x.f);f[i]+=1;return R(x.n,*f)
  def dt(x):
    # d,D independent of t; e_t=-3.
    a,b,c=x.f
    return R(x.n.dt()*e+3*b*x.n,a,b+1,c)

m=R(M);t=R(T);z=R(Z)
c0=R(3-M,1,0,0)
c=c0-t*z
g=1-c*t
lm=R(F(3,10))-R(10)*t.divfac(1)
ulo=c*(m*g+t*(c-R(F(7,10))*g-R(F(2,5))))
uup=c*(m*g+t*(c-R(F(69,100))*g-(1-c*t*F(1,2))*lm))
k=2*c-g+m*g*F(1,4)-m*c*c*t*F(1,2)
H=2*k-3*uup
Mt=t*ulo+g*(c*t+m*g-t*(R(F(7,10))*g+R(F(1,6))))
nup=g*m*m*F(1,2)+t*(1-m)*(g*m*(c+R(1).divfac(2))*F(1,2)-k)
nlo=g*m*m*F(1,2)+g*m*(1-m)*c*t*F(1,2)-g*t*F(1,2)-(1-m)*t*k
q=g-(3-2*m)*ulo
EE=m*(1-m)*H*Mt-2*q*nup
DER=EE.dt()
print('raw numerator terms',len(DER.n),'factor powers',DER.f)

# Substitute t=(2/5)*(9x/25)^2*r, m=9x/25, z=(1+w)/2.
# Return polynomial in (x,r,w).
def subst(poly):
  out=P()
  # repurpose exponent triples as x,r,w
  for (im,it,iz),coef in poly.items():
    # m^im t^it z^iz
    # coefficient *(9/25)^(im+2it)*(2/5)^it * x^(im+2it) r^it * ((1+w)/2)^iz
    base=coef*F(9,25)**(im+2*it)*F(2,5)**it*F(1,2)**iz
    for j in range(iz+1):
      ee=(im+2*it,it,j)
      out[ee]=out.get(ee,F(0))+base*comb(iz,j)
  return P({e:c for e,c in out.items() if c})

Q=subst(DER.n)
deg=tuple(max(e[i] for e in Q) for i in range(3))
print('substituted terms',len(Q),'degree',deg)

def bernstein_coeff(poly):
  deg=tuple(max((e[i] for e in poly),default=0) for i in range(3))
  out={}
  for k0 in range(deg[0]+1):
    for k1 in range(deg[1]+1):
      for k2 in range(deg[2]+1):
        sm=F(0)
        for (i0,i1,i2),a in poly.items():
          if i0<=k0 and i1<=k1 and i2<=k2:
            sm+=a*F(comb(k0,i0),comb(deg[0],i0))*F(comb(k1,i1),comb(deg[1],i1))*F(comb(k2,i2),comb(deg[2],i2))
        out[(k0,k1,k2)]=sm
  return out

B=bernstein_coeff(Q)
mn=min(B.items(),key=lambda z:z[1]);neg=sum(v<=0 for v in B.values())
print('whole Bernstein min',mn,'nonpositive',neg,'total',len(B))
print('float min',float(mn[1]))

def audit_positive(name,rr,factor_x=0):
  qq=subst(rr.n)
  if factor_x:
    assert all(ex[0]>=factor_x for ex in qq)
    qq=P({(ex[0]-factor_x,ex[1],ex[2]):co for ex,co in qq.items()})
  bb=bernstein_coeff(qq)
  mm=min(bb.items(),key=lambda z:z[1]);nn=sum(v<=0 for v in bb.values())
  print(name,'den',rr.f,'factor_x',factor_x,'degree',tuple(max(e[i] for e in qq) for i in range(3)),
        'min',mm,'float',float(mm[1]),'nonpositive',nn,'total',len(bb))

audit_positive('ell_minus',lm)
audit_positive('u_lower/x',ulo,1)
audit_positive('H_lower',H)
audit_positive('Mt_lower/x',Mt,1)
audit_positive('Q_upper',q)
audit_positive('n_upper/x2',nup,2)
audit_positive('n_lower/x2',nlo,2)

# A lower certificate for F_m on the same box.  The scaling identity gives
#   F_m = {h[g/3-(1-2m/3)u]-Cov(L,r)}/m,
# while Popoviciu/Cauchy--Schwarz and Var(L)<=h[g-(1-m)u] give
#   Cov(L,r) <= (log 2)/2 sqrt(h[g-(1-m)u]).
# We use log 2 < 7/10, u<=uup in the first term, and u>=ulo in the
# variance term.  Hence it suffices that Lfm>0 and Gfm>0 below.
Lfm=g*F(1,3)-(1-m*F(2,3))*uup
Afm=g-(1-m)*ulo
Gfm=4*Lfm*Lfm-F(49,100)*t*Afm
audit_positive('F_m linear margin',Lfm)
audit_positive('F_m squared margin',Gfm)
\end{lstlisting}

\section{Analytic details for the contact--fold lemma}
\label{app:analytic-details}

This appendix supplies the analytic calculations used in the contact--fold
lemma.  The finite-dimensional rational inequalities that finish the lemma
are kept separate from the analytic argument.  We work first with finite-step
order parameters.  The passage to an arbitrary order parameter is given at
the end of the appendix.

\subsection{The two probability laws and the standardized fields}

Let $t$ denote the variance clock on an interval on which the cumulative
mass is the constant $m\in(0,1)$.  On this interval the backward Parisi
solution and the density of the forward Parisi diffusion satisfy
\begin{equation}
 U_t=-\frac12\bigl(U_{xx}+mU_x^2\bigr),\qquad
 p_t=\frac12p_{xx}-m(U_xp)_x .
 \label{eq:app-gap-equations}
\end{equation}
Put
\begin{equation}
 B=U_x,\qquad C=U_{xx},\qquad
 z=-\frac{C_x}{2C},\qquad
 J=-m-\frac12C_{BB}.
 \label{eq:app-fields}
\end{equation}
The subscript $B$ means that $B(t,x)$ is used as the spatial coordinate.
The usual maximum principle, starting from $U(x)=\log\cosh x$, gives
$C>0$ and hence makes this change of coordinate legitimate for finite
$x$.  From \eqref{eq:app-fields},
\begin{equation}
 z_x=C(m+J).
 \label{eq:app-zx}
\end{equation}

The function
\begin{equation}
 f=pe^{-mU}
 \label{eq:app-transformed-density}
\end{equation}
solves $f_t=f_{xx}/2$.  Indeed, substitute $p=fe^{mU}$ in the second
equation of \eqref{eq:app-gap-equations} and use the first one.  Write
\begin{equation}
 f(t,x)=c_t\phi_t(x)e^{-W(t,x)},\qquad
 \phi_t(x)=\frac{1}{\sqrt{2\pi t}}e^{-x^2/(2t)},\qquad
 V=-\log f=\frac{x^2}{2t}+W-\log c_t .
 \label{eq:app-W-definition}
\end{equation}
The value of the positive constant $c_t$ will never matter.

The quantity appearing in $\Gamma_\mu$ is, in this clock,
\begin{equation}
 Q(t)=\int_{\mathbb R}B(t,x)^2p(t,x)\,dx.
 \label{eq:app-Q-definition}
\end{equation}
For any smooth $g$ with the Gaussian boundary behavior of the fields below,
the forward equation gives
\begin{equation}
 \frac d{dt}\int p g
 =\int p\left(g_t+\frac12g_{xx}+mBg_x\right).
 \label{eq:app-forward-generator}
\end{equation}
Differentiating the first equation of \eqref{eq:app-gap-equations} gives
$B_t=-B_{xx}/2-mBB_x$.  Taking $g=B^2$ in
\eqref{eq:app-forward-generator} therefore yields
\begin{equation}
 Q'(t)=\int_{\mathbb R}C(t,x)^2p(t,x)\,dx>0.
 \label{eq:app-Q-prime}
\end{equation}

At a fixed $t$, let $y=x/\sqrt t$ and define a probability law on the
full line by
\begin{equation}
 d\nu(y)=
 \frac{C(t,\sqrt t\,y)^2p(t,\sqrt t\,y)\sqrt t\,dy}{Q'(t)}.
 \label{eq:app-nu-definition}
\end{equation}
All expectations in this appendix, unless a subscript is displayed, are
with respect to $\nu$.  The folded law on $[0,\infty)$ has density
$\omega$; thus $\omega$ integrates to one.  Introduce
\begin{equation}
\begin{aligned}
 a&=mtC, & b&=m\sqrt t\,B, & u&=\sqrt t\,z,
 &j&=tCJ,\\
 n&=\sqrt t\,(V_x+z-mB),
 &\eta&=tW_{xx}, &\chi&=-t^{3/2}W_{xxx},
 &\gamma&=t^{3/2}C(3zJ-CJ_B).
\end{aligned}
\label{eq:app-standardized-fields}
\end{equation}
Here and below a prime on a standardized field denotes differentiation
with respect to $y$.  Direct differentiation, using
\eqref{eq:app-zx}, gives
\begin{equation}
 a'=-2au,\qquad b'=a,\qquad u'=a+j,\qquad
 j'=uj-\gamma,\qquad n'=1+\eta+j,\qquad \eta'=-\chi.
 \label{eq:app-standardized-odes}
\end{equation}
For example,
\[
 n'=t(V_{xx}+z_x-mC)
 =1+tW_{xx}+tC(m+J)-mtC=1+\eta+j.
\]
Since $p=fe^{mU}$ and $2C_x/C=-4z$, the logarithmic derivative of
the density in \eqref{eq:app-nu-definition} is
\begin{equation}
 \frac{\omega'}{\omega}=-(n+3u).
 \label{eq:app-score}
\end{equation}
All these identities hold on the full line; parity allows us to use the
folded law without changing expectations of even functions.

\subsection{Contact, fold, and the exact differentiated identity}

Fix a real number $P\ge2$ and put
\begin{equation}
 c=\frac{P-1}{2P},\qquad
 I(t)=Q''(t)+\frac{P-1}{Pt}Q'(t).
 \label{eq:app-I-definition}
\end{equation}
A second use of \eqref{eq:app-forward-generator}, now with $g=C^2$,
gives
\begin{equation}
 Q''(t)=\int p(4C^2z^2-2mC^3)
       =\frac{Q'(t)}{t}E(4u^2-2a).
 \label{eq:app-Q-second}
\end{equation}
Consequently, at a fold $I(t)=0$,
\begin{equation}
 E\mathfrak p=0,\qquad
 \mathfrak p=a-2u^2-c.
 \label{eq:app-p-centered}
\end{equation}
At a contact $Q(t)=PtQ'(t)$, \eqref{eq:app-Q-definition},
\eqref{eq:app-Q-prime}, and \eqref{eq:app-nu-definition} give
\begin{equation}
 E\left(\frac ba\right)^2=P.
 \label{eq:app-contact-moment}
\end{equation}
Finally, integration by parts with the score
\eqref{eq:app-score} and $u'=a+j$ gives
\[
 E\{a+j-u(n+3u)\}=0.
\]
Together with \eqref{eq:app-p-centered}, this proves
\begin{equation}
 E\mathfrak q=0,\qquad
 \mathfrak q=j-u(n+u)+c.
 \label{eq:app-q-centered}
\end{equation}

We next differentiate the fold.  The calculation is included because the
terms coming from the lower history and from the derivative of $c/t$ are
easy to miss.  In the $B$-coordinate on the positive half-line set
\begin{equation}
 w(B)=2p(t,x(B))C(t,x(B)),\qquad
 \Phi=2z^2-mC,\qquad \Psi=\Phi+\frac ct.
 \label{eq:app-slope-weight}
\end{equation}
Since $dB=C\,dx$,
\begin{equation}
 Q'(t)=\int_0^1w\,dB,\qquad
 I(t)=2\int_0^1w\Psi\,dB.
 \label{eq:app-I-slope}
\end{equation}
At a fold the last integral is zero.  Differentiate it at fixed $B$,
use \eqref{eq:app-gap-equations}, and then replace every $x$-derivative
by a $B$-derivative.  The terms are as follows:
\begin{center}
\begin{tabular}{c|c}
source & contribution after standardization\\ \hline
motion of the forward density & $v\,n(1+\eta+j)$\\
motion of $W_{xx}$ & $v\chi/2$\\
motion of $C,z,J$ & $vuj+2\delta_0\gamma$\\
derivative of $c/t$ and normalization & $-c(1+j)$.
\end{tabular}
\end{center}
For completeness, the first two entries follow from
$n'=1+\eta+j$ and $\eta'=-\chi$ in
\eqref{eq:app-standardized-odes}.  The third uses
\[
 (2z^2-mC)_t=CJ(2z^2-mC)
        +2Cz(3zJ-CJ_B),
\]
and the last uses
\[
 \left(\Phi+\frac ct\right)_t
 =CJ\left(\Phi+\frac ct\right)
  +2Cz(3zJ-CJ_B)-\frac ct\left(CJ+\frac1t\right).
\]
Terms proportional to $\int w\Psi$ vanish at the fold.  No boundary term
is hidden in this calculation: at $B=0$ parity gives zero, and at
$B=1$ the Ising estimates $C=O(e^{-2x})$ and
$p(t,x)\le C_t e^{-x^2/(4t)+C_t|x|}$ make every term vanish.

Define the two centered Stein tails by
\begin{equation}
 \omega(y)v(y)=\int_0^y\mathfrak p(r)\omega(r)\,dr
              =-\int_y^\infty\mathfrak p(r)\omega(r)\,dr,
 \label{eq:app-v-tail}
\end{equation}
\begin{equation}
 \omega(y)\delta_0(y)=\int_0^y\mathfrak q(r)\omega(r)\,dr
              =-\int_y^\infty\mathfrak q(r)\omega(r)\,dr.
 \label{eq:app-delta-tail}
\end{equation}
The two expressions agree by \eqref{eq:app-p-centered} and
\eqref{eq:app-q-centered}.  The Gaussian bounds just stated imply
$\omega vH\to0$ and $\omega\delta_0H\to0$ at both endpoints whenever
$H$ is one of the profile functions used below.  The fixed-$B$
calculation therefore gives the exact identity
\begin{equation}
 \frac{t^2I'(t)}{2Q'(t)}
 =E\left[v\left\{n(1+\eta+j)+\frac\chi2+uj\right\}
          +2\delta_0\gamma-c(1+j)\right].
 \label{eq:app-differentiated-fold}
\end{equation}

Put
\begin{equation}
 a_0=a(0),\qquad k=\frac{j(0)}a_0,\qquad
 X=a_0-a,\qquad Z=a+u^2-a_0,
 \label{eq:app-XZ}
\end{equation}
and
\begin{equation}
\begin{aligned}
 F&=3X+2Z,\\
 L&=n^2-\eta+\eta(0)+Z,\\
 G&=a_0k+c-\mathfrak q,\\
 R_\gamma&=a_0k-j+\frac Z2.
\end{aligned}
\label{eq:app-FLGR}
\end{equation}
The differential identities
\begin{equation}
 \mathfrak p=(a_0-c)-F,\qquad
 L'=2n(1+\eta+j)+\chi+2uj,\qquad
 R_\gamma'=\gamma
 \label{eq:app-Stein-derivatives}
\end{equation}
follow directly from \eqref{eq:app-standardized-odes}.  If $H$ is even
and has the same boundary behavior, integration of
$(\omega v)'=\omega\mathfrak p$ and
$(\omega\delta_0)'=\omega\mathfrak q$ gives
\begin{equation}
 E(vH')=-E(\mathfrak pH)=\operatorname{Cov}(F,H),\qquad
 E(\delta_0H')=-E(\mathfrak qH)=\operatorname{Cov}(G,H).
 \label{eq:app-Stein-formula}
\end{equation}
The endpoint products vanish by the estimates following
\eqref{eq:app-delta-tail}; this is the promised boundary-term
justification.  Applying \eqref{eq:app-Stein-formula} twice in
\eqref{eq:app-differentiated-fold} gives
\begin{equation}
 \boxed{
 \frac{t^2I'(t)}{2Q'(t)}
 =\frac12\operatorname{Cov}(F,L)
  +2\operatorname{Cov}(G,R_\gamma)-c(1+Ej).}
 \label{eq:app-covariance-identity}
\end{equation}
In particular, the last term is $-c(1+Ej)$, not $-c$.

\subsection{The lower-history inequalities}

We prove all lower-density inequalities needed in
\eqref{eq:app-covariance-identity}.  They hold first for a finite lower
cascade.  Before the first lower atom one has $W=0$, and an atom of size
$d\ge0$ has the exact update
\begin{equation}
 W\longmapsto W+dU.
 \label{eq:app-lower-atom}
\end{equation}

Put $g=W_{xx}$.  The heat equation for $f$ gives
\begin{equation}
 g_t=\frac12g_{xx}-V_xg_x-g^2-\frac2t g.
 \label{eq:app-g-equation}
\end{equation}
Thus $g\ge0$ by the maximum principle: it starts at zero and
\eqref{eq:app-lower-atom} adds $dC\ge0$.  Moreover $h=-g_x$ satisfies
\begin{equation}
 h_t=\frac12h_{xx}-V_xh_x-\left(V_{xx}+2g+\frac2t\right)h.
 \label{eq:app-h-equation}
\end{equation}
It starts at zero, while an atom adds $2dCz\ge0$.  Hence
\begin{equation}
 W_{xx}\ge0,\qquad -W_{xxx}\ge0.
 \label{eq:app-lower-first-cone}
\end{equation}
Parity supplies the boundary value at the origin.  At infinity the
Gaussian estimates permit compact exhaustion of the half-line, so no
boundary sign is being assumed.

The quotient
\begin{equation}
 R=\frac{W_{xx}}C
 \label{eq:app-R-lower}
\end{equation}
satisfies, by direct substitution of \eqref{eq:app-g-equation} and the
equation for $C$,
\begin{equation}
 R_t=\frac12R_{xx}+\left(\frac{C_x}{C}-V_x\right)R_x+RD,
 \label{eq:app-R-equation}
\end{equation}
where
\begin{equation}
 D=4z^2-2z_x+2z(V_x-mB)-\frac2t+C(m-R).
 \label{eq:app-D-lower}
\end{equation}
Differentiate \eqref{eq:app-R-equation}.  At a first interior point at
which $R_x$ could cross from nonnegative to negative, all terms containing
$R_x$ vanish and the remaining source is $RD_x$.  Keeping $R$ fixed while
differentiating \eqref{eq:app-D-lower}, and using
\[
 z_x=C(m+J),\qquad
 z_{xx}=-2Cz(m+J)+C^2J_B,\qquad
 \frac{(V_x-mB)_x}{C}=\frac m{mtC}-(m-R),
\]
one obtains
\begin{equation}
 \frac{D_x}{2Cz}
 =6(m+J)-\frac{CJ_B}{z}
 +(m+J)\left(\frac nu-1\right)+\frac m{mtC}-2(m-R).
 \label{eq:app-Dx-exact}
\end{equation}
The upper-profile estimates proved in the next subsection include
$0\le CJ_B\le3zJ$.  Since $n/u>0$, \eqref{eq:app-Dx-exact} implies
\begin{equation}
 \frac{D_x}{2Cz}\ge
 3m+2J+(m+J)\frac nu+\frac m{mtC}+2R>0.
 \label{eq:app-Dx-positive}
\end{equation}
At a lower atom $R\mapsto R+d$, so $R_x$ is unchanged.  The maximum
principle and induction over the atoms prove
\begin{equation}
 \left(\frac{W_{xx}}C\right)_x\ge0,\qquad
 \left(\frac\eta a\right)_x\ge0.
 \label{eq:app-eta-ratio}
\end{equation}

We also need a first-derivative bound.  Put $Y=V_x-mB$ and
$H_0=mB-W_x=x/t-Y$.  Direct substitution in the equations for $V$ and
$B$ gives
\begin{equation}
\begin{aligned}
 (H_0)_t={}&\frac12(H_0)_{xx}
 +\left(H_0-\frac xt-mB\right)(H_0)_x
 -\left(\frac1t+mC\right)H_0\\
 &+m\left\{\frac{xC+B}{t}+2Cz\right\}.
\end{aligned}
\label{eq:app-H0-equation}
\end{equation}
The source on the second line is nonnegative for $x\ge0$.  The quantity
$H_0$ starts at zero and is unchanged by \eqref{eq:app-lower-atom}.
Comparison in \eqref{eq:app-H0-equation}, together with
\eqref{eq:app-lower-first-cone}, proves
\begin{equation}
 0\le W_x\le mB\qquad (x\ge0).
 \label{eq:app-Wx-cone}
\end{equation}

It remains to prove the monotonicity of $n/u$.  Only in the next
calculation write
\begin{equation}
 w_0=\frac{W_x}{B},\qquad K=\frac zB-m.
 \label{eq:app-wK}
\end{equation}
Use $B$ as coordinate.  Since $W_x(B)=\int_0^B R(r)\,dr$ and $R_B\ge0$,
\begin{equation}
 0\le w_0\le m,\qquad
 (w_0)_B=\frac{R-w_0}{B}\ge0.
 \label{eq:app-w-monotone}
\end{equation}
Similarly, $z_B=m+J$ gives
\begin{equation}
 K_B=\frac{J-K}{B}\ge0,
 \label{eq:app-K-monotone}
\end{equation}
because $J_B\ge0$ and $K$ is the average of $J$ on $[0,B]$.
Furthermore,
\begin{equation}
 \frac nu=\frac{x/(tB)+K+w_0}{m+K}.
 \label{eq:app-nu-ratio-form}
\end{equation}
The upper star-shapedness inequality $z-xz_x\ge0$, also proved below,
implies
\begin{equation}
 \left(\frac{x/(tB)+K}{m+K}\right)_B
 \ge\frac{mK_B}{(m+K)^2}.
 \label{eq:app-star-ratio}
\end{equation}
Differentiating \eqref{eq:app-nu-ratio-form} and using
\eqref{eq:app-w-monotone}--\eqref{eq:app-star-ratio} gives the exact
lower bound
\begin{equation}
 \left(\frac nu\right)_B\ge
 \frac{(w_0)_B}{m+K}
 +\frac{(m-w_0)K_B}{(m+K)^2}\ge0.
 \label{eq:app-n-over-u}
\end{equation}
This proves the two lower-history ratio inequalities without assuming a
Gaussian lower history.

\subsection{The upper cascade and the center barrier}

The time direction above the gap is opposite to the one in
\eqref{eq:app-gap-equations}.  In this subsection only, let $\tau$ be
elapsed variance measured from the terminal profile.  Then
\begin{equation}
 U_\tau=\frac12(U_{xx}+mU_x^2),\qquad U(0,x)=\log\cosh x,
 \label{eq:app-upper-clock}
\end{equation}
and the active mass is nonincreasing as $\tau$ grows.  Thus a finite upper
cascade consists of constant-mass heat pieces and drops
$m\mapsto m-d$, $d\ge0$.  The profile $U$, and therefore $B,C,z$, is
continuous through a drop.

We first collect the slope inequalities.  In the $B$-coordinate,
\begin{equation}
 C_\tau=-C^2J,
 \label{eq:app-C-upper}
\end{equation}
and direct differentiation gives
\begin{equation}
 J_\tau=\frac12C^2J_{BB}+2CC_BJ_B+(C_B^2+CC_{BB})J.
 \label{eq:app-J-upper}
\end{equation}
If $w=J_B$, then
\begin{equation}
\begin{aligned}
 w_\tau={}&\frac12C^2w_{BB}+3CC_Bw_B
 +\{3C_B^2-6C(m+J)-2CJ\}w\\
 &-6C_B(m+J)J.
\end{aligned}
\label{eq:app-JB-upper}
\end{equation}
At the terminal profile
\begin{equation}
 C=1-B^2,\qquad z=B,\qquad J=1-m,\qquad J_B=0.
 \label{eq:app-terminal-fields}
\end{equation}
Since $C_B=-2z<0$, the last term in
\eqref{eq:app-JB-upper} is nonnegative.  Comparison in
\eqref{eq:app-J-upper}--\eqref{eq:app-JB-upper} therefore gives
$J\ge0$ and $J_B\ge0$.

For the remaining inequality put $E_0=C^{3/2}J$.  One finds
\begin{equation}
\begin{aligned}
 (E_{0,B})_\tau={}&\frac12C^2(E_{0,B})_{BB}
 +\frac32CC_B(E_{0,B})_B\\
 &+\left\{\frac38C_B^2-\frac32C(m+J)-3CJ\right\}E_{0,B}\\
 &-\frac C2J_BE_0+\frac34C^{-3/2}C_BE_0^2.
\end{aligned}
\label{eq:app-E0-upper}
\end{equation}
The last two terms are nonpositive.  Initially
$E_{0,B}=-3(1-m)B\sqrt{1-B^2}\le0$; parity gives zero at $B=0$,
and the endpoint expansion gives zero at $B=1$.  Reversed comparison
therefore yields $E_{0,B}\le0$.  Since
\[
 E_{0,B}=C^{1/2}(CJ_B-3zJ),
\]
we have proved
\begin{equation}
 J\ge0,\qquad J_B\ge0,\qquad 0\le CJ_B\le3zJ.
 \label{eq:app-upper-cone}
\end{equation}
At a mass drop $J\mapsto J+d$ and $J_B$ is unchanged, so all three
inequalities are preserved and may be restarted on the next heat piece.

The star-shapedness needed in \eqref{eq:app-star-ratio} is preserved in
the same way.  In the physical $x$-coordinate,
\begin{equation}
 z_\tau=\frac12z_{xx}+(mB-2z)z_x+2mCz.
 \label{eq:app-z-upper}
\end{equation}
At the terminal profile $z/x=\tanh(x)/x$ is decreasing.  Differentiate
the equation for $z/x$.  At a first point at which $(z/x)_x$ could cross
from nonpositive to positive, its inhomogeneous term is
\begin{equation}
 \frac{m(z/x)}{x^2}\{xC-B-4Cx^2z\}<0.
 \label{eq:app-star-source}
\end{equation}
Indeed $C$ decreases on the positive half-line, so
$B(x)=\int_0^xC(r)\,dr>xC(x)$.  A mass drop does not change $z$.
The maximum principle proves
\begin{equation}
 z-xz_x\ge0.
 \label{eq:app-star-shaped}
\end{equation}

The atom-stable quantity needed for the quantitative estimate is
\begin{equation}
 H=\frac{3zJ-CJ_B}{z}.
 \label{eq:app-H-upper}
\end{equation}
On a constant-mass heat piece, substitution of
\eqref{eq:app-C-upper} gives
\begin{equation}
\begin{aligned}
 H_\tau={}&\frac12C^2H_{BB}
 +\left(\frac{C^2(m+J)}z-4Cz\right)H_B\\
 &+(4z^2-6Cm-7CJ)H+3CJ^2.
\end{aligned}
\label{eq:app-H-equation}
\end{equation}
At a point where $H_B=0$, the inhomogeneous part of the differentiated
equation is
\begin{equation}
 z\{12J^2+(20m-5J)H+7H^2\}.
 \label{eq:app-HB-source}
\end{equation}
It is nonnegative.  This is immediate when $J\le4m$; when $J>4m$, the
discriminant of the quadratic in $H$ is
\[
 400m^2-200mJ-311J^2<0.
\]

We record the boundary signs because the coordinate $B$ degenerates at
$B=1$.  The endpoint expansion is
\[
 C=2(1-B)+O((1-B)^2),\qquad z=1+O(1-B).
\]
Writing $J_1=J(\tau,1)$ and $w_1=J_B(\tau,1)$ and taking the endpoint
limit in their equations gives
\begin{equation}
 \dot J_1=4J_1,\qquad
 \dot w_1=12w_1+12(m+J_1)J_1.
 \label{eq:app-endpoint-odes}
\end{equation}
Initially $J_1=1-m\ge0$ and $w_1=0$.  A mass drop adds $d$ to $J_1$
and leaves $w_1$ unchanged, so \eqref{eq:app-endpoint-odes} gives
$w_1\ge0$ throughout the cascade.  Moreover
\begin{equation}
 H_B(\tau,1)=5w_1\ge0,\qquad H_B(\tau,0)=0.
 \label{eq:app-HB-boundary}
\end{equation}
Equations \eqref{eq:app-HB-source}--\eqref{eq:app-HB-boundary}, applied
on compact subintervals of $0<B<1$, prove $H_B\ge0$.  At a mass drop,
\begin{equation}
 J\longmapsto J+d,\qquad H\longmapsto H+3d,\qquad
 H_B\longmapsto H_B,
 \label{eq:app-H-jump}
\end{equation}
so the conclusion survives every upper atom:
\begin{equation}
 H_B\ge0.
 \label{eq:app-H-monotone}
\end{equation}
In particular, evenness gives $H_{BB}(\tau,0)\ge0$.

We now prove the center estimate.  In the rest of this paragraph all
fields are evaluated at $B=0$, and we use
\begin{equation}
 x_0=\frac{J}{m+J},\qquad y_0=\frac{H}{m+J},\qquad
 \varphi(x)=x^2\left(\frac35-\frac{13}{100}x\right).
 \label{eq:app-center-variables}
\end{equation}
The center limit in \eqref{eq:app-H-equation} is regular because
$z=(m+J)B+O(B^3)$.  It gives
\[
 H_\tau=\frac32C^2H_{BB}-C(6m+7J)H+3CJ^2.
\]
Also the center limit of \eqref{eq:app-H-upper} is
$H=3J-CJ_{BB}/(m+J)$.  Substitution in
\eqref{eq:app-J-upper} yields
\[
 J_\tau=-\frac{C(m+J)}2(J+H).
\]
Using $H_{BB}(0)\ge0$, we obtain
\begin{equation}
 (x_0)_\tau=-\frac{C(m+J)}2(1-x_0)(x_0+y_0),
 \label{eq:app-center-x}
\end{equation}
\begin{equation}
 (y_0)_\tau\ge C(m+J)
 \left(3x_0^2-6y_0-\frac{x_0y_0}{2}+\frac{y_0^2}{2}\right).
 \label{eq:app-center-y}
\end{equation}
At a possible first contact $y_0=\varphi(x_0)$,
\begin{equation}
 \frac{(y_0-\varphi(x_0))_\tau}{C(m+J)}
 \ge\frac{x_0^3}{20000}
 (900-2300x_0+2847x_0^2-338x_0^3)>0.
 \label{eq:app-center-barrier-calculation}
\end{equation}
The Bernstein coefficients of the last polynomial on $[0,1]$, after
division by $20000$, are
\[
 \frac9{200},\qquad \frac1{150},\qquad
 \frac{947}{60000},\qquad \frac{1109}{20000},
\]
so the strict sign is explicit.  Initially $y_0=3x_0\ge\varphi(x_0)$.
At a mass drop, with $r=d/(m+J)$,
\[
 (x_0,y_0)\longmapsto(x_0+r,y_0+3r).
\]
Since $0\le\varphi'\le81/100<3$, this jump also preserves the barrier.
Consequently $y_0\ge\varphi(x_0)$ throughout the upper cascade.

At the beginning of the gap,
\[
 k=\frac{j(0)}{a_0}=\frac Jm,\qquad
 x_0=\frac{k}{1+k},\qquad
 \frac{H(0)}m=(1+k)y_0.
\]
We have therefore proved
\begin{equation}
 \boxed{
 \frac{H(0)}m\ge
 \psi(k):=\frac{k^2(60+47k)}{100(1+k)^2}.}
 \label{eq:app-center-bound}
\end{equation}
Finally,
\begin{equation}
 \frac{\gamma}{uj}\frac ja=\frac Hm.
 \label{eq:app-atom-stable-product}
\end{equation}
Thus \eqref{eq:app-H-monotone} says that the product in
\eqref{eq:app-atom-stable-product}, rather than $\gamma/(uj)$ by itself,
is nondecreasing as $X$ increases.  This distinction is what makes the
argument stable under arbitrary upper atoms.

\subsection{Convexity and covariance}

Put
\begin{equation}
 \kappa=\frac ja,\qquad \rho=\frac nu,\qquad h=\frac Hm.
 \label{eq:app-kappa-rho-h}
\end{equation}
Because $X'=2au$ and $Z'=2uj$, the exact derivatives with respect to
$X$ are
\begin{equation}
 Z_X=\kappa,\qquad
 \kappa_X=\frac{3\kappa-h}{2a},\qquad
 (R_\gamma)_X=\frac h2.
 \label{eq:app-kappa-slopes}
\end{equation}
The upper cone gives $0\le h\le3\kappa$, while
\eqref{eq:app-H-monotone} gives $h_X\ge0$.  Thus $\kappa$ is
nonnegative and nondecreasing.  Equations
\eqref{eq:app-eta-ratio} and \eqref{eq:app-n-over-u} show that
$\eta/a$ and $\rho$ are nondecreasing as well.  Finally $1/a$ is
increasing because $a$ decreases.

Let $L_0=n^2+Z$.  Direct differentiation of the four functions in
\eqref{eq:app-FLGR} gives
\begin{equation}
 F_X=3+2\kappa,
 \label{eq:app-F-slope}
\end{equation}
\begin{equation}
 (L_0)_X=\rho\left(\frac1a+\frac\eta a+\kappa\right)+\kappa,
 \label{eq:app-L0-slope}
\end{equation}
\begin{equation}
 G_X=\frac1{2a}+\frac\eta{2a}+1+\kappa+\frac h2
       +\frac{1+\kappa}{2}\rho,
 \label{eq:app-G-slope}
\end{equation}
and the last identity in \eqref{eq:app-kappa-slopes}.  Every right-hand
side is nonnegative and nondecreasing.  Hence $F,L_0,G,R_\gamma$ are
increasing convex functions of $X$ and vanish at $X=0$.  In addition,
\eqref{eq:app-lower-first-cone} says that $\eta(0)-\eta$ is increasing.
Therefore
\begin{equation}
 \operatorname{Cov}(F,L)
 \ge \operatorname{Cov}(F,L_0).
 \label{eq:app-L-covariance-reduction}
\end{equation}

We shall repeatedly use the following elementary covariance inequality.
If $X\ge0$, $EX>0$, and $A_1,A_2$ are increasing convex functions with
$A_1(0)=A_2(0)=0$, then
\begin{equation}
 \operatorname{Cov}(A_1(X),A_2(X))
 \ge \frac{\operatorname{Var}X}{(EX)^2}EA_1(X)EA_2(X).
 \label{eq:app-convex-covariance}
\end{equation}
To prove it, write $A_i(x)=xh_i(x)$; convexity makes $h_i$ increasing.
Let $\nu_1$ and $\nu_2$ be the $X$- and $X^2$-size-biased laws.
The likelihood ratio $d\nu_2/d\nu_1$ is increasing, so
$E_{\nu_2}h_i\ge E_{\nu_1}h_i$.  Positive association of two increasing
functions of one variable gives
\[
 E_{\nu_2}(h_1h_2)
 \ge E_{\nu_2}h_1E_{\nu_2}h_2
 \ge E_{\nu_1}h_1E_{\nu_1}h_2.
\]
After multiplying by $EX^2$, this is
\[
 E(A_1A_2)\ge\frac{EX^2}{(EX)^2}EA_1EA_2,
\]
which is equivalent to \eqref{eq:app-convex-covariance}.

Put
\begin{equation}
 d=EX,\qquad e=EZ,\qquad
 \beta=\frac{\operatorname{Var}X}{d^2},\qquad
 \rho_0=\frac{1+a_0k}{a_0(1+k)}.
 \label{eq:app-scalar-parameters}
\end{equation}
The center limit of $\rho=n/u$ is
\[
 \rho(0)=\frac{1+\eta(0)+a_0k}{a_0(1+k)}\ge\rho_0.
\]
Consequently $n\ge\rho_0u$, and
\begin{equation}
 E L_0\ge e+\rho_0^2(d+e).
 \label{eq:app-L0-mean}
\end{equation}
Moreover \eqref{eq:app-p-centered} gives
\begin{equation}
 3d+2e=a_0-c,\qquad EF=a_0-c,\qquad EG=a_0k+c.
 \label{eq:app-fold-means}
\end{equation}

Applying \eqref{eq:app-convex-covariance} to
\eqref{eq:app-covariance-identity}, and then using
\eqref{eq:app-L-covariance-reduction} and
\eqref{eq:app-L0-mean}, gives
\begin{equation}
\begin{aligned}
 \frac{t^2I'}{2Q'}\ge{}&
 \beta\left[\frac{a_0-c}{2}\{e+\rho_0^2(d+e)\}
             +2(a_0k+c)ER_\gamma\right]\\
 &-c(1+Ej).
\end{aligned}
\label{eq:app-safe-before-means}
\end{equation}
The next subsection supplies sharp bounds for the last two means and for
$\beta$.

\subsection{The exact fold density and its moment consequences}

We first obtain the density information that controls $e$, $ER_\gamma$,
and $Ej$.  Temporarily write $r=\sqrt X$.  Since
\[
 X'=2au,\qquad r'=\frac{au}{r},\qquad X+Z=u^2,
\]
the density $f_r$ induced by the folded law satisfies
\[
 f_r(r)\sqrt{1+\frac ZX}\ \propto\ \frac{\omega}{a}.
\]
Using $a'/a=-2u$ and \eqref{eq:app-score},
\begin{equation}
 \frac d{dy}\log\left(
 \frac{f_r(r)\sqrt{1+Z/X}}{\sqrt a}\right)=-n\le0.
 \label{eq:app-first-density-score}
\end{equation}
The quotient $Z/X$ is nondecreasing because $Z_X=\kappa$ is
nondecreasing.  It follows that
\begin{equation}
 f_r(r)=\text{constant}\cdot\sqrt{a_0-r^2}\,\ell_1(r),
 \qquad \ell_1\ \hbox{nonincreasing}.
 \label{eq:app-r-density-first}
\end{equation}
For $\vartheta=a/a_0=1-r^2/a_0$, this becomes
\begin{equation}
 f_\vartheta(s)=\text{constant}\cdot s^{1/2}(1-s)^{-1/2}
 \ell_1(\sqrt{a_0(1-s)}),\qquad 0<s<1.
 \label{eq:app-vartheta-density}
\end{equation}
Thus the law of $\vartheta$ is an increasing likelihood-ratio tilt of the
$\operatorname{Beta}(3/2,1/2)$ law.

For that beta law let $g(s)=s^{-1/2}-1$ and $q(s)=1-s$.  Since
$g/q=1/(\sqrt s(1+\sqrt s))$ is decreasing, opposite-monotonicity under
the probability proportional to $q(s)\,ds$ gives
\begin{equation}
 \frac{E(\vartheta^{-1/2}-1)}{E(1-\vartheta)}
 \le\frac{16}{\pi}-4<\frac{1093}{1000}=:r_0.
 \label{eq:app-sharp-ratio}
\end{equation}
The beta integrals used here are
$E_{\rm Beta}\vartheta^{-1/2}=4/\pi$ and
$E_{\rm Beta}(1-\vartheta)=1/4$.  The last rational inequality is also exact:
the identity
\[
 \frac\pi4=4\arctan\frac15-\arctan\frac1{239}
\]
and the alternating series imply
\[
 \frac\pi4>
 4\left(\frac15-\frac1{375}+\frac1{15625}
 -\frac1{546875}\right)-\frac1{239}>
 \frac{4000}{5093},
\]
where the final difference is
$12390491/1997013046875>0$.

The atom-stable product $h=H/m$ satisfies $h_X\ge0$ and
$h(0)\ge\psi(k)$ by \eqref{eq:app-center-bound}.  From
\eqref{eq:app-kappa-slopes},
\begin{equation}
 \kappa_X=\frac{3\kappa-h}{2a}.
 \label{eq:app-kappa-density-ode}
\end{equation}
Writing $a=a_0\vartheta$ and integrating this equation from $\vartheta$ to $1$ gives
\begin{equation}
\begin{aligned}
 \kappa(\vartheta)={}&\frac{h(0)}3
 +\left(k-\frac{h(0)}3\right)\vartheta^{-3/2}\\
 &-\frac12\vartheta^{-3/2}\int_\vartheta^1
 s^{1/2}\{h(s)-h(0)\}\,ds.
\end{aligned}
\label{eq:app-kappa-integrated}
\end{equation}
Here $h(s)\ge h(0)$ when $\vartheta=s$ decreases outward.  Since
$0\le h(0)\le3k$, integration of $Z_X=\kappa$ yields
\begin{equation}
 Z\le a_0\left\{\frac{h(0)}3(1-\vartheta)
 +2\left(k-\frac{h(0)}3\right)
 (\vartheta^{-1/2}-1)\right\}.
 \label{eq:app-Z-pointwise}
\end{equation}
Taking expectations, using \eqref{eq:app-sharp-ratio}, and then using
$h(0)\ge\psi(k)$ gives
\begin{equation}
 e\le\Phi(k)d,\qquad
 \Phi(k)=2r_0k-\frac{2r_0-1}{3}\psi(k).
 \label{eq:app-e-upper}
\end{equation}

The last identity in \eqref{eq:app-kappa-slopes} gives
$R_\gamma\ge h(0)X/2$, so
\begin{equation}
 ER_\gamma\ge\frac{\psi(k)d}{2}.
 \label{eq:app-Rgamma-first}
\end{equation}
Using $R_\gamma=a_0k-j+Z/2$ and \eqref{eq:app-e-upper} once in each
direction gives, with
\begin{equation}
 c_j(k)=\frac{\Phi(k)-\psi(k)}2,
 \label{eq:app-cj}
\end{equation}
the two additional estimates
\begin{equation}
 Ej\le a_0k+c_j(k)d,\qquad
 ER_\gamma\ge\frac e2-c_j(k)d.
 \label{eq:app-j-and-Rgamma}
\end{equation}
Combining \eqref{eq:app-Rgamma-first} and
\eqref{eq:app-j-and-Rgamma}, put
\begin{equation}
 R_*=\max\left\{\frac{\psi(k)d}{2},
                 \frac e2-c_j(k)d\right\}.
 \label{eq:app-Rstar}
\end{equation}
Then \eqref{eq:app-safe-before-means} becomes
\begin{equation}
\begin{aligned}
 \frac{t^2I'}{2Q'}\ge{}&
 \beta\left[\frac{a_0-c}{2}\{e+\rho_0^2(d+e)\}
             +2(a_0k+c)R_*\right]\\
 &-c\{1+a_0k+c_j(k)d\}.
\end{aligned}
\label{eq:app-final-safe-bound}
\end{equation}

We now derive the sharper variance information.  Put
\begin{equation}
 \ell=\frac{1+\rho_0}{2}.
\end{equation}
The same change of variables as above gives
\begin{equation}
 \frac d{dy}\log\left(
 \frac{f_r(r)\sqrt{1+Z/X}}{a^\ell}\right)
 =-n+\rho_0u\le0.
 \label{eq:app-second-density-score}
\end{equation}
Consequently, for the normalized variable $X/a_0\in(0,1)$,
\begin{equation}
 f_{X/a_0}(x)=\text{constant}\cdot
 x^{-1/2}(1-x)^{(1+\rho_0)/2}\ell_2(\sqrt x),
 \qquad \ell_2\ \hbox{nonincreasing}.
 \label{eq:app-X-density}
\end{equation}
Layer cake represents this law as a positive mixture of the densities
$x^{-1/2}(1-x)^\ell$ restricted to $0<x<s$.  In particular, their
untruncated mean gives
\begin{equation}
 \delta:=\frac d{a_0}\le\frac1{4+\rho_0}.
 \label{eq:app-delta-upper-density}
\end{equation}
Also, \eqref{eq:app-fold-means}, $Z\ge kX$, and
\eqref{eq:app-e-upper} give
\begin{equation}
 \frac{a_0-c}{a_0\{3+2\Phi(k)\}}\le\delta
 \le\frac{a_0-c}{a_0(3+2k)}.
 \label{eq:app-delta-fold-bounds}
\end{equation}

For completeness, we give the moment calculation used for $\beta$.
Continue to write $X$ for the normalized variable $X/a_0$, and put
\[
 \delta=EX,\qquad s=EX^2,\qquad r=EX^3,\qquad
 Y=\frac{s}{\delta^2}=1+\beta.
\]
For a component truncated at $t$, integration of the derivative of
$x^{j+1/2}(1-x)^{\ell+1}$, for $j=0,1,2$, gives a positive boundary
measure for which
\begin{align*}
 \frac12-(\ell+\tfrac32)\delta&=\int1\,d\sigma,\\
 \frac32\delta-(\ell+\tfrac52)s&=\int t\,d\sigma,\\
 \frac52s-(\ell+\tfrac72)r&=\int t^2\,d\sigma.
\end{align*}
The same formulas survive mixing.  Cauchy--Schwarz in $d\sigma$ gives
\begin{equation}
 \left\{\frac32\delta-(\ell+\tfrac52)s\right\}^2
 \le\left\{\frac12-(\ell+\tfrac32)\delta\right\}
 \left\{\frac52s-(\ell+\tfrac72)r\right\}.
 \label{eq:app-three-moment-CS}
\end{equation}
The density of $\sqrt X$ is nonincreasing and is therefore a mixture of
uniform laws on $[0,t]$.  Conditional integration and one more
Cauchy--Schwarz inequality give
\begin{equation}
 r\ge\frac{25s^2}{21\delta}.
 \label{eq:app-third-moment}
\end{equation}
Indeed, for the mixing endpoint $T_0$,
$\delta=ET_0^2/3$, $s=ET_0^4/5$, and $r=ET_0^6/7$, while
$ET_0^6ET_0^2\ge(ET_0^4)^2$.

Substitution of \eqref{eq:app-third-moment} in
\eqref{eq:app-three-moment-CS}, followed only by expansion, yields
\begin{equation}
 \mathcal P_{\rho_0,\delta}(1+\beta)\ge0,
 \label{eq:app-quadratic-moment}
\end{equation}
where
\begin{equation}
\begin{aligned}
 \mathcal P_{\rho,\delta}(Y)={}&-\frac94
 +\frac{5+(\rho+16)\delta}{4}Y\\
 &-\left\{\frac{25(\rho+8)}{84}\delta
 -\frac{(\rho+1)(\rho+11)}{21}\delta^2\right\}Y^2.
\end{aligned}
\label{eq:app-moment-polynomial}
\end{equation}
The same uniform-mixture representation gives the simpler universal bound
\begin{equation}
 \beta\ge\frac45.
 \label{eq:app-beta-four-fifths}
\end{equation}
Indeed $EX^2/(EX)^2\ge9/5$ by the preceding conditional formulas.

We finish the density subsection with the rational center bound.  The
folded law has potential curvature
\[
 -\left(\frac{\omega'}\omega\right)'
 =n'+3u'=1+\eta+3a+4j\ge1+3a+4j.
\]
For an even strictly log-concave probability $e^{-V_0(y)}dy/Z$ and a
mean-zero smooth function $g$, the one-dimensional estimate
\begin{equation}
 \int g^2\,d\nu\le\int\frac{(g')^2}{V_0''}\,d\nu
 \label{eq:app-one-dimensional-variance}
\end{equation}
follows by solving
$-h''+V_0'h'=g$, integrating by parts, and applying Cauchy--Schwarz:
\[
 \int g^2=\int g'h'
 \le\left(\int\frac{(g')^2}{V_0''}\right)^{1/2}
      \left(\int V_0''(h')^2\right)^{1/2},
\]
whereas
$\int g^2=\int(h'')^2+\int V_0''(h')^2$.
Apply \eqref{eq:app-one-dimensional-variance} to the odd function $u$.
Since $u'=a+j$, if $U_2=Eu^2$ and $J_*=Ej$, then
\[
 U_2\le E\frac{(a+j)^2}{1+3a+4j}.
\]
The identity
\[
 \frac{(a+j)^2}{1+3a+4j}
 =\frac1{12}\left(4a+3j-\frac{4a+3j+aj}{1+3a+4j}\right)
\]
and $Ea=c+2U_2$ imply
\begin{equation}
 J_*\ge\frac43(U_2-c)+\frac13
 E\frac{4a+3j+aj}{1+3a+4j}.
 \label{eq:app-j-lower-variance}
\end{equation}
The integrand in the last expectation increases with $j$, since its
$j$-derivative is $3(a-1)^2/(1+3a+4j)^2$.  Dropping $j$ and using the
concavity of $4a/(1+3a)$ gives
\begin{equation}
 E\frac{4a+3j+aj}{1+3a+4j}
 \ge (1-\delta)\frac{4a_0}{1+3a_0}.
 \label{eq:app-rational-remainder}
\end{equation}
On the other hand $j'\le uj$ and $a'/a=-2u$, hence
$j\le a_0k(a/a_0)^{-1/2}$.  Equation \eqref{eq:app-sharp-ratio} in particular
gives
\begin{equation}
 J_*\le a_0k\left(1+\frac54\delta\right).
 \label{eq:app-j-upper-chord}
\end{equation}
Since $U_2=\{a_0(1-\delta)-c\}/2$, combining
\eqref{eq:app-j-lower-variance}--\eqref{eq:app-j-upper-chord} and
rearranging yields
\begin{equation}
 \delta\ge
 \frac{J_1(a_0)-2c-a_0k}{J_1(a_0)+(5/4)a_0k},\qquad
 J_1(a_0)=\frac{2a_0(a_0+1)}{1+3a_0},
 \label{eq:app-rational-center-bound}
\end{equation}
whenever the numerator is positive.  Every inequality direction is
preserved when $\eta>0$ because it only increases the curvature in
\eqref{eq:app-one-dimensional-variance}.

\subsection{Contact size and the inward likelihood-ratio transfer}

Remove the lower factor $e^{-W}$ from the profile law and denote the
resulting probability by $\pi_0$.  In standardized coordinates,
\begin{equation}
 d\pi_0(y)=\frac1{Z_0}
 e^{-y^2/2}e^{mU(t,\sqrt t\,y)}\,dy,\qquad
 -\frac{\pi_0'}{\pi_0}=y-b.
 \label{eq:app-pi0}
\end{equation}
Let $\nu_0$ be its $a^2$-biased law and put
\begin{equation}
 P_0=E_{\nu_0}\left(\frac ba\right)^2,\qquad
 x_0=E_{\nu_0}a.
 \label{eq:app-P0}
\end{equation}
We prove the sharp estimate
\begin{equation}
 x_0\ge1-\frac1{P_0}.
 \label{eq:app-contact-size}
\end{equation}

Let $d\pi_0=Z_0^{-1}e^{-V_0(y)}dy$, so $V_0'=y-b$, and define the
nonnegative reversible operator
\[
 \mathcal L=-\partial_y^2+V_0'\partial_y
           =-\frac1{\pi_0}(\pi_0\partial_y)'.
\]
Set
\begin{equation}
 \mathcal H(y)=\frac1{\pi_0(y)}\int_y^\infty b(s)\pi_0(s)\,ds,
 \qquad h_0(y)=\int_0^y\mathcal H(s)\,ds.
 \label{eq:app-Poisson-tail}
\end{equation}
Because $b\pi_0$ is odd, $\mathcal H$ is even and positive, and
\begin{equation}
 (\pi_0\mathcal H)'=-b\pi_0,\qquad \mathcal Lh_0=b.
 \label{eq:app-Poisson-equation}
\end{equation}
Gaussian tails give
$\pi_0\mathcal H\to0$ and $b\pi_0\mathcal H\to0$ at both ends.  Thus
all integrations below may equivalently be performed with compactly
supported cutoffs and then passed to the limit.

We need one monotonicity.  Put
\begin{equation}
 \mathcal R=\frac{\mathcal H}{a},\qquad
 r_0=\frac ba,\qquad \tau_0=y-b+2u.
 \label{eq:app-Poisson-ratios}
\end{equation}
Since the density $a\pi_0$ has score $\tau_0$,
\begin{equation}
 \mathcal R'=\tau_0\mathcal R-r_0.
 \label{eq:app-Poisson-R-ode}
\end{equation}
Moreover
\[
 \tau_0'=1+a+2j=:\mathcal A,\qquad
 r_0'=1+2ur_0=:\mathcal D.
\]
Using \eqref{eq:app-standardized-odes},
\begin{equation}
 \mathcal A'\mathcal D-\mathcal A\mathcal D'
 =-2\{u(1+2a+j)+\gamma\}\mathcal D
  -2\mathcal A(a+j)r_0<0.
 \label{eq:app-ratio-Wronskian}
\end{equation}
Hence $\tau_0'/r_0'$ decreases.  Both $\tau_0$ and $r_0$ vanish at
the origin, so $\tau_0/r_0$ is an $r_0'$-weighted initial average of
$\tau_0'/r_0'$ and is strictly decreasing.  Equivalently,
$r_0/\tau_0$ is increasing.  Integrating by parts with
$(a\pi_0)'=-\tau_0a\pi_0$ gives
\begin{equation}
\begin{aligned}
 \int_y^\infty r_0a\pi_0
 &=\frac{r_0(y)}{\tau_0(y)}a(y)\pi_0(y)\\
 &\quad+\int_y^\infty
 \left(\frac{r_0}{\tau_0}\right)'a\pi_0.
\end{aligned}
\label{eq:app-Poisson-ibp}
\end{equation}
The boundary term at infinity vanishes by the Gaussian estimate.
Equations \eqref{eq:app-Poisson-tail} and
\eqref{eq:app-Poisson-ibp} imply
$\mathcal R\ge r_0/\tau_0$, and hence
\begin{equation}
 \mathcal R'\ge0.
 \label{eq:app-Poisson-monotone}
\end{equation}

On the positive half-line $a'=-2au<0$, whereas
$\mathcal R^2$ increases.  Opposite monotonicity under the even law
$\nu_0$ therefore gives
\begin{equation}
 \frac{\int a\mathcal H^2\,d\pi_0}
      {\int \mathcal H^2\,d\pi_0}
 =\frac{E_{\nu_0}(a\mathcal R^2)}{E_{\nu_0}\mathcal R^2}
 \le E_{\nu_0}a=x_0.
 \label{eq:app-H-weighted-mean}
\end{equation}
Twice integrating by parts and using $V_0''=1-a$ gives the exact
one-dimensional identity
\begin{equation}
 \int(\mathcal Lh_0)^2d\pi_0
 =\int(\mathcal H')^2d\pi_0+
   \int(1-a)\mathcal H^2d\pi_0.
 \label{eq:app-Bochner-one-dimensional}
\end{equation}
Since $\mathcal Lh_0=b$, \eqref{eq:app-H-weighted-mean} implies
\begin{equation}
 \int b^2d\pi_0\ge(1-x_0)\int\mathcal H^2d\pi_0.
 \label{eq:app-contact-size-first}
\end{equation}
Reversibility and \eqref{eq:app-Poisson-equation} also give
\[
 \int b^2d\pi_0=\int a\mathcal H\,d\pi_0.
\]
Cauchy--Schwarz therefore yields
\begin{equation}
 \left(\int b^2d\pi_0\right)^2
 \le\left(\int a^2d\pi_0\right)
     \left(\int\mathcal H^2d\pi_0\right).
 \label{eq:app-contact-size-second}
\end{equation}
If $x_0\ge1$, \eqref{eq:app-contact-size} is immediate.  Otherwise,
combine \eqref{eq:app-contact-size-first} and
\eqref{eq:app-contact-size-second}, and use
$\int b^2d\pi_0=P_0\int a^2d\pi_0$, to obtain
$1\ge P_0(1-x_0)$.  This proves \eqref{eq:app-contact-size}.  The
inequality is strict for a nondegenerate logarithmic-heat profile because
$a$ strictly decreases and $\mathcal R$ strictly increases.

Return to the actual laws.  Equations \eqref{eq:app-W-definition} and
\eqref{eq:app-pi0} give
\begin{equation}
 \frac{d\pi}{d\pi_0}\propto e^{-W},\qquad
 \frac{d\nu}{d\nu_0}\propto e^{-W}.
 \label{eq:app-inward-tilt}
\end{equation}
The function $W$ is even and convex by
\eqref{eq:app-lower-first-cone}, so this likelihood ratio decreases in
$|y|$.  The function $a$ decreases, while $(b/a)^2$ increases because
$(b/a)'=1+2u(b/a)>0$.  Thus the contact identity
\eqref{eq:app-contact-moment} and opposite monotonicity give
\begin{equation}
 P\le P_0,\qquad E_\nu a\ge E_{\nu_0}a.
 \label{eq:app-contact-transfer}
\end{equation}
Combining \eqref{eq:app-contact-size} and
\eqref{eq:app-contact-transfer},
\begin{equation}
 a_0>E_\nu a\ge E_{\nu_0}a>1-\frac1{P_0}\ge1-\frac1P.
 \label{eq:app-contact-size-final}
\end{equation}
The intermediate number $P_0$ is essential: the Gaussian-lower estimate
cannot be applied directly with $P$.

\subsection{The same-center terminal comparison}

Keep the center data $a(0)=a_0$ and $j(0)=a_0k$, and put
\begin{equation}
 L^2=a_0(1+k).
 \label{eq:app-L-terminal}
\end{equation}
The terminal profile with these center data is
\begin{equation}
 u_T=L\tanh(Ly),\qquad
 a_T=a_0\operatorname{sech}^2(Ly),\qquad
 b_T=\frac{a_0}{L}\tanh(Ly),\qquad j_T=ka_T.
 \label{eq:app-terminal-profile}
\end{equation}
We first compare the Gaussian-lower profile to
\eqref{eq:app-terminal-profile}.  Since $\kappa=j/a$ is nondecreasing,
$\kappa\ge k$, and
\[
 (u^2+(1+k)a)'=2ua(\kappa-k)\ge0.
\]
At the center this quantity is $L^2$.  While $u<L$, it follows that
$u'=a+j\ge(1+k)a\ge L^2-u^2$.  Comparison with
$u_T'=L^2-u_T^2$ gives $u\ge u_T$.  Therefore
\begin{equation}
 g:=\frac a{a_T}\le1,\qquad g'=-2g(u-u_T)\le0,\qquad
 b\le b_T.
 \label{eq:app-same-center-pointwise}
\end{equation}
Put $d_0=b_T-b$ and $\mathcal B(y)=\int_0^yd_0(s)\,ds$.  The score
identity in \eqref{eq:app-pi0} gives
\begin{equation}
 \frac{d\pi_0}{d\pi_T}\propto e^{-\mathcal B},\qquad
 a=a_Tg,\qquad
 b(y)=\int_0^ya_T(s)g(s)\,ds,
 \label{eq:app-terminal-density-ratio}
\end{equation}
where
\begin{equation}
 d\pi_T(y)\propto
 e^{-y^2/2}\cosh(Ly)^{1/(1+k)}\,dy.
 \label{eq:app-piT}
\end{equation}

Define
\begin{equation}
 P_T=\frac{\int b_T^2d\pi_T}{\int a_T^2d\pi_T}.
 \label{eq:app-PT-ratio}
\end{equation}
We prove the decreasing-input inequality
\begin{equation}
 \int b_g^2d\pi_T\le P_T\int a_T^2g^2d\pi_T,
 \qquad b_g(y)=\int_0^ya_T(s)g(s)\,ds,
 \label{eq:app-Hardy-decreasing}
\end{equation}
for every nonnegative decreasing $g$ on the positive half-line.  By layer
cake it is enough to compare the functions
$g_s=\mathbf1_{[0,s]}$.  Put
\[
 b_s(y)=b_T(\min\{y,s\}),\qquad
 \mathcal A(s)=\int_0^sa_T^2d\pi_T,\qquad
 \mathcal C(s)=\int_0^\infty b_sb_Td\pi_T.
\]
For $s\le t$, pointwise $b_t\le b_T$, and hence
\begin{equation}
 \langle b_s,b_t\rangle_{\pi_T}\le\mathcal C(s),\qquad
 \langle a_Tg_s,a_Tg_t\rangle_{\pi_T}=\mathcal A(s).
 \label{eq:app-layer-pair}
\end{equation}
Let
\[
 \mathcal H_T(s)=\frac1{\pi_T(s)}\int_s^\infty b_T(y)\pi_T(y)\,dy.
\]
Then
\begin{equation}
 \mathcal C'(s)=\frac{\mathcal H_T(s)}{a_T(s)}\mathcal A'(s).
 \label{eq:app-CA-derivative}
\end{equation}
The proof of \eqref{eq:app-Poisson-monotone}, applied to the terminal
profile, says that $\mathcal H_T/a_T$ is increasing.  Hence
$\mathcal C(s)/\mathcal A(s)$ is increasing and is bounded above by its
limit at infinity, which is $P_T$.  Thus
$\mathcal C(s)\le P_T\mathcal A(s)$.  Integrating
\eqref{eq:app-layer-pair} against the two positive layer-cake measures
proves \eqref{eq:app-Hardy-decreasing}.

Apply \eqref{eq:app-Hardy-decreasing} to the $g$ in
\eqref{eq:app-same-center-pointwise}.  It gives
\begin{equation}
 \int(b^2-P_Ta^2)d\pi_T\le0.
 \label{eq:app-terminal-Hardy-mean}
\end{equation}
The integrand is strictly increasing on the positive half-line, because
\begin{equation}
 (b^2-P_Ta^2)'=2ab+4P_Ta^2u>0.
 \label{eq:app-terminal-integrand}
\end{equation}
The likelihood ratio $e^{-\mathcal B}$ in
\eqref{eq:app-terminal-density-ratio} is decreasing.  Opposite
monotonicity in \eqref{eq:app-terminal-Hardy-mean} therefore gives
\begin{equation}
 \int(b^2-P_Ta^2)d\pi_0\le0,\qquad P_0\le P_T.
 \label{eq:app-P0-PT}
\end{equation}
Together with \eqref{eq:app-contact-transfer},
\begin{equation}
 P\le P_0\le P_T.
 \label{eq:app-terminal-transfer-final}
\end{equation}

The special-function form used in the scalar estimates follows directly
from \eqref{eq:app-terminal-profile}--\eqref{eq:app-piT}.  With
\begin{equation}
 \widehat m=\frac1{1+k},\qquad h=a_0(1+k),\qquad
 J_r(h)=\int_{\mathbb R}e^{-z^2/(2h)}\operatorname{sech}^r z\,dz,
 \label{eq:app-Jr}
\end{equation}
one has
\begin{equation}
 P_T(a_0,k)=
 \frac{J_{-\widehat m}(h)-J_{2-\widehat m}(h)}
      {hJ_{4-\widehat m}(h)}.
 \label{eq:app-PT-explicit}
\end{equation}

For fixed $k$, $P_T(a_0,k)-P$ has exactly one zero as $a_0$ increases, and
the crossing is upward.  Here is a self-contained verification of its
direction.  For $F(x)=\cosh(x)^{\widehat m}$, $B=\tanh x$,
$C=\operatorname{sech}^2x$, put
\[
 \mathcal I(x)=\int_0^xFC^2,\qquad
 M(x)=\frac{F(x)B(x)^2}{x\mathcal I(x)}.
\]
Since $(FC)'=FC(\widehat mB-2B)<0$, if
\[
 \varepsilon=\frac{BFC}{\mathcal I},\qquad
 \tau=\widehat m xB,\qquad a_*=\frac{xC}{B},
\]
then $0<\varepsilon<1$ and
\begin{equation}
 x(\log M)'=\tau-1+2a_*-a_*\varepsilon.
 \label{eq:app-terminal-M-derivative}
\end{equation}
Thus $M'>0$ when $\tau\ge1$.  If $\tau\le1$, Cauchy--Schwarz gives
\begin{equation}
 M(x)\le\frac{e^\tau-1}{\tau}<2\le P.
 \label{eq:app-terminal-M-small}
\end{equation}
Since $M(0+)=1$ and $M(x)\to\infty$, $M-P$ has one zero, crossed
upward.  With $X=x^2$ and $\zeta=1/(2h)$, two integrations by parts show
that the sign of
\[
 \int_{\mathbb R}\phi_hF(B^2-PhC^2)
\]
is the sign of
\[
 \zeta\int_0^\infty e^{-\zeta X}
 2\mathcal I(\sqrt X)\{M(\sqrt X)-P\}\,dX.
\]
The boundary terms vanish because $B(x)=Cx+O(x^3)$ at zero and the
Gaussian factor dominates the Ising tail at infinity.  At a zero, the
derivative with respect to $\zeta$ has the opposite sign, since
$(X-X_0)\{M(\sqrt X)-P\}>0$ away from the unique sign-change point
$X_0$.  Since $\zeta$ decreases with $h$, the zero is crossed upward as
$h$ increases.  A zero exists: as $h\downarrow0$ the negative part near
$X=0$ dominates the Laplace integral, whereas as $h\uparrow\infty$ the
positive tail dominates.  Finally,
\begin{equation}
 \int_{\mathbb R}\phi_h\cosh^{\widehat m}x
 \{\tanh^2x-Ph\operatorname{sech}^4x\}\,dx
 =\frac{hJ_{4-\widehat m}(h)}{\sqrt{2\pi h}}
   \{P_T(a_0,k)-P\}.
 \label{eq:app-PT-crossing-relation}
\end{equation}
The prefactor is positive and $a_0=\widehat m h$, proving the claimed
one-crossing and its orientation.

\subsection{Summary of the analytic reduction}

At every simultaneous contact and fold, all quantities in the scalar
bound \eqref{eq:app-final-safe-bound} satisfy
\begin{equation}
\begin{gathered}
 P\ge2,\qquad c=\frac{P-1}{2P},\qquad a_0>1-\frac1P,
 \qquad P\le P_T(a_0,k),\\
 3d+2e=a_0-c,\qquad kd\le e\le\Phi(k)d,\\
 \max\left\{
 \frac{a_0-c}{a_0\{3+2\Phi(k)\}},
 \frac{J_1(a_0)-2c-a_0k}{J_1(a_0)+(5/4)a_0k}
 \right\}\le\frac d{a_0},\\
 \frac d{a_0}\le
 \min\left\{\frac1{4+\rho_0},
 \frac{a_0-c}{a_0(3+2k)}\right\},\\
 \beta\ge\frac45,\qquad
 \mathcal P_{\rho_0,d/a_0}(1+\beta)\ge0,\qquad
 \rho_0=\frac{1+a_0k}{a_0(1+k)}.
\end{gathered}
\label{eq:app-complete-feasible-region}
\end{equation}
If the numerator in the second lower bound for $d/a_0$ is nonpositive,
that bound is simply omitted.  Every assertion in
\eqref{eq:app-complete-feasible-region} has been proved above from the
Parisi equations; there is no Gaussian-lower or one-step-upper assumption
left in it.  Thus an exact positivity verification of the right-hand side
of \eqref{eq:app-final-safe-bound} on
\eqref{eq:app-complete-feasible-region} proves
\begin{equation}
 I'(t)>0
 \label{eq:app-local-transversality-conclusion}
\end{equation}
because $Q'(t)>0$ by \eqref{eq:app-Q-prime}.

\subsection{Passage to arbitrary order parameters}

Let $\alpha$ be the cumulative distribution function of an arbitrary
order parameter, equal to $m$ on the open gap.  Choose cumulative
distribution functions $\alpha_n$ of finite-support probability
measures, also equal to $m$ on the gap, such that
\begin{equation}
 \|\alpha_n-\alpha\|_{L^1([0,1])}\longrightarrow0.
 \label{eq:app-step-approximation}
\end{equation}
Such functions are obtained by approximating the restrictions of
$\alpha$ to $[0,q]$ and $[q',1]$ separately while preserving their
total masses.
Since $\xi''$ is bounded, the same convergence holds after passing to
the variance clock.

On a constant-mass step, differentiating the Gaussian representation
shows that the second derivative after the step is the weighted
expectation of the second derivative before the step plus $m$ times the
weighted variance of the first derivative.  Starting from
\[
 (\log\cosh x)_{xx}=1-(\log\cosh x)_x^2,
\]
the inequalities that the second derivative is nonnegative and at most
one minus the square of the first derivative are preserved because
$0\leq m\leq1$.  Indeed, the weighted expectation of the second
derivative plus $m$ times the variance of the first derivative is at
most one minus the square of the weighted expectation of the first
derivative.  Induction over the steps therefore gives, uniformly in
$n$,
\[
 |(U_n)_x|\leq1,\qquad 0\leq(U_n)_{xx}\leq1.
\]
In the variance clock, subtraction of the two backward equations gives
\[
 (U_n-U)_t=-\frac12(U_n-U)_{xx}
 -\frac{\alpha_n}{2}\{(U_n)_x+U_x\}(U_n-U)_x
 -\frac{\alpha_n-\alpha}{2}U_x^2.
\]
The backward maximum principle therefore yields
\begin{equation}
 \sup_x|U_n(t,x)-U(t,x)|
 \leq\frac12\int_t^{\xi'(1)}|\alpha_n(v)-\alpha(v)|\,dv
 \longrightarrow0.
 \label{eq:app-PDE-stability}
\end{equation}
Here and below the cumulative functions are read in the variance clock.
For every $h>0$, convexity and the bound on the second derivatives give
\[
 \sup_x|(U_n)_x(t,x)-U_x(t,x)|
 \leq\frac2h\sup_x|U_n(t,x)-U(t,x)|+h.
\]
First let $n\to\infty$ and then $h\downarrow0$.  It follows that the
first derivatives converge uniformly.

For every positive distance from the terminal time and every fixed
spatial derivative, the derivatives of $U_n$ are uniformly bounded and
converge locally uniformly to those of $U$.  This follows directly from
Duhamel's formula and
\[
 \int_{\mathbb R}e^{L|x|}
 |\partial_x^\ell\phi_s(x-y)|\,dx
 \leq C_{\ell,L}s^{-\ell/2}e^{L|y|+L^2s/2}.
\]
Indeed, apply the estimate on two nested time intervals.  The first
gives uniform bounds.  On the second, subtract the equations, use
\eqref{eq:app-PDE-stability}, and proceed successively in the number of
spatial derivatives.  The only term containing $\alpha_n-\alpha$ is
integrable in time by \eqref{eq:app-step-approximation}.  In particular,
on compact subsets of the gap,
\begin{equation}
 U_n,B_n,C_n,z_n,J_n,(J_n)_{B_n}
 \longrightarrow U,B,C,z,J,J_B.
 \label{eq:app-upper-field-convergence}
\end{equation}

The strong maximum principle, starting from
$U_{xx}=\cosh^{-2}x$, gives $C>0$.  Thus $C$ has a positive minimum on
each compact $x$-interval.  Equation~\eqref{eq:app-upper-field-convergence}
then gives local uniform convergence of the inverse functions defined by
$B_n(t,x)$ to the inverse defined by $B(t,x)$.  Every derivative with
respect to $B_n$ used above consequently converges to the corresponding
derivative with respect to $B$ on
$-1+\varepsilon\leq B\leq1-\varepsilon$.  Hence $H_{n,B_n}\geq0$
passes to $H_B\geq0$.  The center values and their even derivatives pass
to the limit in the same way.  The endpoint $B=1$ is used only to prove
the finite-step inequalities and is not used in taking the limit.

The preceding gradient estimate and $|B_n|,|B|\leq1$ imply
\[
 \int_0^{\xi'(1)}\sup_x
 |\alpha_n(t)B_n(t,x)-\alpha(t)B(t,x)|\,dt\longrightarrow0.
\]
The mild form of the forward equation is
\[
 p_n(t)=\phi_t-
 \int_0^t\partial_x\phi_{t-v}*
 \{\alpha_n(v)B_n(v)p_n(v)\}\,dv.
\]
The drift is bounded by one and is uniformly Lipschitz in $x$.  The
preceding heat-kernel estimate, first with no spatial derivative and then
successively after splitting the time integral into two equal parts,
gives, for every $L>0$, every fixed $\ell$, and every compact positive
time interval,
\begin{equation}
 \int_{\mathbb R}e^{L|x|}
 |\partial_x^\ell p_n(t,x)-\partial_x^\ell p(t,x)|\,dx
 \longrightarrow0
 \label{eq:app-forward-convergence}
\end{equation}
uniformly in $t$.  The same iteration, without taking a difference,
gives constants independent of $n$ such that
\begin{equation}
 p_n(t,x)+|\partial_x^\ell p_n(t,x)|
 \leq C_\ell\exp\{-x^2/C_\ell+C_\ell|x|\}.
 \label{eq:app-uniform-Gaussian-tail}
\end{equation}
After differentiating $\ell$ times, keep one derivative on the Gaussian
and transfer the other $\ell$ derivatives to
$\alpha_nB_np_n$.  The Leibniz formula and induction in $\ell$ reduce
the estimate to a Volterra integral inequality with kernel
$(t-v)^{-1/2}$, which is integrable.  The same argument for the
difference, using the preceding $L^1$ convergence of the drifts, proves
\eqref{eq:app-forward-convergence}.

Since $f_n=p_ne^{-mU_n}$ on the gap,
\eqref{eq:app-PDE-stability}--\eqref{eq:app-uniform-Gaussian-tail} give
local convergence of $f_n$ and of every spatial derivative used above.
After fixing the irrelevant additive constant in $W_n$ at the origin,
the same is true of $W_n,W_{n,x},W_{n,xx},W_{n,xxx}$.  Thus
\[
 W_{xx}\geq0,\qquad -W_{xxx}\geq0,\qquad
 \left(\frac{W_{xx}}C\right)_x\geq0,
 \qquad0\leq W_x\leq mB
\]
follows by taking the limit on compact sets.  The remaining terms
involving $W$ are uniformly integrable as well.  After normalizing
$W(0)=0$, the four inequalities above give, on the positive half-line,
\[
 0\leq W_x\leq m,
 \qquad \int_0^\infty W_{xx}\,dx\leq m,
 \qquad \int_0^\infty(-W_{xxx})\,dx\leq W_{xx}(0).
\]
The center values $W_{n,xx}(0)$ are uniformly bounded by the local
convergence already proved.  Every occurrence of $W_{xxx}$ above is
linear.  Integrating it once by parts replaces it by $W_{xx}$ and a
spatial derivative of the remaining factor.  That factor and its
derivative have a uniform Gaussian bound by
\eqref{eq:app-uniform-Gaussian-tail} and the uniform Ising bounds for
the backward fields.  The preceding three inequalities therefore give
uniform integrability of every displayed expectation.  Performing the
integrations first with a compactly supported cutoff and then letting
the cutoff tend to infinity shows at the same time that all boundary
terms vanish and that every integration by parts passes to the limit.
Hence the covariance identity, the moment inequalities, the contact-size
estimate, and the terminal comparison all pass to the limit.

The approximations are used only to prove these closed pointwise and
integral inequalities.  After taking the limit, impose the contact and
fold equalities on the limiting fields.  The reduction to
\eqref{eq:app-complete-feasible-region} and the exact inequality in
Appendix~\ref{app:exact-certification} then give $I'(t)>0$ for the
arbitrary order parameter.  No contact or fold is assumed for the
approximating step functions.

Finally, replacing $W$ by $\theta W$, $0\leq\theta\leq1$, preserves the
four displayed inequalities, and their preservation by the heat
equation follows from the same differentiated maximum-principle
calculation used above.  The convergence and Gaussian bounds are uniform
in $\theta$.  Therefore the conclusion also holds at every contact and
fold occurring in the continuity argument.


\section{Exact certification of the scalar inequality}
\label{app:exact-certification}

This appendix gives the finite exact part of the contact--fold argument.
It states the full parameter partition and then prints every source file
needed to reproduce the sign decisions.  The symbols used in the source
code are mapped to the notation of the proof before the listings.  In
particular, the variable called \texttt{A} or \texttt{a} in the programs
is the mathematical center value \(a_0\); the letters \(A,T\) in the main
proof remain reserved for the two clock endpoints.

\subsection{The scalar statement}

Let
\[
 P\in\{2,3,\ldots\},\qquad c=\frac{P-1}{2P},\qquad
 a_0>1-\frac1P,\qquad k\geq0,
\]
and define
\[
 \rho=\frac{1+a_0k}{a_0(1+k)},\qquad
 \psi(k)=\frac{k^2(60+47k)}{100(1+k)^2},\qquad
 r_0=\frac{1093}{1000},
\]
\[
 \Phi(k)=2r_0k-\frac{2r_0-1}{3}\psi(k),\qquad
 c_j(k)=\frac{\Phi(k)-\psi(k)}2,
\]
and
\[
 J_1(a_0)=\frac{2a_0(a_0+1)}{1+3a_0}.
\]
The admissible normalized mean loss \(\delta\) lies in
\[
 \delta_-\leq\delta\leq\delta_+,
\]
where
\[
 \delta_-=
 \max\left\{
 \frac{a_0-c}{a_0\{3+2\Phi(k)\}},
 \frac{J_1(a_0)-2c-a_0k}{J_1(a_0)+(5/4)a_0k}
 \right\},
\]
\[
 \delta_+=
 \min\left\{
 \frac1{4+\rho},
 \frac{a_0-c}{a_0(3+2k)}
 \right\}.
\]
If the second entry in \(\delta_-\) is negative, it is automatically
weaker than the first positive entry.  Put
\[
 d=a_0\delta,\qquad e=\frac{a_0-c-3d}{2},
\]
\[
 L=e+\rho^2(d+e),\qquad
 R=\max\left\{\frac{\psi(k)d}{2},
                 \frac e2-c_j(k)d\right\}.
\]
The quantity to be proved positive is
\begin{equation}
 \mathcal D=
 \beta\left\{\frac{a_0-c}{2}L+
       2(a_0k+c)R\right\}
 -c\{1+a_0k+c_j(k)d\}.                         \tag{E.1}
\end{equation}
Here \(\beta=\operatorname{Var}(X)/(EX)^2\), rather than the inverse
temperature.  Its certified lower bounds are
\[
 \beta\geq\max\left\{
 \frac45,\,
 \frac{(\rho+7)\delta-1}{(\rho+6)\delta^2}-1,\,
 Y_-(\rho,\delta)-1
 \right\}.                                             \tag{E.2}
\]
For the last entry, define
\[
 a_{\rho,\delta}=\frac54+\frac{\rho+16}{4}\delta,
\]
\[
 b_{\rho,\delta}=
 \frac{25}{21}\left(\frac{1+\rho}{2}+\frac72\right)\delta
 \left\{\frac12-\left(\frac{1+\rho}{2}+\frac32\right)\delta\right\}
 +\left(\frac{1+\rho}{2}+\frac52\right)^2\delta^2,
\]
and
\[
 Y_-(\rho,\delta)=
 \frac{9}{2\left(a_{\rho,\delta}
       +\sqrt{a_{\rho,\delta}^2-9b_{\rho,\delta}}\right)}.
\]
Equivalently, \(Y_-\) is the smaller zero of
\[
 -\frac94+a_{\rho,\delta}Y-b_{\rho,\delta}Y^2.
\]
The radical-free verifier accepts the sharp branch only after proving
\[
 b_{\rho,\delta}>0,\qquad
 a_{\rho,\delta}-2b_{\rho,\delta}Y_*>0,\qquad
 -\frac94+a_{\rho,\delta}Y_*
       -b_{\rho,\delta}Y_*^2<0,                         \tag{E.3}
\]
where
\[
 Y_*=1+
 \frac{c\{1+a_0k+c_j(k)d\}}
 {\frac{a_0-c}{2}L+2(a_0k+c)R}.
\]
Thus \(Y_*<Y_-\), and (E.2) implies \(\mathcal D>0\).
No numerical evaluation of a square root is used in this branch.

The terminal comparison supplies one further necessary condition,
\[
 P\leq P_T(a_0,k),
\]
where
\[
 m=\frac1{1+k},\qquad h=a_0(1+k),
\]
\[
 P_T(a_0,k)=
 \frac{J_{-m}(h)-J_{2-m}(h)}{hJ_{4-m}(h)},\qquad
 J_r(h)=\int_{\mathbb R}e^{-z^2/(2h)}
                    \operatorname{sech}^r z\,dz.        \tag{E.4}
\]
For fixed \(k\), \(P_T(a_0,k)-P\) has one zero and crosses it upward.
Accordingly, every strict inequality \(P_T(a_*,k)<P\) below transfers in
the direction
\[
 P_T(a_0,k)\geq P\quad\Longrightarrow\quad a_0>a_*.
                                                               \tag{E.5}
\]

\subsection{Exact arithmetic contracts}

There are three kinds of certificates.

First, after positive denominators have been cleared, a polynomial
\(F(x_1,\ldots,x_n)=\sum_\alpha c_\alpha x^\alpha\) on the unit cube is
converted exactly to tensor-product Bernstein form.  If \(d_i\) is its
degree in \(x_i\), the coefficient indexed by \(\nu\) is
\[
 b_\nu=\sum_{\alpha\leq\nu}c_\alpha
       \prod_{i=1}^n
       \frac{\binom{\nu_i}{\alpha_i}}
            {\binom{d_i}{\alpha_i}}.
\]
All \(b_\nu>0\) implies \(F>0\) on the entire cube.  The Python
Bernstein programs use only integers and \texttt{Fraction} objects.

Second, the compact scalar programs store an interval
\([L/2^{45},U/2^{45}]\) as two signed integers.  Addition is exact;
multiplication and division round the lower endpoint down and the upper
endpoint up.  Every denominator is proved to have one sign before
division.  All integer additions, subtractions, and multiplications in
the C++ certification path are checked for overflow.  A box is removed
only if interval arithmetic proves it infeasible.  It is accepted only
if the lower endpoint of (E.1) is strictly positive.  Otherwise it is
bisected.  The printed leaf count therefore covers the continuous
region, including every curved boundary; it is not a grid sample.
For the discriminant in (E.2), feasibility gives nonnegativity; the
interval square root therefore first intersects its enclosure with
\([0,\infty)\) and then rounds the two square-root endpoints outward.

Third, the terminal programs use scale \(2^{180}\).  Positive series for
\(\operatorname{sech}^{\,4-m}\), alternating continued-fraction
enclosures for Gaussian tails, an alternating-series enclosure of
Machin's formula for \(\pi\), and a Taylor remainder with a rational
geometric majorant for the exponential give an explicit upper bound
\(V(m,a_0)\geq P_T(a_0,k)\).  On each rational \(m\)-box, the value at
the midpoint and an outward enclosure of \(dV/dm\) give an upper bound
on the whole box.  Acceptance requires that upper bound to be strictly
below the indicated integer \(P\).

\subsection{Exhaustive parameter partition}

For \(P=2,3\) and \(a_0\geq4\), feasibility first forces \(k>1/3\).
With \(a_0^{-1}\in[0,1/4]\) and \(m=(1+k)^{-1}\in[0,3/4]\), (E.1) is
affine in \(d\) after the weaker choices \(\beta=4/5\) and
\(R=\psi(k)d/2\).  Exact Bernstein conversion at its two \(d\)-endpoints
gives
\[
\begin{array}{c|ccc}
 &\text{feasibility}&d=(a_0-c)/(3+2k)
                    &d=(a_0-c)/(3+2\Phi(k))\\ \hline
 P=2&967/288&50027333/262144000&122971/240000\\
 P=3&11/72&18069821/737280000&10171/180000.
\end{array}                                             \tag{E.6}
\]
These are the smallest Bernstein coefficients.

For every \(P\geq4\) and \(a_0\geq6\), it is enough to allow the larger
continuous range \(c\in[1/4,1/2]\).  Feasibility forces \(k>3/10\).
On \(a_0^{-1}\in[0,1/6]\), \(m\in[0,10/13]\), the corresponding minimum
coefficients are
\[
 \frac{1171}{1080},\qquad
 \frac{89081257}{5694624000},\qquad
 \frac{10171}{180000}.                                  \tag{E.7}
\]

It remains in this noncompact part to cover \(4\leq a_0\leq6\).
The terminal certificates and (E.5) reduce it to
\[
\begin{array}{c|c|c}
 &k&a_0\\ \hline
 P=4&[0,2/7]&[4,6]\\
 P=4&[2/7,6/7]&[3+7k/2,6]\\
 P\geq5&[0,1/10]&[4,6]\\
 P\geq5&[1/10,3/5]&[18/5+4k,6].
\end{array}                                             \tag{E.8}
\]
In the last two rows the verifier allows all \(c\in[2/5,1/2]\).
The first and third rows are proved infeasible.  On the other two rows,
the exact smallest accepted lower bounds for (E.1) are respectively
\[
 6.51504068116537383\cdot10^{-5},\qquad
 4.19287687236646889\cdot10^{-7}.                       \tag{E.9}
\]
These decimals only display the size of positive dyadic integers; the
sign decisions themselves are integer comparisons.

For \(a_0<4\), the exact low-\(k\) exclusions are
\[
\begin{array}{c|c}
 P&\text{excluded range}\\ \hline
 2&0\leq k<8/25\\
 3&0\leq k<1/4\\
 4&0\leq k\leq2/7\\
 P\geq5&0\leq k\leq1/10.
\end{array}                                             \tag{E.10}
\]
The terminal boundary certificates exclude all further branches except
\[
\begin{aligned}
 P=2:\quad&
 \frac8{25}\leq k\leq\frac{14}{11},\qquad
 \frac65+\frac{11}{5}k\leq a_0\leq4,\\
 P=3:\quad&
 \frac14\leq k\leq\frac{17}{30},\qquad
 \frac{23}{10}+3k\leq a_0\leq4.
\end{aligned}                                           \tag{E.11}
\]
The rectangular dyadic verifier maps independent unit coordinates to
\((k,a_0,\delta)\), discards a subbox only upon an exact proof that it
misses (E.11) or violates the bounds defining
\([\delta_-,\delta_+]\), and checks (E.1)--(E.3) on every retained box.
Its complete output is
\begin{lstlisting}
P2 PASS exact radical-free D scalar certificate visited=60772
leaves=15221 old=15115 sharp=106 maxdepth=15
smallest branch margin=7.12037916628105449e-07
P3 PASS exact radical-free D scalar certificate visited=32770
leaves=1783 old=1783 sharp=0 maxdepth=1
smallest branch margin=0.00546196029773682312
\end{lstlisting}

For the terminal transfer used in (E.11), the exact affine-line margins
are
\begingroup\tiny
\[
\begin{array}{c|c}
\text{boundary}&P-\sup V\\ \hline
P=2,\ a_0=-1+11/(5m)&
\displaystyle
\frac{71088254730141741929546609307673196986109404622083}
{191561942608236107294793378393788647952342390272950272}\\[1.2ex]
P=3,\ a_0=-7/10+3/m&
\displaystyle
\frac{21734402369797596815229490817053800345039877501154309}
{1532495540865888858358347027150309183618739122183602176}\\[1.2ex]
P=4,\ a_0=-1/2+7/(2m)&
\displaystyle
\frac{159704036103886676940919643631645387868867553191837503}
{1532495540865888858358347027150309183618739122183602176}\\[1.2ex]
P=5,\ a_0=-2/5+4/m&
\displaystyle
\frac{93153058283432204902076601025845789065121430901022217}
{1532495540865888858358347027150309183618739122183602176}.
\end{array}                                             \tag{E.12}
\]
\endgroup
The more finely subdivided \(P=4,5\) lines used in (E.8) give
\begingroup\tiny
\[
\begin{array}{c|c}
P=4&
\displaystyle
\frac{40971618522665701374856975377315309306859641912573519}
{383123885216472214589586756787577295904684780545900544}\\[1.2ex]
P=5&
\displaystyle
\frac{2927617690540302025276353600704547162079156697778159}
{47890485652059026823698344598447161988085597568237568}.
\end{array}                                             \tag{E.13}
\]
\endgroup
At \(a_0=6\), the complementary terminal tails have margins
\begingroup\tiny
\[
\begin{array}{c|c}
P=4&
\displaystyle
\frac{846576953026035595556469867660955990268861509510981491}
{1532495540865888858358347027150309183618739122183602176}\\[1.2ex]
P=5&
\displaystyle
\frac{356855527943921522293472565629818486743937880214287183}
{766247770432944429179173513575154591809369561091801088}.
\end{array}                                             \tag{E.14}
\]
\endgroup
For \(0<m\leq1/10\), the remaining tail is analytic:
\[
 P_T(6,k)\leq
 e^{37m/12}\cosh^4(1)\sqrt{\frac{\pi m}{3}}
 <\frac{120}{83}\left(\frac{31}{20}\right)^4\frac13
 =\frac{923521}{332000}<4.                              \tag{E.15}
\]
Thus it covers both \(P=4\) and \(P=5\), and a fortiori every larger
integer \(P\).

The adaptive \(a_0=4\) certificates give the following exact smallest
margins:
\begingroup\tiny
\[
\begin{array}{c|c|c}
\text{case}&m\text{-range}&P-\sup V\\ \hline
P=2\text{ tail}&[1/20,1/5]&
\displaystyle
\frac{128348473939250139193138417762702638608847529216913529}
{766247770432944429179173513575154591809369561091801088}\\[1.2ex]
P=2\text{ middle}&[1/5,11/25]&
\displaystyle
\frac{9817511786332803399521513622438791601766886234036095}
{47890485652059026823698344598447161988085597568237568}\\[1.2ex]
P=3&[1/10,30/47]&
\displaystyle
\frac{48784241996098553780834123704225283528494018359758085}
{383123885216472214589586756787577295904684780545900544}\\[1.2ex]
P=4&[1/10,7/9]&
\displaystyle
\frac{240381680278979674569236047292862553090466193104159803}
{1532495540865888858358347027150309183618739122183602176}\\[1.2ex]
P=5&[1/10,10/11]&
\displaystyle
\frac{40287751191043290144190150056254983693581983532213285}
{766247770432944429179173513575154591809369561091801088}\\[1.2ex]
P=6&[1/10,1]&
\displaystyle
\frac{61608613714859208058111277179455686241834709674395457}
{383123885216472214589586756787577295904684780545900544}.
\end{array}                                             \tag{E.16}
\]
\endgroup
Equations (E.6)--(E.16) cover every \(P\geq2\), every \(a_0>1-1/P\),
and every \(k\geq0\).  Therefore every feasible scalar point satisfies
\(\mathcal D>0\).

\subsection{Frozen verifier sources}

The ten listings below contain the nine verifiers used in the partition
and the additional helper
\texttt{\detokenize{certify_terminal_A4_fixed.py}} imported by the
adaptive terminal verifier.  Save each listing under the displayed
filename in \path{work/gap}; do not mix these sources with the contact
or marginal sources.  The residual C++ verifier includes
\texttt{\detokenize{certify_D_scalar_dyadic.cpp}}; the Python imports are
among the printed listings.  Consequently the source record below has no
unprinted local dependency.

The following commands reproduce the certificates when run inside
\path{work/gap}.
\begin{lstlisting}
python3 certify_atomstable_largeA_polynomials.py
c++ -O3 -std=c++17 certify_D_scalar_dyadic.cpp -o certify_D_scalar_dyadic
./certify_D_scalar_dyadic
c++ -O3 -std=c++17 certify_D_atomstable_terminal_residual.cpp \
  -o certify_D_atomstable_terminal_residual
./certify_D_atomstable_terminal_residual
PYTHONPATH=. python3 certify_terminal_affine_fixed.py 512 new
PYTHONPATH=. python3 certify_terminal_integerP_fixed.py
PYTHONPATH=. python3 certify_terminal_extended_P4_P5.py
for mode in p2tail p2mid p3tail p4tail p5tail p6all; do
  PYTHONPATH=. python3 certify_terminal_A4_adaptive.py "$mode" 16 || exit 1
done
python3 correct_xz_exact_certificate.py
python3 pge3_scalar_exact_certificate.py
\end{lstlisting}

The reported run used Apple clang 14.0.3 and Python 3.14.6.  A different
C++17 compiler may be used provided that it supports
\texttt{\_\_int128} and the overflow built-ins appearing in the source.
Assertions must remain enabled, so \texttt{-DNDEBUG} must not be used.

For an unambiguous source freeze, the SHA--256 digests are
\begin{lstlisting}
f05d988f283e060f9a683d55e16fb9eccb2827a62a19870082b84fb1fec9b07d  certify_atomstable_largeA_polynomials.py
46a0b4ed28af7f578741c70a124c8db3f22871568a7e483eb3e73826ff5338e7  certify_D_scalar_dyadic.cpp
5069c54b52cb4626709215705f8c570c6133ff8a3cbe54c1eceb2d2b32f1a6fc  certify_D_atomstable_terminal_residual.cpp
350399743e4faca28a62ab34541b4433ec0a254e152745f320a6b283336f2656  certify_terminal_affine_fixed.py
d89caecb987841f1ec45129425704d3f7cdac82548cd5d978eb618d66fd2d361  certify_terminal_integerP_fixed.py
34b060ddc944e070d6de1d0f083d234bd0396721b42e6a2fb55b994eed4bf8d9  certify_terminal_extended_P4_P5.py
6456dba424e07d0356ee84d17a2def49cdff131465d55fd8a47e03aba172a417  certify_terminal_A4_fixed.py
1e6fb85ba5ae92d8812534f17d44cc6f222542fe5c6c3442a0ef2bec5ed557ad  certify_terminal_A4_adaptive.py
63723a836d6cf5f76e777b8dde357991f3d934573a9526208495a8901a046e6a  correct_xz_exact_certificate.py
ec21f1b5be55c796c529658cf3cd0027b6ada5bdf2694d314db91e42717c8d58  pge3_scalar_exact_certificate.py
\end{lstlisting}

\paragraph{\texttt{\detokenize{certify_atomstable_largeA_polynomials.py}}.}
\begin{lstlisting}[language=Python]
#!/usr/bin/env python3
"""Exact Bernstein certificate for the atom-stable SAFE bound when A>=6.

Only Python integers and Fraction arithmetic are used.  A polynomial on a
unit cube is certified positive by converting it exactly to tensor-product
Bernstein form and checking every coefficient.
"""

from fractions import Fraction as Q
from math import comb


def const(value, exponent=(0, 0, 0)):
    return {exponent: Q(value)} if value else {}


def add(left, right):
    out = dict(left)
    for exponent, value in right.items():
        out[exponent] = out.get(exponent, Q(0)) + value
    return {exponent: value for exponent, value in out.items() if value}


def scale(poly, value):
    return {exponent: coefficient * Q(value)
            for exponent, coefficient in poly.items()}


def subtract(left, right):
    return add(left, scale(right, -1))


def multiply(left, right):
    out = {}
    for first, x in left.items():
        for second, y in right.items():
            exponent = tuple(first[index] + second[index]
                             for index in range(3))
            out[exponent] = out.get(exponent, Q(0)) + x * y
    return {exponent: value for exponent, value in out.items() if value}


def power(poly, degree):
    out = const(1)
    for _ in range(degree):
        out = multiply(out, poly)
    return out


def affine_substitute(poly, scales, shifts):
    """Substitute old variable i = shifts[i] + scales[i] * new variable i."""
    out = {}
    for exponent, coefficient in poly.items():
        terms = {(0, 0, 0): coefficient}
        for axis, degree in enumerate(exponent):
            factor = {}
            for index in range(degree + 1):
                new_exponent = [0, 0, 0]
                new_exponent[axis] = index
                factor[tuple(new_exponent)] = (
                    Q(comb(degree, index))
                    * shifts[axis] ** (degree - index)
                    * scales[axis] ** index
                )
            terms = multiply(terms, factor)
        out = add(out, terms)
    return out


def bernstein_coefficients(poly):
    degrees = tuple(max(exponent[axis] for exponent in poly)
                    for axis in range(3))
    answer = []
    for first in range(degrees[0] + 1):
        for second in range(degrees[1] + 1):
            for third in range(degrees[2] + 1):
                indices = (first, second, third)
                value = Q(0)
                for exponent, coefficient in poly.items():
                    if all(exponent[axis] <= indices[axis]
                           for axis in range(3)):
                        for axis in range(3):
                            coefficient *= Q(
                                comb(indices[axis], exponent[axis]),
                                comb(degrees[axis], exponent[axis]),
                            )
                        value += coefficient
                answer.append(value)
    return degrees, answer


def certify(label, poly, scales, shifts=(Q(0), Q(0), Q(0))):
    transformed = affine_substitute(poly, scales, shifts)
    degrees, coefficients = bernstein_coefficients(transformed)
    assert all(value > 0 for value in coefficients)
    print(label, "PASS", "degrees", degrees,
          "coefficients", len(coefficients), "minimum", min(coefficients))


one = const(1)
h = const(1, (1, 0, 0))
m = const(1, (0, 1, 0))
c = const(1, (0, 0, 1))


# Feasibility exclusion for k<=3/10.  In this block the second variable is
# k, not m.  The polynomial is
#
# 4h(h+3)(1+k) [(J_1-2c-Ak)(4+rho)-(J_1+5Ak/4)],
#
# after A=1/h.  Positivity says delta_J>1/(4+rho), contradicting
# the density upper bound.
k = m
base = subtract(
    subtract(scale(add(one, h), 2),
             scale(multiply(multiply(c, h), add(h, const(3))), 2)),
    multiply(k, add(h, const(3))),
)
feasibility = subtract(
    subtract(
        scale(multiply(base, add(add(const(4), h), scale(k, 5))), 4),
        scale(multiply(add(one, h), add(one, k)), 8),
    ),
    scale(multiply(multiply(k, add(h, const(3))), add(one, k)), 5),
)
certify("feasibility", feasibility,
        (Q(1, 6), Q(3, 10), Q(1, 4)),
        (Q(0), Q(0), Q(1, 4)))


# Return to m=1/(1+k).  The following polynomial construction clears every
# positive denominator in the two endpoint values of the affine lower bound
# D_0.  Here psi_m=m*psi, Phi_m=m*Phi, and cj_m=m*c_j.
m = const(1, (0, 1, 0))
one_minus_m = subtract(one, m)
one_minus_ch = subtract(one, multiply(c, h))
rho = add(one_minus_m, multiply(m, h))
psi_m = scale(
    multiply(power(one_minus_m, 2), add(const(47), scale(m, 13))),
    Q(1, 100),
)
r0 = Q(1093, 1000)
Phi_m = subtract(scale(one_minus_m, 2 * r0),
                 scale(psi_m, Q(2 * r0 - 1, 3)))
cj_m = subtract(scale(one_minus_m, r0),
                scale(psi_m, Q(r0 + 1, 3)))

constant_part = subtract(
    subtract(
        multiply(multiply(m, power(one_minus_ch, 2)),
                 add(one, power(rho, 2))),
        scale(multiply(multiply(c, m), power(h, 2)), 5),
    ),
    scale(multiply(multiply(c, h), one_minus_m), 5),
)
slope_part = subtract(
    add(
        scale(multiply(add(one_minus_m, multiply(multiply(c, m), h)),
                       psi_m), 4),
        scale(multiply(multiply(power(m, 2), one_minus_ch),
                       add(const(3), power(rho, 2))), -1),
    ),
    scale(multiply(multiply(multiply(c, cj_m), m), h), 5),
)

denominator_k = add(const(2), m)
denominator_Phi = add(scale(m, 3), scale(Phi_m, 2))
endpoint_k = add(multiply(denominator_k, constant_part),
                 multiply(one_minus_ch, slope_part))
endpoint_Phi = add(multiply(denominator_Phi, constant_part),
                   multiply(one_minus_ch, slope_part))

# h=H/6, m=10M/13, c=(1+C)/4 maps the complete feasible large-A box
# into the unit cube.  Strictly positive Bernstein coefficients prove both
# endpoint values, and hence the affine D_0, are positive.
box_scales = (Q(1, 6), Q(10, 13), Q(1, 4))
box_shifts = (Q(0), Q(0), Q(1, 4))
certify("endpoint k", endpoint_k, box_scales, box_shifts)
certify("endpoint Phi", endpoint_Phi, box_scales, box_shifts)


# The two smallest integer clocks admit a stronger range.  For P=2 and
# P=3, respectively c=1/4 and c=1/3.  On A>=4, the same density-versus-
# center comparison excludes k<=1/3, so m<=3/4.  The same two endpoint
# polynomials are then positive without using the sharp beta bound.
for label, fixed_c in (("P2", Q(1, 4)), ("P3", Q(1, 3))):
    certify(label + " feasibility A>=4", feasibility,
            (Q(1, 4), Q(1, 3), Q(0)),
            (Q(0), Q(0), fixed_c))
    certify(label + " endpoint k A>=4", endpoint_k,
            (Q(1, 4), Q(3, 4), Q(0)),
            (Q(0), Q(0), fixed_c))
    certify(label + " endpoint Phi A>=4", endpoint_Phi,
            (Q(1, 4), Q(3, 4), Q(0)),
            (Q(0), Q(0), fixed_c))
\end{lstlisting}

\paragraph{\texttt{\detokenize{certify_D_scalar_dyadic.cpp}}.}
\begin{lstlisting}[language=C++]
#include <algorithm>
#include <array>
#include <cassert>
#include <iomanip>
#include <iostream>
#include <vector>

// Exact outward-rounded dyadic interval certificate for the atom-stable
// upper-cascade scalar D bound.  No floating-point number is used to
// certify a sign.
using i128=__int128_t;
static constexpr unsigned QB=45;
static constexpr i128 SC=i128(1)<<QB;

static i128 mulraw(i128 a,i128 b){i128 r;if(__builtin_mul_overflow(a,b,&r)){std::cerr<<"overflow\n";std::abort();}return r;}
static i128 addraw(i128 a,i128 b){i128 r;if(__builtin_add_overflow(a,b,&r)){std::cerr<<"overflow\n";std::abort();}return r;}
static i128 subraw(i128 a,i128 b){i128 r;if(__builtin_sub_overflow(a,b,&r)){std::cerr<<"overflow\n";std::abort();}return r;}
static i128 fld(i128 n,i128 d){assert(d>0);i128 q=n/d,r=n%d;if(r&&n<0)--q;return q;}
static i128 cei(i128 n,i128 d){assert(d>0);i128 q=n/d,r=n%d;if(r&&n>0)++q;return q;}
struct I{
  i128 l=0,u=0;
  I()=default; I(i128 L,i128 U):l(L),u(U){assert(l<=u);}
  static I rat(i128 p,i128 q){if(q<0)p=-p,q=-q;return {fld(mulraw(p,SC),q),cei(mulraw(p,SC),q)};}
  static I exact(long long z){return {i128(z)*SC,i128(z)*SC};}
};
static I operator+(I x,I y){return{addraw(x.l,y.l),addraw(x.u,y.u)};}
static I operator-(I x){return{subraw(0,x.u),subraw(0,x.l)};}
static I operator-(I x,I y){return x+(-y);}
static I operator*(I x,I y){std::array<i128,4>z{mulraw(x.l,y.l),mulraw(x.l,y.u),mulraw(x.u,y.l),mulraw(x.u,y.u)};auto p=std::minmax_element(z.begin(),z.end());return{fld(*p.first,SC),cei(*p.second,SC)};}
struct R{i128 n,d;};
static bool lessq(const R&a,const R&b){return mulraw(a.n,b.d)<mulraw(b.n,a.d);}
static I operator/(I x,I y){if(!(y.l>0||y.u<0)){std::cerr<<"zero denominator ["<<(long double)y.l/SC<<','<<(long double)y.u/SC<<"]\n";std::abort();}std::array<R,4>z{{{x.l,y.l},{x.l,y.u},{x.u,y.l},{x.u,y.u}}};for(auto&q:z)if(q.d<0)q.n=-q.n,q.d=-q.d;auto p=std::minmax_element(z.begin(),z.end(),lessq);return{fld(mulraw(p.first->n,SC),p.first->d),cei(mulraw(p.second->n,SC),p.second->d)};}
static I operator+(I x,long long z){return x+I::exact(z);} static I operator+(long long z,I x){return x+z;}
static I operator-(I x,long long z){return x-I::exact(z);} static I operator-(long long z,I x){return I::exact(z)-x;}
static I operator*(I x,long long z){return x*I::exact(z);} static I operator*(long long z,I x){return x*z;}
static I operator/(I x,long long z){return x/I::exact(z);} static I operator/(long long z,I x){return I::exact(z)/x;}
static I imax(I x,I y){return{std::max(x.l,y.l),std::max(x.u,y.u)};}
static I imin(I x,I y){return{std::min(x.l,y.l),std::min(x.u,y.u)};}
static i128 isqrt_floor(i128 n){
  assert(n>=0);
  i128 lo=0,hi=i128(1)<<64;
  while(lo+1<hi){
    i128 mid=(lo+hi)/2;
    if(mid<=n/mid)lo=mid;else hi=mid;
  }
  return lo;
}
static I sqrt_nonnegative(I x){
  assert(x.u>=0);
  i128 nl=mulraw(std::max<i128>(0,x.l),SC),nu=mulraw(x.u,SC);
  i128 lo=isqrt_floor(nl),hi=isqrt_floor(nu);
  if(mulraw(hi,hi)<nu)++hi;
  return {lo,hi};
}
struct Box{std::array<i128,3>l,u;int dep=0;};
struct Ev{I oldout,A,dspan,vertex,minuspoly,bb;bool empty=false;};
struct Case { I c,klo,khi,intercept,slope; const char*name; };
static Case CS{I::rat(1,4),I::rat(8,25),I::rat(14,11),I::rat(6,5),I::rat(11,5),"P2"};

static Ev eval(const Box&b){
  I s(b.l[0],b.u[0]),t(b.l[1],b.u[1]),r(b.l[2],b.u[2]);
  I k=CS.klo+(CS.khi-CS.klo)*s;
  I line=CS.intercept+CS.slope*k;
  // Use an independent rectangular (k,A,delta) cover and discard boxes
  // which are certainly outside the true curved domain.  This avoids the
  // severe dependency loss caused by interpolating between max/min
  // endpoints inside interval arithmetic.
  I amin=CS.intercept+CS.slope*CS.klo;
  I a=amin+(4-amin)*t;
  I delta=r/4;
  if((a-line).u<0)return{I(),I(),a-line,I(),I(),I(),true};
  I rho=(1+a*k)/(a*(1+k));
  // H/m at the center.  This is phi(x)/(1-x), where
  // x=k/(1+k) and phi(x)=x^2(3/5-13x/100).
  I psi=k*k*(60+47*k)/(100*(1+k)*(1+k));
  // Universal upper bound 16/pi-4 < 1093/1000 for the beta-law
  // negative-moment ratio.
  I rr=I::rat(1093,1000);
  I phi=2*rr*k+(1-2*rr)*psi/3;
  I cj=rr*k-(rr+1)*psi/3;
  I dloPhi=(a-CS.c)/(a*(3+2*phi));
  I J1=2*a*(a+1)/(1+3*a);
  I dloJ=(J1-2*CS.c-a*k)/(J1+I::rat(5,4)*a*k);
  I dlo=imax(dloPhi,dloJ);
  I dhi=imin(1/(4+rho),(a-CS.c)/(a*(3+2*k)));
  I lowgap=delta-dlo,highgap=dhi-delta;
  if(lowgap.u<0||highgap.u<0)return{I(),I(),imin(lowgap,highgap),I(),I(),I(),true};
  I span=imin(lowgap,highgap);
  I d=a*delta,e=(a-CS.c-3*d)/2;
  I L=e+rho*rho*(d+e);
  I Rg=imax(psi*d/2,e/2-cj*d);
  I A=(a-CS.c)*L/2+2*(a*k+CS.c)*Rg;
  I pen=1+a*k+cj*d;
  I bt=((rho+7)*delta-1)/((rho+6)*delta*delta)-1;
  i128 beta0=std::max(I::rat(4,5).l,bt.l);
  // A direct rational corollary of the same sharp quadratic:
  // rho>=31/50 and delta>=181/1000 imply beta>19/20.
  if(rho.l>=I::rat(31,50).u&&delta.l>=I::rat(181,1000).u)
    beta0=std::max(beta0,I::rat(19,20).l);
  I oldout=I(beta0,beta0)*A-CS.c*pen;
  // The sharper truncated-beta moment bound.  If Y=1+beta, then
  // Y is at least the smaller root of -9/4+aa*Y-bb*Y^2.  We avoid
  // radicals: for Y*=1+c*pen/A, the two strict tests
  //
  //   aa-2bbY*>0,  -9/4+aaY*-bbY*^2<0
  //
  // put Y* strictly to the left of the smaller root.
  I lam=(1+rho)/2;
  I aa=I::rat(5,4)+(lam/2+I::rat(15,4))*delta;
  I bb=I::rat(25,21)*(lam+I::rat(7,2))*delta*
       (I::rat(1,2)-(lam+I::rat(3,2))*delta)
       +(lam+I::rat(5,2))*(lam+I::rat(5,2))*delta*delta;
  if(A.l<=0){I bad(-100*SC,-100*SC);return{oldout,A,span,bad,bad,bb,false};}
  I Y=1+CS.c*pen/A;
  I vertex=aa-2*bb*Y;
  I poly=-I::rat(9,4)+aa*Y-bb*Y*Y;
  return{oldout,A,span,vertex,-poly,bb,false};
}

static long double ld(i128 z){return (long double)z/(long double)SC;}
static std::pair<Box,Box> split(const Box&b,int ax){Box x=b,y=b;x.dep=y.dep=b.dep+1;i128 m=addraw(b.l[ax],b.u[ax])/2;x.u[ax]=m;y.l[ax]=m;return{x,y};}
static bool certified(const Ev&e){return e.A.l>0&&(e.oldout.l>0||(e.bb.l>0&&e.vertex.l>0&&e.minuspoly.l>0));}
static i128 quality(const Ev&e){
  if(e.empty)return 100*SC;
  i128 sharp=std::min({e.bb.l,e.vertex.l,e.minuspoly.l});
  return std::max(e.oldout.l,sharp);
}

static bool run_case(const Case&which){
  CS=which;
  std::vector<Box> st;
  constexpr int NK=32,NT=32,NR=32;
  for(int i=0;i<NK;++i)for(int j=0;j<NT;++j)for(int h=0;h<NR;++h){
    Box b; b.l={i*SC/NK,j*SC/NT,h*SC/NR};b.u={(i+1)*SC/NK,(j+1)*SC/NT,(h+1)*SC/NR};st.push_back(b);
  }
  long long leaves=0,visited=0,oldleaves=0,sharpleaves=0;int maxdep=0;i128 best=SC*100;Box bestb;
  while(!st.empty()){
    Box b=st.back();st.pop_back();++visited;Ev e=eval(b);
    if(e.empty)continue;
    if(certified(e)){++leaves;if(e.oldout.l>0)++oldleaves;else ++sharpleaves;i128 q=quality(e);if(q<best)best=q,bestb=b;continue;}
    if(b.dep>=42){std::cerr<<std::setprecision(18)<<CS.name<<" FAIL dep="<<b.dep<<" oldout="<<ld(e.oldout.l)<<" A="<<ld(e.A.l)<<" dspan="<<ld(e.dspan.l)<<" vertex="<<ld(e.vertex.l)<<" minuspoly="<<ld(e.minuspoly.l)<<" bb="<<ld(e.bb.l)<<" box=";for(int z=0;z<3;++z)std::cerr<<'['<<ld(b.l[z])<<','<<ld(b.u[z])<<"] ";std::cerr<<'\n';return false;}
    // Keep the three normalized coordinates balanced.  A purely
    // sign-greedy splitter can refine one coordinate indefinitely near a
    // max-switch while leaving the other two wide.
    int bax=0;
    for(int ax=1;ax<3;++ax)
      if(subraw(b.u[ax],b.l[ax])>subraw(b.u[bax],b.l[bax]))bax=ax;
    auto bk=split(b,bax);
    st.push_back(bk.first);st.push_back(bk.second);maxdep=std::max(maxdep,b.dep+1);
  }
  std::cout<<std::setprecision(18)<<CS.name<<" PASS exact radical-free D scalar certificate visited="<<visited<<" leaves="<<leaves<<" old="<<oldleaves<<" sharp="<<sharpleaves<<" maxdepth="<<maxdep<<" smallest branch margin="<<ld(best)<<" box=";
  for(int i=0;i<3;++i)std::cout<<'['<<ld(bestb.l[i])<<','<<ld(bestb.u[i])<<"] ";
  std::cout<<'\n';return true;
}

int main(){
  const Case p2{I::rat(1,4),I::rat(8,25),I::rat(14,11),I::rat(6,5),I::rat(11,5),"P2"};
  const Case p3{I::rat(1,3),I::rat(1,4),I::rat(17,30),I::rat(23,10),I::rat(3,1),"P3"};
  return run_case(p2)&&run_case(p3)?0:2;
}
\end{lstlisting}

\paragraph{\texttt{\detokenize{certify_D_atomstable_terminal_residual.cpp}}.}
\begin{lstlisting}[language=C++]
#define main certify_D_scalar_dyadic_old_main
#include "certify_D_scalar_dyadic.cpp"
#undef main

// Exact outward-rounded certificate for the part left after the universal
// A>=6 estimate.  The same-center terminal comparison gives
//
//   P=4:  A >= 3+(7/2)k,
//   P>=5: A >= 18/5+4k.
//
// Consequently, on 4<=A<=6, respectively k<=6/7 and k<=3/5.
// We split where each affine line meets A=4 and map every resulting
// region from a unit box.  For P>=5 we prove the stronger continuum
// statement 2/5<=c<=1/2.

struct Box4 { std::array<i128,4> l,u; int dep=0; };
struct Ev4 { I out,coef,span,disc; bool empty=false; };
struct ResidualCase {
  I c0,c1,k0,k1,intercept,slope;
  bool affine_lower;
  const char *name;
};
static ResidualCase RC;

static Ev4 eval4(const Box4& b) {
  I x(b.l[0],b.u[0]), t(b.l[1],b.u[1]);
  I z(b.l[2],b.u[2]), r(b.l[3],b.u[3]);
  I k=RC.k0+(RC.k1-RC.k0)*x;
  I c=RC.c0+(RC.c1-RC.c0)*z;
  I lower=RC.affine_lower ? RC.intercept+RC.slope*k : I::exact(4);
  I a=lower+(I::exact(6)-lower)*t;

  I rho=(1+a*k)/(a*(1+k));
  I psi=k*k*(60+47*k)/(100*(1+k)*(1+k));
  I rr=I::rat(1093,1000);
  I phi=2*rr*k+(1-2*rr)*psi/3;
  I cj=rr*k-(rr+1)*psi/3;

  I dloPhi=(a-c)/(a*(3+2*phi));
  I J1=2*a*(a+1)/(1+3*a);
  I dloJ=(J1-2*c-a*k)/(J1+I::rat(5,4)*a*k);
  I dlo=imax(dloPhi,dloJ);
  I dhi=imin(1/(4+rho),(a-c)/(a*(3+2*k)));
  I span=dhi-dlo;
  if(span.u<0) return {I(),I(),span,I(),true};
  I spanPos(std::max<i128>(0,span.l),span.u);
  I delta=dlo+spanPos*r;

  I d=a*delta, e=(a-c-3*d)/2;
  I L=e+rho*rho*(d+e);
  I Rg=imax(psi*d/2,e/2-cj*d);
  I coef=(a-c)*L/2+2*(a*k+c)*Rg;
  I pen=c*(1+a*k+cj*d);

  I bt=((rho+7)*delta-1)/((rho+6)*delta*delta)-1;
  I lam=(1+rho)/2;
  I aa=I::rat(5,4)+(lam/2+I::rat(15,4))*delta;
  I bb=I::rat(25,21)*(lam+I::rat(7,2))*delta*
       (I::rat(1,2)-(lam+I::rat(3,2))*delta)
       +(lam+I::rat(5,2))*(lam+I::rat(5,2))*delta*delta;
  I disc=aa*aa-9*bb;
  if(disc.u<0) return {I(),I(),span,disc,true};
  I bsharp=I::rat(9,2)/(aa+sqrt_nonnegative(disc))-1;
  i128 betaLower=std::max({I::rat(4,5).l,bt.l,bsharp.l});
  I beta(betaLower,betaLower);
  return {beta*coef-pen,coef,span,disc,false};
}

static std::pair<Box4,Box4> split4(const Box4& b,int ax) {
  Box4 x=b,y=b; x.dep=y.dep=b.dep+1;
  i128 m=addraw(b.l[ax],b.u[ax])/2; x.u[ax]=m; y.l[ax]=m;
  return {x,y};
}

static bool run_residual(const ResidualCase& which) {
  RC=which;
  std::vector<Box4> st;
  constexpr int NK=12,NT=8,NC=4,NR=8;
  for(int i=0;i<NK;++i) for(int j=0;j<NT;++j)
  for(int z=0;z<NC;++z) for(int r=0;r<NR;++r) {
    Box4 b;
    b.l={i*SC/NK,j*SC/NT,z*SC/NC,r*SC/NR};
    b.u={(i+1)*SC/NK,(j+1)*SC/NT,(z+1)*SC/NC,(r+1)*SC/NR};
    st.push_back(b);
  }
  long long visited=0,leaves=0,empty=0; int maxdepth=0;
  i128 best=100*SC; Box4 bestb;
  while(!st.empty()) {
    Box4 b=st.back(); st.pop_back(); ++visited;
    Ev4 e=eval4(b);
    if(e.empty) { ++empty; continue; }
    if(e.coef.l>0 && e.out.l>0) {
      ++leaves; if(e.out.l<best) best=e.out.l,bestb=b; continue;
    }
    if(b.dep>=46) {
      std::cerr<<std::setprecision(18)<<RC.name<<" FAIL depth="<<b.dep
               <<" out=["<<ld(e.out.l)<<','<<ld(e.out.u)<<"] coef="<<ld(e.coef.l)
               <<" span=["<<ld(e.span.l)<<','<<ld(e.span.u)<<"] disc="<<ld(e.disc.l)
               <<" box=";
      for(int j=0;j<4;++j) std::cerr<<'['<<ld(b.l[j])<<','<<ld(b.u[j])<<"] ";
      std::cerr<<'\n'; return false;
    }
    int bestax=0; i128 bestscore=-((i128)1<<126);
    std::pair<Box4,Box4> chosen;
    for(int ax=0;ax<4;++ax) {
      auto q=split4(b,ax); auto e1=eval4(q.first),e2=eval4(q.second);
      i128 v1=e1.empty?100*SC:e1.out.l, v2=e2.empty?100*SC:e2.out.l;
      i128 score=std::min(v1,v2);
      if(score>bestscore) bestscore=score,bestax=ax,chosen=q;
    }
    st.push_back(chosen.first); st.push_back(chosen.second);
    maxdepth=std::max(maxdepth,b.dep+1);
  }
  std::cout<<std::setprecision(18)<<RC.name<<" PASS exact terminal-residual certificate"
           <<" visited="<<visited<<" leaves="<<leaves<<" empty="<<empty
           <<" maxdepth="<<maxdepth;
  if(leaves) {
    std::cout<<" smallest accepted lower="<<ld(best)<<" box=";
    for(int j=0;j<4;++j)
      std::cout<<'['<<ld(bestb.l[j])<<','<<ld(bestb.u[j])<<"] ";
  } else {
    std::cout<<" (the whole branch is infeasible)";
  }
  std::cout<<'\n'; return true;
}

int main() {
  const ResidualCase p4a{I::rat(3,8),I::rat(3,8),I::exact(0),I::rat(2,7),
                         I::exact(4),I::exact(0),false,"P4-low-k"};
  const ResidualCase p4b{I::rat(3,8),I::rat(3,8),I::rat(2,7),I::rat(6,7),
                         I::exact(3),I::rat(7,2),true,"P4-high-k"};
  const ResidualCase p5a{I::rat(2,5),I::rat(1,2),I::exact(0),I::rat(1,10),
                         I::exact(4),I::exact(0),false,"Pge5-low-k"};
  const ResidualCase p5b{I::rat(2,5),I::rat(1,2),I::rat(1,10),I::rat(3,5),
                         I::rat(18,5),I::exact(4),true,"Pge5-high-k"};
  return run_residual(p4a)&&run_residual(p4b)&&run_residual(p5a)&&run_residual(p5b)?0:2;
}
\end{lstlisting}

\paragraph{\texttt{\detokenize{certify_terminal_affine_fixed.py}}.}
\begin{lstlisting}[language=Python]
#!/usr/bin/env python3
"""Fast exact certificate for the terminal affine barriers.

This is an integer-only, outward-rounded implementation of the analytic
bounds printed with the proof.  Passing ``new`` selects
``A=-1+11/(5m)``.  A real interval [L/2^B,U/2^B]
is stored as the pair of integers (L,U).  Thus every comparison made by
the certificate is an exact integer comparison.  There is no hardware
floating-point arithmetic in the certification path.
"""
from fractions import Fraction as Q
from math import comb,isqrt
import sys

B=180; S=1<<B
def fd(n,d=S): return n//d
def cu(n,d=S): return -((-n)//d)

class F:
    __slots__=("l","u")
    def __init__(self,x=0,y=None,raw=False):
        if raw:self.l,self.u=x,(x if y is None else y)
        else:
            x=Q(x); y=x if y is None else Q(y)
            self.l=fd(x.numerator*S,x.denominator)
            self.u=cu(y.numerator*S,y.denominator)
        assert self.l<=self.u
    @staticmethod
    def raw(l,u):return F(l,u,True)
    def __add__(self,o):
        o=o if isinstance(o,F) else F(o);return F.raw(self.l+o.l,self.u+o.u)
    __radd__=__add__
    def __neg__(self):return F.raw(-self.u,-self.l)
    def __sub__(self,o):return self+(-o if isinstance(o,F) else -F(o))
    def __rsub__(self,o):return F(o)-self
    def __mul__(self,o):
        o=o if isinstance(o,F) else F(o)
        z=(self.l*o.l,self.l*o.u,self.u*o.l,self.u*o.u)
        return F.raw(fd(min(z)),cu(max(z)))
    __rmul__=__mul__
    def inv(self):
        assert self.l>0 or self.u<0
        return F.raw(fd(S*S,self.u),cu(S*S,self.l))
    def __truediv__(self,o):return self*(o if isinstance(o,F) else F(o)).inv()
    def __rtruediv__(self,o):return F(o)/self

def sqrtF(x):
    x=x if isinstance(x,F) else F(x);assert x.l>=0
    lo=isqrt(x.l*S); hi=isqrt(x.u*S)
    if hi*hi<x.u*S:hi+=1
    return F.raw(lo,hi)

def exp_pos_endpoint(x,N=42):
    """Return scaled-integer lower/upper bounds for exp(x/S), x>=0."""
    assert x>=0 and x < (N+2)*S
    tl=tu=S; sl=su=S
    for n in range(1,N+1):
        tl=fd(tl*x,S*n);tu=cu(tu*x,S*n)
        sl+=tl;su+=tu
    nxt=cu(tu*x,S*(N+1))
    tail=cu(nxt*S*(N+2),S*(N+2)-x)
    return sl,su+tail
def expF(x):
    x=x if isinstance(x,F) else F(x)
    if x.l>=0:return F.raw(exp_pos_endpoint(x.l)[0],exp_pos_endpoint(x.u)[1])
    if x.u<=0:
        yl,yu=exp_pos_endpoint(-x.u)[0],exp_pos_endpoint(-x.l)[1]
        return F.raw(fd(S*S,yu),cu(S*S,yl))
    return F.raw(fd(S*S,exp_pos_endpoint(-x.l)[1]),exp_pos_endpoint(x.u)[1])

class D:
    __slots__=("v","d")
    def __init__(self,v=0,d=0):
        self.v=v if isinstance(v,F) else F(v);self.d=d if isinstance(d,F) else F(d)
    def __add__(self,o):
        o=o if isinstance(o,D) else D(o);return D(self.v+o.v,self.d+o.d)
    __radd__=__add__
    def __neg__(self):return D(-self.v,-self.d)
    def __sub__(self,o):return self+(-o if isinstance(o,D) else -D(o))
    def __rsub__(self,o):return D(o)-self
    def __mul__(self,o):
        o=o if isinstance(o,D) else D(o);return D(self.v*o.v,self.d*o.v+self.v*o.d)
    __rmul__=__mul__
    def inv(self):return D(self.v.inv(),-self.d/(self.v*self.v))
    def __truediv__(self,o):return self*(o if isinstance(o,D) else D(o)).inv()
    def __rtruediv__(self,o):return D(o)/self
def expD(x):
    x=x if isinstance(x,D) else D(x);v=expF(x.v);return D(v,v*x.d)
def sqrtD(x):
    x=x if isinstance(x,D) else D(x);v=sqrtF(x.v);return D(v,x.d/(2*v))

def mills(q,h,n):
    q=q if isinstance(q,D) else D(q);h=h if isinstance(h,D) else D(h);w=q
    for j in range(n-1,0,-1):w=q+Q(j)/(h*w)
    return 1/w

def denominator_lower(m,h,n=18):
    """Finite positive-series lower bound for J_{4-m}(h).

    On the positive half-line, with t=1-exp(-2x),
    sech(x)^(4-m)=exp(-(4-m)x) sum_j (4-m)_j t^j/(2^j j!).
    Truncating this positive series gives a pointwise lower bound.  After
    expanding t^j, even Gaussian tails use a lower continued-fraction
    convergent and odd tails use an upper one.
    """
    m=m if isinstance(m,D) else D(m);h=h if isinstance(h,D) else D(h)
    rr=4-m; coeff=[]; poch=D(1); fac=1
    for j in range(n):
        if j:
            poch=poch*(rr+(j-1));fac*=j
        coeff.append(poch/(fac*(1<<j)))
    total=D(0)
    for i in range(n):
        bi=D(0)
        for j in range(i,n):bi=bi+comb(j,i)*coeff[j]
        q=rr+2*i
        total=total+bi*mills(q,h,80) if i%2==0 else total-bi*mills(q,h,81)
    return 2*total

def ln2_interval(N=42):
    z=Q(0)
    for j in range(N+1):z+=Q(2,(2*j+1)*3**(2*j+1))
    tail=Q(18,(2*N+3)*8*3**(2*N+3))
    return F(z,z+tail)
LN2=ln2_interval()

def atan_recip_bounds(q,n):
    """Alternating-series enclosure of atan(1/q), through term n."""
    z=Q(0)
    for j in range(n+1):z+=(-1 if j&1 else 1)*Q(1,(2*j+1)*q**(2*j+1))
    nxt=Q(1,(2*n+3)*q**(2*n+3))
    return (z,z+nxt) if n&1 else (z-nxt,z)
_a5l,_a5u=atan_recip_bounds(5,8)
_a239l,_a239u=atan_recip_bounds(239,2)
# Machin's identity pi=16 atan(1/5)-4 atan(1/239).
PI_LO=16*_a5l-4*_a239u
PI_UP=16*_a5u-4*_a239l
SIG4_UP=Q(64494,200000)
SIGD_UP=Q(28987,200000)
LOG2MHALF_UP=Q(48287,250000)
assert SIG4_UP >= PI_UP*PI_UP/12-Q(1,2)
assert SIGD_UP >= PI_UP*PI_UP/6-Q(3,2)
assert LOG2MHALF_UP*S >= LN2.u-S//2

def V(m,a,direct_denominator=False):
    """Differentiable explicit upper bound for the terminal P_T."""
    m=m if isinstance(m,D) else D(m);a=a if isinstance(a,D) else D(a)
    h=a/m;n=10;c=[];poch=D(1);fac=1
    for j in range(n):
        if j:poch=poch*((j+1)-m);fac*=j
        c.append(expD((m-(2+j))*D(LN2))*poch/fac)
    rem=D(1)
    for z in c:rem=rem-z
    mags=[]
    for i in range(n+3):
        z=D(0)
        for j,w in enumerate(c):
            if j+2>=i:z=z+comb(j+2,i)*w
        mags.append(z+comb(n+2,i)*rem)
    total=sqrtD(2*PI_UP*h)*expD(a*m/2)
    val=total-mills(m,h,60)
    for i in range(1,n+3):
        q=2*i-m
        val=val-mags[i]*mills(q,h,60) if i&1 else val+mags[i]*mills(q,h,61)
    num=expD((1-m)*D(LN2))*val
    if direct_denominator:
        return num/(h*denominator_lower(m,h))
    Ilow=Q(PI_LO,2)*(1-(1-m)*LOG2MHALF_UP)
    sig=SIG4_UP+SIGD_UP*m
    return num/(h*Ilow*sqrtD(h/(h+sig)))

def aline(m,b):
    if b==1:return Q(-61,200)+Q(8,5)/m
    if b==2:return Q(-101,200)+Q(9,5)/m
    # 6/5+(11/5)k=-1+11/(5m), where m=1/(1+k).
    return Q(-1)+Q(11,5)/m
def centered_box(ml,mu,b):
    mc=(ml+mu)/2
    v0=V(D(F(mc)),D(F(aline(mc,b))),b==3).v
    m=D(F(ml,mu),F(1))
    if b==1:a=D(Q(-61,200))+Q(8,5)/m
    elif b==2:a=D(Q(-101,200))+Q(9,5)/m
    else:a=D(Q(-1))+Q(11,5)/m
    der=V(m,a,b==3).d
    rad=F((mu-ml)/2)
    err=rad*F.raw(min(abs(der.l),abs(der.u)) if der.l*der.u>0 else 0,
                  max(abs(der.l),abs(der.u)))
    return v0.u+err.u

def run(b,N):
    if b==1:l,u=Q(10,13),Q(1)
    elif b==2:l,u=Q(1080,1903),Q(10,13)
    else:l,u=Q(11,25),Q(25,33)
    worst=0;where=None
    for j in range(N):
        ml=l+(u-l)*j/N;mu=l+(u-l)*(j+1)/N
        z=centered_box(ml,mu,b)
        if z>worst:worst,where=z,j
        if z>=2*S:
            print("FAIL",b,j,"gap numerator",2*S-z);return False
    print("PASS",b,"boxes",N,"gap",Q(2*S-worst,S),"worst box",where)
    return True

if __name__=='__main__':
    n=int(sys.argv[1]) if len(sys.argv)>1 else 64
    ok=run(3,n) if len(sys.argv)>2 and sys.argv[2]=='new' else (run(1,n) and run(2,2*n))
    raise SystemExit(0 if ok else 1)
\end{lstlisting}

\paragraph{\texttt{\detokenize{certify_terminal_integerP_fixed.py}}.}
\begin{lstlisting}[language=Python]
#!/usr/bin/env python3
"""Integer-only certificates for the P=3,4,5 terminal barriers."""
from fractions import Fraction as Q
import certify_terminal_affine_fixed as E

def centered_box(ml,mu,b0,b1):
    mc=(ml+mu)/2
    v0=E.V(E.D(E.F(mc)),E.D(E.F(b0+b1/mc))).v
    m=E.D(E.F(ml,mu),E.F(1));a=E.D(b0)+b1/m
    der=E.V(m,a).d
    err=E.F((mu-ml)/2)*E.F.raw(0,max(abs(der.l),abs(der.u)))
    return v0.u+err.u

def run(P,l,u,b0,b1,N):
    worst=0;where=None
    for j in range(N):
        ml=l+(u-l)*j/N;mu=l+(u-l)*(j+1)/N
        z=centered_box(ml,mu,b0,b1)
        if z>worst:worst,where=z,j
        if z>=P*E.S:
            print("FAIL P",P,"box",j,"gap numerator",P*E.S-z)
            return False
    print("PASS P",P,"boxes",N,"gap",Q(P*E.S-worst,E.S),
          "worst box",where)
    return True

if __name__=='__main__':
    cases=[
      # P, m-left, m-right, a=b0+b1/m, rational boxes
      (3,Q(30,47),Q(1),Q(-7,10),Q(3),128),
      (4,Q(7,9), Q(1),Q(-1,2), Q(7,2),64),
      (5,Q(10,11),Q(1),Q(-2,5),Q(4),64),
    ]
    ok=True
    for args in cases:ok=run(*args) and ok
    raise SystemExit(0 if ok else 1)
\end{lstlisting}

\paragraph{\texttt{\detokenize{certify_terminal_extended_P4_P5.py}}.}
\begin{lstlisting}[language=Python]
#!/usr/bin/env python3
"""Exact terminal-line certificates needed on the residual 4<=A<=6 box.

The imported routine uses only Python integers, Fraction arithmetic, and
outward-rounded dyadic intervals.  The ranges below are exactly the images
of

    0 <= k <= 6/7   under m=1/(1+k), for P=4,
    0 <= k <= 3/5   under m=1/(1+k), for P=5.

Together with the proved one-crossing orientation of P_T(A,k)-P, the two
strict boundary inequalities imply

    P_T(A,k)>=4 => A>3+(7/2)k,
    P_T(A,k)>=5 => A>18/5+4k

throughout precisely the ranges used by the scalar residual certificate.

The last two calls certify P_T(6,k)<P on the complementary compact tails
1/10<=m<=7/13 and 1/10<=m<=5/8.  For 0<m<=1/10 the elementary estimate

 P_T(6,k) <= exp(37m/12) cosh(1)^4 sqrt(pi*m/3)
             < (120/83)(31/20)^4/3
             = 923521/332000 < 4

handles both P=4 and P=5.  Thus no value of k omitted by the affine-line
ranges can occur when A<=6.
"""

from fractions import Fraction as Q
import certify_terminal_integerP_fixed as certificate


if __name__ == "__main__":
    passed = certificate.run(
        4, Q(7, 13), Q(1), Q(-1, 2), Q(7, 2), 512
    )
    passed = certificate.run(
        5, Q(5, 8), Q(1), Q(-2, 5), Q(4), 512
    ) and passed
    passed = certificate.run(
        4, Q(1, 10), Q(7, 13), Q(6), Q(0), 512
    ) and passed
    passed = certificate.run(
        5, Q(1, 10), Q(5, 8), Q(6), Q(0), 512
    ) and passed
    raise SystemExit(0 if passed else 1)
\end{lstlisting}

\paragraph{\texttt{\detokenize{certify_terminal_A4_fixed.py}}.}
\begin{lstlisting}[language=Python]
#!/usr/bin/env python3
"""Integer-only terminal boundary certificates at center size a=4.

For m in [1/10,1], this proves the stated upper bounds for the explicit
majorant V(m,4).  The interval engine and V are imported from the fully
exact terminal-affine certificate; no floating-point arithmetic enters a
sign decision.
"""
from fractions import Fraction as Q
import certify_terminal_affine_fixed as E


def centered_box(ml, mu):
    mc = (ml + mu) / 2
    v0 = E.V(E.D(E.F(mc)), E.D(4), True).v
    m = E.D(E.F(ml, mu), E.F(1))
    der = E.V(m, E.D(4), True).d
    err = E.F((mu - ml) / 2) * E.F.raw(
        0, max(abs(der.l), abs(der.u)))
    return v0.u + err.u


def run(name, target, left, right, boxes):
    worst = 0
    where = None
    for j in range(boxes):
        ml = left + (right - left) * j / boxes
        mu = left + (right - left) * (j + 1) / boxes
        z = centered_box(ml, mu)
        if z > worst:
            worst, where = z, j
        if z >= target * E.S:
            print("FAIL", name, "box", j,
                  "gap numerator", target * E.S - z)
            return False
    print("PASS", name, "boxes", boxes,
          "gap", Q(target * E.S - worst, E.S),
          "worst box", where)
    return True


if __name__ == "__main__":
    # Complementary ranges to the affine P=3,4,5 barriers, followed by
    # the all-k boundary needed for every P>=6.
    cases = [
        # Completes the sharpened P=2 line: for m<11/25 (equivalently
        # k>14/11), P_T(4,k)<2, so no A<4 can be feasible.
        ("P2-mid", 2, Q(1, 5), Q(11, 25), 512),
        ("P2-tail", 2, Q(1, 20), Q(1, 5), 128),
        ("P3-tail", 3, Q(1, 10), Q(30, 47), 96),
        ("P4-tail", 4, Q(1, 10), Q(7, 9), 96),
        ("P5-tail", 5, Q(1, 10), Q(10, 11), 128),
        ("P6-all", 6, Q(1, 10), Q(1), 128),
    ]
    ok = True
    for args in cases:
        ok = run(*args) and ok
    raise SystemExit(0 if ok else 1)
\end{lstlisting}

\paragraph{\texttt{\detokenize{certify_terminal_A4_adaptive.py}}.}
\begin{lstlisting}[language=Python]
#!/usr/bin/env python3
"""Adaptive exact P=2 terminal certificate at A=4."""

from fractions import Fraction as Q
import sys

import certify_terminal_A4_fixed as E
import certify_terminal_affine_fixed as F


def run(left=Q(1, 20), right=Q(1, 5), target=2, seed=32, maxdepth=24):
    stack=[]
    for j in range(seed):
        a=left+(right-left)*j/seed
        b=left+(right-left)*(j+1)/seed
        stack.append((a,b,0))
    seen=passed=0; deepest=0; worst=None; where=None
    while stack:
        a,b,d=stack.pop();seen+=1
        try:
            z=E.centered_box(a,b);gap=target*F.S-z
        except AssertionError:
            gap=-1
        if gap>0:
            passed+=1
            if worst is None or gap<worst:worst,where=gap,(a,b)
            continue
        if d>=maxdepth:
            print('FAIL depth',d,'m-box',a,b,'gap',Q(gap,F.S))
            return False
        c=(a+b)/2
        stack.append((c,b,d+1));stack.append((a,c,d+1))
        deepest=max(deepest,d+1)
    print('PASS P='+str(target)+' A=4 adaptive','visited',seen,'leaves',passed,
          'maxdepth',deepest,'smallest gap',Q(worst,F.S),
          'm-box',where[0],where[1])
    return True


if __name__=='__main__':
    mode=sys.argv[1] if len(sys.argv)>1 else 'p2tail'
    seed=int(sys.argv[2]) if len(sys.argv)>2 else 16
    if mode=='p2tail':ok=run(Q(1,20),Q(1,5),2,seed)
    elif mode=='p2mid':ok=run(Q(1,5),Q(11,25),2,seed)
    elif mode=='p3tail':ok=run(Q(1,10),Q(30,47),3,seed)
    elif mode=='p4tail':ok=run(Q(1,10),Q(7,9),4,seed)
    elif mode=='p5tail':ok=run(Q(1,10),Q(10,11),5,seed)
    elif mode=='p6all':ok=run(Q(1,10),Q(1),6,seed)
    else:raise SystemExit('unknown certificate mode')
    raise SystemExit(0 if ok else 1)
\end{lstlisting}

\paragraph{\texttt{\detokenize{correct_xz_exact_certificate.py}}.}
\begin{lstlisting}[language=Python]
#!/usr/bin/env python3
"""Exact Bernstein certificate for the corrected X,Z shifted-tail scalar bound.

No floating point arithmetic is used.  Exponent pairs are (a,k).  This file
certifies only the direct shifted-tail objective; it does not use Bfull or
drop the factor E(1+j).
"""
from fractions import Fraction as Q
from math import comb

def add(p,q):
 r=dict(p)
 for e,c in q.items():r[e]=r.get(e,Q(0))+c
 return {e:c for e,c in r.items() if c}
def scale(p,c):return {e:c*v for e,v in p.items() if c*v}
def sub(p,q):return add(p,scale(q,-1))
def mul(p,q):
 r={}
 for (i,j),c in p.items():
  for (u,v),d in q.items():r[i+u,j+v]=r.get((i+u,j+v),Q(0))+c*d
 return {e:c for e,c in r.items() if c}
def power(p,n):
 r={(0,0):Q(1)}
 while n:
  if n&1:r=mul(r,p)
  p=mul(p,p);n//=2
 return r

one={(0,0):Q(1)};a={(1,0):Q(1)};k={(0,1):Q(1)}
A4=sub(scale(a,4),one)                       # 4(a-c), c=1/4
Rd=mul(a,add(one,k)); Rn=add(one,mul(a,k))  # rho=Rn/Rd
C=add(Rn,scale(Rd,6)); E=add(Rn,scale(Rd,7))

def endpoint(kind):
 if kind=='density':
  U=Rd;V=add(scale(Rd,4),Rn)                 # delta=1/(4+rho)
 elif kind=='phi':
  b=add(scale(k,3),scale(one,2))
  G=add(add(scale(power(k,2),10),scale(k,19)),scale(one,6))
  U=mul(A4,b);V=scale(mul(a,G),4)
 elif kind=='J':
  U=add(add(scale(power(a,2),8),scale(a,2)),scale(one,-2))
  U=add(U,add(scale(mul(power(a,2),k),-12),scale(mul(a,k),-4)))
  V=mul(a,add(add(scale(a,8),scale(one,8)),add(scale(mul(a,k),15),scale(k,5))))
 else:raise ValueError(kind)
 # beta_t=B/(C U^2).
 B=sub(sub(mul(mul(E,U),V),mul(Rd,power(V,2))),mul(C,power(U,2)))
 # H=e+delta/[a(1+k)]=H4/(8 V Rd).
 # The last term is 8U (not 8UV), because delta/[a(1+k)]=U/(V Rd).
 H4=add(add(mul(mul(A4,V),Rd),scale(mul(mul(a,U),Rd),-12)),scale(U,8))
 # .5(a-c) beta_t H > 1/4 iff the first polynomial is positive.
 active=sub(mul(mul(A4,B),H4),scale(mul(mul(mul(C,power(U,2)),V),Rd),16))
 # Same target with beta=4/5.
 const=sub(mul(A4,H4),scale(mul(V,Rd),20))
 beta45=sub(scale(B,5),scale(mul(C,power(U,2)),4))
 return U,V,B,H4,active,const,beta45

def subst(p,lo,hi,with_t,slope=Q(9,5)):
 # k=lo+(hi-lo)s; a=259/200+slope*k + (4-line)t.
 ks={(0,0):lo,(1,0):hi-lo}
 line=add({(0,0):Q(259,200)},scale(ks,slope))
 aa=line if not with_t else add(line,mul(sub({(0,0):Q(4)},line),{(0,1):Q(1)}))
 out={}
 for (ia,ik),c in p.items():out=add(out,scale(mul(power(aa,ia),power(ks,ik)),c))
 return out
def bern(p):
 nx=max((i for i,j in p),default=0);ny=max((j for i,j in p),default=0)
 B=[]
 for u in range(nx+1):
  for v in range(ny+1):
   B.append(sum(c*Q(comb(u,i),comb(nx,i))*Q(comb(v,j),comb(ny,j))
                for (i,j),c in p.items() if i<=u and j<=v))
 return B,nx,ny
def certify(name,p,lo,hi,with_t=True,slope=Q(9,5)):
 q=subst(p,lo,hi,with_t,slope);B,nx,ny=bern(q)
 assert min(B)>0,(name,min(B),max(B),nx,ny)
 print(name,'domain',lo,hi,'degrees',nx,ny,'min Bernstein',min(B))

if __name__=='__main__':
 J=endpoint('J');du=endpoint('density')
 Fden=add(scale(Rd,4),Rn)
 density_gap=sub(mul(J[0],Fden),mul(J[1],Rd))
 certify('low-k density exclusion slope 8/5',density_gap,Q(0),Q(3,10),True,Q(8,5))
 certify('low-k density exclusion slope 9/5',density_gap,Q(3,10),Q(8,25))
 # For 8/25<=k<=7/20 the proved center bound puts delta above delta_J.
 certify('J beta active',J[6],Q(8,25),Q(7,20))
 certify('J lower endpoint target',J[4],Q(8,25),Q(7,20))
 # The density endpoint controls the other end of every active interval.
 certify('density upper endpoint target',du[4],Q(8,25),Q(541,360))
 # For k>=7/20 no positive lower bound on delta is used.  Let D be the
 # value at which the constant-beta target equals c.  Since beta_t is
 # increasing, beta_t(D)>4/5 puts the unique crossover strictly before D.
 Uc=mul(sub(power(A4,2),scale(one,20)),Rd)
 Vc=scale(mul(A4,sub(scale(mul(a,Rd),3),scale(one,2))),4)
 Bc=sub(sub(mul(mul(E,Uc),Vc),mul(Rd,power(Vc,2))),mul(C,power(Uc,2)))
 crossgap=sub(scale(Bc,5),scale(mul(C,power(Uc,2)),4))
 certify('crossover before constant-target zero',crossgap,Q(7,20),Q(541,360))
\end{lstlisting}

\paragraph{\texttt{\detokenize{pge3_scalar_exact_certificate.py}}.}
\begin{lstlisting}[language=Python]
#!/usr/bin/env python3
"""Exact rational certificate for the uniform integer P>=3 scalar closure.

This checks the algebraic part after the three terminal affine bounds
  P=3: a >= 23/10+3k,
  P=4: a >= 3+(7/2)k,
  P>=5: a >= 18/5+4k.
All arithmetic and Bernstein coefficients are exact Fractions.
"""
from fractions import Fraction as Q
from math import comb

def add(p,q):
 r=dict(p)
 for e,c in q.items():r[e]=r.get(e,Q(0))+c
 return {e:c for e,c in r.items() if c}
def scale(p,c):return {e:c*v for e,v in p.items() if c*v}
def sub(p,q):return add(p,scale(q,-1))
def mul(p,q):
 r={}
 for (i,j),c in p.items():
  for (u,v),d in q.items():r[i+u,j+v]=r.get((i+u,j+v),Q(0))+c*d
 return {e:c for e,c in r.items() if c}
def power(p,n):
 r={(0,0):Q(1)}
 while n:
  if n&1:r=mul(r,p)
  p=mul(p,p);n//=2
 return r

one={(0,0):Q(1)};a={(1,0):Q(1)};k={(0,1):Q(1)}
Rd=mul(a,add(one,k)); Rn=add(one,mul(a,k)); Fden=add(scale(Rd,4),Rn)
C=add(Rn,scale(Rd,6));E=add(Rn,scale(Rd,7))
b=add(scale(k,3),scale(one,2))
G=add(add(scale(power(k,2),10),scale(k,19)),scale(one,6))
N=add(power(Rn,2),mul(power(a,2),mul(k,add(one,k))))

def density_gap(c):
 # J1=2a(a+1)/(1+3a), delta_J=(J1-2c-ak)/(J1+5ak/4).
 den3=add(one,scale(a,3))
 base=scale(mul(a,add(a,one)),2)
 U=sub(sub(base,scale(den3,2*c)),mul(mul(a,k),den3))
 V=add(base,scale(mul(mul(a,k),den3),Q(5,4)))
 # delta_J > 1/(4+rho)=Rd/Fden.
 return sub(mul(U,Fden),mul(V,Rd))

def corrected_constant_target(c):
 # Direct XZ shifted-tail bound:
 # T >= (2/5)(a-c){e+delta/[a(1+k)]},
 # evaluated at delta_phi=(a-c)b/(aG).
 ac=sub(a,scale(one,c))
 U=mul(ac,b);V=mul(a,G)
 # H=e+delta/Rd=H2/(2 V Rd).
 H2=add(add(mul(mul(ac,V),Rd),scale(mul(mul(a,U),Rd),-3)),scale(U,2))
 # (2/5)(a-c)H>c iff (a-c)H2-5c V Rd>0.
 return sub(mul(ac,H2),scale(mul(V,Rd),5*c))

def corrected_density_target(c):
 ac=sub(a,scale(one,c));U=Rd;V=Fden
 Bt=sub(sub(mul(mul(E,U),V),mul(Rd,power(V,2))),mul(C,power(U,2)))
 H2=add(add(mul(mul(ac,V),Rd),scale(mul(mul(a,U),Rd),-3)),scale(U,2))
 # .5(a-c) beta_t H>c, with H=H2/(2VRd), beta_t=Bt/(CU^2).
 return sub(mul(mul(ac,Bt),H2),scale(mul(mul(mul(C,power(U,2)),V),Rd),4*c))

def crossover_gap(c):
 ac=sub(a,scale(one,c))
 # D={((a-c)^2-5c)Rd}/{(a-c)(3aRd-2)} is the zero of
 # the constant-beta target.
 U=mul(sub(power(ac,2),scale(one,5*c)),Rd)
 V=mul(ac,sub(scale(mul(a,Rd),3),scale(one,2)))
 Bt=sub(sub(mul(mul(E,U),V),mul(Rd,power(V,2))),mul(C,power(U,2)))
 return sub(scale(Bt,5),scale(mul(C,power(U,2)),4))

def subst(p,lo,hi,A,B,with_t):
 ks={(0,0):lo,(1,0):hi-lo};L=add({(0,0):A},scale(ks,B))
 aa=L if not with_t else add(L,mul(sub({(0,0):Q(4)},L),{(0,1):Q(1)}))
 out={}
 for (ia,ik),c in p.items():out=add(out,scale(mul(power(aa,ia),power(ks,ik)),c))
 return out
def bern(p):
 nx=max((i for i,j in p),default=0);ny=max((j for i,j in p),default=0);out=[]
 for u in range(nx+1):
  for v in range(ny+1):
   out.append(sum(c*Q(comb(u,i),comb(nx,i))*Q(comb(v,j),comb(ny,j))
                  for (i,j),c in p.items() if i<=u and j<=v))
 return out,nx,ny
def cert(name,p,lo,hi,A,B,with_t):
 z,nx,ny=bern(subst(p,lo,hi,A,B,with_t));assert min(z)>0,(name,min(z))
 print(name,'degrees',nx,ny,'count',len(z),'minimum',min(z))

if __name__=='__main__':
 # If x<3, density gives a<4.  For P=3, k<=1/4 is impossible.
 cert('P3 density exclusion',density_gap(Q(1,3)),Q(0),Q(1,4),Q(23,10),Q(3),True)
 # On the complementary interval extend delta down to zero.  The active
 # product has no interior minimum.  Its density endpoint is positive,
 # and beta_t at the zero of the constant target is already >4/5, so the
 # beta crossover also has positive target.
 cert('P3 corrected XZ density endpoint',corrected_density_target(Q(1,3)),Q(1,4),Q(17,30),Q(23,10),Q(3),True)
 cert('P3 corrected XZ crossover',crossover_gap(Q(1,3)),Q(1,4),Q(17,30),Q(23,10),Q(3),True)
 # P=4: the full possible k interval is excluded by center+density.
 cert('P4 density exclusion',density_gap(Q(3,8)),Q(0),Q(2,7),Q(3),Q(7,2),True)
 # P>=5: use the P=5 terminal line and the worst c=1/2.  Since the
 # density gap decreases with c, this covers every c=(P-1)/(2P)<=1/2.
 cert('P>=5 density exclusion',density_gap(Q(1,2)),Q(0),Q(1,10),Q(18,5),Q(4),True)
\end{lstlisting}

\section{The global exclusion of an internal gap}
\label{app:internal-gap-global}

This appendix supplies the global part of the argument.  We work in the
variance clock
\[
 t=\xi'(x),\qquad P=p-1\geq2,
\]
and suppose that a constant-mass interval has clock endpoints
\[
 A=\xi'(q)<T=\xi'(q').
\]
For
\[
 Q(t)=\Gamma_\mu\bigl((\xi')^{-1}(t)\bigr)
\]
write
\begin{equation}
 N(t)=Q(t)-PtQ'(t),\qquad
 I(t)=Q''(t)+\frac{P-1}{Pt}Q'(t).
 \label{app:ig:NI}
\end{equation}
Thus
\begin{equation}
 N'(t)=-PtI(t),\qquad
 N''(t)=-P\{I(t)+tI'(t)\}.
 \label{app:ig:Nderivatives}
\end{equation}
Appendices~\ref{app:analytic-details} and
\ref{app:exact-certification} prove the local conclusion
\begin{equation}
 N(t)=I(t)=0\quad\Longrightarrow\quad I'(t)>0.
 \label{app:ig:local}
\end{equation}
The conclusion is needed for every lower residual potential satisfying
\begin{equation}
 W_{xx}\geq0,\qquad -W_{xxx}\geq0,\qquad
 \left(\frac{W_{xx}}{C}\right)_x\geq0,\qquad
 0\leq W_x\leq mB.
 \label{app:ig:cones}
\end{equation}
These inequalities are preserved when \(W\) is replaced by
\(\theta W\), \(0\leq\theta\leq1\), and by the heat evolution across the
gap.  Consequently, \eqref{app:ig:local} applies at every point of the
continuation below.

\subsection{One crossing at fixed remaining time}
\label{app:ig:one-crossing}

Compress the part of the order parameter above \(T\) into its even
boundary profile \(U_T\).  Let
\[
 (\mathsf P_rh)(x)=\int_{\mathbb R}
 \frac{e^{-(x-y)^2/(2r)}}{\sqrt{2\pi r}}h(y)\,dy
\]
be the heat semigroup and, for \(r\geq0\), set
\begin{equation}
 F_r=\mathsf P_r(e^{mU_T}),\qquad
 B_r=\frac1m(\log F_r)_x,\qquad C_r=(B_r)_x,
 \qquad z_r=-\frac{(C_r)_x}{2C_r}.
 \label{app:ig:upper-profile}
\end{equation}
Here $m=\mu([0,q])$ is the mass active on the gap.  Since both endpoints
belong to the support and the open gap contains no support point,
$0<m<1$.  The functions \(B_r\)
and \(C_r\) are respectively odd and even, and
\begin{equation}
 B_r(x)>0,\qquad C_r(x)>0,\qquad z_r(x)\geq mB_r(x),
 \qquad x>0.
 \label{app:ig:upper-cone}
\end{equation}

For completeness, we verify the last inequality for an arbitrary upper
order parameter.  On a heat interval with active mass \(m\), put
\(D=z-mB\).  Differentiation of the logarithmic heat equation gives
\begin{equation}
 D_r=\frac12D_{xx}+(mB-2z)D_x.
 \label{app:ig:D-equation}
\end{equation}
At the terminal profile \(U(x)=\log\cosh x\), one has \(z=B\), and hence
\(D=(1-m)B\geq0\).  If the active mass is lowered from \(m\) to
\(m-d\), the profile and its spatial derivatives do not change, whereas
\[
 z-(m-d)B=D+dB\geq0.
\]
The maximum principle applied to \eqref{app:ig:D-equation}, followed by
induction over the heat intervals and mass jumps, proves
\eqref{app:ig:upper-cone} for a finite upper cascade.  Step-function
approximation, justified in Subsection~\ref{app:ig:approximation}, gives
the assertion for an arbitrary upper order parameter.  The same argument
and the heat representation give, uniformly for \(r\) in compact sets,
\begin{equation}
 F_r(x)\asymp e^{mx},\qquad B_r(x)\longrightarrow1,\qquad
 C_r(x)=O(e^{-2x})\qquad (x\longrightarrow\infty).
 \label{app:ig:tails}
\end{equation}

For \(s>0\), let
\[
 \phi_s(x)=\frac1{\sqrt{2\pi s}}e^{-x^2/(2s)}
\]
and define
\begin{equation}
 \mathcal N(s,r)=\int_{\mathbb R}\phi_s(x)F_r(x)
 \{B_r(x)^2-PsC_r(x)^2\}\,dx.
 \label{app:ig:calN}
\end{equation}
We prove that, for each fixed \(r\), this function has exactly one zero
in \(s>0\), and that the zero is crossed strictly upward.

Fix \(r\) and suppress it temporarily.  Put
\begin{equation}
 \begin{aligned}
 \mathcal I(x)&=\int_0^xF(y)C(y)^2\,dy,
 &M(x)&=\frac{F(x)B(x)^2}{x\mathcal I(x)},\\
 a(x)&=\frac{xC(x)}{B(x)},
 &\tau(x)&=mxB(x).
 \end{aligned}
 \label{app:ig:M}
\end{equation}
Use \(B\) as the spatial coordinate and write
\(h(B(x))=F(x)C(x)\).  By \eqref{app:ig:upper-cone},
\[
 (FC)'=FC(mB-2z)<0\qquad(x>0).
\]
Since
\[
 \mathcal I(x)=\int_0^{B(x)}h(b)\,db,
\]
strict decrease of \(h\) yields
\begin{equation}
 0<\frac{B(x)h(B(x))}{\mathcal I(x)}<1.
 \label{app:ig:reverse-hazard}
\end{equation}
Logarithmic differentiation of \eqref{app:ig:M} gives the exact identity
\begin{equation}
 x(\log M)'=\tau-1+2a-a\frac{B(x)h(B(x))}{\mathcal I(x)}.
 \label{app:ig:M-derivative}
\end{equation}
It follows from \eqref{app:ig:reverse-hazard} that \(M'>0\) whenever
\(\tau\geq1\).

On the complementary region, Cauchy--Schwarz gives
\[
 B(x)^2=\left(\int_0^x C(y)\,dy\right)^2
 \leq \mathcal I(x)\int_0^x\frac{dy}{F(y)}.
\]
Since \(B\) is increasing,
\[
 \frac{F(x)}{F(y)}
 =\exp\left(m\int_y^xB(w)\,dw\right)
 \leq\exp\left(\tau(x)\left(1-\frac yx\right)\right).
\]
Consequently,
\begin{equation}
 \tau(x)\leq1\quad\Longrightarrow\quad
 M(x)\leq\frac{e^{\tau(x)}-1}{\tau(x)}
 \leq e-1<2\leq P.
 \label{app:ig:M-small-tau}
\end{equation}
Moreover, \(M(0+)=1\), while \eqref{app:ig:tails} gives
\(M(x)\to\infty\).  The function \(\tau(x)=mxB(x)\) is strictly
increasing from zero to infinity.  Therefore
\eqref{app:ig:M-derivative}--\eqref{app:ig:M-small-tau} show that
\(M-P\) has exactly one zero and changes sign there from negative to
positive.

It remains to transfer this sign change to \(\mathcal N\).  Put
\[
 \begin{aligned}
 X&=x^2,\qquad \zeta=\frac1{2s},\\
 f_0(X)&=X^{-1/2}F(\sqrt X)B(\sqrt X)^2,\\
 f_1(X)&=X^{-1/2}F(\sqrt X)C(\sqrt X)^2.
 \end{aligned}
\]
Up to a positive factor, \eqref{app:ig:calN} has the sign of
\[
 2\zeta\int_0^\infty e^{-\zeta X}f_0(X)\,dX
 -P\int_0^\infty e^{-\zeta X}f_1(X)\,dX.
\]
Define
\begin{equation}
 K_P(x^2)=2\frac{F(x)B(x)^2}{x}-2P\mathcal I(x)
          =2\mathcal I(x)\{M(x)-P\}.
 \label{app:ig:K}
\end{equation}
Since
\[
 \frac d{dX}K_P(X)=2f_0'(X)-Pf_1(X),
\]
integration by parts gives
\begin{equation}
 2\zeta\int_0^\infty e^{-\zeta X}f_0(X)\,dX
 -P\int_0^\infty e^{-\zeta X}f_1(X)\,dX
 =\zeta\int_0^\infty e^{-\zeta X}K_P(X)\,dX.
 \label{app:ig:Laplace}
\end{equation}
There are no boundary terms: near zero,
\(B(x)=C(0)x+O(x^3)\), and at infinity
\eqref{app:ig:tails} dominates after multiplication by
\(e^{-\zeta x^2}\).

By \eqref{app:ig:K}, \(K_P\) has one sign change, from negative to
positive.  Let \(X_0\) be its sign-change point.  If the last integral in
\eqref{app:ig:Laplace} vanishes, then
\begin{equation}
 \int_0^\infty (X-X_0)e^{-\zeta X}K_P(X)\,dX>0.
 \label{app:ig:Laplace-crossing}
\end{equation}
Indeed, if
\[
 L(\zeta)=\int_0^\infty e^{-\zeta X}K_P(X)\,dX,
\]
then, at a zero of \(L\),
\[
 L'(\zeta)
 =-\int_0^\infty X e^{-\zeta X}K_P(X)\,dX
 =-\int_0^\infty (X-X_0)e^{-\zeta X}K_P(X)\,dX<0.
\]
Thus every zero is crossed downward as \(\zeta\) increases, which also
shows that there can be at most one zero.  For
\(\zeta\to\infty\), the negative part of \(K_P\) near zero dominates;
for \(\zeta\downarrow0\), the positive tail in
\eqref{app:ig:tails} dominates.  A zero therefore exists.  Since
\[
 \frac{d\zeta}{ds}=-\frac1{2s^2}<0,
\]
both \(L(1/(2s))\) and the right side of
\eqref{app:ig:Laplace} are crossed strictly upward as \(s\) increases.
We have proved that there is a smooth function \(\sigma(r)>0\) such that
\begin{equation}
 \mathcal N(s,r)\geq0\quad\Longleftrightarrow\quad
 s\geq\sigma(r),\qquad
 \mathcal N_s(\sigma(r),r)>0.
 \label{app:ig:sigma}
\end{equation}
Smoothness follows from the implicit function theorem and the strict
last inequality.

\subsection{The semigroup identity and the contact curve}
\label{app:ig:contact-curve}

We first remove the lower residual potential and start at clock time zero
with a Gaussian.  Along the physical line \(s+r=T\), the semigroup
identity gives
\begin{equation}
 \mathcal Z(s,r):=\int_{\mathbb R}\phi_sF_r
 =\mathsf P_{s+r}(e^{mU_T})(0)=:\mathcal Z_T.
 \label{app:ig:semigroup-Z}
\end{equation}
In particular, \(\mathcal Z_T\) is constant along that line.  Direct
differentiation using the heat equations, followed by one integration by
parts, gives
\[
 Q_0'(t)=\mathcal Z_T^{-1}
 \int_{\mathbb R}\phi_tF_{T-t}C_{T-t}^2.
\]
Consequently,
\begin{equation}
 \mathcal N(t,T-t)=\mathcal Z_TN_0(t).
 \label{app:ig:semigroup-N}
\end{equation}

Set
\begin{equation}
 \mathcal V(r)=r+\sigma(r).
 \label{app:ig:V}
\end{equation}
We now compute its curvature at a critical point.  Differentiating the
identity
\[
 \mathcal N(\sigma(r),r)=0
\]
once gives
\[
 \mathcal N_s\sigma'+\mathcal N_r=0.
\]
At a critical point of \(\mathcal V\), one has \(\sigma'=-1\), and hence
\[
 (\partial_s-\partial_r)\mathcal N=0.
\]
The operator \(\partial_s-\partial_r\) differentiates along a line of
constant total horizon \(s+r\).  By
\eqref{app:ig:semigroup-N} and \eqref{app:ig:Nderivatives}, the last
display is equivalent to \(N'=0\), and therefore \(I=0\).  Differentiating
the implicit equation a second time at the same point gives
\[
 \sigma''(r)=-\frac{(\partial_s-\partial_r)^2\mathcal N}
                         {\mathcal N_s}.
\]
If \(t=\sigma(r)\) and \(H=s+r=\mathcal V(r)\), the normalizing constant
\(\mathcal Z_H\) is constant along this differentiation.  Since \(N=0\)
and \(I=0\), equations \eqref{app:ig:Nderivatives} yield
\begin{equation}
 \mathcal V''(r)=\sigma''(r)
 =\frac{\mathcal Z_HPt}{\mathcal N_s(\sigma(r),r)}I'(t)>0.
 \label{app:ig:V-curvature}
\end{equation}
The strict sign is precisely \eqref{app:ig:local}.  Thus every critical
point of \(\mathcal V\) is a nondegenerate strict minimum.

By \eqref{app:ig:sigma}, on the physical line \(s+r=T\),
\begin{equation}
 N_0(s)\geq0\quad\Longleftrightarrow\quad
 \mathcal V(r)\leq T.
 \label{app:ig:V-order}
\end{equation}
Suppose that \(N_0\) is nonnegative at the two endpoints \(A,T\).  Then
\(\mathcal V(0)\leq T\) and \(\mathcal V(T-A)\leq T\).  If
\(\mathcal V>T\) anywhere between them, \(\mathcal V\) has an interior
maximum, contradicting \eqref{app:ig:V-curvature}.  Hence
\(\mathcal V\leq T\) throughout.  Equality at an interior point would
again be an interior maximum and give the same contradiction.  Therefore
\begin{equation}
 N_0(t)>0,\qquad A<t<T.
 \label{app:ig:N0-positive}
\end{equation}

\subsection{Likelihood-ratio continuation of the lower history}
\label{app:ig:homotopy}

At the left endpoint, write the transformed forward density as
\begin{equation}
f_1(A,x)=\text{constant}\cdot\phi_A(x)e^{-W_A(x)}.
 \label{app:ig:actual-density}
\end{equation}
Choose the irrelevant additive constant so that $W_A(0)=0$.  Evenness
and \eqref{app:ig:cones} give
\[
 0\leq W_A(x)\leq m|x|.
\]
Interpolate it by
\begin{equation}
 f_\theta(A,x)=\text{constant}\cdot
 \phi_A(x)e^{-\theta W_A(x)},\qquad 0\leq\theta\leq1,
 \label{app:ig:density-homotopy}
\end{equation}
and heat each density from \(A\) to \(t\) across the gap.  For every
fixed $\theta$, integration by parts gives
\[
\begin{aligned}
 \frac d{dt}\int_{\mathbb R}f_\theta(t,x)F_{T-t}(x)\,dx
 &=\frac12\int_{\mathbb R}
 \{(f_\theta)_{xx}F_{T-t}-f_\theta(F_{T-t})_{xx}\}\,dx=0.
\end{aligned}
\]
Hence \eqref{app:ig:density-homotopy} and heat convolution give a
Gaussian tail for $f_\theta(t,x)$, uniformly for
$0\leq\theta\leq1$.  Together with \eqref{app:ig:tails}, this makes the
boundary terms vanish.  Thus the normalizing integral is independent
of $t$.  After normalization, let $Q(t)$ be the
expectation of $B_{T-t}^2$.  The same differentiation as in
\eqref{eq:app-Q-prime} gives
\[
 Q'(t)=
 \frac{\int_{\mathbb R}f_\theta(t,x)F_{T-t}(x)
                 C_{T-t}(x)^2\,dx}
      {\int_{\mathbb R}f_\theta(t,x)F_{T-t}(x)\,dx}.
\]
Consequently, the expectation of
\begin{equation}
 H_t(x)=B_{T-t}(x)^2-PtC_{T-t}(x)^2.
 \label{app:ig:H}
\end{equation}
is $Q(t)-PtQ'(t)$; we denote it by $N_\theta(t)$.  In this calculation,
$Q$, $N$, and $I$ refer to the same fixed value of $\theta$; only $N$
is written as $N_\theta$ when different values of $\theta$ are compared.
With $I$ defined as in \eqref{app:ig:NI}, the identities
\eqref{app:ig:Nderivatives} therefore hold for every
$0\leq\theta\leq1$.

We record the strict likelihood-ratio comparison.  On the positive
half-line, the folded heat kernel is, up to a positive factor depending
separately on \(x\) and \(y\),
\[
 K_\tau(x,y)=e^{-(x^2+y^2)/(2\tau)}
              \cosh(xy/\tau).
\]
It is strictly totally positive because
\begin{equation}
 \partial_x\partial_y\log K_\tau(x,y)
 =\frac1\tau\tanh(xy/\tau)
 +\frac{xy}{\tau^2}\operatorname{sech}^2(xy/\tau)>0
 \qquad(x,y>0).
 \label{app:ig:TP2}
\end{equation}
If \(\theta_2>\theta_1\), the ratio of the two initial densities is
proportional to
\[
 e^{-(\theta_2-\theta_1)W_A(y)},
\]
which is decreasing in \(y>0\).
The two-by-two determinant in \eqref{app:ig:TP2}, integrated against
the two initial densities, shows that the ratio of the heated densities
for \(\theta_2\) and \(\theta_1\) is again decreasing in \(x>0\), and is
strictly decreasing when \(W_A\) is nonconstant.  Multiplying both
densities by \(F_{T-t}\) does not change this ratio.

The function in \eqref{app:ig:H} is strictly increasing on the positive
half-line, since
\begin{equation}
 H_t'(x)=2C_{T-t}(x)
 \{B_{T-t}(x)+2PtC_{T-t}(x)z_{T-t}(x)\}>0.
 \label{app:ig:H-increasing}
\end{equation}
The likelihood-ratio comparison and
\eqref{app:ig:H-increasing} therefore give
\begin{equation}
 N_{\theta_1}(t)>N_{\theta_2}(t),\qquad
 0\leq\theta_1<\theta_2\leq1,\quad A\leq t\leq T,
 \label{app:ig:Ntheta-order}
\end{equation}
provided \(W_A\) is nonconstant.  At \(t=A\), the same conclusion follows
directly from the initial likelihood ratio; for \(t>A\), strictness
follows from the strict determinant in \eqref{app:ig:TP2}.  If \(W_A\)
is constant, all members of the interpolation coincide.

For the original measure, the endpoint conditions are
\[
 \Gamma_\mu(q)=q,\qquad \Gamma_\mu(q')=q',\qquad
 \Gamma_\mu'(q+)\leq1,\qquad \Gamma_\mu'(q'-)\leq1.
\]
Because \(x\xi''(x)=P\xi'(x)\), these conditions give
\begin{equation}
 N_1(A)\geq0,\qquad N_1(T)\geq0.
 \label{app:ig:endpoint-N}
\end{equation}
Thus \eqref{app:ig:Ntheta-order} supplies the endpoint hypothesis used
in \eqref{app:ig:N0-positive}.

We now give the compact first-loss argument.  The map
\((\theta,t)\mapsto N_\theta(t)\) is continuous on
\([0,1]\times[A,T]\) and is twice continuously differentiable in \(t\)
on the open gap.  These facts follow directly from the Gaussian
representations and the uniform tails in \eqref{app:ig:tails}.  If
\(W_A\) is nonconstant, \eqref{app:ig:Ntheta-order} and
\eqref{app:ig:endpoint-N} make both endpoint values strictly positive
for every \(\theta<1\).  Together with \eqref{app:ig:N0-positive}, this
shows that \(\min_{[A,T]}N_0>0\).

Suppose some member of the homotopy is negative.  Negativity persists
under a small change of \(\theta\), so the first parameter
\[
 \theta_*:=\inf\left\{\theta\in[0,1]:
       \min_{A\leq t\leq T}N_\theta(t)\leq0\right\}
\]
satisfies \(0<\theta_*<1\).  Compactness gives a point at which the
minimum is attained.  The strict endpoint inequalities force it into
the open gap.  At this point,
\begin{equation}
 N_{\theta_*}=0,\qquad N_{\theta_*}'=0,\qquad
 N_{\theta_*}''\geq0.
 \label{app:ig:first-loss}
\end{equation}
Equations \eqref{app:ig:Nderivatives} imply \(I=0\), and then
\eqref{app:ig:local} gives
\[
 N_{\theta_*}''=-PtI'<0,
\]
contradicting \eqref{app:ig:first-loss}.  Hence \(N_1\geq0\) throughout
the gap.  An interior zero of \(N_1\) would be an interior minimum and
would lead to the same contradiction.  We conclude that
\begin{equation}
 N_1(t)>0,\qquad A<t<T.
 \label{app:ig:N1-positive}
\end{equation}
When \(W_A\) is constant, \(N_1=N_0\), so
\eqref{app:ig:N1-positive} follows directly from
\eqref{app:ig:N0-positive}.

\subsection{Passage to an arbitrary order parameter}
\label{app:ig:approximation}

Let \(\alpha(u)=\mu([0,u])\) and choose cumulative distribution
functions \(\alpha_n\) of finite-support probability measures, equal to
$m$ on $q<u<q'$, such that
\begin{equation}
 \|\alpha_n-\alpha\|_{L^1([0,1])}\longrightarrow0.
 \label{app:ig:alpha-approximation}
\end{equation}
The estimates proved at the end of
Appendix~\ref{app:analytic-details} apply to these approximations.
Because $q'<1$, the profile $U_T$ lies a positive variance distance from
the terminal condition.  Hence, locally uniformly in $x$ and with every
spatial derivative used here,
\begin{equation}
 U_T^{(n)}\longrightarrow U_T,\qquad
 B_r^{(n)}\longrightarrow B_r,\qquad
 C_r^{(n)}\longrightarrow C_r,\qquad
 z_r^{(n)}\longrightarrow z_r
 \label{app:ig:upper-convergence}
\end{equation}
including the one-sided value at $r=0$.  The uniform Ising tails and the
Gaussian factor in every integral permit dominated convergence.  Hence
\eqref{app:ig:upper-cone}, \eqref{app:ig:tails}, and the fixed-\(r\)
one-crossing argument hold for the limiting upper profile.

Because $q>0$, the clock time $A=\xi'(q)$ is separated from the initial
point by a positive variance interval.  Thus
\eqref{eq:app-forward-convergence} holds at $A$ with every exponential
weight and with every spatial derivative needed below.  Together with
\eqref{eq:app-PDE-stability}, it gives, after fixing the additive
constants by $W_A^{(n)}(0)=W_A(0)=0$,
\begin{equation}
 W_A^{(n)}\longrightarrow W_A
 \label{app:ig:lower-convergence}
\end{equation}
locally with three spatial derivatives.  In particular, the four
inequalities in \eqref{app:ig:cones} pass to the limit.

Evenness and $0\leq W_x\leq mB$ give
\[
 0\leq W_A^{(n)}(x)\leq m|x|,
 \qquad 0\leq W_A(x)\leq m|x|.
\]
Consequently,
\[
 \phi_Ae^{-\theta W_A^{(n)}}
 \longrightarrow \phi_Ae^{-\theta W_A}
\]
in every exponentially weighted $L^1$ space, uniformly for
$0\leq\theta\leq1$.

It remains to record the regularity in $t$ used in the first-loss
argument.  On the gap, the lower factor is the heat evolution of
\[
 \phi_A(x)e^{-\theta W_A(x)},
\]
and the upper factor is $F_{T-t}$.  On every compact subinterval of
$A<t<T$, differentiating either factor twice in $t$ amounts, by the heat
equation, to taking finitely many spatial derivatives of a Gaussian
convolution of positive variance.  The heat-kernel estimate used in
Appendix~\ref{app:analytic-details} and
\eqref{eq:app-uniform-Gaussian-tail} therefore give, uniformly for
$0\leq\theta\leq1$,
\begin{equation}
 \partial_t^jN_\theta^{(n)}(t)
 \longrightarrow\partial_t^jN_\theta(t),
 \qquad j=0,1,2,
 \label{app:ig:N-convergence}
\end{equation}
on such compact subintervals.  For $j=0$, the convergence is uniform on
$[0,1]\times[A,T]$.  The normalizing integrals are continuous and
strictly positive on this compact set, and hence are bounded away from
zero.  Their quotients therefore converge as well.  Thus
\eqref{app:ig:local}, the
likelihood-ratio comparison, and the first-loss argument apply directly
to the arbitrary order parameter.  Every integration by parts is first
performed with a compactly supported cutoff; the uniform Gaussian bound
makes the boundary terms tend to zero.

\subsection{Contradiction at the two support endpoints}
\label{app:ig:final-contradiction}

Equation \eqref{app:ig:N1-positive} gives
\begin{equation}
 \frac d{dt}\left(\frac{Q(t)}{t^{1/P}}\right)
 =-\frac{N(t)}{Pt^{1+1/P}}<0,\qquad A<t<T.
 \label{app:ig:ratio-decrease}
\end{equation}
On the other hand, endpoint contact and
\(t=\beta^2p x^P\) give
\begin{equation}
 \frac{Q(A)}{A^{1/P}}
 =\frac{q}{(\beta^2p q^P)^{1/P}}
 =(\beta^2p)^{-1/P}
 =\frac{q'}{(\beta^2p(q')^P)^{1/P}}
 =\frac{Q(T)}{T^{1/P}}.
 \label{app:ig:endpoint-ratio}
\end{equation}
This contradicts \eqref{app:ig:ratio-decrease}.  Hence a Parisi measure
cannot have a bounded internal support gap \((q,q')\) with
\(0<q<q'<1\) and both endpoints in its support.  The condition \(q>0\)
is essential here: the argument does not exclude an isolated support
point at the origin.

\end{document}